\documentclass[12pt]{article}
\usepackage[T1]{fontenc}
\usepackage[utf8]{inputenc}
\usepackage{lmodern}
\usepackage[margin=1in]{geometry}
\usepackage{amsmath,amssymb,amsthm,mathtools,bm}
\usepackage{mathrsfs}
\usepackage{setspace}
\usepackage{booktabs,array,enumitem,microtype}
\usepackage{needspace}
\usepackage{placeins}
\usepackage{graphicx}
\usepackage{comment}
\usepackage{xcolor}
\usepackage[authoryear]{natbib}
\usepackage[colorlinks=true,linkcolor=blue,citecolor=blue,urlcolor=blue]{hyperref}
\usepackage[nameinlink,noabbrev]{cleveref}
\hypersetup{pdftitle={Sharp Asymptotic Theory of Maximum Likelihood Estimation for Gaussian Processes with a Matern Kernel},pdfsubject={Unknown smoothness, sub-Gaussian measurement errors, minimax risk, and uncertainty quantification}}
\allowdisplaybreaks[2]
\makeatletter
\renewcommand{\l@subsection}{\@dottedtocline{2}{1.5em}{2.8em}}
\renewcommand{\l@subsubsection}{\@dottedtocline{3}{4.3em}{3.2em}}
\makeatother
\setlist[enumerate]{leftmargin=*,itemsep=3pt,topsep=5pt}
\setlist[itemize]{leftmargin=*,itemsep=3pt,topsep=5pt}
\numberwithin{equation}{section}
\newtheorem{theorem}{Theorem}[section]

\newtheorem{corollary}[theorem]{Corollary}
\newtheorem{lemma}[theorem]{Lemma}
\theoremstyle{definition}
\newtheorem{definition}{Definition}
\newtheorem{assumption}{Assumption}

\crefname{assumption}{assumption}{assumptions}
\Crefname{assumption}{Assumption}{Assumptions}
\newcommand{\R}{\mathbb R}

\newcommand{\Z}{\mathbb Z}
\newcommand{\E}{\mathbb E}
\newcommand{\PP}{\mathbb P}
\newcommand{\T}{\top}

\newcommand{\tr}{\operatorname{tr}}
\newcommand{\diag}{\operatorname{diag}}

\newcommand{\Cov}{\operatorname{Cov}}
\newcommand{\Var}{\operatorname{Var}}
\newcommand{\op}{\mathrm{op}}

\newcommand{\Frob}{\mathrm F}
\newcommand{\KL}{\operatorname{KL}}
\newcommand{\TV}{\operatorname{TV}}
\newcommand{\Aff}{\operatorname{Aff}}
\newcommand{\argmax}{\operatorname*{arg\,max}}
\newcommand{\calI}{\mathcal I}

\newcommand{\norm}[1]{\left\lVert#1\right\rVert}

\DeclareMathOperator{\interior}{int}

\title{Sharp Asymptotic Theory of Maximum Likelihood Estimation for Gaussian Processes with a Mat\'ern Kernel}
\author{Ameer Qaqish$^1$,~Sudipto Banerjee$^2$ and Didong Li$^1$\\
$^1$Department of Biostatistics, University of North Carolina at Chapel Hill\\
$^2$Department of Biostatistics, University of California, Los Angeles}
\date{}
\begin{document}
\maketitle

\begin{abstract}
Gaussian processes (GPs) provide a flexible framework for modeling complex dependence structures through a covariance function, or kernel, whose parameters determine the variability, correlation, and smoothness of the underlying process. The Mat\'ern kernel is among the most popular choices due to its flexibility in modeling different degrees of smoothness. Although the identifiability of Mat\'ern kernel parameters under fixed-domain asymptotics is completely understood, identifiability alone does not establish the consistency or asymptotic behavior of the maximum likelihood estimators (MLEs) routinely used in practice. In this paper, we address this gap by developing a sharp asymptotic theory for the joint MLE of Mat\'ern kernel parameters. Since parameter identifiability depends critically on the dimension of the observation domain, we characterize the estimable parameters and establish their convergence rates across different dimensions. We establish joint asymptotic normality and minimax optimality of the convergence rates. In addition, we establish asymptotically valid confidence regions for the estimable parameters. Our results accommodate irregular sampling designs and extend to sub-Gaussian measurement errors. Extensive simulation studies illustrate the theoretical results. 
\end{abstract}

{\singlespacing
\tableofcontents
}
\section{Introduction}\label{sec:intro}

Gaussian processes (GPs) provide a flexible framework for modeling
complex dependence structures and have been extensively used in
spatial statistics, machine learning, Bayesian optimization, and
scientific applications
\citep{williams2006gaussian,stein1999interpolation,garnett2023bayesian}.
A fundamental component of a GP model is its covariance function,
also known as the kernel, which determines the dependence between
observations and the regularity of the underlying random function.
In practice, the kernel is typically chosen from a parametric family,
and its parameters are estimated by maximizing the Gaussian likelihood.
These estimates play an important role in prediction, uncertainty
quantification, and statistical inference.

The Mat\'ern kernel is a particularly popular choice among covariance kernels because it models the smoothness of the underlying process~\citep{stein1999interpolation}. Here, we consider the Mat\'ern kernel of the form
\begin{equation}\label{eq:matern-intro}
K_{\sigma^2,l,\nu}(x,x')
=
\sigma^2\frac{2^{1-\nu}}{\Gamma(\nu)}
\left(\frac{\|x-x'\|}{l}\right)^\nu
\mathsf K_\nu\left(\frac{\|x-x'\|}{l}\right),~~x,x'\in\mathbb{R}^d,
\end{equation}
where $\mathsf K_\nu$ is the modified Bessel function of the second kind, $\sigma^2$ is the spatial variance, $l$ is the length-scale, and $\nu$ is the smoothness parameter. As $\|x-x'\| \to 0^{+}$, the Mat\'ern covariance function converges to $\sigma^2$. Unlike kernels that fix smoothness, the Mat\'ern family can model random functions with different levels of regularity by varying $\nu$. Observations at $n$ locations are often modeled as a partial realization of a latent GP 
with i.i.d. measurement errors of variance $\tau^2$, commonly referred to as the nugget. The covariance between observations is
\[
\Cov(Y_{n,i},Y_{n,j}\mid x_{n,1:n})
=K_{\sigma^2,l,\nu}(x_{n,i},x_{n,j})
   +\tau^2\mathbf 1_{\{i=j\}}.
\]
Thus repeated observations at the same location share the
latent field value but have independent measurement errors. The parameter vector $(\sigma^2,l,\nu,\tau^2)$ is estimated jointly by maximizing the Gaussian likelihood.

Despite the widespread use of the Mat\'ern kernel, statistical
inference for its parameters remains challenging, particularly
under fixed-domain asymptotics, where the observation locations
become increasingly dense within a fixed bounded region.
Unlike in standard estimation problems based on independent
observations, increasing the sample size in this setting provides
more information about a single realization of the underlying GP.
Consequently, the ability to estimate kernel parameters depends
on the covariance structure and the dimension of the observation
domain. A fundamental question is whether the parameters are
identifiable, or more precisely, microergodic, under fixed-domain
asymptotics. For the Mat\'ern kernel with fixed smoothness,
\citet{zhang2004} showed that, when the dimension $d\leq3$,
the spatial variance and lengthscale cannot be consistently
estimated separately, whereas 
$c=\sigma^2l^{-2\nu}$ is microergodic.
For $d>4$, \citet{anderes2010} established that all parameters are identifiable, and the critical case $d=4$
has been addressed more recently by \citet{pillai2026}, who showed that all parameters are identifiable.

However, identifiability does not establish an estimator's consistency or asymptotic distribution. Existing likelihood analyses for the Mat\'ern kernel impose restrictive assumptions: \citet{kaufman2013} studied noiseless joint estimation of spatial variance and lengthscale with known smoothness; \citet{chen2000measurement,tang2021} included a positive nugget under restricted covariance parameters and observation designs; and \citet{petit2025} analyzed joint amplitude and smoothness estimation in a noiseless periodized model. Consequently, the joint maximum likelihood estimators (MLEs) of spatial variance, lengthscale, smoothness, and nugget variance remain uncharacterized---leaving open their consistency for all identifiable parameters, convergence and optimal rates, and asymptotic normality.

In this paper, we address these questions by developing a sharp
asymptotic theory for the joint MLE of Mat\'ern kernel parameters.
Our results establish consistency, convergence rates, and joint
asymptotic normality for all identifiable parameters across
different dimensions. In particular, when $d\leq3$, we establish
the asymptotic behavior of the estimators of the microergodic
parameter $c$, smoothness $\nu$, and nugget variance $\tau^2$,
while allowing all four parameters to be
optimized jointly.
When $d\geq4$, all four parameters are consistently estimable.
Write $\nu_0$ for the true smoothness, we show that the
estimators of $\nu$, $c$, and $\tau^2$ converge at rates 
\[
n^{-d/(4\nu_0+2d)},\qquad
(\log n)n^{-d/(4\nu_0+2d)},\qquad
n^{-1/2},
\]
respectively. For $d\geq4$, the spatial variance and lengthscale
converge at rate $(\log n)^{-1/2}$ when $d=4$, and at rate
$n^{-(d-4)/(4\nu_0+2d)}$ when $d>4$.
We further establish matching minimax lower bounds, showing
that these convergence rates cannot be improved in general.
Thus, the MLE achieves the minimax optimal rates for all
identifiable parameters.

Beyond consistency and asymptotic normality, we establish
convergence of moments and justify confidence intervals based
on estimated information matrices. Our analysis allows
sub-Gaussian measurement errors, which need not
follow a Gaussian distribution. In this setting, the Gaussian
likelihood remains the estimation criterion, and we show how
the asymptotic covariance is adjusted to account for
non-Gaussian errors. The theoretical results accommodate
deterministic and random observation designs, including
irregular sampling locations under suitable coverage conditions.
Numerical experiments investigate the finite-sample behavior
of the estimators and their agreement with the asymptotic theory.

The primary challenge is characterizing the parameter information contained in noisy, densely sampled observations of a Mat\'ern GP. Because the Mat\'ern spectral density decays polynomially, recovering covariance parameter information requires different tools than for the RBF kernel~\citep{qaqishli2026rbf}. Using localized Fourier coefficients and Sobolev sampling arguments to quantify this information, we identify explicit coordinates in which the Fisher information is $\asymp I$ and establish that the MLE error in these coordinates is $O_P(1)$, yielding consistency. A quadratic log-likelihood approximation in these coordinates and a central limit theorem for the normalized score then yield asymptotic normality of the MLEs. Our results hold for any choice of global MLE, including when the maximizer is non-unique.

The remainder of the paper is organized as follows.
\Cref{sec:back} introduces the GP model and reviews related
work on Mat\'ern parameter inference.
\Cref{sec:main} presents the main theoretical results,
including consistency, convergence rates, asymptotic normality,
minimax optimality, and confidence intervals.
\Cref{sec:simulation} reports simulation studies, and
\Cref{sec:discussion} discusses limitations and future
directions. All proofs and additional simulation details
are provided in the Appendix.

\section{Background and related work}\label{sec:back}

\subsection{Background}\label{sec:gp-background}\
\begin{definition}[Gaussian process]
A process $\{F(x):x\in D\}$ is a Gaussian process, written $F\sim\mathrm{GP}(m,K)$, if every finite vector of its evaluations is multivariate normal, with mean entries $m(x_i)$ and covariance entries
$K(x_i,x_j)$.
\end{definition}

Throughout this paper, we consider a zero-mean GP
$F\sim\mathrm{GP}(0,K)$ on a fixed bounded domain
$D\subset\mathbb R^d$ containing a cube with nonempty interior,
where $K$ is the Mat\'ern kernel defined in
\eqref{eq:matern-intro}. We observe
\begin{equation}\label{eq:model}
Y_{n,i}=F(x_{n,i})+\varepsilon_{n,i},
\qquad i=1,\ldots,n,
\end{equation}
at locations $x_{n,1},\ldots,x_{n,n}\in D$, where
$\varepsilon_{n,i}=\tau\xi_{n,i}$,
$\xi_{n,i}\stackrel{\mathrm{iid}}{\sim}P$,
$\mathbb E_P[\xi]=0$, and $\mathbb E_P[\xi^2]=1$.
The measurement errors are independent of $F$ and of the
observation locations. When the locations are random, they are
sampled independently of $F$. We allow non-Gaussian errors satisfying the conditions in \Cref{ass:noise}. The error distribution $P$ may also vary with $n$.

Let $Y_n=(Y_{n,1},\ldots,Y_{n,n})^\top$.
Conditional on the observation locations, its covariance matrix is 
\begin{equation*}
C_n(\phi)=\sigma^2 \mathsf{R}_n(l,\nu)+\tau^2I_n,
\qquad
\phi=(\sigma^2,l,\nu,\tau^2).
\end{equation*}
Here $[\mathsf R_n(l,\nu)]_{ij}
=K_{1,l,\nu}(x_{n,i},x_{n,j})$ is the Mat\'ern correlation
between sites $i$ and $j$.
When the measurement errors are Gaussian,
$Y_n\mid x_{n,1:n}\sim N_n(0,C_n(\phi))$.
For non-Gaussian errors, $Y_n$ is generally not Gaussian,
but its covariance remains $C_n(\phi)$.
We estimate the kernel parameters by maximizing the
Gaussian likelihood based on this covariance matrix (see \Cref{sec:setup}).

We focus on fixed-domain asymptotics, where $D$ remains
fixed as the number of observations $n$ increases.
The sampling locations may be irregular and need not
form a regular grid. Repeated locations are also permitted,
with independent measurement errors for each observation.
The precise sampling conditions are introduced in
\Cref{sec:setup}. Our analysis makes use of the spectral representation of the Mat\'ern kernel. Throughout, we adopt the Fourier transform convention 
\begin{equation*}
\widehat K(\omega)
=
\int_{\mathbb R^d}
K(x)e^{-2\pi i\omega^\top x}\,dx,
\qquad
K(x)
=
\int_{\mathbb R^d}
\widehat{K}(\omega)e^{2\pi i\omega^\top x}\,d\omega.
\end{equation*}
For the Mat\'ern covariance, regarded as a function of the lag,
$f=\widehat K$ is its spectral density.
The Mat\'ern spectral density is given
in \eqref{eq:spectrum}.

\subsection{Related work on Mat\'ern parameter inference}
\label{sec:related}

A fundamental inferential challenge under fixed-domain asymptotics is determining the identifiable parameters in the Mat\'ern family. This question concerns the equivalence and orthogonality of the GP laws induced by different covariance parameters. A parameter is 
identifiable if different values of the parameter imply orthogonality of the GP laws. For known smoothness $\nu$ and $d\le3$, \citet{zhang2004} showed that two Mat\'ern GPs induce equivalent Gaussian measures if and only if they have the same microergodic parameter $c=\sigma^2l^{-2\nu}$. Consequently, the spatial variance $\sigma^2$ and lengthscale $l$ cannot be consistently estimated separately. In contrast, \citet{anderes2010} established that these parameters can be distinguished when $d>4$. More recently, \citet{pillai2026} resolved the critical case $d=4$ by establishing orthogonality of Gaussian measures with different lengthscales, even when their microergodic parameters agree. Together, these results fully characterize the dimension-dependent identifiability of the Mat\'ern covariance parameters. 

However, identifiability alone does not establish consistency of a particular estimator. Several studies have constructed consistent estimators of Mat\'ern parameters, under various assumptions. For example, for noiseless observations, \citet{anderes2010} developed estimators based on higher-order increments to separate spatial variance and lengthscale in dimensions greater than four. For noiseless observations,
\citet{loh2021fixed} constructed estimators of microergodic parameters, including smoothness, using higher-order quadratic variations under several irregular sampling designs. Extending this approach to observations with a nugget, \citet{loh2023estimating} established consistency and convergence rate bounds for estimators of the smoothness, microergodic
parameter, and nugget variance. These results demonstrate that consistent estimation is possible in several settings, including when smoothness and nugget variance are unknown. However, the estimators are specifically constructed using quadratic variations rather than maximizing the Gaussian likelihood, so their results do not characterize the MLE.

A separate line of research has investigated MLEs under fixed-domain asymptotics. With known smoothness and no nugget, \citet{zhang2004} established consistency of the MLE of the microergodic parameter when the lengthscale is fixed. \citet{du2009fixed} further established its asymptotic normality for observations on a bounded interval.  \citet{kaufman2013} extended these results by establishing consistency and asymptotic normality of the MLE of the microergodic parameter when the spatial variance and lengthscale are jointly estimated. \citet{chen2000measurement} studied the one-dimensional Ornstein--Uhlenbeck process, corresponding to a Mat\'ern kernel with $\nu=1/2$, under equally spaced observations. They established consistency and asymptotic normality, revealing that the signal parameter converges at a slower rate than the nugget variance. More generally, \citet{tang2021} studied maximum likelihood
estimation of the microergodic parameter and nugget variance
with known smoothness and a fixed lengthscale, under additional
sampling and eigenvalue assumptions. Although these studies provide important results for
likelihood-based inference, they assume known smoothness and do not establish the asymptotic behavior of the unrestricted joint MLE of all four Mat\'ern parameters with a positive nugget.

Other studies have considered likelihood estimation in
related but different settings. \citet{petit2025}
investigated joint estimation of amplitude and smoothness for a noiseless periodized Mat\'ern model, while \citet{karvonen2023} studied smoothness estimation in GP interpolation with deterministic response functions. These results provide additional insights into smoothness estimation, but their observation models differ from the noisy Euclidean Mat\'ern model considered here. 

Existing literature establishes the dimension-dependent identifiability of Mat\'ern parameters, consistency of specially constructed estimators, and MLE asymptotics in restricted settings. However, a general asymptotic theory for the joint MLE of the spatial variance, lengthscale, smoothness, and nugget variance remains unavailable---a critical gap since maximum likelihood is the standard approach for estimating GP kernel parameters, as in \texttt{GPML}~\citep{rasmussen2010gaussian} and \texttt{scikit-learn}~\citep{pedregosa2011scikit}. What remains open is whether the joint MLE is consistent for every identifiable parameter, what rates and limiting distribution it has, and whether any estimator can achieve faster rates uniformly over the parameter class. We address these questions in \Cref{sec:main}.

\section{Theory}\label{sec:main}

\subsection{Model and assumptions}\label{sec:setup}
For the theoretical analysis, it is convenient to reparameterize the Mat\'ern kernel through its spectral density. In particular, we introduce the log spectral amplitude $\eta$,
the squared inverse lengthscale $a$, and the nugget variance $t$, and write $\theta=(\eta,\nu,a,t)$, where
 \begin{equation}\label{eq:chart}
 (\eta,\nu,a,t)
 =\left(\log\{A_d(\nu)\sigma^2l^{-2\nu}\},\ \nu,\ l^{-2},\ \tau^2\right),\qquad A_d(\nu)=\frac{(4\pi)^{d/2}\Gamma(\nu+d/2)}{\Gamma(\nu)}.
\end{equation}
These coordinates give a simple form for the Mat\'ern spectral density:
\begin{equation}\label{eq:spectrum}
f_\theta(\omega)=e^\eta(a+4\pi^2\norm\omega^2)^{-\beta},
 \qquad \beta=\nu+d/2.
\end{equation}
Under the observation model \eqref{eq:model}, let $S_n(\eta,\nu,a)$ be the covariance matrix of the underlying Mat\'ern GP evaluated at the observed locations. Conditional on these locations, the
covariance of $Y_n$ is $C_n(\theta)=S_n(\eta,\nu,a)+tI_n$ and $C_{n,0}=C_n(\theta_0)$. We estimate $\theta$ by maximizing the Gaussian likelihood over a parameter space $\Theta$. Ignoring the additive constant, the log-likelihood is
\begin{equation}\label{eq:likelihood}
 \ell_n(\theta)=-\frac12\left\{\log\det C_n(\theta)
                      +Y_n^\top C_n(\theta)^{-1}Y_n\right\}.
\end{equation}
Throughout, we use MLEs to mean
any measurable choice
$\widehat\theta_n\in\argmax_{\theta\in\Theta}\ell_n(\theta)$,
where $\ell_n$ is the Gaussian log-likelihood in \eqref{eq:likelihood}. This is the exact log-likelihood when the measurement errors are
Gaussian. For non-Gaussian errors, we maximize the same Gaussian
criterion, whose covariance model remains correctly specified.
The parameter space $\Theta$ and compact set of true parameters
$\mathcal K$ are specified below.

Since the identifiable parameters depend on the domain dimension $d$ (see \Cref{sec:related}), we define
\[
 \vartheta=
 \begin{cases}
 (\eta,\nu,t),&d\le3,\\
 (\eta,\nu,a,t),&d\ge4,
 \end{cases}
 \qquad k_d=\dim(\vartheta)=\begin{cases}
 3,&d\le3,\\
 4,&d\ge4.\end{cases}
\]

\begin{assumption}[Parameter space and interior truths]\label{ass:param}
The likelihood is maximized over a fixed compact set $\Theta\subset\R\times(0,\infty)^3$ in coordinates $\theta=(\eta,\nu,a,t)$. Let $\mathcal K\subset\Theta$ be a fixed compact set of true parameters. For $d\ge4$, assume $\mathcal K\subset\operatorname{int}(\Theta)$. For $d\le3$, there is an open neighborhood $U$ of the projection of $\mathcal K$ onto the $\vartheta$ coordinates such that the feasible range values at each $\vartheta\in U$ form an interval $I(\vartheta)\coloneqq \{a:(\vartheta,a)\in\Theta\} =[\underline a(\vartheta),\overline a(\vartheta)]$ with $0<\underline a(\vartheta)<\overline a(\vartheta)$. Here $(\vartheta,a)$ denotes the corresponding permutation of $(\eta,\nu,a,t)$, and both endpoint functions are Lipschitz on $U$. The true range may lie at either endpoint.
\end{assumption}

The log-likelihood is continuous on the compact set $\Theta$, so a global maximizer exists. In some cases, there may be multiple maximizers. Our results apply to any measurable choice of global maximizer. The likelihood may be optimized in any one-to-one coordinates, provided the same parameter set is transformed with them. For example, a fixed box in $(\sigma^2,l,\nu,\tau^2)$ satisfies Assumption~\ref{ass:param} provided the true parameters are interior when $d\ge4$, or the true identifiable parameters are interior to the set of identifiable values allowed by the box when $d\le3$. When $d\le3$, the true range may lie at an endpoint of its
feasible interval.

We consider two sampling conditions. The first requires the
observation locations to cover a fixed interior cube at a
polynomial rate and gives preliminary polynomial consistency
of the amplitude, smoothness, and nugget estimators. The second
requires the fastest possible deterministic coverage rate on an interior cube or
independent random sampling with a positive density lower bound on an interior
cube. Under the second assumption, we obtain sharp estimation rates
and asymptotic normality. Neither condition requires a minimum
separation between observation locations.

\begin{assumption}[Polynomial coverage]\label{ass:coverage}
There exist a fixed cube $Q\subset D$ and constants $C>0$ and $\gamma>0$ such that, for all sufficiently large $n$,
\begin{equation*}
 h_n(Q)\coloneqq \sup_{x\in Q}\min_{1\le i\le n}\norm{x-x_{n,i}}
 \le Cn^{-\gamma}.
\end{equation*}
\end{assumption}
Since the closed balls at $x_{n,i}$ with radius $h_n(Q)$ cover $Q$, necessarily $\gamma\le1/d$. The following assumption requires the fastest possible convergence rate for deterministic designs and also accommodates independent random sampling.

\begin{assumption}[Sampling for sharp rates]\label{ass:sharp}
Either the deterministic sites satisfy $h_n(Q)\le Cn^{-1/d}$ on a
fixed interior cube, or the sites are iid in $D$ with a density bounded
below by a positive constant on such a cube, independently of the
field and errors.
\end{assumption}

For random designs, $h_n(Q)=O\{(\log n/n)^{1/d}\}$ almost surely.
This bound includes a logarithmic factor relative to the deterministic coverage condition, but the same sharp estimation rates hold.

\begin{assumption}[Uniform sub-Gaussian measurement errors]\label{ass:noise}
The standardized measurement errors are iid with law $P\in\mathcal P$, independently of the underlying GP and observation locations. We assume that
\begin{equation*}
 \E_P\xi=0,\qquad \E_P\xi^2=1,\qquad
 \sup_{P\in\mathcal P}\E_Pe^{\xi^2/K_\varepsilon^2}\le2,\qquad
 \inf_{P\in\mathcal P}(\E_P\xi^4-1)\ge v_*
\end{equation*}
for fixed $K_\varepsilon<\infty$ and $v_*>0$. Uniform assertions allow the common law to vary with $n$ and the truth. No symmetry or density is required.
\end{assumption}

The sub-Gaussian condition provides uniform control of the error moments, while the last condition bounds $\Var_P(\xi^2)=\E_P\xi^4-1$ away from zero, preventing degeneracy of the asymptotic variance of the nugget estimator.
Let $C_i=\partial_{\vartheta_i}C_n(\theta)|_{\theta=\theta_0}$, holding $a=a_0$ fixed when $d\le3$. Conditional on the observation
locations, define
\begin{equation}\label{eq:info}
 \calI_n=\E_0\!\left[
       -\nabla_\vartheta^2\ell_n(\theta_0)\mid x_{n,1:n}\right],\quad
 [\calI_n]_{ij}=\frac12\operatorname{tr}(C_{n,0}^{-1}C_iC_{n,0}^{-1}C_j),\quad 1\le i,j\le k_d.
\end{equation}
Here $\E_0$ is expectation under the true response law. The expected negative Hessian depends only on the covariance matrix of $Y_n$, so $\calI_n$ is independent of the error fourth moment and equals the Fisher information in the Gaussian experiment. For non-Gaussian errors, however, it need not equal the covariance matrix of the score. When $d\leq 3$, $\calI_n$ is the $(\eta,\nu,t)$ block of the full information matrix at $a_0$.

\subsection{MLE asymptotics and convergence rates}\label{sec:mle}
We first introduce the following sequences, which determine the convergence rates of the MLEs:
\[
R_n=n^{1/(2\nu_0+d)},~
L_n=\log R_n,~
r_n=R_n^{d/2}=n^{d/(4\nu_0+2d)},~
b_n=
\begin{cases}
(\log n)^{1/2},&d=4,\\
R_n^{(d-4)/2}=n^{(d-4)/(4\nu_0+2d)},&d>4.
\end{cases}
\]
Unlike the nugget and range, the spectral amplitude and smoothness have strongly coupled estimation errors. To account for this dependence, we introduce a coordinate transformation from $(\eta - \eta_0, \nu - \nu_0)$ into the centered amplitude and smoothness errors
\begin{equation}\label{eq:shear}
 u_n(\theta)=\eta-\eta_0-2L_n(\nu-\nu_0),\qquad v=\nu-\nu_0.
\end{equation}
Subtracting $2L_n(\nu-\nu_0)$ removes the component of the amplitude
error that grows logarithmically relative to the smoothness error. We define the normalized estimation errors of the identifiable parameters as

\begin{equation*}
 h_n(\theta)=
 \begin{cases}
 (r_nu_n,r_nv,\sqrt n(t-t_0))^\T,&d\le3,\\
 (r_nu_n,r_nv,b_n(a-a_0),\sqrt n(t-t_0))^\T,&d\ge4.
 \end{cases}
\end{equation*}
Note that
$h_n(\theta)=E_n^{-1}(\vartheta-\vartheta_0)$, where

\begin{equation}\label{eq:E}
E_n=\begin{cases} \begin{pmatrix}r_n^{-1}&2L_nr_n^{-1}&0\\0&r_n^{-1}&0\\0&0&n^{-1/2}\end{pmatrix} & d\leq 3,\\
\begin{pmatrix}r_n^{-1}&2L_nr_n^{-1}&0&0\\0&r_n^{-1}&0&0\\0&0&b_n^{-1}&0\\0&0&0&n^{-1/2}\end{pmatrix} &d\geq 4.
\end{cases}
\end{equation}
The matrix $E_n$ depends on the true parameter and was chosen so that the Gaussian information in $h$ coordinates satisfies $E_n^{\top}\mathcal{I}_nE_n \asymp I$. Theorem~\ref{thm:normalized} shows that the normalized MLE $\widehat h_n\coloneqq h_n(\widehat\theta_n)$ is $O_P(1)$ in these coordinates. To characterize the asymptotic distribution of $\hat{h}_n$, we define
\begin{equation}\label{eq:normalizedscoreinfo}J_n=E_n^\T\calI_nE_n,
~\Delta_n=E_n^\T\nabla_\vartheta\ell_n(\theta_0),
 ~K_n=\Cov_{\theta_0,P}(\Delta_n\mid x_{n,1:n}),
~V_n=J_n^{-1}K_nJ_n^{-1}.
\end{equation}
Here $J_n$ is the Gaussian information in $h$ coordinates, $\Delta_n$ is the score in $h$-coordinates at $h=0$, $K_n$ is the covariance matrix of this score, and $V_n$ is the asymptotic covariance matrix of the normalized estimator $\hat{h}_n$. For Gaussian measurement errors, the model is correctly specified, so $K_n=J_n$ and $V_n=J_n^{-1}$. For non-Gaussian errors, $K_n$ and $J_n$ need not be equal.

Our results are uniform over $P\in\mathcal P$ and $\theta_0\in\mathcal K$ under compactness and local feasibility in Assumption~\ref{ass:param}. In low dimensions, uniformity includes true ranges at either endpoint of $I(\vartheta_0)$. Our first result establishes asymptotic normality of the normalized MLE $\widehat h_n$ and its asymptotic covariance matrix.

\Needspace{13\baselineskip}
\begin{theorem}[$\hat{h}_n$ asymptotics]\label{thm:normalized}
Under \Cref{ass:param,ass:sharp,ass:noise}, let $\widehat\theta_n$ be
any measurable global maximizer of the Gaussian criterion over $\Theta$.
Write $\widehat h_n=h_n(\widehat\theta_n)$ and
$\kappa_4(P,t_0)=t_0^2(E_P\xi^4-3)$.
Then the following conclusions hold uniformly over
$\theta_0\in\mathcal K$ and $P\in\mathcal P$.
\begin{enumerate}[label=(\roman*),leftmargin=*,itemsep=8pt]
\item \textbf{Expansion and normal approximation.}
\begin{equation}\label{eq:Jbounds}
 cI_{k_d}\preceq J_n,K_n,V_n\preceq CI_{k_d},
\end{equation}
and
\begin{equation}\label{eq:MLEexpansion}
 \widehat h_n=J_n^{-1}\Delta_n+o_{\PP}(1),\quad
 K_n^{-1/2}\Delta_n\Rightarrow N_{k_d}(0,I_{k_d}),\quad
 V_n^{-1/2}\widehat h_n\Rightarrow N_{k_d}(0,I_{k_d}).
\end{equation}
Convergence of laws is measured in bounded Lipschitz distance.

\item \textbf{Covariance structure.}
Write $e_t$ for the last coordinate vector, the covariance has the
following structure. Let $J_{n,\mathrm{sig}}$ be the principal block
of $J_n$ for the signal coordinates $(u,v)$ when $d\le3$ and
$(u,v,a)$ when $d\ge4$.
\begin{subequations}
\begin{equation}\label{eq:kurtosis-correction}
 V_n=J_n^{-1}+\kappa_4(P,t_0)e_te_t^\T+o(1).
\end{equation}
Equivalently,
\begin{equation}\label{eq:covariance-blocks}
 V_n=\begin{pmatrix}
       J_{n,\mathrm{sig}}^{-1} & 0\\
       0&t_0^2(E_P\xi^4-1)
     \end{pmatrix}+o(1).
\end{equation}
At $d=4$, the signal block further separates range from centered
amplitude and smoothness:
\begin{equation}\label{eq:critical-covariance-block}
J_{n,\mathrm{sig}}^{-1}
 =\begin{pmatrix}J_{n,(u,v)}^{-1}&0\\0&J_{n,aa}^{-1}\end{pmatrix}+o(1).
\end{equation}
Here $J_{n,(u,v)}$ is the principal two-by-two block for $(u,v)$,
and $J_{n,aa}$ is the normalized range entry. All $o(1)$ terms in
these covariance formulas vanish in operator norm, uniformly.
\end{subequations}

\item \textbf{Uniform tails and moments.}
There are uniform constants $c,C,\delta>0$ such that
\begin{equation}\label{eq:tails}
 \sup_{\theta_0\in\mathcal K,\ P\in\mathcal P}
 P_{\theta_0,P}(\|\widehat h_n\|>H)
 \le Ce^{-cH^2}+Ce^{-cn^\delta},\qquad H\ge0.
\end{equation}
For every fixed $\lambda>0$ and $0<\alpha<\min\{2,2\delta\}$,
\begin{equation}\label{eq:allmoments}
 \sup_{n\ge n_0}\sup_{\theta_0\in\mathcal K,\ P\in\mathcal P}
 E_{\theta_0,P}\exp\{\lambda\|\widehat h_n\|^\alpha\}<\infty.
\end{equation}
In particular, every fixed polynomial moment is uniformly bounded.
For every fixed continuous function $F$ satisfying
$|F(x)|\le C_F e^{\lambda\|x\|^\alpha}$ for some such
$\lambda,\alpha$ and finite $C_F$, and $Z\sim N_{k_d}(0,I)$,
\begin{equation}\label{eq:momentcomparison}
 \sup_{\theta_0\in\mathcal K,\ P\in\mathcal P}
 \left|E_{\theta_0,P}F(\widehat h_n)-E_ZF(V_n^{1/2}Z)\right|
 \longrightarrow0.
\end{equation}
Consequently $E\widehat h_n\to0$ and
$E(\widehat h_n\widehat h_n^\top)-V_n\to0$ uniformly.
For every fixed continuous $\varphi$ satisfying the same growth condition,
\begin{equation}\label{eq:normalizedmoments}
 \sup_{\theta_0,P}\left|E\varphi(V_n^{-1/2}\widehat h_n)
                              -E\varphi(Z)\right|\longrightarrow0.
\end{equation}

\item \textbf{Deterministic and random locations.}
For fixed designs these statements concern the response law at the
given locations. For random designs, on one event of location probability
one, all statements hold conditionally on the realized sequence of
locations, uniformly over $\mathcal K\times\mathcal P$, eventually with
deterministic constants. In \eqref{eq:momentcomparison}, $V_n$ then
uses those locations and $E_Z$ averages only over $Z$.
The stochastic expansion, standardized normal approximation, tail and
moment bounds, and standardized moment convergence also hold under the
joint law of locations and responses. For the joint-law assertions, if $J_n$ is singular, define
$J_n^{-1}\Delta_n$, $K_n^{-1/2}\Delta_n$, and
$V_n^{-1/2}\widehat h_n$ to be zero.
\end{enumerate}
\end{theorem}

To first order, the signal estimators have the same covariance as under Gaussian measurement errors. The error distribution enters through the nugget MLE's asymptotic variance, which depends on $\Var_P(\xi^2)$. By \eqref{eq:covariance-blocks}, the nugget MLE satisfies
\begin{equation}\label{eq:lognugget}
 \frac{\sqrt n(\widehat t_n-t_0)}{t_0\sqrt{E_P\xi^4-1}}
      \Rightarrow N(0,1),\qquad
 \frac{\sqrt n(\log\widehat t_n-\log t_0)}
      {\sqrt{E_P\xi^4-1}}\Rightarrow N(0,1).
\end{equation}
These are also the limits of the Gaussian likelihood variance
estimator $n^{-1}\sum_{i=1}^n\varepsilon_{n,i}^2$ and its
logarithm if the errors
$\varepsilon_{n,i}=\sqrt{t_0}\xi_i$, $\xi_i\sim P$,
were observed directly.

We now use \Cref{thm:normalized} and the delta method to obtain the convergence rates of the MLEs in various common parameterizations of the kernel.

\begin{corollary}[Uniform estimation rates]\label{cor:rates}
Under Theorem~\ref{thm:normalized}, let $g$ be a parameter with a
rate listed in Table~\ref{tab:rates}, and denote that rate by
$\rho_{g,n}(\nu_0)$. Then
\[
 g(\widehat\theta_n)-g(\theta_0)
   =O_{P,\mathcal K\times\mathcal P}\{\rho_{g,n}(\nu_0)\},
 \qquad
 \left\{\E_{\theta_0,P}
   |g(\widehat\theta_n)-g(\theta_0)|^2\right\}^{1/2}
   \asymp\rho_{g,n}(\nu_0),
\]
uniformly over $\theta_0\in\mathcal K$ and $P\in\mathcal P$.
The centered amplitude error
$\widehat\eta_n-\eta_0-2L_n(\widehat\nu_n-\nu_0)$ has rate
$r_n^{-1}$ in both senses. The rates hold for deterministic locations
and, for random locations, both conditionally for almost every
location sequence and under the joint law, as in
Theorem~\ref{thm:normalized}.
\end{corollary}

\begin{table}[!htbp]\centering\small
\caption{Uniform stochastic and root mean squared error rates over
$\mathcal K\times\mathcal P$. Each positive parameter and its logarithm
have the same rate. ``NA'' denotes lack of 
consistent estimation over the full family.}\label{tab:rates}
\begin{tabular}{@{}lccc@{}}\toprule
Parameter & $d\le3$ & $d=4$ & $d>4$\\\midrule
$\nu,\ \log\nu$ & $r_n^{-1}$ & $r_n^{-1}$ & $r_n^{-1}$\\
$\eta,\ c,\ \log c$ & $(\log n)/r_n$ & $(\log n)/r_n$ & $(\log n)/r_n$\\
$a,\ \log a,\ l,\ \log l$ & NA & $(\log n)^{-1/2}$ & $b_n^{-1}$\\
$\sigma^2,\ \log\sigma^2$ & NA & $(\log n)^{-1/2}$ & $b_n^{-1}$\\
$\tau^2,\ \log\tau^2$ & $n^{-1/2}$ & $n^{-1/2}$ & $n^{-1/2}$\\\bottomrule
\end{tabular}
\end{table}
\FloatBarrier

For $d\ge4$, the faster $(\log n)/r_n$ rate for
$c=\sigma^2l^{-2\nu}$, compared with the $b_n^{-1}$ rates for
$\sigma^2$ and $l$, results from cancellation between the variance
and lengthscale errors.

\subsection{Minimax optimality}\label{sec:minimax}

Next, we show that the rates in Corollary~\ref{cor:rates} are minimax optimal.
For any parameter $g$ with a rate listed in Table~\ref{tab:rates},
write $\rho_{g,n}(\nu)$ for that rate. Here the worst-case risk
takes the supremum over both the covariance parameter $\theta$
and the error law $P$.

\begin{theorem}[Minimax squared-error rates]\label{thm:minimax}
Assume \Cref{ass:param,ass:sharp,ass:noise}. In addition, suppose that
$N(0,1)\in\mathcal P$. Let
$\mathcal K\Subset\interior(\Theta)$ be a compact rectangle with
nonempty interior and $\bar\nu=\max_{\theta\in\mathcal K}\nu$.
For every coordinate $g$ specified above,
\begin{equation*}
 \inf_{T_n}\sup_{\theta\in\mathcal K,P\in\mathcal P}
       E_{\theta,P}|T_n-g(\theta)|^2
 \asymp\sup_{\theta\in\mathcal K,P\in\mathcal P}
       E_{\theta,P}|g(\widehat\theta_n)-g(\theta)|^2
 \asymp\rho_{g,n}(\bar\nu)^2.
\end{equation*}
At each truth the global MLE has mean squared error of order
$\rho_{g,n}(\nu)^2$, uniformly on the compact set and over the noise
class. The lower bound also holds in shrinking neighborhoods: for
any fixed interior truth $\theta_0$ and fixed $H>0$, let
$\mathcal U_{n,H}(\theta_0)$ consist of the parameters $\theta\in\Theta$ satisfying
$\|h_n(\theta)\|\le H$, with $a=a_0$ additionally fixed when $d\le3$.
Then
\begin{equation}\label{eq:localminimax}
 \liminf_{n\to\infty}\inf_{T_n}
 \sup_{\theta\in\mathcal U_{n,H}(\theta_0),\ P\in\mathcal P}
 \frac{E_{\theta,P}|T_n-g(\theta)|^2}
      {\rho_{g,n}(\nu_0)^2}>0.
\end{equation}
Thus even within each such neighborhood, the worst-case mean squared
error cannot have a smaller order. All comparisons are up to fixed
multiplicative constants.

For random locations, these risk comparisons hold conditionally for almost every location sequence and under the joint law of locations and responses. In the latter case, the infimum is over estimators based on both the locations and responses.
\end{theorem}
The extra Gaussian-member assumption is used for the minimax lower
bound. The minimax result
identifies the optimal order of mean squared error, but it does not identify
the smallest attainable leading variance for each fixed error law.

\subsection{Asymptotic normality in common parameterizations}\label{sec:parameterizations}

Corollary~\ref{cor:rates} gives the coordinatewise estimation rates.
We now translate the normal approximation in Theorem~\ref{thm:normalized}
into the common parameterizations in Table~\ref{tab:parameterizations}. Write $\phi$
for a full vector in the table. When $d\ge4$, all its coordinates
are identifiable and $\zeta=\phi$. When $d\le3$, the indicated
three-coordinate block is $\zeta$ and the remaining coordinate
$\rho$ is its nuisance. A dash means the full vector does not contain
such an identifiable coordinate block. In particular, the
low-dimensional variance--lengthscale vectors must be reparameterized to expose
$c=\sigma^2l^{-2\nu}$ before applying the result below.
Every admissible vector has the form $\zeta=g(\vartheta)$ for some function $g$.

\begin{table}[htbp]\centering\small
\caption{Common full parameterizations and their identifiable block
and nuisance coordinate in dimensions at most three. All four
coordinates are identifiable when $d\ge4$.}
\label{tab:parameterizations}
\begin{tabular}{@{}lll@{}}\toprule
Parameterization & Full vector $\phi$ & $\zeta\ ;\ \rho$ for $d\le3$\\\midrule
Spectral & $(\eta,\nu,a,t)$ & $(\eta,\nu,t)\ ;\ a$\\
Spectral, log scale & $(\eta,\log\nu,\log a,\log t)$
  & $(\eta,\log\nu,\log t)\ ;\ \log a$\\
Microergodic & $(c,\nu,l,\tau^2)$ & $(c,\nu,\tau^2)\ ;\ l$\\
Microergodic, log scale & $(\log c,\log\nu,\log l,\log\tau^2)$
  & $(\log c,\log\nu,\log\tau^2)\ ;\ \log l$\\
Variance--lengthscale & $(\sigma^2,l,\nu,\tau^2)$ & ---\\
\shortstack[l]{Variance--lengthscale, log scale} & $(\log\sigma^2,\log l,\log\nu,\log\tau^2)$ & ---\\\bottomrule
\end{tabular}
\end{table}

\FloatBarrier

\begin{corollary}[Uniform inference in common parameterizations]\label{cor:log}
Under Theorem~\ref{thm:normalized}, choose an identifiable vector
$\zeta=g(\vartheta)$ from Table~\ref{tab:parameterizations}.
Define the linearized covariance by
\begin{equation*}
 \Sigma_n^\zeta=Dg(\vartheta_0)E_nV_nE_n^\top Dg(\vartheta_0)^\top.
\end{equation*}
Let $\calI_n^\phi$ be the full Gaussian Fisher information in the
chosen parameterization at the truth and actual locations, and let
$\calI_n^\zeta$ be its identifiable block (the whole matrix when
$d\ge4$). Then
\begin{equation}\label{eq:reporting-information}
 \calI_n^\zeta=Dg(\vartheta_0)^{-\top}\calI_n Dg(\vartheta_0)^{-1},
\end{equation}
and the comparison with inverse Gaussian information is
\begin{equation}\label{eq:reporting-information-correction}
 \Sigma_n^\zeta =(\calI_n^\zeta)^{-1}
       +\frac{\kappa_4(P,t_0)}n q_\zeta(t_0)^2e_te_t^\top+Dg(\vartheta_0)E_n\,o(1)\,E_n^\top Dg(\vartheta_0)^\top.
\end{equation}
Here $q_\zeta(t)=1$ if the last coordinate is $t=\tau^2$, and
$q_\zeta(t)=t^{-1}$ if it is $\log t$; $e_t$ selects that coordinate.
The matrix $o(1)$ tends to zero in operator norm uniformly.
For Gaussian errors, $\Sigma_n^\zeta=(\calI_n^\zeta)^{-1}$ exactly.
In low dimension, $(\calI_n^\zeta)^{-1}$ in
\eqref{eq:reporting-information-correction} may also be replaced by
$[(\calI_n^\phi)^{-1}]_{\zeta\zeta}$, with the same form of
vanishing remainder.

For every listed identifiable coordinate,
\begin{equation*}
 \frac{\widehat\zeta_{n,j}-\zeta_{0,j}}
      {\sqrt{[\Sigma_n^\zeta]_{jj}}}\Rightarrow N(0,1)
\end{equation*}
uniformly over $\mathcal K\times\mathcal P$. These standardized errors
satisfy the exponential tail bound of Theorem~\ref{thm:normalized},
with different constants and a positive exponent if needed. Convergence of expectations
also holds for every fixed continuous function of polynomial growth.

For every listed identifiable vector other than the original-scale
variance--lengthscale vector when $4\le d\le8$, the joint conclusion also holds:
\begin{equation}\label{eq:joint-reporting-normal}
  (\Sigma_n^\zeta)^{-1/2}(\widehat\zeta_n-\zeta_0)\Rightarrow N_{k_d}(0,I).
\end{equation}
The joint standardized errors have the same tail and polynomial-moment
convergence properties. Conditional and unconditional random-design
assertions have the meanings in Theorem~\ref{thm:normalized}.
\end{corollary}

For $4\le d\le8$, joint normal convergence on the original
variance--lengthscale coordinates cannot be asserted in general. Appendix~\ref{app:raw-public-failure} proves its failure when $\nu_0\ne1/2$, even with Gaussian errors. The joint conclusion holds on that scale when $d>8$.

\subsection{Consistency under polynomial coverage}
\label{sec:preliminary-consistency}
Consistency of the amplitude, smoothness, and nugget MLEs requires less coverage
than the sharp rates and normal approximation in
Theorem~\ref{thm:normalized}. The following result applies whenever
the sites cover an interior cube at any fixed polynomial rate.

\begin{theorem}[Uniform consistency under polynomial coverage]\label{thm:preliminary}
Assume \Cref{ass:param,ass:coverage,ass:noise}, with deterministic sites,
and let $\widehat\theta_n$ be any measurable global maximizer in
\eqref{eq:likelihood}. There are constants $c,C,\delta,\delta_1>0$
such that, for all sufficiently large $n$,
\begin{equation}\label{eq:preliminary-global-bound}
 \sup_{\theta_0\in\mathcal K,\ P\in\mathcal P}
 P_{\theta_0,P}\!\left\{
   \frac{|\widehat\eta_n-\eta_0|}{\log n}
   +|\widehat\nu_n-\nu_0|+|\widehat t_n-t_0|
       >Cn^{-\delta}\right\}
 \le Ce^{-cn^{\delta_1}}.
\end{equation}
Consequently $(\widehat\eta_n,\widehat\nu_n,\widehat t_n)$ is
uniformly consistent in every dimension, with
\[
 \widehat\eta_n-\eta_0
       =O_{P,\mathcal K\times\mathcal P}(n^{-\delta}\log n),
 \qquad
 |\widehat\nu_n-\nu_0|+|\widehat t_n-t_0|
       =O_{P,\mathcal K\times\mathcal P}(n^{-\delta}).
\]
These upper rates also hold in every fixed $L^p$ norm,
$1\le p<\infty$, uniformly over $\mathcal K\times\mathcal P$.
Replacing $\eta$ by $c$ or $\log c$, $\nu$ by $\log\nu$, and
$t=\tau^2$ by $\log\tau^2$ gives the same respective rates.

For the iid design alternative in Assumption~\ref{ass:sharp}, these
conclusions also hold conditionally for almost every sequence of
locations, eventually with deterministic constants, and under the
joint law of locations and responses.
\end{theorem}

\subsection{Practical uncertainty quantification}\label{sec:uncertainty}
Corollary~\ref{cor:log} gives marginal normal approximations
and, for the specified parameterizations, joint normal
approximations with covariance $\Sigma_n^\zeta$.
This covariance depends on the unknown truth and noise law.
We now construct a data-based estimate $\widehat\Sigma_n^\zeta$
for confidence intervals and joint confidence regions. 

Choose a full parameterization $\phi$ from
Table~\ref{tab:parameterizations} and fit the model in those coordinates
over the image of the same parameter set $\Theta$. For $d\ge4$, take
$\zeta=\phi$. For $d\le3$, use a row with an identifiable coordinate
block $\zeta$ and let $\rho$ be the remaining coordinate already in
that row: respectively $a,\log a,l$, or $\log l$. Thus, after a
permutation, $\phi=(\zeta,\rho)$. Variance and lengthscale are not separately consistently estimable in
the low-dimensional variance--lengthscale parameterizations, so
those rows require reparameterization before applying this construction.

Let $\widehat{\mathcal A}_n^\phi$ be either the full Gaussian Fisher
information or the full negative Hessian of the Gaussian criterion,
computed in the chosen full coordinates at the global fit. The Fisher
information uses the trace formula in \eqref{eq:info}, with derivatives
taken in $\phi$ and evaluated at $\widehat\phi_n$. Estimate the error
fourth cumulant with the corrected residual estimator
$\widetilde\kappa_{4,n}$ of Lemma~\ref{lem:kurtosis}. Use $q_\zeta$ and
$e_t$ as defined in Corollary~\ref{cor:log}. The covariance estimate is
\begin{equation}\label{eq:estimatedlogcov}
 \widehat\Sigma_n^\zeta
 =\left[\{\widehat{\mathcal A}_n^\phi\}^{-1}\right]_{\zeta\zeta}
 +\frac{\widetilde\kappa_{4,n}}n
       q_\zeta(\widehat t_n)^2e_te_t^\top.
\end{equation}
For $d\ge4$ the selected block is the entire inverse. For $d\le3$,
the full matrix is inverted before its identifiable block is selected. The added
cumulant term accounts for non-Gaussian errors, where the nugget MLE asymptotic variance depends on the fourth cumulant of the measurement errors.

\begin{corollary}[Uniform feasible Wald inference]\label{cor:uncertainty}
Under the assumptions of Theorem~\ref{thm:normalized}, use the
construction above. The matrix $\widehat\Sigma_n^\zeta$ is positive definite
with probability tending to one, uniformly, and
\begin{equation*}
 \frac{[\widehat\Sigma_n^\zeta]_{jj}}{[\Sigma_n^\zeta]_{jj}}
       \longrightarrow_P1.
\end{equation*}
Consequently, for every coordinate and fixed $0<\alpha<1$,
\begin{equation*}
 \sup_{\theta_0\in\mathcal K,\ P\in\mathcal P}
 \left|P_{\theta_0,P}\!\left\{
 \zeta_{0,j}\in\left[\widehat\zeta_{n,j}
 \mathbin{\pm}z_{1-\alpha/2}\sqrt{[\widehat\Sigma_n^\zeta]_{jj}}\right]
 \right\}-(1-\alpha)\right|\longrightarrow0.
\end{equation*}
For the vectors covered by the joint conclusion
\eqref{eq:joint-reporting-normal},
\begin{equation}\label{eq:studentizednormal}
 (\widehat\Sigma_n^\zeta)^{-1/2}(\widehat\zeta_n-\zeta_0)
       \Rightarrow N_{k_d}(0,I_{k_d}).
\end{equation}
The corresponding quadratic statistic tends to $\chi^2_{k_d}$.
For the original-scale variance--lengthscale vector in dimensions four through eight,
\eqref{eq:estimatedlogcov} gives valid marginal Wald intervals; joint
regions are obtained by transforming a region from any parameterization
with the joint normal limit. Conditional and unconditional random-design
assertions have the meanings in Theorem~\ref{thm:normalized}.
\end{corollary}

If the covariance estimate is undefined, nonfinite, or not
positive definite, it may be replaced by any fixed positive
definite matrix; the probability of using this fallback tends to zero uniformly. In low dimension the same asymptotic conclusions hold if the
identifiable curvature block is inverted directly. In our
dimension-three simulations, however, this gives appreciably poorer
coverage for smoothness and the microergodic parameter
(Section~\ref{sec:sim-results} and Appendix~\ref{apdx:sim-range}).
For a positive definite full Fisher information matrix at the same
parameter and locations, block inversion gives 
\[
 \left[\begin{pmatrix}A&B\\B^\top&D\end{pmatrix}^{-1}
       \right]_{\zeta\zeta}
 =(A-BD^{-1}B^\top)^{-1}\succeq A^{-1},
\]    
where $A$ is the identifiable information block and $D$ is the
nuisance block. Direct block inversion treats the nuisance as known;
full inversion accounts for its coupling with the identifiable
parameters. Appendix~\ref{app:full-inverse-block} derives the identity
and explains the qualification needed for an observed Hessian.

\subsection{Known smoothness}\label{sec:known}
Smoothness is often fixed at a prespecified value in applications.
Here we consider the case where that value equals the true
smoothness. Knowing the smoothness improves the rate for $\eta$ from $(\log n)/r_n$ to $r_n^{-1}$ and leaves the rates for $a$ and $t$ unchanged.
\begin{corollary}[Known smoothness]\label{cor:known}
Under Assumptions~\ref{ass:param}, \ref{ass:sharp}, and
\ref{ass:noise}, suppose smoothness is fixed at its true value. The same conclusions hold for every identifiable vector in Table~\ref{tab:parameterizations}, with its smoothness coordinate removed. For this restricted fit, remove the smoothness row and column from both the expected negative Hessian and the score covariance before forming the sandwich covariance. The restrictions on joint inference in the original variance--lengthscale coordinates remain as stated in Corollary~\ref{cor:log}. The rates for $\eta,c,\log c$ improve to $r_n^{-1}$. The tail, moment, and full inverse-Hessian conclusions also hold for the restricted fit. The minimax lower bounds additionally assume $N(0,1)\in\mathcal P$. For the global comparison, assume also that, in the free coordinates,
$\mathcal K$ is a compact rectangle with nonempty interior contained in the relative interior of the fixed-smoothness parameter set. For the local comparison, take the neighborhoods in Theorem~\ref{thm:minimax} within the slice $\nu=\nu_0$. The other rates and the fourth-cumulant
correction are unchanged. The estimator is the global maximizer over the parameter
set with smoothness fixed.
\end{corollary}

\subsection{Proof heuristic}\label{sec:outline}

We give a brief heuristic for the proof of \Cref{thm:normalized} with Gaussian measurement
errors and deterministic locations. A detailed proof outline appears in
Appendix~\ref{app:proof-outline}.

Let $Q$ have side length $L$, and choose
$\chi\in C_c^\infty(Q)$, with $0\le\chi\le1$, equal to $1$
on a cube $Q_0\Subset Q$ of side length $L/m$, where $m\ge3$
is a fixed integer. Consider the localized Fourier coefficients
\[
 X_\lambda
 =L^{-d/2}\int_Q\chi(x)F(x)e^{-2\pi i\lambda^\top x}\,dx,
 \qquad \lambda\in\mathscr L\coloneqq (m/L)\mathbb Z^d.
\]
These are Fourier series coefficients of $\chi F$ on $Q$,
retaining every $m$th frequency in each coordinate.
Localized Fourier coefficients also play a central role in
\citet[Sections~5--7]{pillai2026}.

Since $F$ is Gaussian, the coefficients are jointly zero-mean complex Gaussian. Lemma~\ref{lem:windows}(vi) gives
\[
 \begin{aligned}
 v_\theta(\lambda)
 \coloneqq \E_\theta|X_\lambda|^2
 &=L^{-d}\int_{\mathbb R^d}
       f_\theta(\omega)|\widehat\chi(\lambda-\omega)|^2\,d\omega =c_\chi f_\theta(\lambda)
       \{1+O(\langle\lambda\rangle^{-2})\},
 \end{aligned}
\]
where $c_\chi=L^{-d}\int_Q\chi(x)^2\,dx$ and $\langle z\rangle=(1+\|z\|^2)^{1/2}$.
Moreover,
\[
 \nabla_\theta\log v_\theta(\lambda)
 =
 \begin{pmatrix}
  1\\
  -\log m_a(\lambda)+O(\langle\lambda\rangle^{-2})\\
  -\beta/m_a(\lambda)+O(\langle\lambda\rangle^{-4})\\
  0
 \end{pmatrix},
 \qquad \theta=(\eta,\nu,a,t),
\]
where $\beta=\nu+d/2$ and $m_a(\lambda)=a+4\pi^2\|\lambda\|^2$. Hence, the coefficient variances and their parameter sensitivities are closely approximated by those of the spectral density.
Their correlations also decay rapidly with frequency separation.
For every $M>0$,
\[
 \begin{aligned}
 \frac{|\E_\theta[X_\lambda\overline{X_\mu}]|}
      {\sqrt{v_\theta(\lambda)v_\theta(\mu)}}
 &\le C_M\langle\lambda-\mu\rangle^{-M} \qquad \text{and}\qquad 
 \frac{|\E_\theta[X_\lambda X_\mu]|}
      {\sqrt{v_\theta(\lambda)v_\theta(\mu)}}
 \le C_M\langle\lambda+\mu\rangle^{-M}.
 \end{aligned}
\]
Because $X_{-\lambda}=\overline{X_\lambda}$, we retain only
frequencies with $\lambda_1>0$. The resulting coefficients
are weakly correlated, with variances approximately
$c_\chi f_\theta(\lambda)$.

The coefficients are not observed directly, and measurement
error makes recovery harder at high frequencies, where their
variances are small. Let $\widehat X_{\lambda,n}=\E_0(X_\lambda\mid Y_n)$ be the best linear predictor under the truth. With
$R_n=n^{1/(2\nu_0+d)}$, Lemma~\ref{lem:sampling} gives
\[
 \frac{\{\E_0|\widehat X_{\lambda,n}-X_\lambda|^2\}^{1/2}}
      {\sqrt{v_{\theta_0}(\lambda)}}
 \le Cn^{-1/2}\langle\lambda\rangle^{\nu_0+d/2}
 \le C\epsilon^{\nu_0+d/2},
 \qquad \|\lambda\|\le\epsilon R_n,
\]
for each fixed $\epsilon>0$ and all sufficiently large $n$.
Thus frequencies below a small multiple of $R_n$ can be
recovered with arbitrarily small relative root mean squared error.

For the signal parameters, this suggests the idealized model
\[
 X_\lambda^{\mathrm{id}}
 \stackrel{\mathrm{ind}}{\sim}
 N_{\mathbb C}(0,f_\theta(\lambda)),
 \qquad
 \lambda\in\mathscr L,\quad
 \lambda_1>0,\quad \|\lambda\|\le\epsilon R_n,
\]
after dividing the original coefficients by $\sqrt{c_\chi}$. Here $N_{\mathbb C}(0,v)$ denotes a complex Gaussian with independent $N(0,v/2)$ real and imaginary parts. Its likelihood suggests the signal-parameter rates and asymptotic normality in Theorem~\ref{thm:normalized}. In the proof, energies of the predicted coefficients provide the Fisher information lower bound (Lemma~\ref{lem:information}), while spectral calculations give matching upper bounds (Lemma~\ref{lem:upper}). Localization and a quadratic expansion of the likelihood then yield the normal approximation.

Finally, the normalized nugget error is asymptotically
independent of the identifiable signal-parameter errors.
Lemma~\ref{lem:noise-diagonals} and the likelihood expansion give
\[
 \sqrt n(\widehat t_n-t_0)
 =\frac1{\sqrt n}\sum_{i=1}^n(\varepsilon_{n,i}^2-t_0)
   +o_P(1)
 \ \Rightarrow\ N(0,2t_0^2).
\]
Thus, to first order, the nugget estimator behaves as if the
measurement errors themselves were observed.

\section{Simulation}\label{sec:simulation}
We evaluate how closely the asymptotic theory in \Cref{thm:normalized} describes empirical MLEs across dimensions $d=3,4,5$ (spanning distinct identifiability and convergence-rate regimes) and smoothnesses $\nu_0=1/2, 3/2$. For each latent-field realization, we generate datasets under four measurement-error laws differing in fourth moments and shape. We also check coverage of the feasible Wald intervals (\Cref{cor:uncertainty}) with and without the fourth-cumulant correction~\eqref{eq:estimatedlogcov}.

\subsection{Simulation setup}\label{sec:sim-setup}

\noindent\textbf{Model and designs.}
We simulate from \eqref{eq:model} with the Mat\'ern covariance \eqref{eq:matern-intro} on the unit cube $D=[0,1]^d$, for $d\in\{3,4,5\}$ and $\nu_0\in\{1/2,3/2\}$, with $(\sigma_0^2,l_0,\tau_0^2)=(1,\ 0.4,\ 0.25)$. In the parameterization \eqref{eq:matern-intro} (without a factor $\sqrt{2\nu}$), $\nu_0=1/2$ gives the exponential covariance $\sigma^2e^{-r/l}$ and $\nu_0=3/2$ gives $\sigma^2(1+r/l)e^{-r/l}$. The true values were chosen, using the exact Fisher information alone and before any Monte Carlo run, so that the asymptotic regime is approached within the computationally accessible range of $n$; shorter lengthscales or smaller nuggets require much larger $n$, especially for $\nu_0=1/2$ (Appendix~\ref{apdx:sim-truth}). The locations form a regular midpoint grid with $G$ points per axis,
\[
 x\in\{(i-\tfrac12)/G:\ i=1,\ldots,G\}^d,
\]
so that $n=G^d$. This design satisfies \Cref{ass:sharp}. We use
\begin{align*}
 G&\in\{5,7,10,14,22,32,40,46\}, &&d=3\quad(125\le n\le 97{,}336),\\
 G&\in\{3,4,6,8,10,13,18\}, &&d=4\quad(81\le n\le 104{,}976),\\
 G&\in\{3,4,5,6,8,10\}, &&d=5\quad(243\le n\le 100,000).
\end{align*}
For every $(d,\nu_0,n)$ we generate $R=400$ independent GP realizations. For each realization, we generate four data sets using independent measurement errors from the four distributions described below.

\noindent\textbf{Measurement error laws.}
The standardized errors $\xi_i=\varepsilon_i/\tau_0$ are iid from one of four laws satisfying \Cref{ass:noise}: (i) the standard normal law, with $E\xi^4=3$; (ii) the uniform law on $[-\sqrt3,\sqrt3]$, with $E\xi^4=9/5$, which has light tails; (iii) the scale mixture $0.9\,N(0,s^2)+0.1\,N(0,9s^2)$ with $s^2=5/9$, with $E\xi^4=25/3$, which is symmetric with a large fourth moment; and (iv) the centered and scaled Bernoulli variable $(B-0.1)/0.3$ with $B\sim\mathrm{Bernoulli}(0.1)$, with $E\xi^4=73/9$, which is skewed (skewness $8/3$) and discrete. By \eqref{eq:lognugget}, the limiting standard deviation of $\sqrt n(\log\widehat\tau^2-\log\tau_0^2)$ is $\sqrt{E\xi^4-1}$, which equals $1.41$, $0.89$, $2.71$, and $2.67$ for laws (i)--(iv). Laws (iii)
and (iv) have nearly the same fourth moment but very different shapes, so the theory predicts nearly the same nugget behavior for both. Within a replicate, all four laws share one realization of the latent field, and the errors for each law are drawn independently. Differences between laws within a replicate are therefore due to the errors alone, while for each law the replicates remain independent.

\noindent\textbf{Likelihood evaluation and MLE.}
To evaluate the exact Gaussian likelihood for sample sizes up to approximately $10^5$, we exploit the symmetries of the regular grids to block-diagonalize the covariance matrix. This allows us to evaluate the likelihood, its derivatives, and the Fisher information without factorizing the full covariance matrix. Appendix~\ref{apdx:sim-lik} describes the likelihood computation, and Appendix~\ref{apdx:sim-generation} describes data generation. The numerical fits target the global maximizer over the box
\[
 \Theta=\{\eta\in[-15,40],\ \nu\in[0.1,4],\ a\in[0.25,2500],\ \tau^2\in[10^{-3},10]\}.
\] The fits use L-BFGS-B with the exact gradient, initialized and preconditioned at the truth. The implementation searches an enlarged, sample-size-dependent set, but all $67{,}200$ retained estimates satisfy the displayed bounds. We restarted the $143$ fits whose
projected-gradient residual exceeded $10^{-3}$ from five starting points and retained the largest criterion value. All final residuals are below $10^{-3}$. Appendix~\ref{apdx:sim-mle} gives the restart rule, constraint audit, and additional multistart checks; these
diagnostics do not certify global maximization. Boundary solutions occur only at small $n$, mostly with $\widehat\nu$ at its upper bound or $\widehat\tau^2$ at its lower bound. Under each of the four error laws, at most $0.5\%$ of the estimates lie on a reported parameter bound for $n\ge10^4$, and none for $n\ge2.5\times10^4$ (\Cref{tab:sim-boundary}).

\noindent\textbf{Error summaries.}
For each estimable parameter we report, as functions of $n$: 
\begin{enumerate}[label=(\alph*)]
\item the root mean squared error (RMSE) over the $R=400$ replicates;
\item the robust standard deviation $\mathrm{IQR}/1.349$ of the replicate
errors, which is calibrated to the standard deviation of a normal distribution and is
insensitive to a few extreme estimates;

\item the information-based standard deviation at the true parameters and the observed locations. Write the log-likelihood in coordinates $w=(\eta,\nu,\log a,\log t)^\top$. Let $\calI_n^{\mathrm{full}}$ be the four-parameter Gaussian information at $w_0$, and let $\Omega_n$ be the covariance of the score $\nabla_w\ell_n(w_0)$ under the true error law, at the same locations. If $J_\zeta(w_0)$ is the Jacobian of the reported vector $\zeta$ with respect to $w$, the covariance benchmark from the first-order MLE approximation is
\[
 B_n^\zeta
 = J_\zeta(w_0)(\calI_n^{\mathrm{full}})^{-1}
   \Omega_n(\calI_n^{\mathrm{full}})^{-1}J_\zeta(w_0)^\top .
\]
For coordinate $j$, we report $\sqrt{[B_n^\zeta]_{jj}}$.
Under Gaussian errors, $\Omega_n=\calI_n^{\mathrm{full}}$.
For $d=3$, we also compute the benchmark with the range treated as known: before forming the sandwich, we retain the $(\eta,\nu,\log t)$ rows and columns of both matrices and the corresponding columns of $J_\zeta(w_0)$. Both benchmarks are compared with the same four-parameter fits.
Appendix~\ref{apdx:sim-information} gives the computation.

\item the empirical coverage of nominal $95\%$ Wald intervals
$\widehat\zeta_{n,j}\pm1.96\,[\widehat\Sigma_n^\zeta]_{jj}^{1/2}$.
We use the fourth-cumulant correction in
\eqref{eq:estimatedlogcov}, with $\widetilde\kappa_{4,n}$ estimated
as in \Cref{lem:kurtosis}. The Fisher and observed-Hessian versions
are computed as described below. For comparison, we also report
uncorrected intervals, setting $\widetilde\kappa_{4,n}=0$.
\end{enumerate}

Bootstrap $95\%$ intervals for (a) and (b) are obtained by resampling the replicates. Appendix~\ref{apdx:sim-range} compares the two benchmarks for $d=3$ in (c). The parameters shown differ by dimension. For $d=3$ we report
$\nu$, $c=\sigma^2l^{-2\nu}$ and $\tau^2$, and show the distribution of the estimates of $\sigma^2$ and $l$, which are not consistent. For $d=4,5$ we report $\nu,\sigma^2,l,\tau^2$, all of which are consistent.

\noindent\textbf{Standard error comparisons.}
We compare standard errors based on the Fisher information and the observed negative Hessian, evaluated at the fitted parameters or at the truth. Both are computed in $w=(\eta,\nu,\log a,\log t)$, inverted, and mapped to the reported parameters by the delta method, with the fourth-cumulant correction. For the Hessian, this need not equal inversion after direct differentiation in the reported coordinates at the truth or at a fit that is on the boundary of the parameter space. For $d=3$, we also compare direct inversion of the identifiable block with selection of that block from the full inverse. The versions evaluated at the truth are oracle benchmarks. Appendix~\ref{apdx:sim-mle} specifies the convention.

\noindent\textbf{Log parametrization.}
As in \Cref{cor:log}, we report the results on the log scale. The log transformation preserves the rates and provides a natural parametrization for positive parameters. Results on the original scale give the same marginal rate comparisons and are given in Appendix~\ref{apdx:sim-original}. 

\begin{figure}[!t]
\centering
\includegraphics[width=\textwidth]{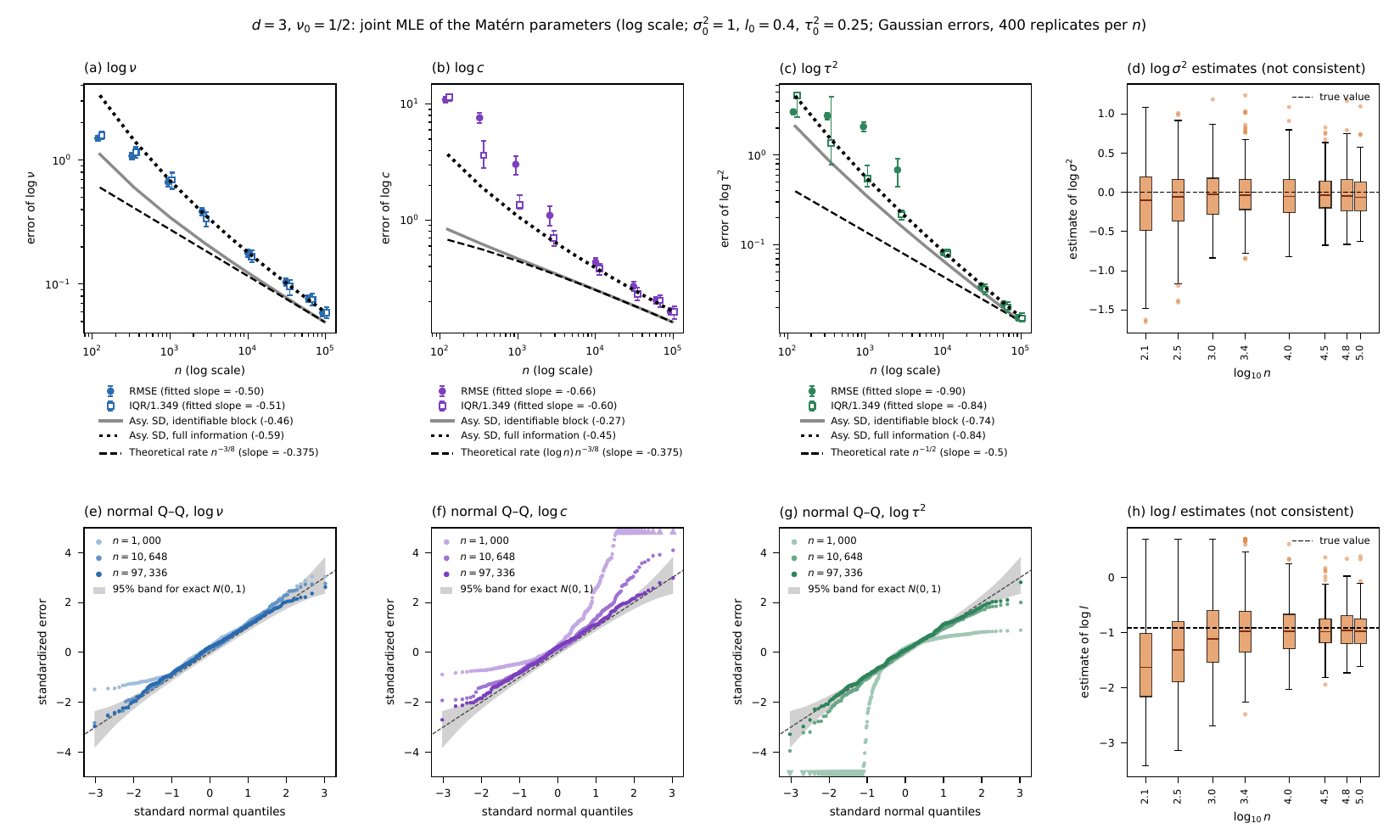}
\caption{Simulation results for $d=3$, $\nu_0=1/2$ on the log scale, with $(\sigma_0^2,l_0,\tau_0^2)=(1,0.4,0.25)$, Gaussian errors, and $400$ replicates per $n$. (a)--(c): RMSE (filled circles) and robust standard deviation $\mathrm{IQR}/1.349$ (open squares), both with bootstrap $95\%$ intervals; information-based standard deviation from the identifiable block of the Fisher information as in \Cref{thm:normalized} (grey line) and from the full information, accounting for estimation of the range (dotted line); dashed lines show the rates in Table~\ref{tab:rates}, written in $n$. (e)--(g): normal Q--Q plots of the errors standardized by the full-information standard deviation at $n\approx10^3,10^4,10^5$, with a pointwise \(95\%\) band under exact normality; triangles mark points outside the frame. (d), (h): box plots of $\log\widehat\sigma^2$ and $\log\widehat l$, which are not consistent for $d\le3$; dashed lines show the true values.}
\label{fig:sim-d3-nu05}
\end{figure}

\begin{figure}[!th]
\centering
\includegraphics[width=\textwidth]{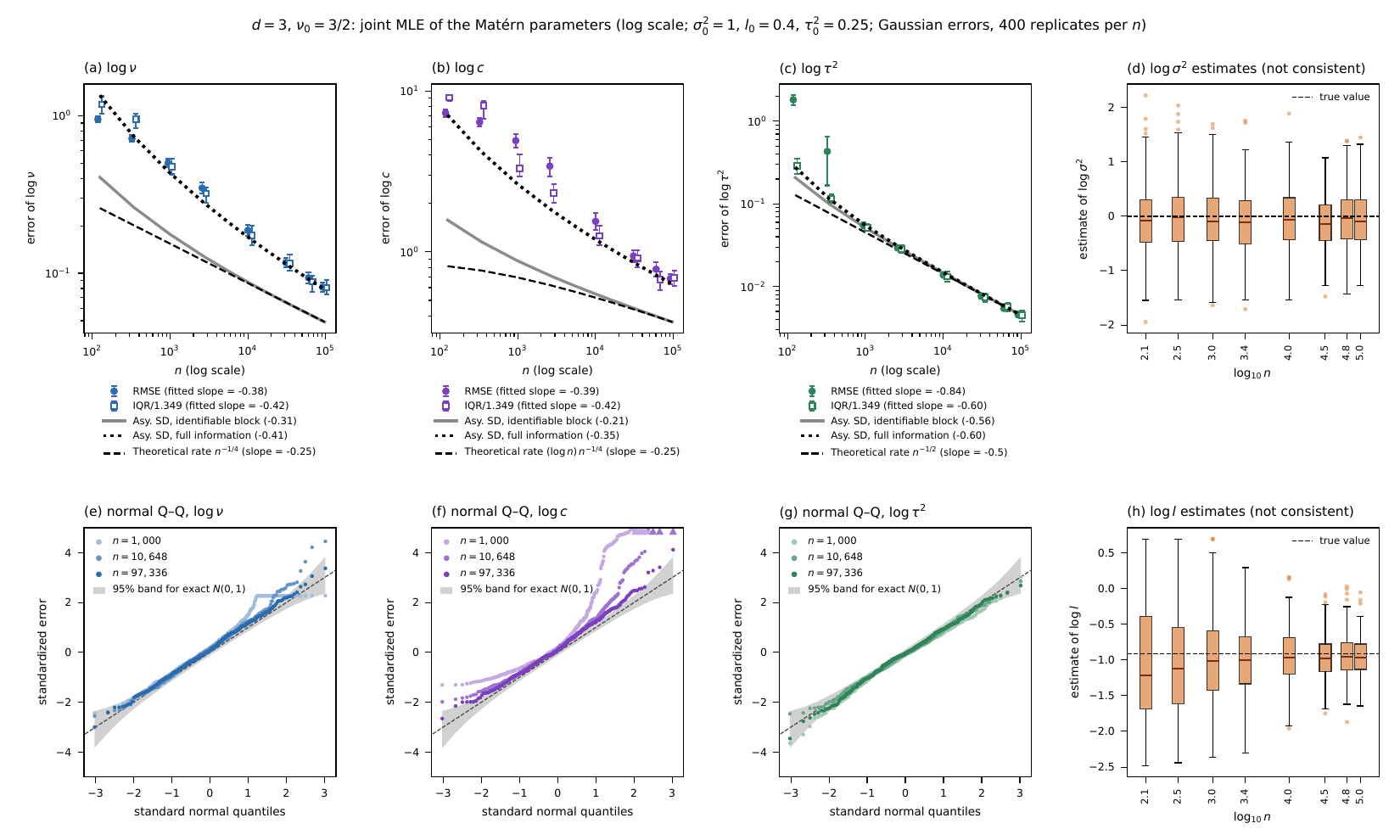}
\caption{Simulation results for $d=3$, $\nu_0=3/2$, with the layout as in \Cref{fig:sim-d3-nu05}.}
\label{fig:sim-d3-nu15}
\end{figure}

\begin{figure}[th]
\centering
\includegraphics[width=\textwidth]{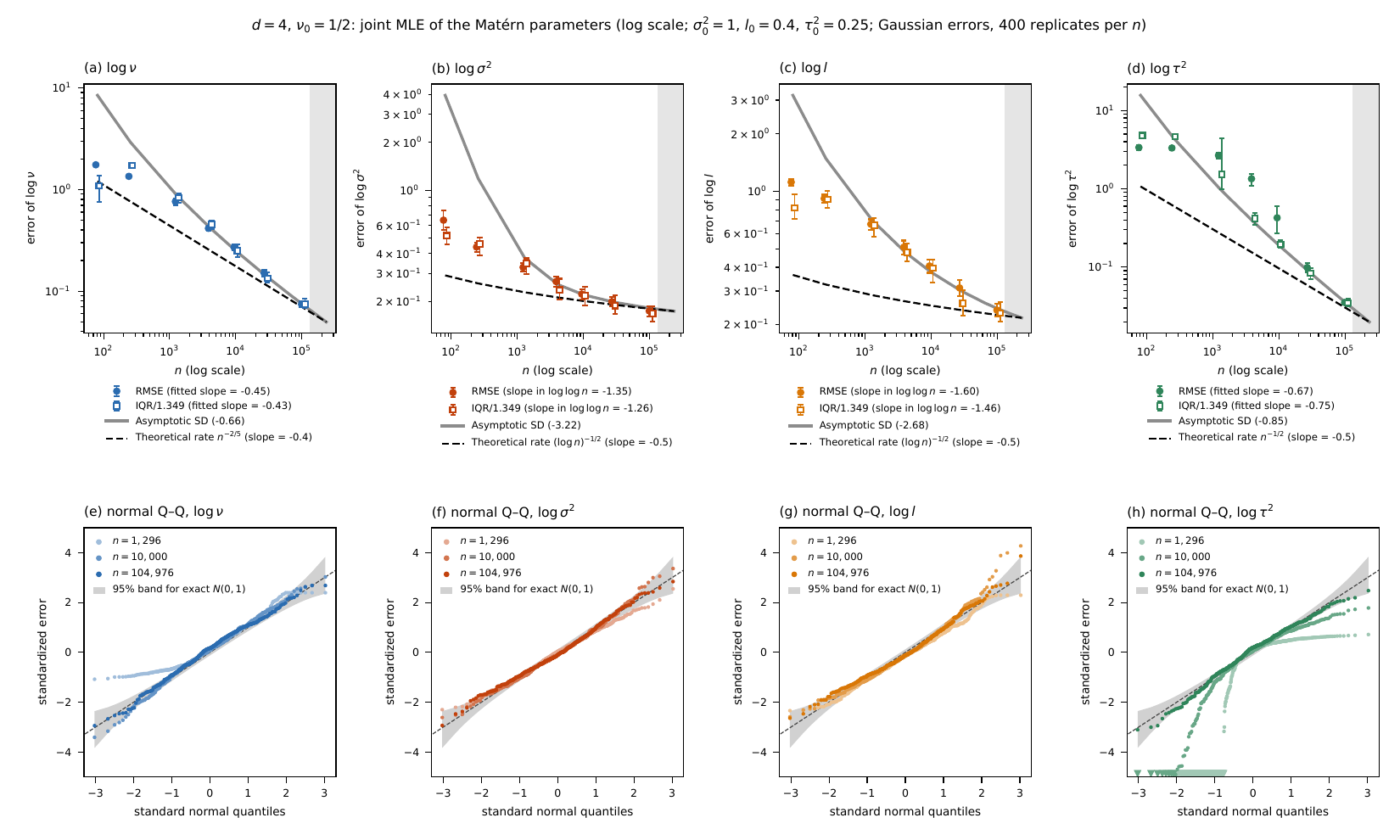}
\caption{Simulation results for $d=4$, $\nu_0=1/2$, with Gaussian errors. Panels (a)--(d) show estimation errors for
$\log\nu$, $\log\sigma^2$, $\log l$, and $\log\tau^2$;
panels (e)--(h) show the corresponding normal Q--Q plots.
Symbols and reference bands are as in Figure~\ref{fig:sim-d3-nu05}. The grey region marks sample sizes beyond the Monte Carlo range, where only the Fisher information is evaluated.}
\label{fig:sim-d4-nu05}
\end{figure}

\begin{figure}[th]
\centering
\includegraphics[width=\textwidth]{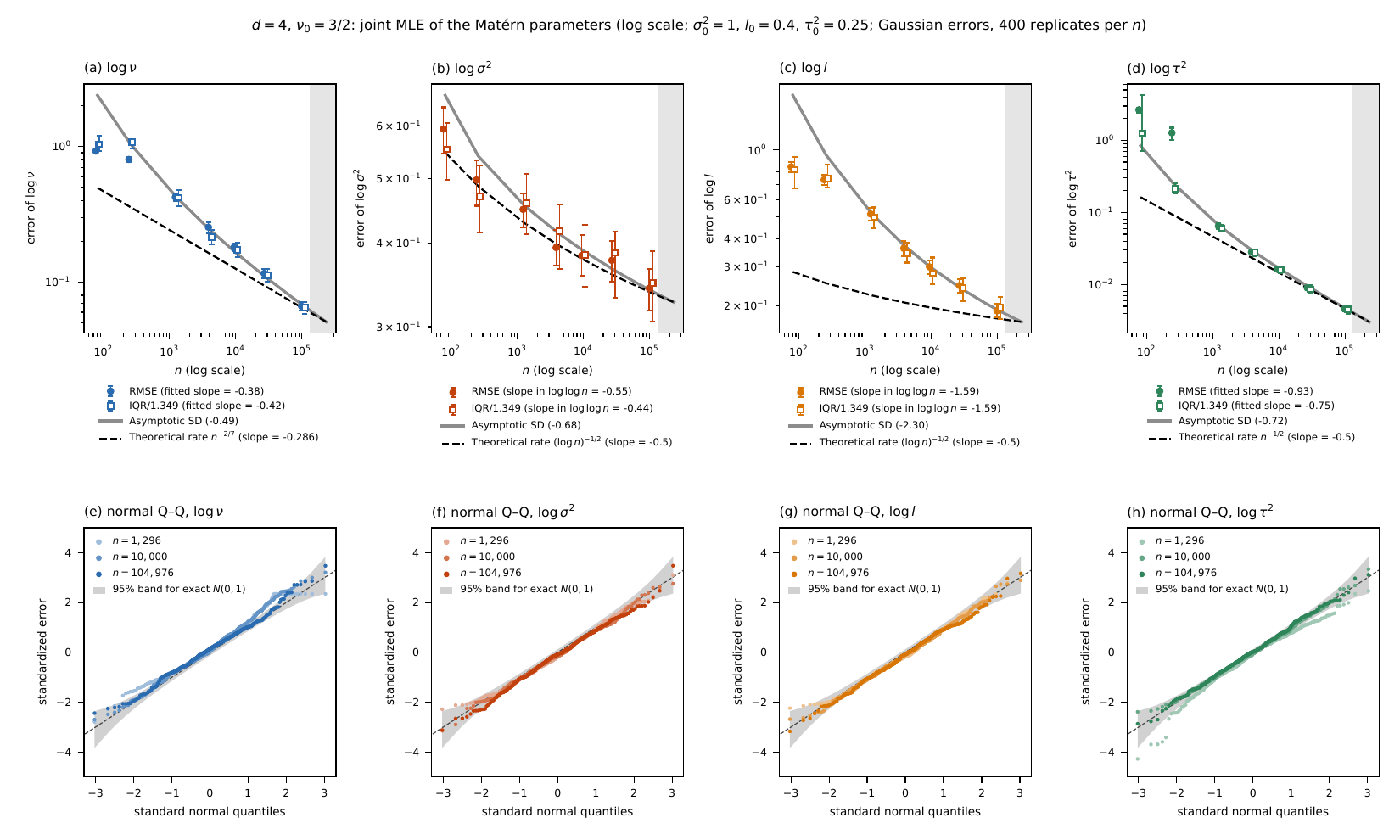}
\caption{Simulation results for $d=4$, $\nu_0=3/2$, with the layout and
graphical elements as in \Cref{fig:sim-d4-nu05}.}
\label{fig:sim-d4-nu15}
\end{figure}

\begin{figure}[th]
\centering
\includegraphics[width=\textwidth]{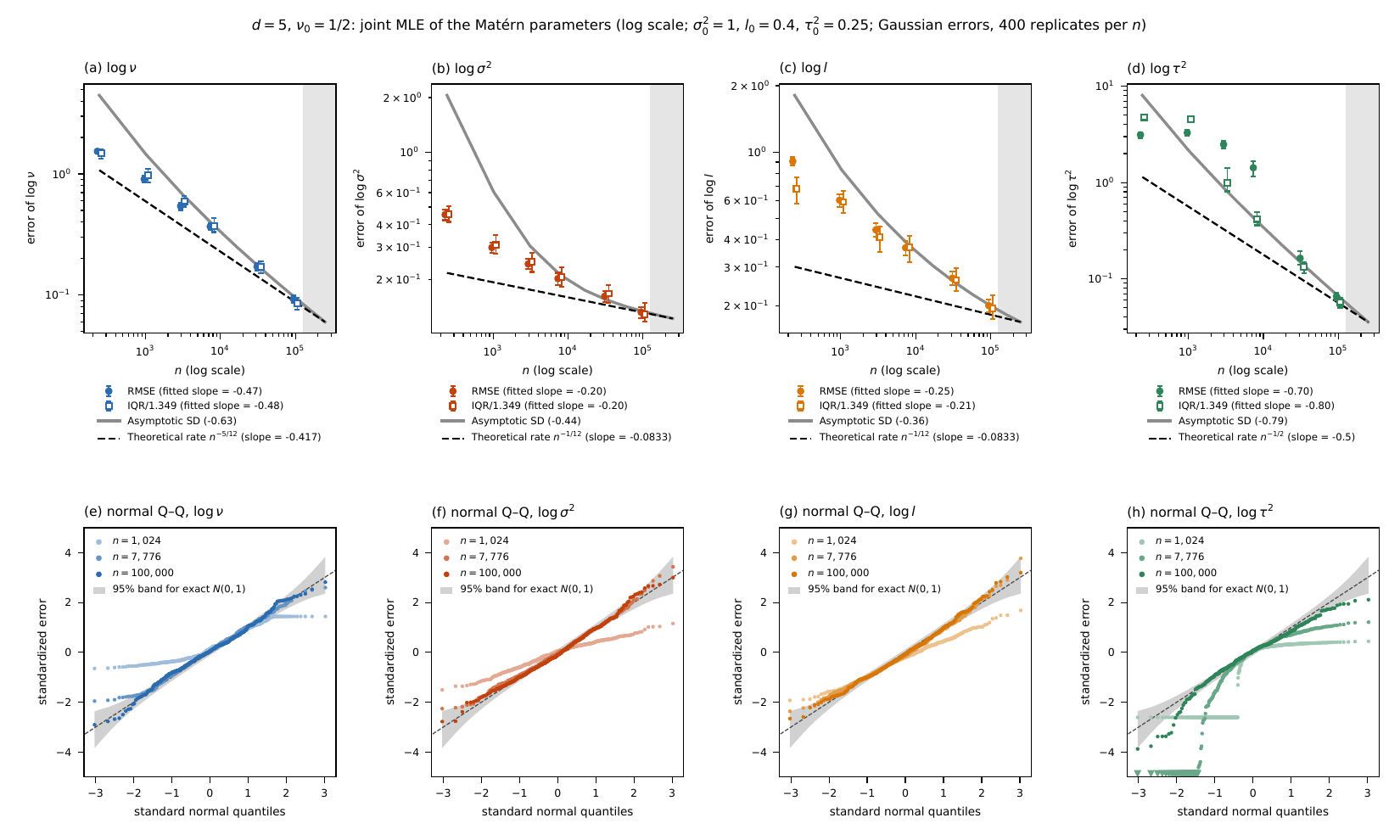}
\caption{Simulation results for $d=5$, $\nu_0=1/2$, with the layout as in \Cref{fig:sim-d4-nu05}.}
\label{fig:sim-d5-nu05}
\end{figure}

\begin{figure}[th]
\centering
\includegraphics[width=\textwidth]{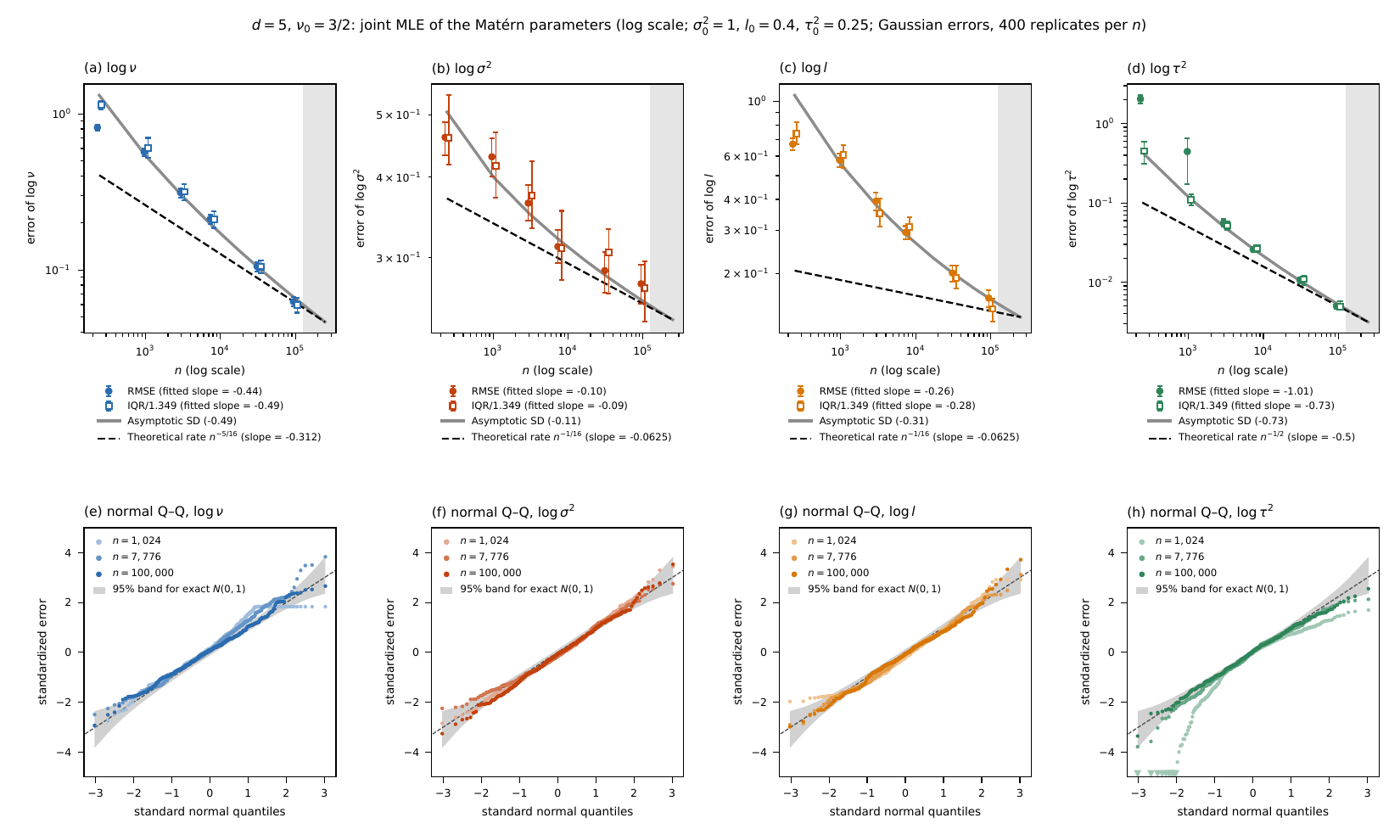}
\caption{Simulation results for $d=5$, $\nu_0=3/2$, with the layout as in \Cref{fig:sim-d4-nu05}.}
\label{fig:sim-d5-nu15}
\end{figure}

\begin{figure}[th]
\centering
\includegraphics[width=\textwidth]{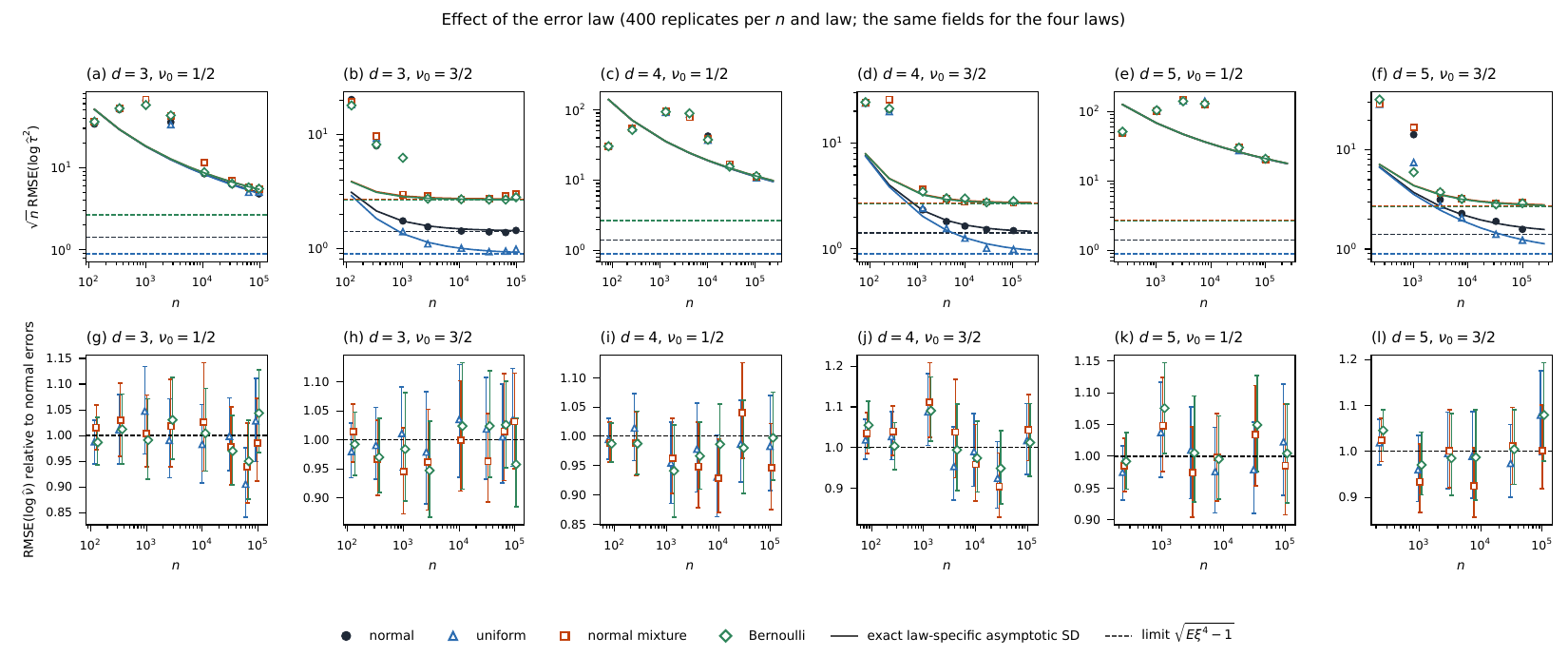}
\caption{Effect of the error law, $400$ replicates per $n$ and law. (a)--(f): $\sqrt n\,\mathrm{RMSE}(\log\widehat\tau^2)$ for normal (circles), uniform (triangles), normal-mixture (squares), and Bernoulli (diamonds) errors, one panel per $(d,\nu_0)$. The solid lines show $\sqrt n$ times the law-specific information-based standard deviation of $\log\widehat\tau^2$; the dashed lines show the corresponding limits $\sqrt{\E_P\xi^4-1}$. (g)--(l): ratio of the RMSE of
$\log\widehat\nu$ under each non-Gaussian law to that under Gaussian errors, computed from the same field realizations; the theory predicts one. The intervals in (g)--(l) are paired percentile-bootstrap $95\%$ intervals, obtained by resampling the same replicate indices across error laws.}
\label{fig:sim-laws}
\end{figure}

\subsection{Simulation results}\label{sec:sim-results}

\noindent\textbf{Overview of simulation results.}
\Cref{fig:sim-d3-nu05,fig:sim-d3-nu15,fig:sim-d4-nu05,fig:sim-d4-nu15,fig:sim-d5-nu05,fig:sim-d5-nu15}
present the estimation errors and normal Q--Q plots under Gaussian measurement errors across different dimensions and smoothness
parameters. The estimation errors are compared with the information-based standard deviations and theoretical convergence rates. For $d=3$, we additionally examine the estimates of the
spatial variance and lengthscale, which are not identifiable. \Cref{fig:sim-laws} compares the estimation errors under four different measurement error distributions.
\Cref{tab:sim-coverage,tab:sim-joint} report the empirical coverage of marginal confidence intervals and joint confidence regions, respectively. Additional comparisons of convergence rates and
standard error estimators are provided in the Appendix.

\noindent\textbf{Estimation errors and convergence rates.}
For Gaussian measurement errors, the Monte Carlo errors agree closely with the information-based standard deviations once $n$ is moderately large, as predicted by \Cref{cor:log}. For $d=4,5$ and both values of $\nu_0$, the RMSEs of all four log-parameters are within $14\%$ of their design-specific information-based standard deviations for every simulated $n\ge2.5\times10^4$, and within $7\%$ at the largest sample size. For $\nu_0=3/2$ they are within $11\%$ from $n\approx8\times10^3$, whereas for $\nu_0=1/2$ the nugget needs larger $n$ (see the nugget paragraph below). For $d=3$ and $\nu_0=1/2$, the
RMSEs of $\log\widehat\nu$, $\log\widehat c$, and $\log\widehat\tau^2$ are within $15\%$ of the full-information standard deviation for $n\ge10^4$, and within $4\%$ at $n=97{,}336$. For $d=3$ and $\nu_0=3/2$, where the smoothness rate $n^{-1/4}$ is the slowest in the
study, the RMSE exceeds the full-information standard deviation at $n=97{,}336$ by $3\%$ for $\log \nu$, $8\%$ for $\log c$, and $0.3\%$ for $\log \tau^2$. At small $n$,
the RMSE can be much smaller than the information-based standard deviation for $d\ge4$; there the Fisher standard deviation of $\log \nu$ or $\log \tau^2$ is comparable to the width of the corresponding log-parameter interval, and the estimates are truncated by its boundary.

The comparisons in Appendix~\ref{apdx:sim-rates} show that the finite-sample slopes can differ from their theoretical limits. Even the information-based standard deviations have not reached their limiting slopes over the simulated sample sizes. Moreover, the empirical slopes are affected by boundary estimates and large estimation errors at small sample sizes. Consequently, these slopes should not be interpreted as direct estimates of the asymptotic exponents.

\noindent\textbf{Dimension-dependent identifiability.} 
For $d=3$, the estimators of smoothness, the microergodic parameter $c$, and the nugget variance become more concentrated as $n$ increases. In contrast, the distributions of the estimated spatial variance and lengthscale retain substantial variation even at the largest sample size, consistent with their nonidentifiability under fixed-domain asymptotics. For $d=4$ and $d=5$ both parameters are consistent, and their errors track the information-based standard deviation, but they decrease slowly: from $n\approx10^3$ to $n\approx10^5$, the RMSE of $\log\widehat\sigma^2$ decreases by a factor of only $1.3$--$2.3$, compared with factors between $6.4$ and $10.3$ for $\log\widehat\nu$. This slow decrease is consistent with the qualitative distinction between the rates in \Cref{cor:rates}.

\noindent\textbf{Smoothness and amplitude.}
The smoothness is estimated at the polynomial rate $n^{-d/(4\nu_0+2d)}$, which is faster for rougher fields: at $n\approx10^5$ the RMSE of $\log\widehat\nu$ is between $0.06$ and $0.09$ for $d\ge4$. For $d=3$ the amplitude $c$ is estimated jointly with $\nu$, and its error carries
the additional factor $\log n$ in Table~\ref{tab:rates}: at the largest $n$, the RMSE of $\log\widehat c$ is $2.8$ ($\nu_0=1/2$) and $8.3$ ($\nu_0=3/2$) times that of $\log\widehat\nu$, close to the ratios $2.8$ and $8.0$ of the full-information standard deviations.

\noindent\textbf{Nugget variance.}
For $\nu_0=3/2$, the nugget errors are approaching the parametric limit. With Gaussian errors, at the largest $n$,
$\sqrt n\,\mathrm{RMSE}(\log\widehat\tau^2)$ equals $1.44$,
$1.49$, and $1.58$ for $d=3,4,5$, compared with the limit
$\sqrt2\approx1.41$ in \eqref{eq:lognugget} and with $\sqrt n$ times the information-based standard deviation, $1.43$, $1.49$, and $1.66$. For the rough exponential
field, $\nu_0=1/2$, the nugget error again matches its information-based standard deviation, but the latter is still several times larger than $\sqrt{2/n}$ at the largest simulated $n$ (by factors of about $3.5$, $7.9$, and $14.8$ for $d=3,4,5$). The proportion of signal-dominated directions decreases at rate $n^{-2\nu_0/(2\nu_0+d)}$, which is particularly slow for rough fields. This helps explain why the information-based
standard deviation of the nugget remains substantially larger than its limiting value over the simulated sample sizes. The same effect causes the heavy lower tail of $\log\widehat\tau^2$ for $d=3$, $\nu_0=1/2$ at $n=10^3$, where the fit occasionally explains the noise by a rougher signal (Appendix~\ref{apdx:sim-nugget}).

\noindent\textbf{Non-Gaussian measurement errors.}
\Cref{fig:sim-laws} compares the four error distributions. At the largest sample sizes, the RMSE of $\log\widehat\tau^2$ agrees with its law-specific information-based standard deviation to within $11\%$ across all settings. For $\nu_0=3/2$, the differences
between error distributions are substantial and agree with the dependence of the asymptotic variance on the fourth moment. In particular, the normal-mixture and Bernoulli errors produce similar estimation errors despite their different distributional shapes. For $\nu_0=1/2$, the differences between error laws are much smaller because the nugget remains strongly confounded with the rough signal at the simulated sample sizes.

In contrast, the estimation of the signal parameters
is largely unaffected by the measurement error
distribution. Their information-based standard
deviations differ from the Gaussian benchmarks
by less than $0.1\%$, and their empirical RMSEs
differ by at most $8\%$ at the largest sample sizes.
These results agree with the asymptotic covariance
structure established in \Cref{thm:normalized}.

\noindent\textbf{Asymptotic normality.}
The Q--Q plots in
\Cref{fig:sim-d3-nu05,fig:sim-d3-nu15,fig:sim-d4-nu05,fig:sim-d4-nu15,fig:sim-d5-nu05,fig:sim-d5-nu15} show that the standardized estimation errors generally approach the normal reference as the sample size increases. For $d=4,5$, the Kolmogorov--Smirnov distances from the standard normal distribution range from $0.03$ to $0.09$ at $n\approx10^5$. For $d=3$, the corresponding distances range from $0.04$ to $0.09$ at $n=97{,}336$. Deviations from normality are more pronounced at smaller sample sizes, particularly when parameter estimates approach the boundary of $\Theta$ or when the underlying process is relatively rough.

\noindent\textbf{Coverage of the Wald intervals.}
\Cref{tab:sim-coverage} reports the empirical coverage of nominal $95\%$ confidence intervals at the largest sample sizes. For $\nu_0=3/2$, the uncorrected nugget intervals exhibit substantial coverage errors under non-Gaussian measurement errors. In particular, the coverage is only $66\%$--$74\%$ for the normal-mixture and Bernoulli errors, whereas it is approximately $99\%$--$100\%$ for uniform errors. These findings agree with the coverage predicted by the information-based standard deviations.

The fourth-cumulant correction substantially improves the coverage. Across all $24$ combinations of dimension, smoothness, and error distribution, the corrected nugget intervals achieve coverage between $92.0\%$ and $96.8\%$, using either the Fisher information or observed Hessian. The signal-parameter intervals achieve coverage between $91\%$ and $97\%$ without any correction. Replacing the estimated fourth cumulant with its true value changes
nugget coverage by at most $0.5$ percentage points, corresponding to two of the $400$ replicates. The paired Monte Carlo standard error at that comparison is $0.35$ percentage points.

\begin{table}[!ht]
\centering\small
\caption{Coverage (\%) of nominal $95\%$ Wald intervals at the largest $n$ ($97{,}336$, $104{,}976$, and $10^5$ for $d=3,4,5$), $400$ replicates. ``Signal'': range over the identifiable signal log-parameters of the interval from \eqref{eq:estimatedlogcov} with Fisher information at the fit. The four columns for $\log\tau^2$: predicted normal-approximation coverage using information-based standard deviations for the uncorrected interval ($\widetilde\kappa_{4,n}=0$); its Monte Carlo coverage; and
Monte Carlo coverage of the corrected interval of \Cref{cor:uncertainty} with Fisher information and observed
Hessian. The Monte Carlo standard error is about $1.1$ percentage points at $95\%$. Every reported marginal variance was finite and positive, so all $400$ replicates contribute to each estimate.}
\label{tab:sim-coverage}
\vspace{3pt}
\begin{tabular}{lllrrrrr}
\toprule
 & & & & \multicolumn{4}{c}{$\log\tau^2$}\\
\cmidrule(l){5-8}
 & & Error law & Signal & \shortstack{Normal\\approximation} & Uncorrected & Fisher & Hessian\\
\midrule
$d=3$ & $\nu_0=1/2$ & Normal & $94.8$--$95.5$ & $95.0$ & $95.8$ & $95.8$ & $95.5$\\
 & & Uniform & $95.2$ & $95.6$ & $96.8$ & $96.5$ & $96.8$\\
 & & Mixture & $97.0$ & $92.5$ & $92.8$ & $95.8$ & $96.2$\\
 & & Bernoulli & $94.8$--$95.0$ & $92.6$ & $92.8$ & $94.0$ & $93.8$\\
 & $\nu_0=3/2$ & Normal & $96.0$ & $95.0$ & $93.8$ & $93.8$ & $93.8$\\
 & & Uniform & $92.8$--$93.5$ & $99.8$ & $99.8$ & $95.0$ & $95.0$\\
 & & Mixture & $96.0$--$96.2$ & $69.9$ & $66.8$ & $92.0$ & $92.0$\\
 & & Bernoulli & $95.8$--$96.0$ & $70.6$ & $66.5$ & $96.2$ & $96.2$\\
\midrule
$d=4$ & $\nu_0=1/2$ & Normal & $94.5$--$95.5$ & $95.0$ & $94.0$ & $94.0$ & $94.5$\\
 & & Uniform & $94.8$ & $95.1$ & $94.8$ & $94.8$ & $94.2$\\
 & & Mixture & $94.5$--$96.2$ & $94.5$ & $96.0$ & $96.2$ & $96.5$\\
 & & Bernoulli & $94.2$--$95.0$ & $94.5$ & $95.8$ & $95.8$ & $96.5$\\
 & $\nu_0=3/2$ & Normal & $92.2$--$94.8$ & $95.0$ & $95.8$ & $95.8$ & $95.8$\\
 & & Uniform & $93.0$--$95.5$ & $99.6$ & $99.8$ & $95.5$ & $95.5$\\
 & & Mixture & $93.5$--$95.0$ & $71.2$ & $72.2$ & $95.8$ & $95.8$\\
 & & Bernoulli & $93.5$--$96.8$ & $71.9$ & $67.8$ & $95.0$ & $95.0$\\
\midrule
$d=5$ & $\nu_0=1/2$ & Normal & $93.5$--$94.8$ & $95.0$ & $94.5$ & $94.5$ & $94.8$\\
 & & Uniform & $94.2$--$95.2$ & $95.0$ & $93.8$ & $93.8$ & $93.8$\\
 & & Mixture & $93.8$--$95.2$ & $94.9$ & $95.2$ & $95.2$ & $95.2$\\
 & & Bernoulli & $93.0$--$95.0$ & $94.9$ & $95.8$ & $96.0$ & $96.0$\\
 & $\nu_0=3/2$ & Normal & $91.2$--$96.8$ & $95.0$ & $96.5$ & $96.5$ & $96.5$\\
 & & Uniform & $92.8$--$94.2$ & $99.1$ & $98.8$ & $96.2$ & $96.2$\\
 & & Mixture & $92.8$--$96.0$ & $74.7$ & $72.8$ & $93.8$ & $93.8$\\
 & & Bernoulli & $91.0$--$94.0$ & $75.4$ & $74.0$ & $95.2$ & $95.2$\\
\bottomrule
\end{tabular}
\end{table}

\noindent\textbf{Estimated standard errors.}
We also compare different estimators of the
asymptotic standard errors. For $d=3$,
inverting only the identifiable information
block substantially underestimates the standard errors of smoothness and the microergodic parameter, leading to undercoverage even at the largest
sample size. In contrast, extracting the
identifiable block from the inverse of the full information matrix provides coverage close to the nominal level. The Fisher information and observed Hessian evaluated at the MLE give similar results for sufficiently large $n$, although the Fisher information is more stable at smaller sample sizes. Detailed comparisons are reported in \Cref{tab:sim-se} and \Cref{fig:sim-se,fig:sim-se-d4,fig:sim-se-d5}.

\noindent\textbf{Joint confidence regions.} Joint coverage is appreciably less accurate than marginal coverage
(Table~\ref{tab:sim-joint}). For $d=5$, $\nu_0=1/2$, the nominal $95\%$ log-scale region covers only $76\%$--$82\%$ at $n=10^5$. To understand this shortfall, Table~\ref{tab:sim-joint-steps} compares five calculations on the same $400$ Gaussian data sets.

Write $\theta=(\eta,\nu,a,t)$, $s_n=\nabla_\theta\ell_n(\theta_0)$, and $I_0=\calI_n(\theta_0)$. The first-order approximation to the MLE error is $I_0^{-1}s_n$, whose covariance under Gaussian errors is exactly $I_0^{-1}$. The reported log-parameter vector is $\zeta=(\log\nu,\log\sigma^2,\log l,\log\tau^2)$. For each row of the table, we calculate $v^\top V^{-1}v$ using the error vector $v$ and standardizing covariance $V$ listed there, and report the percentage below $\chi^2_{4,0.95}$. The first row is a score diagnostic; the remaining rows give coverage of Wald regions.

\begin{table}[htbp]
\centering\small
\caption{Successive comparisons for $d=5$, $\nu_0=1/2$,
$n=10^5$, under Gaussian errors. All rows use the same $400$
data sets and the same $\chi^2_{4,0.95}$ cutoff.
``At the truth'' means that both the information and the
Jacobian used for the delta method are evaluated at $\theta_0$.
In row (5), both are evaluated at the MLE.}
\label{tab:sim-joint-steps}
\setlength{\tabcolsep}{4pt}\begin{tabular}{@{}p{0.35\linewidth}p{0.45\linewidth}r@{}}
\toprule
Error vector $v$ & Standardizing covariance $V$
& Below cutoff (\%)\\
\midrule
(1) First-order approximation: $I_0^{-1}s_n$
& Inverse information at the truth, $I_0^{-1}$
& $95.00$\\[5pt]

(2) Actual MLE error: $\widehat\theta_n-\theta_0$
& The same inverse information, $I_0^{-1}$
& $86.00$\\[5pt]

(3) Actual log-parameter error: $\widehat\zeta_n-\zeta_0$
& Inverse information at the truth, transformed to $\zeta$
by the delta method; no cumulant correction
& $82.25$\\[5pt]

(4) The same log-parameter error
& The covariance in (3), with the estimated cumulant
correction, using $t_0$
& $82.00$\\[5pt]

(5) The same log-parameter error
& Fitted information and delta method, with the estimated
cumulant correction, using $\widehat t_n$
& $79.00$\\
\bottomrule
\end{tabular}
\end{table}

Rows (1) and (2) use the same inverse information at the truth. Replacing the first-order approximation by the actual MLE error lowers the proportion below the common $\chi^2_4$ cutoff from $95\%$ to $86\%$. Passing from (2) to (3) changes the parameterization and applies the delta method, while still using information at the truth; this lowers coverage by a further $3.75$ points. The estimated cumulant correction in (4) lowers coverage by $0.25$ percentage points. Finally, replacing the truth by the fitted parameters in the covariance calculation lowers coverage from $82\%$ to $79\%$. The paired Monte Carlo standard errors of the differences between (1) and (2), and between (4) and (5), are both about $1.6$ percentage points.

These comparisons point primarily to the first-order MLE expansion remainder, alongside effects from reporting coordinates and covariance estimation; Appendix~\ref{apdx:sim-joint-diagnostic} provides additional sample-size diagnostics and comparisons. Thus, accurate finite-sample joint confidence regions may require substantially larger sample sizes than marginal inference.

\begin{table}[ht]\centering\small
\caption{Coverage (\%) of nominal $95\%$ joint regions for the estimable log-parameter vector at the largest sample size, using the full inverse Fisher information at the fit and the estimated cumulant correction. Ranges span the four error laws, and both smoothness values in the first row. Each estimate uses $400$ replicates; its Monte Carlo standard error is about $1.2$--$2.1$ percentage points over the displayed range. Every covariance estimate was positive definite, so all $400$ replicates contribute to each coverage estimate.}
\label{tab:sim-joint}
\begin{tabular}{@{}llrr@{}}\toprule
Dimension & Smoothness & $n$ & Joint coverage\\\midrule
$d=3$ & $1/2,\ 3/2$ & $97{,}336$ & $89$--$94$\\
 \hline
$d=4$ & $1/2$ & $104{,}976$ & $85$--$89$\\
 & $3/2$ & $104{,}976$ & $89.5$--$92.5$\\
 \hline
$d=5$ & $1/2$ & $100{,}000$ & $76$--$82$\\
 & $3/2$ & $100{,}000$ & $88.5$--$91$\\\bottomrule
\end{tabular}
\end{table}
For $d=4,5$, joint ellipsoids formed directly by a linear delta method in $(\sigma^2,l,\nu,\tau^2)$ perform worse, with $58\%$--$81\%$ coverage at the largest $n$. The generally valid construction in Corollary~\ref{cor:uncertainty} instead transforms a log-scale region back to the original units, preserving coverage exactly via a one-to-one transformation. Appendix~\ref{app:raw-public-failure} explains why direct linear ellipsoids need not share this asymptotic validity.

\section{Discussion and future work}\label{sec:discussion}

We establish a sharp fixed-domain asymptotic theory for MLEs of Mat\'ern parameters, proving asymptotic normality and minimax-optimal convergence rates---the first comprehensive characterization of the joint MLE with unknown smoothness and nugget variance. We also construct asymptotically valid confidence intervals and extend the theory to non-Gaussian measurement errors. Simulations show close agreement with the information-based
standard deviations and near-nominal marginal coverage at
large sample sizes. The limiting rate exponents can be
approached slowly, and joint confidence regions exhibit
appreciable undercoverage over the simulated sample sizes. Irregular and random designs are covered by our theory. There are, however, several limitations of our study that motivate important directions for future work.

First, for general locations, dense evaluation of the exact Gaussian likelihood requires $O(n^3)$ computation, which is prohibitive for large datasets. Our simulations exploit additional grid symmetries that are unavailable for general irregular designs. There are scalable GP methods, such as inducing-point, nearest-neighbor, and Vecchia approximations \citep{vecchia1988,banerjee2008,titsias2009,datta2016,liu2020gaussian}, which reduce the computational cost of estimating the parameters. For noiseless Mat\'ern observations with known smoothness and a fixed working range, \citet{zhang2024vecchia} studied estimation
of the microergodic parameter by Vecchia likelihood. An open question for the model studied here is whether a scalable approximation preserves the rates and joint limit of the exact MLE when both smoothness and the nugget variance are unknown. Understanding how the asymptotic behavior of scalable GP estimators compares to that of the MLE would clarify the computational–statistical tradeoff in large-scale GP models.

Next, we assume a zero mean, whereas applications model it via regression covariates or an unknown trend estimated jointly with covariance parameters. Because these nuisance parameters affect likelihood-based inference, extending our consistency, rate, and asymptotic distribution results to unknown means---especially with many covariates or a nonparametric trend---is important.

Finally, we focus on stationary, isotropic Mat\'ern kernels with a single length-scale, whereas spatial dependence varies across directions or locations. Anisotropic and nonstationary models introduce extra parameters and richer dependence, raising fixed-domain identifiability and estimation questions; extending our asymptotic theory and optimal rates to them broadens applicability in environmental sciences.

\section*{Disclosure on generative AI and AI-assisted technologies}

ChatGPT (OpenAI) assisted with developing and checking mathematical arguments and with drafting and revising the manuscript. Claude Code (Anthropic) assisted with implementing, debugging, and validating the simulation code and producing numerical summaries and figures. The authors reviewed and revised the resulting material and remain fully accountable for
the accuracy and integrity of the paper.

\section*{Code availability}

The simulation code and instructions for reproducing the
numerical results are available at
\url{https://github.com/Ameer-eng/euclidean-matern-mle-simulations}.

\bibliographystyle{apalike}
\bibliography{ref}

\appendix
\section*{Appendix}
Appendix~\ref{app:main} states the lemmas used in the paper. Appendix~\ref{app:proof-main} proves Theorem~\ref{thm:normalized}, then explains its main constructions and proves the remaining main-body results. Appendix~\ref{app:lemmas} proves every lemma in Appendix~\ref{app:main}. Appendix~\ref{app:auxiliary} gives auxiliary calculations and interpretations. Appendix~\ref{apdx:sim} gives additional simulation details and results.

\section{Lemmas}\label{app:main}
All auxiliary results are stated here and proved in Appendix~\ref{app:lemmas}, where the proof of Lemma A.$j$ appears in Section C.$j$, $j=1,\cdots,27$. The statements are grouped by their role: probabilistic and spectral foundations, likelihood control, feasible inference, and supporting design and analytic results. Appendix~\ref{app:proof-main} proves the main-body results in the order they are stated. Probabilistic statements are uniform over the error class wherever observations, rather than only latent field quantities, are involved.

\subsection{Notation for the proofs}\label{sec:proofnotation}
\noindent\textbf{Uniformity, feasible sets, and an enclosing rectangle.} Throughout the appendices, $\Theta$ and the compact truth set $\mathcal K$ satisfy Assumption~\ref{ass:param}. For $d\ge4$, compactness gives $\operatorname{dist}(\mathcal K,\Theta^c)>0$. For $d\le3$, put $\mathcal K_\vartheta=\{\vartheta:(\vartheta,a)\in\mathcal K\}$. Choose $r_*>0$ such that the closed $r_*$-neighborhood of $\mathcal K_\vartheta$ lies in $U$. The feasible range interval is $I(\vartheta)=[\underline a(\vartheta),\overline a(\vartheta)]$ there. Write $L_*$ for a common Lipschitz constant of its two endpoints. The interval endpoints, its positive length on this neighborhood, and the positive model parameters have uniform bounds by compactness. A chosen global maximizer is understood to be a measurable function of the responses and, when random, the locations. The stronger rectangle condition on $\mathcal K$ in Theorem~\ref{thm:minimax} is imposed only for that risk comparison;
it does not require $\Theta$ to be a rectangle.

For analytic bounds, fix a compact enclosing rectangle
\[
 \mathcal R=[\eta_-,\eta_+]\times[\nu_-,\nu_+]
             \times[a_-,a_+]\times[t_-,t_+],
 \qquad \Theta\subset\operatorname{int}(\mathcal R),
\]
with $0<\nu_-<\nu_+$, $0<a_-<a_+$, and $0<t_-<t_+$. The Mat\'ern covariance and Gaussian criterion are defined on this whole rectangle, including parameter values outside the optimization set. Set $I_*=[a_-,a_+]$. Bounds involving the parameter endpoints or a ``parameter rectangle'' below refer to $\mathcal R$. Uniform spectral, derivative, and covariance estimates are proved on this enclosing family and then restricted to $\Theta$. No maximization over $\mathcal R$ is required.

In particular, let $\mathcal N_n^\square$ denote the set in \eqref{eq:localset} with $\Theta$ replaced by $\mathcal R$. Its affine image in normalized coordinates is compact and convex, even if $\Theta$ is not. Segments, meshes, and likelihood-process
suprema used to bound the fit can therefore be taken on this larger set. The actual preliminary set remains $\mathcal N_n=\Theta\cap\mathcal N_n^\square$.
For each fixed $H$, the ambient low-dimensional local set is $\Gamma_H=\{h:\|h\|\le H\}\times I_*$; the feasible interval at $h$ is instead
\[
 I_n(h)=I(\vartheta_0+E_nh).
\]
All $h$ in this fixed ball have nonempty feasible intervals for sufficiently large $n$, uniformly, because $H\|E_n\|_{\op}<r_*$. Their endpoints obey
\begin{equation}\label{eq:feasible-endpoint-motion}
 \sup_{\|h\|\le H}
 \max\bigl\{|\underline a(\vartheta_0+E_nh)-\underline a(\vartheta_0)|,
            |\overline a(\vartheta_0+E_nh)-\overline a(\vartheta_0)|\bigr\}
 \le L_*H\|E_n\|_{\op}=o(1).
\end{equation}

For example, a fixed box in $(s,l,\nu,t)$, with $s=\sigma^2$, has
\[
 \begin{split}
 \underline a(\vartheta)
 &=\max\left\{l_+^{-2},
       \left(\frac{e^\eta}{A_d(\nu)s_+}\right)^{1/\nu}\right\},\\
 \overline a(\vartheta)
 &=\min\left\{l_-^{-2},
       \left(\frac{e^\eta}{A_d(\nu)s_-}\right)^{1/\nu}\right\}.
 \end{split}
\]
Indeed $s=e^\eta\{A_d(\nu)a^\nu\}^{-1}$, so the bounds on $s$ give the second lower and upper bounds on $a$, and the lengthscale bounds give the first ones. The maximum and minimum of these smooth functions are locally Lipschitz. Thus the assumption holds when the projected truths stay inside the feasible $(\eta,\nu,t)$ region and the two displayed endpoints remain distinct. More explicitly, for a box with bounds $s_-<s_+$, $l_-<l_+$, $\nu_-<\nu_+$, and $t_-<t_+$, these low-dimensional conditions require the
compact truth set to satisfy
\[
 \nu_-<\nu_0<\nu_+,\qquad t_-<t_0<t_+,
\]
and
\[
 \log A_d(\nu_0)+\log s_- -2\nu_0\log l_+
 <\eta_0<
 \log A_d(\nu_0)+\log s_+ -2\nu_0\log l_-.
\]
Here the endpoints are those of this example box, rather than of the enclosing rectangle $\mathcal R$. The strict inequalities give uniform positive margins on a compact truth set. They require interiority of the identifiable projection, not interiority of $a_0$ in $I(\vartheta_0)$. For $d\ge4$, all four variance--lengthscale coordinates must be interior to their box, since the coordinate transformation
is a smooth bijection. The same observations apply to boxes in the corresponding log coordinates. An arbitrary compact set in low dimension need not have this interval structure; no such unrestricted generalization is asserted.

At fixed locations, probabilities and expectations concern the response law. Uniformity over deterministic designs is understood with a common coverage cube, constants $C,\gamma$, and an integer $n_0$ such that $h_n(Q)\le C n^{-\gamma}$ for every $n\ge n_0$. For the deterministic alternative in Assumption~\ref{ass:sharp},
take $\gamma=1/d$. For random locations, the same quantities are conditional on the realized locations unless the joint law is explicitly specified. In particular, the matrices $J_n,K_n,V_n$ are deterministic in the conditional experiment. Uniform convergence in distribution means convergence in bounded Lipschitz distance, uniformly over $\mathcal K\times\mathcal P$, where
\[
 d_{\mathrm{BL}}(\mu,\nu)
 =\sup_{\substack{\|f\|_\infty\le1\\
                  |f(x)-f(y)|\le\|x-y\|}}
       \left|\int f\,d\mu-\int f\,d\nu\right|.
\]
For probability laws on the finite-dimensional Euclidean spaces considered here, convergence in distribution is equivalent to convergence in $d_{\mathrm{BL}}$; see
\citet[Section~1.12, Theorem~1.12.4]{vandervaart2023weak}.

\noindent\textbf{Spectral density and likelihood.} The signal spectral density is \eqref{eq:spectrum}. Write $s=\sigma^2$, the inverse parameter map is
\begin{equation}\label{eq:inverse}
 c=e^\eta/A_d(\nu),\qquad
 s=e^\eta\{A_d(\nu)a^\nu\}^{-1},\qquad l=a^{-1/2}.
\end{equation}

The criterion is defined in \eqref{eq:likelihood}; the global fit maximizes it over $\Theta$ as specified in Section~\ref{sec:setup}.

\noindent\textbf{Normalizations.}
The subscript ``id'' denotes differentiation with respect to $\vartheta$, holding $a$ fixed in low dimension. Put
\begin{equation}\label{eq:rates}
\kappa_0=2\nu_0+d,\quad q_0=d/\kappa_0,\quad
R_n=n^{1/\kappa_0},\quad M_n=R_n^d,\quad
r_n=\sqrt{M_n},\quad L_n=\log R_n.
\end{equation}
For $d\ge4$, set
\begin{equation}\label{eq:arate}
b_n=\begin{cases}\sqrt{\log n},&d=4,\\R_n^{(d-4)/2},&d>4.\end{cases}
\end{equation}
The letter $b_n$ here denotes the range normalization, not a kernel parameter. For the centered increments in \eqref{eq:shear}, we write $u,v$ when $n$ is understood. The local parameter map $E_n$ is given in \eqref{eq:E}. Its inverse is
\[
 E_n^{-1}=\begin{cases}
 \begin{pmatrix}
 r_n&-2L_nr_n&0\\0&r_n&0\\0&0&\sqrt n
 \end{pmatrix},&d\le3,\\[3mm]
 \begin{pmatrix}
 r_n&-2L_nr_n&0&0\\0&r_n&0&0\\
 0&0&b_n&0\\0&0&0&\sqrt n
 \end{pmatrix},&d\ge4.
 \end{cases}
\]
Thus $h_n(\theta)=E_n^{-1}(\vartheta-\vartheta_0)$ and $\vartheta=\vartheta_0+E_nh$ in both dimension regimes, in the coordinate orders already specified. These matrices normalize errors at the truth. The theoretical conclusions concern global maximizers independently of the coordinates used to compute them.

\noindent\textbf{Derivatives in normalized coordinates.} For $d\ge4$, write $\theta(h)=\theta_0+E_nh$ and let $\Phi$ be a differentiable scalar function of $\theta$. The chain rule gives, component by component,
\[
 \partial_{h_j}\Phi(\theta(h))
 =\sum_i\partial_{\theta_i}\Phi(\theta(h))
                  \frac{\partial\theta_i}{\partial h_j}
 =\sum_i(E_n)_{ij}\partial_{\theta_i}\Phi(\theta(h)).
\]
Thus, with row Jacobians and column gradients respectively,
\[
 \frac{\partial\Phi}{\partial h}
       =\frac{\partial\Phi}{\partial\theta}E_n,
 \qquad
 \nabla_h\Phi=E_n^\T\nabla_\theta\Phi,
 \qquad \nabla_h=E_n^\T\nabla_\theta.
\]
All functions on the right are evaluated at $\theta(h)$. For $d\le3$, replace $\theta$ by $\vartheta=(\eta,\nu,t)$ and hold $a$ fixed. In particular, the columns of $E_n$ give
\[
 \partial_{h_u}=r_n^{-1}\partial_\eta,\qquad
 \partial_{h_v}=r_n^{-1}(\partial_\nu+2L_n\partial_\eta),\qquad
 \partial_{h_a}=b_n^{-1}\partial_a,\qquad
 \partial_{h_t}=n^{-1/2}\partial_t,
\]
with the range coordinate omitted when $d\le3$.

For a signal direction $h$, we suppress its zero nugget component. Its relative spectral derivative at the truth is
\[
 \begin{aligned}
 g_h(\omega)
 &\coloneqq \frac{\left.\nabla_z f_{\theta_0+E_nz}(\omega)
                                \right|_{z=0}^{\!\T}h}{f_0(\omega)}\\
 &=\left.\nabla_z\log f_{\theta_0+E_nz}(\omega)
                                \right|_{z=0}^{\!\T}h.
 \end{aligned}
\]
Here $z$ denotes the local coordinate at which the derivative is taken, while $h$ denotes the direction. In low dimension the parameter in this display is $(\vartheta_0+E_nz,a_0)$. We abbreviate the derivative as $g_h(\omega)=(\nabla_h\log f_\theta(\omega)|_{\theta=\theta_0})^\T h$; evaluation at the truth is essential in this identity.

These derivatives explain the normalization. Write $m_0(\omega)=a_0+4\pi^2\norm\omega^2$. Since $\log f_\theta=\eta-(\nu+d/2)\log m_a$,
\[
 \partial_\eta\log f_\theta=1,\qquad
 \partial_\nu\log f_\theta=-\log m_a,\qquad
 \partial_a\log f_\theta=-\frac{\nu+d/2}{m_a}.
\]
At the informative frequency scale $\omega=R_nz$,
\[
 \begin{aligned}
 (\partial_\nu+2L_n\partial_\eta)
          \log f_\theta(R_nz)\big|_{\theta_0}
 &=-\log\{a_0+4\pi^2R_n^2\norm z^2\}+2\log R_n\\
 &=-\log\{a_0/R_n^2+4\pi^2\norm z^2\}.
 \end{aligned}
\]
Thus the $2L_n$ entry of $E_n$ cancels the growing constant in the smoothness derivative. There are order $R_n^d$ informative frequencies, which suggests the amplitude and smoothness scale $r_n=R_n^{d/2}$. The squared range derivative instead accumulates like $\int_1^{R_n}r^{d-5}\,dr$, giving $b_n^2=\log n$ at $d=4$ and $b_n^2=R_n^{d-4}$ at $d>4$, up to constants. The nugget uses the $\sqrt n$ scale. The explicit formula for $g_h$ is recorded in \eqref{eq:signaltangent}.

The expected curvature is defined in \eqref{eq:info}, and $J_n=E_n^\top\calI_nE_n$ as in \eqref{eq:normalizedscoreinfo}. The same trace formula with all four coordinate derivatives defines the full Gaussian information. Corollary~\ref{cor:uncertainty} also justifies covariance estimation by the identifiable block of the inverse full expected curvature or observed Hessian, with the cumulant correction. The normalized error $\widehat h_n$, score $\Delta_n$, and standardized error $V_n^{-1/2}\widehat h_n$ are defined in Section~\ref{sec:mle}.

In bounds, $c,C>0$ denote constants whose values may change from line to line. They may depend on fixed geometry, model bounds, and sampling constants, and on the common noise-class bounds, but not on $n$, the individual law $P$, or the varying parameter in a stated uniform bound. A subscript, as in $C_H$ or $C_{A,B}$, records an additional
dependence when it matters. We write $a_n\asymp b_n$ for the pair of inequalities $cb_n\le a_n\le Cb_n$. Norms with subscript $\Frob$ are Frobenius norms; $\op$ denotes the operator norm. For a random quantity $Z_n$, the notation $Z_n=o_{L^p}(1)$ means
$\sup_{\theta_0\in\mathcal K,\,P\in\mathcal P}\E_{\theta_0,P}\|Z_n\|^p\to0$. When the locations are held fixed, this expectation is conditional on those locations. We use this notation with each fixed $1\le p<\infty$; constants may then also depend on $p$. These probability-space norms are distinct from the function-space norms defined below.

\noindent\textbf{Log coordinates.}
For calculations we also use the log chart
\begin{equation}\label{eq:logcharts}
 \lambda=
 \begin{cases}
 (\log c,\log\nu,\log t),&d\le3,\\
 (\log c,\log\nu,\log a,\log t),&d\ge4.
 \end{cases}
\end{equation}
Throughout the coordinate-change arguments, write $\lambda=g(\vartheta)$ for the selected forward map and use $Dg(\vartheta)$ for its Jacobian. For this chart, with
$B_0=\partial_\nu\log A_d(\nu_0)$,
\[
 Dg(\vartheta_0)=
 \begin{cases}
 \begin{pmatrix}1&-B_0&0\\0&\nu_0^{-1}&0\\0&0&t_0^{-1}\end{pmatrix},&d\le3,\\[3mm]
 \begin{pmatrix}1&-B_0&0&0\\0&\nu_0^{-1}&0&0\\0&0&a_0^{-1}&0\\0&0&0&t_0^{-1}\end{pmatrix},&d\ge4.
 \end{cases}
\]
The information in these coordinates is
\[
 \mathcal I_n^\lambda
 =Dg(\vartheta_0)^{-\T}\mathcal I_n
       Dg(\vartheta_0)^{-1},
 \qquad
 (Dg(\vartheta_0)E_n)^\T\mathcal I_n^\lambda
       Dg(\vartheta_0)E_n=J_n.
\]
For $d\le3$, this is the microergodic log vector in Table~\ref{tab:parameterizations}. For $d\ge4$, we also use the variance--lengthscale log vector $(\log s,\log l,\log\nu,\log t)$, where $s=\sigma^2$. Its forward Jacobian is displayed in \eqref{eq:rawJacobian}.

\noindent\textbf{Localized Fourier coefficients and their predictors.}
Let $K_0$ denote the true signal covariance function, with spectral density $f_0=f_{\theta_0}$. In particular $C_{n,0}$ also includes the nugget, whereas $K_0$ does not.

Our basic latent statistics will be localized Fourier coefficients. Choose the working cube $Q$ with faces parallel to the coordinate axes and with closure strictly inside the cube on which the coverage or density assumption holds. Choose it small enough for the interior sampling inequality in Lemma~\ref{lem:sampling}; this choice is made before fixing any Fourier coefficient or cutoff. Let its side length be $L$, and choose a nonzero real cutoff $\chi\in C_c^\infty(Q)$, with the more specific choice in
Lemma~\ref{lem:windows}. Fix the integer $m\ge3$ used in that lemma and use the single lattice $\mathscr L\coloneqq (m/L)\Z^d$ throughout. For $\lambda\in\mathscr L$, put
\begin{equation}\label{eq:localizedfourier}
 X_\lambda=L^{-d/2}\int_Q\chi(x)F(x)
                         e^{-2\pi i\lambda^\T x}\,dx,
 \qquad v_\theta(\lambda)=\E_\theta|X_\lambda|^2.
\end{equation}
These are Fourier series coefficients of $\chi F$ on $Q$, retaining every $m$th frequency in each coordinate. This fixed selection will give a uniform lower bound for the band spaces in Lemma~\ref{lem:windows}. The cutoff makes the construction local to the observation region; the law of $F$ remains the Euclidean stationary law. Since $F$ and $\chi$ are real, $X_{-\lambda}=\overline{X_\lambda}$. We may either keep one coefficient from each pair or work with its real and imaginary parts.

The coefficients are latent. Preliminary localization uses nearest-cell quadrature, uniformly over the covariance family. The sharp curvature bound uses their best linear predictions
\begin{equation}\label{eq:fourierpredictor}
 \begin{split}
 \widehat X_{\lambda,n}
 &=d_{\lambda,n}C_{n,0}^{-1}Y_n,\\
 [d_{\lambda,n}]_i
 &=\E_0[X_\lambda Y_{n,i}]
   =L^{-d/2}\int_Q\chi(x)K_0(x-x_{n,i})
                   e^{-2\pi i\lambda^\T x}\,dx.
 \end{split}
\end{equation}
Equivalently, these are the Fourier coefficients of the predicted field after applying the same cutoff:
\[
 \widehat X_{\lambda,n}
 =L^{-d/2}\int_Q\chi(x)\widehat F_0(x)
                      e^{-2\pi i\lambda^\T x}\,dx,
 \qquad \widehat F_0(x)=\Cov_0(F(x),Y_n)C_{n,0}^{-1}Y_n.
\]
All prediction weights here are fixed at $\theta_0$, including when their distributions or mean derivatives are evaluated under another parameter. They are auxiliary statistics for proving an information bound, not additional operations required of the MLE. The prediction weights and their Fourier responses will suffice for the analytic calculations. We use the Fourier convention
\[
 \widehat g(\omega)=\int_{\R^d}g(x)e^{-2\pi i\omega^\T x}\,dx,
 \qquad
 \check q(x)=\int_{\R^d}q(\omega)e^{2\pi i\omega^\T x}\,d\omega.
\]
Inverse transforms are interpreted distributionally when necessary. Write $m_0(\omega)=a_0+4\pi^2\norm\omega^2$. Normalize each coefficient by its latent standard deviation at the truth:
\begin{equation}\label{eq:fourierweights}
 \begin{alignedat}{2}
 U_\lambda&=X_\lambda/\sqrt{v_0(\lambda)},&\qquad
 \widehat U_{\lambda,n}
       &=\widehat X_{\lambda,n}/\sqrt{v_0(\lambda)},\\
 b_{\lambda,n}&=
       \frac{d_{\lambda,n}C_{n,0}^{-1}}{\sqrt{v_0(\lambda)}},
 &\qquad
 \widehat U_{\lambda,n}&=\sum_{i=1}^n b_{\lambda,n,i}Y_{n,i},\\
 \ell_\lambda(\omega)&=
       \frac{L^{-d/2}\widehat\chi(\omega+\lambda)}
            {\sqrt{v_0(\lambda)}},&\qquad
 p_{\lambda,n}(\omega)
       &=\sum_{i=1}^n b_{\lambda,n,i}e^{-2\pi i\omega^\T x_{n,i}}.
 \end{alignedat}
\end{equation}
Here $b_{\lambda,n}$ is a row of complex prediction weights. For the whole-space integrals below, take $F$ to be a centered stationary Gaussian field on $\R^d$ with spectral density $f_\theta$, independent of the errors and of any random locations. Its restriction to $D$ has the covariance in \eqref{eq:matern-intro} and hence the
same Gaussian law as the field in \eqref{eq:model}. For compactly supported test functions, write $Z_g=\int_{\R^d}g(x)F(x)\,dx$; also set $Z_{\delta_x}=F(x)$ and extend linearly to finite sums of point masses and test functions. The negative Sobolev isometry in the proof of Lemma~\ref{lem:sampling} gives the general extension of this notation. The inverse transforms of the two Fourier weights are
\[
 \check\ell_\lambda(x)
 =\frac{L^{-d/2}\chi(x)e^{-2\pi i\lambda^\T x}}
              {\sqrt{v_0(\lambda)}},\qquad
 \check p_{\lambda,n}=\sum_i b_{\lambda,n,i}\delta_{x_{n,i}}.
\]
Consequently, with $\varepsilon_n=Y_n-F(x_{n,1:n})$,
\[
 U_\lambda=Z_{\check\ell_\lambda},\qquad
 \widehat U_{\lambda,n}
       =Z_{\check p_{\lambda,n}}+b_{\lambda,n}\varepsilon_n.
\]
Thus $\ell_\lambda$ and $p_{\lambda,n}$ are the Fourier transforms of the test function and the discrete measure defining these signal statistics. For complex random vectors we use $\Cov_0(V,W)=\E_0[(V-\E_0V)(W-\E_0W)^*]$, where the star denotes conjugate transpose, and write $\Cov_0(V)=\Cov_0(V,V)$. Put
\[
 e_\omega(x)=e^{2\pi i\omega^\T x},\qquad
 e_\omega(x_{n,1:n})=(e_\omega(x_{n,i}))_{i=1}^n,
\]
and use the same convention for $F(x_{n,1:n})$. The spectral covariance identities in Lemma~\ref{lem:spectral-covariance} apply to both the integral $Z_{\check\ell_\lambda}$ and the finite sum $Z_{\check p_{\lambda,n}}$. In particular,
\eqref{eq:fouriernormalizations} records their second moments and the separate contribution of observation noise. These identities require no independence between Fourier coefficients.

\noindent\textbf{Covariance matrices and frequency weights.}
For covariance derivatives, it is convenient to give one name to a matrix of Fourier integrals. Put
\begin{equation*}
 \begin{split}
 a_\omega&=e_\omega(x_{n,1:n}),\\
 \mathcal B_n(q)
 &=C_{n,0}^{-1/2}
       \left\{\int_{\R^d}q(\omega)f_0(\omega)
                          a_\omega a_\omega^*\,d\omega\right\}
                          C_{n,0}^{-1/2}.
 \end{split}
\end{equation*}
Thus $\mathcal B_n(q)$ is an $n\times n$ matrix. For real $q$ the matrix $\mathcal B_n(q)$ is Hermitian; for the even real multipliers used in the covariance derivatives it is real symmetric. In particular,
\[
 P_n\coloneqq \mathcal B_n(1)
   =C_{n,0}^{-1/2}S_{n,0}C_{n,0}^{-1/2}
   =I_n-t_0C_{n,0}^{-1},\qquad 0\preceq P_n\preceq I_n.
\]
The map $q\mapsto\mathcal B_n(q)$ is linear and preserves pointwise order for real multipliers. Hence $|q|\le c$ implies $-cI_n\preceq\mathcal B_n(q)\preceq cI_n$.
For any signal derivative, the spectral covariance formula gives
\[
 C_{n,0}^{-1/2}(\partial_h C_n)_{\theta_0}C_{n,0}^{-1/2}
   =\mathcal B_n\left(\frac{(\partial_h f_\theta)_{\theta_0}}{f_0}\right).
\]
The nonnegative scalar weight
\begin{equation}\label{eq:frequencyweightdefinition}
 w_n(\omega)=f_0(\omega)a_\omega^*C_{n,0}^{-1}a_\omega
\end{equation}
satisfies $\tr\mathcal B_n(q)=\int q w_n$ whenever the integral is absolutely convergent. Lemma~\ref{lem:upper} bounds this weight and the resulting matrix norms.

\noindent\textbf{Function-space norms.}
For functions, $\norm g_{L^2}$ denotes the norm in $L^2(\R^d;\mathbb C)$:
\[
 \norm g_{L^2}^2=\int_{\R^d}|g(x)|^2\,dx.
\]
Functions given only on a cube are extended by zero when this whole-space norm is used. A displayed domain, as in $\norm g_{L^2(Q)}$, specifies a norm on that domain.
We use the Japanese bracket
$\langle x\rangle=(1+\norm x^2)^{1/2}$ and the elementary bounds
\[
 \langle x+y\rangle\le\sqrt2\langle x\rangle\langle y\rangle,
 \qquad
 \frac{\langle x\rangle^q}{\langle y\rangle^q}
 \le2^{|q|/2}\langle x-y\rangle^{|q|},\qquad q\in\R.
\]
The sampling argument also uses Sobolev norms. Our convention is
\[
 \norm F_{H^s(\R^d)}^2
   =\int_{\R^d}\langle\omega\rangle^{2s}|\widehat F(\omega)|^2\,d\omega,
 \qquad
 \norm f_{H^s(\Omega)}
   =\inf_{\substack{F\in H^s(\R^d)\\F|_\Omega=f}}
                          \norm F_{H^s(\R^d)},
\]
where $\Omega$ is an enclosing cube when used below. Thus the domain norm is the smallest whole-space extension norm. This convention fixes the fractional orders in the
sampling estimates without requiring a choice of an equivalent intrinsic norm.

\noindent\textbf{Auxiliary Gaussian calculations.}
At the given locations, write $P_{\theta,P}$ for the actual response law when the standardized errors have law $P$, and put $G_{\theta,n}=N(0,C_n(\theta))$ for the Gaussian law with the same covariance. Also write $G_{0,n}=G_{\theta_0,n}$. Unless explicitly marked otherwise, expectations involving observations are under $P_{\theta,P}$; the noise-law index is suppressed when it is immaterial to a covariance identity. The matrices $\calI_n,J_n$ are Gaussian expected curvatures, not the actual Fisher information of $P_{\theta,P}$. In the proof of Lemma~\ref{lem:information} we work explicitly under $G_{\theta,n}$ to bound these deterministic matrices. Gaussian score and conditional expectation identities in that proof are confined to that auxiliary experiment. Elsewhere Fourier predictors denote best linear predictors; they need not be conditional expectations for the actual non-Gaussian law.

\noindent\textbf{The centered smoothness multiplier.}
The centered smoothness multiplier used below is
\begin{equation}\label{eq:spectralcentering}
 g_n(\omega)=2\log R_n-\log m_0(\omega),\qquad
 m_0(\omega)=a_0+4\pi^2\norm\omega^2.
\end{equation}

\subsection{Probabilistic and spectral foundations}\label{sec:lemma-foundations}

The random terms in the Gaussian criterion and its derivatives are quadratic forms in the observations. Their Frobenius norms determine the size of the fluctuations, while a small operator norm ensures that no single direction dominates. The following lemma supplies the concentration, moment, and normal-approximation results used throughout the likelihood arguments, together with the correction to the covariance caused by non-Gaussian errors.

\begin{lemma}[Quadratic-form bounds and the score limit]
\label{lem:gaussian}
Write $Y=S_0^{1/2}G+\sqrt{t_0}\xi$, with independent $G\sim N(0,I_n)$ and standardized errors satisfying Assumption~\ref{ass:noise}. Put $X=C_0^{-1/2}Y$ and $Q(A)=X^\T AX-\tr A$ for real symmetric matrices $A$. Then
\begin{equation}\label{eq:quadratic-cumulant}
 \Cov\{Q(A),Q(B)\}
 =2\tr(AB)+\kappa_4\sum_i
 (C_0^{-1/2}AC_0^{-1/2})_{ii}
 (C_0^{-1/2}BC_0^{-1/2})_{ii}.
\end{equation}
Uniformly in the noise class, $\min\{2,v_*\}\norm A_\Frob^2\le\Var Q(A)\le C\norm A_\Frob^2$. For every $z\ge0$,
\begin{equation}\label{eq:quadconc}
 P\{|Q(A)|>C(\|A\|_{\rm F}\sqrt z+\|A\|_{\rm op}z)\}
 \le2e^{-z}.
\end{equation}
Consequently $\|Q(A)\|_{L^p}\le C(\sqrt p\|A\|_{\rm F}
+p\|A\|_{\rm op})$ for $p\ge1$. We also use the corresponding exponential-moment bound:
\[
 \log \E e^{tQ(A)}\le C t^2\|A\|_\Frob^2,
 \qquad |t|\|A\|_\op\le c.
\]
The constants depend only on the common sub-Gaussian bound. If $\norm{A_n}_\Frob\le C$, $\norm{A_n}_\op\to0$, and $\Var Q(A_n)\ge c>0$, then $Q(A_n)/\sqrt{\Var Q(A_n)}\Rightarrow N(0,1)$ uniformly whenever these bounds and the operator-norm convergence hold uniformly. The vector version holds by Cram\'er--Wold with uniform covariance lower bounds.

The same upper bounds hold for $Q_\varepsilon(A)=\varepsilon^\T A\varepsilon-t_0\tr A$,
with $\E|Q_\varepsilon(A)|\le C\norm A_\Frob$; no positive lower variance bound is needed for these upper bounds.

\end{lemma}

The localized Fourier coefficients provide a way to read the spectrum from the field on a bounded domain. At a large frequency, their variances and logarithmic derivatives closely follow those of the spectral density, and their correlations decay as the
frequencies separate. The next lemma makes this description precise and shows that nearest-site quadrature recovers the slowly growing frequency bands used for preliminary localization.

\begin{lemma}[Localized Fourier coefficients, recovery, and band spaces]
\label{lem:windows}
Assume \Cref{ass:param}. Let $Q$ be the axis-parallel coverage cube of side $L$, put
$h_n=h_n(Q)$, and fix $0<\zeta<\min(\nu_-,1)$. For a fixed integer $m\ge3$, choose a translate $Q_0$ of $(0,L/m)^d$ with $Q_0\Subset Q$, and a real cutoff $\chi\in C_c^\infty(Q)$ equal to $1$ on $Q_0$. Use the Fourier coefficients \eqref{eq:localizedfourier} on $\mathscr L=(m/L)\Z^d$ for all frequency sets below. This lattice is chosen so that its exponentials are orthogonal on $Q_0$; its fixed spacing changes only constants in frequency counts. Real test functions are extended by zero outside $Q$. Put $Z_g=\int_Qg(x)F(x)\,dx$ and
\[
 \norm g_{H^s}^2=\int_{\R^d}\langle\omega\rangle^{2s}
                         |\widehat g(\omega)|^2\,d\omega,
 \qquad\beta=\nu+d/2.
\]
\begin{enumerate}[label=(\roman*),ref=(\roman*)]
\item\label{it:windows-recovery}
\emph{Integral recovery.} Assign each point of $Q$ to a nearest observation site, breaking ties by index, and denote the cells by $V_{n,i}$. For $g\in L^2(Q)$ define
\[
 w_{n,i}(g)=\int_{V_{n,i}}g(x)\,dx,\qquad
 \widehat Z_g=\sum_iw_{n,i}(g)Y_{n,i}.
\]
These statistics are linear in $g$, and
\begin{equation}\label{eq:intrecovery}
 \sup_{\theta\in\Theta}\E_\theta|\widehat Z_g-Z_g|^2
 \le C(h_n^{2\zeta}+h_n^d)\norm g_{L^2}^2.
\end{equation}
For an $L^2(Q)$-orthonormal family $g_1,\ldots,g_N$, the error vector consequently has covariance at most $C(h_n^{2\zeta}+h_n^d)I_N$, with $C$ independent of $N$. For Fourier coefficients the same construction is
\begin{equation}\label{eq:fourierquadrature}
 \widetilde X_{\lambda,n}
 =L^{-d/2}\sum_i
     \left\{\int_{V_{n,i}}\chi(x)e^{-2\pi i\lambda^\T x}\,dx\right\}
     Y_{n,i}.
\end{equation}
Its weights use only the sites and the test function. The bound \eqref{eq:intrecovery} applies to every real linear combination of its real and imaginary parts, uniformly over $\theta$.

\item\label{it:windows-space}
\emph{Test functions near radius $R$.} Fix a relative bandwidth $\omega_*\in(0,1/2)$ and set
\[
 \Lambda_R^+=\mathscr L\cap B(Re_1,\omega_*R/4),\qquad
 \Lambda_R=\Lambda_R^+\cup(-\Lambda_R^+),
\]
and define the real space
\[
 \mathcal V_R=\left\{L^{-d/2}\chi(x)
        \sum_{\lambda\in\Lambda_R}b_\lambda e^{2\pi i\lambda^\T x}:
                 b_{-\lambda}=\overline{b_\lambda}\right\}.
\]
For all sufficiently large $R$, the clusters are disjoint, $\#\Lambda_R^+\asymp R^d$, and $\dim_{\R}\mathcal V_R=2\#\Lambda_R^+\asymp R^d$. The vector $(\Re X_\lambda,\Im X_\lambda)_{\lambda\in\Lambda_R^+}$ therefore generates exactly the latent integrals $Z_g$, $g\in\mathcal V_R$.

\item\label{it:windows-leak}
\emph{Mass outside the band.} Define
\[
 \mathcal A_R=\{\omega:\ |\norm\omega-R|\le\omega_*R\}.
\]
For any fixed $A,B\ge0$,
\begin{equation}\label{eq:windowleak}
 \int_{\R^d\setminus\mathcal A_R}\langle\omega\rangle^{2B}
             |\widehat g(\omega)|^2\,d\omega
 \le C_{A,B}R^{-A}\norm g_{L^2}^2,
 \qquad g\in\mathcal V_R.
\end{equation}
This is the Fourier leakage estimate: the part of the Fourier mass outside the stated band is smaller than any fixed power of $R^{-1}$, also with the indicated polynomial weight.

\item\label{it:windows-sobolev}
\emph{Sobolev scaling.} For $s$ in any fixed compact interval,
\begin{equation}\label{eq:windowsobolev}
 cR^{2s}\norm g_{L^2}^2\le\norm g_{H^s}^2
       \le CR^{2s}\norm g_{L^2}^2,\qquad g\in\mathcal V_R.
\end{equation}

\item\label{it:windows-matern}
\emph{Variances and ratios.} For $g\in L^2(Q)$, $\Var_\theta(Z_g)=\int f_\theta|\widehat g|^2$. Uniformly in $\Theta$,
\begin{equation}\label{eq:windowvariance}
 cR^{-2\beta}\norm g_{L^2}^2\le\Var_\theta(Z_g)
                   \le CR^{-2\beta}\norm g_{L^2}^2,
 \qquad g\in\mathcal V_R.
\end{equation}
For every fixed $A>0$, all large $R$, and nonzero
$g\in\mathcal V_R$,
\begin{equation}\label{eq:windowratio}
 (1-C_AR^{-A})\inf_{\omega\in\mathcal A_R}\frac{f_\theta}{f_{\theta'}}
 \le\frac{\Var_\theta(Z_g)}{\Var_{\theta'}(Z_g)}
 \le(1+C_AR^{-A})\sup_{\omega\in\mathcal A_R}\frac{f_\theta}{f_{\theta'}}.
\end{equation}
The same bounds hold for every generalized eigenvalue of the latent-integral covariance matrices in any basis of $\mathcal V_R$.

\item\label{it:windows-fourier}
\emph{Covariances and energies of Fourier coefficients.}
For the coefficients in \eqref{eq:localizedfourier},
\begin{equation}\label{eq:fouriervariance}
 v_\theta(\lambda)
 =L^{-d}\int f_\theta(\omega)
                   |\widehat\chi(\lambda-\omega)|^2\,d\omega.
\end{equation}
Put
\[
 c_\chi=L^{-d}\int|\widehat\chi(y)|^2\,dy
       =L^{-d}\norm\chi_{L^2}^2>0,\qquad
 (M_\chi)_{ij}=L^{-d}\int y_i y_j|\widehat\chi(y)|^2\,dy.
\]
Uniformly in $\theta$, as $\norm\lambda\to\infty$,
\begin{equation}\label{eq:fourierlocalvariance}
 \begin{split}
 v_\theta(\lambda)
 &=c_\chi f_\theta(\lambda)
    +\frac12\sum_{i,j}(M_\chi)_{ij}
                   \partial_{\lambda_i\lambda_j}f_\theta(\lambda)
    +O\{f_\theta(\lambda)\langle\lambda\rangle^{-4}\},\\
 \frac{v_\theta(\lambda)}{c_\chi f_\theta(\lambda)}
 &=1+O(\langle\lambda\rangle^{-2}).
 \end{split}
\end{equation}
With $m_a(\lambda)=a+4\pi^2\norm\lambda^2$, the corresponding logarithmic derivatives satisfy
\begin{equation}\label{eq:fourierlocalderivatives}
 \begin{split}
 |\partial_\nu\log v_\theta(\lambda)+\log m_a(\lambda)|
     &\le C\langle\lambda\rangle^{-2},\\
 |\partial_a\log v_\theta(\lambda)+\beta/m_a(\lambda)|
     &\le C\langle\lambda\rangle^{-4}.
 \end{split}
\end{equation}
\begin{samepage}
For each $M>0$, uniformly in $\theta$ and $\lambda,\mu\in\mathscr L$,
\begin{align}
 \frac{|\E_\theta[X_\lambda\overline{X_\mu}]|}
      {\sqrt{v_\theta(\lambda)v_\theta(\mu)}}
 &\le C_M\langle\lambda-\mu\rangle^{-M},\label{eq:fouriercovdecay}\\
 \frac{|\E_\theta[X_\lambda X_\mu]|}
      {\sqrt{v_\theta(\lambda)v_\theta(\mu)}}
 &\le C_M\langle\lambda+\mu\rangle^{-M}.\notag
\end{align}
\end{samepage}
Consequently, for every finite $\Lambda\subset\mathscr L$ and real weights $w_\lambda$,
\begin{equation}\label{eq:fourierenergyvariance}
 \Var_\theta\left\{\sum_{\lambda\in\Lambda}w_\lambda
       \left(\frac{|X_\lambda|^2}{v_\theta(\lambda)}-1\right)\right\}
 \le C\sum_{\lambda\in\Lambda}w_\lambda^2.
\end{equation}
The complex covariance matrix of $(X_\lambda/\sqrt{v_\theta(\lambda)})_{\lambda\in\Lambda}$ is bounded by $CI$, independently of $\Lambda$. The covariance matrix of its real and imaginary parts obeys the same bound. No independence of the Fourier coefficients is asserted or needed.

\item\label{it:windows-envelope}
\emph{Envelope for the normalized Fourier weights.}
For $\ell_\lambda$ in \eqref{eq:fourierweights}, every $B\ge1$, and every finite $\Lambda\subset\{\lambda\in\mathscr L: \norm\lambda\le B\}$,
\[
 \sum_{\lambda\in\Lambda}
        f_0(\omega)|\ell_\lambda(\omega)|^2
 \le C\min\{1,\langle B\rangle^{2\beta_0}f_0(\omega)\},
 \qquad \omega\in\R^d.
\]
In particular, fix $0<\epsilon\le1$ and take $B=\epsilon R_n$. For all sufficiently large $n$, the same sum satisfies
\[
 \sum_{\lambda\in\Lambda}f_0(\omega)|\ell_\lambda(\omega)|^2
 \le C\min\{1,(n/t_0)f_0(\omega)\}.
\]
The constants are independent of $B$, $\epsilon$, and the finite set, uniformly over the true parameter in the compact rectangle.

\item\label{it:windows-iid}
\emph{Random coverage.} Suppose that, for each $n$, the sites $x_{n,1},\ldots,x_{n,n}$ are iid with density at least $p_*>0$ on $Q$, with the same lower bound for all $n$. Then almost surely
\[
 h_n\le C(\log n/n)^{1/d}
\]
for all sufficiently large $n$.
\end{enumerate}
Constants are uniform in $n,R,g$ and the parameters in $\Theta$; they may depend on the fixed geometry, $\zeta$, bandwidth, and displayed exponents. The threshold for $R$ may depend on these quantities as well.
\end{lemma}

The next lemma proves the preliminary consistency for $\eta$, $\nu$, and $t$ in all dimensions. 

\begin{lemma}[Preliminary localization across smoothness classes]
\label{lem:preliminary}
Under \Cref{ass:param,ass:coverage,ass:noise}, for some sufficiently small $\delta>0$,
every global MLE obeys
\begin{equation}\label{eq:prelim}
 |\widehat\nu-\nu_0|+|\widehat\eta-\eta_0|/\log n
       +|\widehat t-t_0|=O_{P,\mathcal K\times\mathcal P}(n^{-\delta}).
\end{equation}
Setting $q_n=n^{-\delta}$, every joint likelihood maximizer lies with probability tending to one in
\begin{equation}\label{eq:localset}
\mathcal N_n=\{\theta\in\Theta:
|u|\le q_n\log n,\ |v|\le q_n,\ |t-t_0|\le q_n\},
\quad u=\eta-\eta_0-2L_n v,\ v=\nu-\nu_0.
\end{equation}
There is no additional localization condition on $a$; it remains subject to the original feasible set $\Theta$. The choices can be made so that, for some $c,C,\delta_1>0$,
\begin{equation}\label{eq:prelimfailure}
 \sup_{\theta_0\in\mathcal K,\ P\in\mathcal P}
 P_{\theta_0,P}\{\text{some global MLE lies outside }\mathcal N_n\} \le Ce^{-cn^{\delta_1}}.
\end{equation}
The constants are uniform when the parameter and coverage constants are fixed.
\end{lemma}

We next prove upper bounds for various matrix norms of matrices that appear in the derivatives of $C_{n,0}^{-1/2} C_n(\theta)C_{n,0}^{-1/2}$ with respect to $h$. These upper bounds hold for arbitrary sites in the bounded region. Their different orders explain the normalizing rates in $E_n$, including the logarithmic range rate at the critical dimension, and are the core calculations used in proving $J_n \preceq CI$ and bounding remainders in the Taylor expansion of the log-likelihood.

\begin{lemma}[Universal spectral and derivative bounds]\label{lem:upper}
Under \Cref{ass:param}, for arbitrary sites in a fixed bounded set,
\begin{equation}\label{eq:effective}
\tr P_n\le CM_n.
\end{equation}
For each fixed real $k>0$, write
$m_0(\omega)=a_0+4\pi^2\norm\omega^2$,
\begin{equation}\label{eq:rangeupper}
\norm{\mathcal B_n\bigl(m_0^{-k}\bigr)}_\Frob\le CL_{k,n},\qquad
L_{k,n}=\begin{cases}
1,&d<4k,\\\sqrt{\log n},&d=4k,\\R_n^{(d-4k)/2},&d>4k.
\end{cases}
\end{equation}
The analogous positive trace bound is
\begin{equation}\label{eq:rangetrace}
\tr(\mathcal B_n\bigl(m_0^{-k}\bigr))\le C
\begin{cases}1,&d<2k,\\\log n,&d=2k,\\R_n^{d-2k},&d>2k.\end{cases}
\end{equation}
For the centered logarithm $g_n$ in \eqref{eq:spectralcentering} and each integer $0\le j\le3$,
\begin{align}
\tr(\mathcal B_n\bigl(g_n^{2j}\bigr))&\le CM_n,\quad 0\le j\le3,\label{eq:logmoments}\\
\norm{\mathcal B_n\bigl(g_n^j\bigr)}_\Frob&\le Cr_n,
\quad
\norm{\mathcal B_n\bigl(g_n^j\bigr)}_\op\le C(\log n)^j,\quad 0\le j\le3.\label{eq:logbounds}
\end{align}
The Frobenius bound at $j=0$ is read as a bound for $P_n$. More explicitly, for fixed integers $j\ge0$, $k\ge1$, and $n\ge2$,
\[
 \norm{\mathcal B_n\bigl(m_0^{-k}g_n^j\bigr)}_\Frob
       \le C_{j,k}L_{k,n}(\log n)^j.
\]
The corresponding parameter derivatives are understood as follows. With $n$ and $\theta_0$ held fixed, put
\[
 D_{\nu,n}=\partial_\nu+2L_n\partial_\eta,\qquad
 m_a(\omega)=a+4\pi^2\norm\omega^2,\qquad
 g_{n,a}=2L_n-\log m_a.
\]
\begin{samepage}
Here $D_{\nu,n}$ is differentiation in $v$ with $u$ held fixed in the coordinates \eqref{eq:shear}. Define the matrix whitened at the true covariance by
\[
 B_{j,k,n}(\theta)
 =C_{n,0}^{-1/2}
       \bigl(D_{\nu,n}^{\,j}\partial_a^k S_n\bigr)(\eta,\nu,a)
       C_{n,0}^{-1/2}.
\]
\end{samepage}
Uniformly over $\theta\in\mathcal N_n$,
\[
 \norm{B_{j,k,n}(\theta)}_\Frob
       \le C_{j,k}L_{k,n}(\log n)^j,\qquad
 \norm{B_{j,k,n}(\theta)}_\op\le C_{j,k}(\log n)^j,
 \qquad k\ge1.
\]
Thus $k$ range derivatives determine the factor $L_{k,n}$, and $j$ smoothness derivatives cost at most $(\log n)^j$. Additional $\eta$ derivatives leave these signal derivatives unchanged. For pure centered derivatives, the sharper bounds are
\[
 \sup_{\theta\in\mathcal N_n}
       \norm{B_{j,0,n}(\theta)}_\Frob\le C_jr_n,\qquad
 \sup_{\theta\in\mathcal N_n}
       \norm{B_{j,0,n}(\theta)}_\op\le C_j(\log n)^j.
\]
Any fixed finite collection of derivative orders is allowed; the constants may depend on those orders, but not on $n$, the sites, or the parameters within the stated compact bounds.
\end{lemma}

To obtain the matching Fisher information lower bound $J_n \succeq cI$, we need to recover Fourier coefficients at frequencies of order $R_n$. The next lemma shows that best linear prediction at the truth gives the required accuracy: for $\|\lambda\|\le\epsilon R_n$, the relative mean squared prediction error is bounded by $C\epsilon^{2\nu_0+d}$ for all sufficiently large $n$. More generally, the lemma bounds the $L^2$ error in the best linear prediction of $Z_g$ for $g \in L^2(Q)$. The proof is via a Sobolev sampling inequality that connects the error to the locations.

\begin{lemma}[Prediction of Fourier coefficients and Sobolev sampling]
\label{lem:sampling}
Under \Cref{ass:param,ass:sharp,ass:noise}, eventually for deterministic designs, and
eventually almost surely for iid designs, the following assertions hold on the fixed working cube $Q$ chosen in Appendix~\ref{sec:proofnotation}. For $Z_g=\int_Qg(x)F(x)\,dx$, let $\widehat Z_g^{\,0}=\Cov_0(Z_g,Y_n)C_{n,0}^{-1}Y_n$. Then
\begin{equation}\label{eq:predictionrisk}
 \E_0|Z_g-\widehat Z_g^{\,0}|^2\le Cn^{-1}\norm g_{L^2}^2,
 \qquad g\in L^2(Q).
\end{equation}
In particular, the predictions in \eqref{eq:fourierpredictor} obey
\begin{equation}\label{eq:fourierpredictionrisk}
 \frac{\E_0|X_\lambda-\widehat X_{\lambda,n}|^2}
      {v_0(\lambda)}
 \le C\frac{(1+\norm\lambda^2)^{\beta_0}}n.
\end{equation}
Write $\widehat U_{\lambda,n} =\widehat X_{\lambda,n}/\sqrt{v_0(\lambda)}$. For every finite $\Lambda\subset\mathscr L$ and real weights,
\begin{equation}\label{eq:predictedenergyvariance}
 \Var_0\left\{\sum_{\lambda\in\Lambda}w_\lambda
       (|\widehat U_{\lambda,n}|^2-\E_0|\widehat U_{\lambda,n}|^2)
        \right\}
 \le C\sum_{\lambda\in\Lambda}w_\lambda^2.
\end{equation}

The scalar sampling inequality underlying these assertions is
\begin{equation}\label{eq:samplingfunctiontarget}
 \norm u_{L^2(Q)}^2
 \le\frac Cn\left\{\norm u_{H^s(\Omega)}^2
                   +\sum_{i=1}^n|u(x_{n,i})|^2\right\},
 \qquad s=\nu_-+d/2>d/2,
\end{equation}
for $u\in H^s(\Omega)$, where $\Omega$ is a fixed enclosing cube. The same design event and constants work throughout the compact parameter rectangle.

For every finite $\Lambda\subset\mathscr L$, the Fourier responses of the predictions also obey
\begin{equation}\label{eq:predictionenvelope}
 \sum_{\lambda\in\Lambda}
       f_0(\omega)|p_{\lambda,n}(\omega)|^2
       \le Cw_n(\omega),\qquad \omega\in\R^d.
\end{equation}
Here $p_{\lambda,n}$ and $w_n$ are defined in \eqref{eq:fourierweights} and \eqref{eq:frequencyweightdefinition}. All predictors and normalizations are fixed at the truth. Constants are independent of the finite frequency set and its cardinality.
\end{lemma}

The next lemma proves that $J_n \succeq cI$. Together with the spectral upper bounds,
this establishes $J_n\asymp I$. We also identify the asymptotic block structure of $J_n$.

\begin{lemma}[Gaussian curvature and its block structure]
\label{lem:information}
Under \Cref{ass:param,ass:sharp}, the normalized information matrices $J_n$ satisfy $cI\preceq J_n\preceq CI$ for all sufficiently large $n$. More precisely, for a normalized signal direction $h=(h_u,h_v)$ when $d\le3$ and $h=(h_u,h_v,h_a)$ when $d\ge4$, put
\begin{equation}\label{eq:signaltangent}
 g_h(\omega)=\frac{h_u+h_vg_n(\omega)}{r_n}
             -\frac{\beta_0h_a}{b_nm_0(\omega)},\qquad
 g_n=2\log R_n-\log m_0,
\end{equation}
omitting the range term when $d\le3$. There is a fixed $c>0$ such that, for all sufficiently large $n$ and all such $h$,
\begin{equation}\label{eq:signalcoercivity}
 \norm{\mathcal B_n(g_h)}_\Frob^2\ge c\norm h^2.
\end{equation}
For iid designs these assertions hold eventually almost surely, conditionally on the design. Write $J_{n,\mathrm{sig}}$ for the signal principal block, we have
\[
 J_n=\begin{pmatrix}J_{n,\mathrm{sig}}&0\\0&(2t_0^2)^{-1}\end{pmatrix}+o(1).
\]
Thus the normalized nugget cross-information tends to zero, and the nugget entry tends to $(2t_0^2)^{-1}$. At $d=4$ there is the additional decomposition, in signal coordinate order $(u,v,a)$,
\[
 J_{n,\mathrm{sig}}
 =\begin{pmatrix}J_{n,(u,v)}&0\\0&J_{n,aa}\end{pmatrix}+o(1).
\]
Here $J_{n,(u,v)}$ is the two-by-two principal block and $J_{n,aa}$ is the normalized range entry. In particular, range cross-information with both centered amplitude and smoothness tends to zero. These matrix remainders vanish uniformly in operator norm. The information upper bounds hold without a sampling assumption.
\end{lemma}

\subsection{Local likelihood and global localization}\label{sec:lemma-likelihood}

The following lemma collects the deterministic covariance estimates needed for the likelihood arguments. Frobenius bounds control information and quadratic-form fluctuations; operator bounds control relative covariance changes and matrix inverses. On the preliminary neighborhood, the covariance separation bounds support localization
of the maximizer. On bounded $h$-sets, the derivative bounds justify the likelihood Taylor expansion. Part (ix) supplies the Gaussian signal bounds used to control the reciprocal of the absolute nuisance curvature at the fit. This control is needed for Wald inference using the identifiable block of the full inverse negative Hessian.

\begin{lemma}[Whitened covariance bounds for the Mat\'ern family]
\label{lem:derivatives}
Assume \Cref{ass:param}, with arbitrary sites in a fixed bounded set. Define $\mathcal N_n$ by \eqref{eq:localset}, with $q_n=n^{-\delta}$ for a fixed $0<\delta<d/(4\nu_++2d)$. The bounds in parts (i)--(vii) hold for these arbitrary sites. Parts \ref{it:covariance-separation}--\ref{it:covariance-signal} additionally assume \Cref{ass:sharp}; for iid sites their assertions hold conditionally on the good designs of Lemma~\ref{lem:randomgood}. All bounds are deterministic once the locations are fixed. They also hold on $\mathcal N_n^\square$ and its normalized image, with the conventions in Appendix~\ref{sec:proofnotation}. In particular, integration along segments takes place in this enclosing set and does not require convexity of the optimization set $\Theta$. Write the covariance whitened at the truth as
\begin{equation*}
 M(\theta)=C_{n,0}^{-1/2}C_n(\theta)C_{n,0}^{-1/2},
 \qquad M(\theta_0)=I_n.
\end{equation*}
We suppress $n$ in $M$. The whitening matrices remain fixed under every differentiation. In the centered coordinates $(u,v,a,t)$, subscripts such as $M_v,M_{va}$ mean derivatives in those coordinates. In normalized coordinates, we use the abbreviations
\[
 \begin{aligned}
 M(h)&=C_{n,0}^{-1/2}C_n(\theta_0+E_nh)C_{n,0}^{-1/2},
       &&d\ge4,\\
 M(h,a)&=C_{n,0}^{-1/2}C_n(\vartheta_0+E_nh,a)C_{n,0}^{-1/2},
       &&d\le3.
 \end{aligned}
\]
Thus $M(0)=I_n$ or $M(0,a_0)=I_n$, respectively. The normalized covariance tangents at the truth and their linear combination are
\[
 B_j=\begin{cases}
       \partial_{h_j}M(0),&d\ge4,\\
       \partial_{h_j}M(0,a_0),&d\le3,
      \end{cases}
 \qquad B(h)=\sum_{j=1}^{k_d}h_jB_j.
\]
By \eqref{eq:info}, $[J_n]_{ij}=\tfrac12\tr(B_iB_j)$. All bounds below are uniform over the preliminary neighborhood $\mathcal N_n$ or the displayed subsets of normalized coordinates: constants are independent of $n$ and the allowed parameters for all sufficiently large $n$, and each displayed rate remains valid after taking the supremum over its stated set. Write $R=M^{-1}$. The bounds are grouped by their index set and derivative order.
\begin{enumerate}[label=(\roman*),ref=(\roman*)]
\item\label{it:covariance-basic} \textbf{Covariance comparison and derivative bounds.}
For $\theta\in\mathcal N_n$,
\[
 cI_n\preceq M(\theta)\preceq CI_n,\quad \|R\|_\op\le C,\qquad
 \norm{\partial_u^j\partial_v^kM(\theta)}_\Frob\le Cr_n,
 \quad 1\le j+k\le3.
\]
The first-derivative operator bounds are
\[
 \norm{M_u}_\op\le C,\qquad
 \norm{M_v}_\op\le C\log n,\qquad
 \norm{M_a}_\op\le C.
\]
For range derivatives,
\[
 \norm{M_a}_\Frob\le
 \begin{cases}C,&d\le3,\\Cb_n,&d\ge4,\end{cases}
\]
and, for $d\ge4$,
\begin{equation*}
 \norm{M_{ua}}_\Frob\le Cb_n,\qquad
 \norm{M_{va}}_\Frob\le Cb_n\log n,\qquad
 \norm{M_{aa}}_\Frob\le CL_{2,n}.
\end{equation*}
Finally, $M_t=C_{n,0}^{-1}$, so $\norm{M_t}_\Frob\le C\sqrt n$ and $\norm{M_t}_\op\le C$; all second derivatives involving $t$ vanish.

\item\label{it:covariance-high-increments} \textbf{Full covariance increments: $d\ge4$.}
Set
\[
 \mathcal H_n=\{h\in\R^4:\theta_0+E_nh\in\mathcal N_n\},
 \qquad
 h=(r_nu,r_nv,b_n(a-a_0),\sqrt n(t-t_0))^\T.
\]
For every $h,h'\in\mathcal H_n$,
\begin{equation}\label{eq:fullcovlip}
 \norm{M(h)-M(h')}_\Frob\le C\norm{h-h'}.
\end{equation}

\item\label{it:covariance-high-remainder} \textbf{Taylor remainder: $d\ge4$.}
Define the remainder by the exact identity
\begin{equation*}
 M(h)=I_n+B(h)+\mathcal R_n(h),
 \qquad B(h)=DM(0)[h].
\end{equation*}
Uniformly for $h\in\mathcal H_n$,
\[
 \norm{\mathcal R_n(h)}_\Frob\le\epsilon_n\norm h,
 \qquad
 \epsilon_n=C\left(n^{-\delta}\log n+\frac{L_{2,n}}{b_n}\right)
       \longrightarrow0.
\]
Here $\delta$ is the exponent defining $\mathcal N_n$.

\item\label{it:covariance-low-increments} \textbf{Range perturbations and local increments: $d\le3$.} The change in the whitened covariance when only range varies is
\[
 D_n(a)=M(0,a)-I_n.
\]
The normalized amplitude and centered smoothness tangents satisfy, uniformly over $a\in[a_-,a_+]$,
\begin{equation}\label{eq:rangetangentbounds}
 \norm{\partial_{h_1}M(0,a)-B_1}_\Frob\le C/r_n,\qquad
 \norm{\partial_{h_2}M(0,a)-B_2}_\Frob\le C\log n/r_n.
\end{equation}
The nugget tangent $\partial_{h_3}M(0,a)=B_3$ is independent of $a$.

In normalized coordinates, let
\[
 \mathcal H_n=\{(h,a):(\vartheta_0+E_nh,a)\in\mathcal N_n\},
 \qquad h=(r_nu,r_nv,\sqrt n(t-t_0))^\T.
\]
Then the exact decomposition
\begin{equation}\label{eq:Mlowexpansion}
 M(h,a)=I_n+D_n(a)+B(h)+\mathcal R_n(h,a)
\end{equation}
satisfies, with
$\epsilon_n=C\{n^{-\delta}\log n+(\log n)/r_n\}\to0$,
\begin{equation}\label{eq:Mlowglobalbounds}
 \sup_a\norm{D_n(a)}_\Frob\le C,\qquad
 \norm{\mathcal R_n(h,a)}_\Frob\le\epsilon_n\norm h,
 \quad (h,a)\in\mathcal H_n.
\end{equation}
For any two points of $\mathcal H_n$,
\begin{equation}\label{eq:Mlowlip}
 \norm{M(h,a)-M(h',a')}_\Frob
       \le C\{\norm{h-h'}+|a-a'|\}.
\end{equation}

\item\label{it:covariance-local} \textbf{Normalized derivatives on bounded local sets.}
Fix $H<\infty$ and put
\[
 \Gamma_H=\begin{cases}
 \{h:\|h\|\le H\},&d\ge4,\\
 \{h:\|h\|\le H\}\times[a_-,a_+],&d\le3.
 \end{cases}
\]
Use $M(h)$ in high dimension and $M(h,a)$ in low dimension. Put $M_i=\partial_{h_i}M$,
$M_{ij}=\partial_{h_i h_j}M$, and $M_{ijk}=\partial_{h_i h_j h_k}M$. Uniformly on these sets,
\begin{equation}\label{eq:studentderiv}
 \begin{aligned}
 \norm{M-I_n}_\Frob+\max_i\norm{M_i}_\Frob&\le C_H,\\
 \max_i\{\norm{M_i-B_i}_\Frob+\norm{M_i}_\op\}
 &\le C_H\begin{cases}
 (\log n)/r_n,&d\le3,\\
 (\log n)/r_n+b_n^{-1},&d\ge4,
 \end{cases}\\
 \max_{i,j}\norm{M_{ij}}_\Frob
       +\max_{i,j,k}\norm{M_{ijk}}_\Frob
 &\le C_H\begin{cases}
 r_n^{-1},&d\le3,\\
 (\log n)/r_n+L_{2,n}/b_n^2,&d\ge4.
 \end{cases}
 \end{aligned}
\end{equation}
For $d\le3$, the nuisance derivative bounds are collected in part~\ref{it:covariance-nuisance}. The rates in the last two lines of \eqref{eq:studentderiv} tend to zero by
\eqref{eq:rates}, \eqref{eq:arate}, and \eqref{eq:rangeupper}.

\item\label{it:covariance-hessian}
\textbf{Matrices in the identifiable Hessian.}
On each fixed local set $\Gamma_H$, write $u=h$ or $u=(h,a)$, let $k$ be the total index dimension, and put $R=M^{-1}$. For any two coordinates $r,s$ of $u$, define
\[
 A_{rs}=RM_{rs}R-RM_rRM_sR-RM_sRM_rR.
\]
For identifiable indices $i,j$,
\[
 \sup_{\Gamma_H}\max_{i,j}\left\{
   \|A_{ij}\|_\Frob+
   \max_{1\le\ell\le k}\|\partial_{u_\ell}A_{ij}\|_\Frob\right\}
 \le C_H\begin{cases}
 (\log n)/r_n,&d\le3,\\
 (\log n)/r_n+b_n^{-1},&d\ge4.
 \end{cases}
\]
The same rate, up to a fixed factor, bounds the Frobenius Lipschitz constant of each $A_{ij}$ on $\Gamma_H$.

\item\label{it:covariance-nuisance}
\textbf{Nuisance derivatives and mixed Hessian matrices: $d\le3$.}
All bounds in this part hold on $\Gamma_H$. Write $\|\cdot\|_1$ for the trace norm and set
\[
 T_{d,n}=\begin{cases}1,&d=1,\\\log n,&d=2,\\R_n,&d=3.\end{cases}
\]
The pure nuisance derivatives satisfy
\[
 \|M_a\|_\Frob\le C_H,\qquad
 \|M_a\|_1\le C_HT_{d,n},\qquad
 \max_{2\le j\le4}\|\partial_a^jM\|_1\le C_H.
\]
The mixed covariance derivatives needed through fourth order obey
\begin{equation}\label{eq:localnuisancederivatives}
 \begin{aligned}
 \max_{i,\,1\le j\le3}\|\partial_a^jM_i\|_\Frob
       &\le C_H\frac{\log n}{r_n},\\
 \max_{i,l,\,1\le j\le2}\|\partial_a^jM_{il}\|_\Frob
       &\le C_H\frac{(\log n)^2}{r_n^2}.
 \end{aligned}
\end{equation}
Derivatives in these displays involving a normalized nugget coordinate are zero. The range information and cross-traces satisfy
\[
 \begin{aligned}
 c&\le\tfrac12\|R^{1/2}M_aR^{1/2}\|_\Frob^2\le C_H,\\
 \max_i|\tr(RM_iRM_a)|
       &\le C_H\frac{\log n}{r_n}T_{d,n}=o(1),\\
 \max_i|\partial_a\tr(RM_iRM_a)|
       &\le C_H\frac{\log n}{r_n}=o(1).
 \end{aligned}
\]
For the mixed Hessian matrices defined in part
\ref{it:covariance-hessian},
\[
 \sup_{\Gamma_H}\max_i\left\{
 \|A_{ia}\|_\Frob+\|\partial_aA_{ia}\|_\Frob
 +\max_\ell\|\partial_{u_\ell}A_{ia}\|_\Frob
 +\max_\ell\|\partial_{u_\ell}\partial_aA_{ia}\|_\Frob
 \right\}\le C_H\frac{\log n}{r_n}.
\]
Thus $A_{ia}$ and $\partial_aA_{ia}$ have vanishing Frobenius norms and Frobenius Lipschitz constants on $\Gamma_H$.

\item\label{it:covariance-separation}
\textbf{Identifiable tangents and covariance separation.}
Under Assumption~\ref{ass:sharp},
\[
 c\|h\|^2\le\|B(h)\|_\Frob^2=2h^\T J_nh\le C\|h\|^2.
\]
On the entire preliminary set $\mathcal H_n$, for all sufficiently large $n$,
\[
 \begin{cases}
 c\|h\|^2\le\|M(h)-I_n\|_\Frob^2\le C\|h\|^2,&d\ge4,\\
 c\|h\|^2-C\le\|M(h,a)-I_n\|_\Frob^2
                         \le C(1+\|h\|^2),&d\le3.
 \end{cases}
\]

\item\label{it:covariance-signal}
\textbf{The Gaussian signal in nuisance derivatives: $d\le3$.}
Assume Assumption~\ref{ass:sharp}. At $h=0$, uniformly over
$a\in I=[a_-,a_+]$, put
\[
 P_0=C_{n,0}^{-1/2}S_{n,0}C_{n,0}^{-1/2},\qquad
 B_a=-RM_aR.
\]
Then
\[
 RP_0R\preceq CP_0,\qquad 0\preceq B_a\preceq CP_0,\qquad
 \lambda_2(P_0^{1/2}B_aP_0^{1/2})\ge c>0.
\]
Here $\lambda_2$ is the second largest eigenvalue. Let $a=\alpha(q)$ be a fixed $C^4$ change of range coordinate on a compact interval, with $|\alpha'|$ bounded away from zero and derivatives through order four bounded, and set
\[
 A_q=RM_{qq}R-2RM_qRM_qR.
\]
For every $\eta>0$, there is $c_\eta>0$ such that
\[
 \|A_q\|_\Frob\ge\eta
 \quad\Longrightarrow\quad
 \|P_0^{1/2}A_qP_0^{1/2}\|_\Frob\ge c_\eta.
\]
The constants may depend on the fixed coordinate change.

\end{enumerate}
\end{lemma}

Preliminary localization places every global maximizer in a shrinking parameter neighborhood, but that neighborhood can still be large in $h$-coordinates, so we need further work to prove that $\hat{h}_n$ is $O_P(1)$. The next lemma uses the covariance separation and the likelihood bounds to further exclude large values of $\|h\|$ and prove that $\hat{h}_n = O_P(1)$.

\begin{lemma}[Uniform tails and moments of every global fit]
\label{lem:localization}\label{lem:tails}
Under \Cref{ass:param,ass:sharp,ass:noise}, every global maximizer satisfies \eqref{eq:tails} and \eqref{eq:allmoments}, uniformly over $\mathcal K\times\mathcal P$. The same statements hold for the restricted known-smoothness fit. These results include maximization over the feasible range intervals $I(\vartheta)$ in low dimension.
\end{lemma}

In $d \leq 3$, the fitted range need not approach the true range, so the local Taylor expansion must work uniformly over its feasible values. The next lemma shows that changing the range does not change the leading identifiable score.

\begin{lemma}[Nuisance-score stability and baseline regularity]\label{lem:nuisance}
For $d\le3$, assume \Cref{ass:param,ass:sharp,ass:noise} and write
\[
 \mathcal L_n(h,a)=\ell_n(\vartheta_0+E_nh,a),\qquad
 s_n(a)=\nabla_h\mathcal L_n(0,a),\qquad \Delta_n=s_n(a_0).
\]
For every fixed $1\le p<\infty$,
\begin{equation}\label{eq:nuisancescorecomparison}
 \sup_{a\in[a_-,a_+]}\|s_n(a)-\Delta_n\|=o_{L^p}(1),
\end{equation}
uniformly over $\mathcal K\times\mathcal P$. Thus varying the unnormalized range nuisance over the enclosing interval $I_*=[a_-,a_+]$ does not change the leading normalized identifiable score, whether or not a particular value is feasible. For the baseline criterion $G_n(a)=\mathcal L_n(0,a)$, we also have
\begin{equation}\label{eq:nuisance-baseline-moments}
 \max_{1\le j\le3}
 \left\|\sup_{a\in I_*}|G_n^{(j)}(a)|\right\|_{L^p}\le C_p.
\end{equation}
The same assertions hold after restricting smoothness to its true value. This is the Mat\'ern application of the general score-comparison Lemma~\ref{lem:score-comparison}.
\end{lemma}

For Gaussian error laws, Markov's inequality applied to the square root of a likelihood ratio converts an affinity upper bound into an upper bound on the probability that an alternative yields larger likelihood than the truth. The following lemma allows us to extend this probability bound to sub-Gaussian error laws.

\begin{lemma}[A Gaussian-affinity bound for the Gaussian criterion]
\label{lem:affinitytails}
Let $Y=F+\varepsilon$, where $F\sim N(0,S_0)$ is independent of the errors in Assumption~\ref{ass:noise}, and $C_0=S_0+t_0I$. For every $C\succeq t_-I$, write
\[
 L_C=-\tfrac12\log\det(CC_0^{-1})
       -\tfrac12Y^\top(C^{-1}-C_0^{-1})Y,
 \qquad A_C=-\log\Aff\{N(0,C),N(0,C_0)\}.
\]
There are uniform constants $c,C_*>0$ such that\begin{equation}\label{eq:likelihoodaffinity}
 P(L_C\ge-z)\le C_*\exp(z/2-cA_C),\qquad z\ge0.
\end{equation}
This bounds the Gaussian criterion under the actual response law using a deterministic Gaussian covariance affinity.
\end{lemma}

The expected curvature and the score covariance can differ under non-Gaussian errors. The fourth-cumulant contribution from Lemma~\ref{lem:gaussian} depends on the diagonals of the quadratic-form matrices. Here those contributions vanish for normalized signal coordinates and remain for the nugget. The representation by centered
squared errors explains the nugget variance and the form of the covariance correction used for inference.

\begin{lemma}[Score covariance and nugget representation]
\label{lem:noise-diagonals}
Assume \Cref{ass:param,ass:sharp,ass:noise}. Let $D_j=\partial_{h_j}C_n(\theta_0)$ and
$A_j=C_{n,0}^{-1}D_jC_{n,0}^{-1}$.
For every normalized signal coordinate,
\begin{equation}\label{eq:diagonal-negligible}
 \sum_i(A_j)_{ii}^2=o(1).
\end{equation}
For the nugget, $A_t=C_{n,0}^{-2}/\sqrt n$ and
$\sum_i(A_t)_{ii}^2=t_0^{-4}+o(1)$.
Consequently, uniformly,
\begin{equation}\label{eq:score-covariance-correction}
 K_n=J_n+\frac{\kappa_4}{4t_0^4}e_te_t^\T+o(1),
 \qquad V_n=J_n^{-1}+\kappa_4e_te_t^\T+o(1).
\end{equation}
The covariance $V_n$ has the signal--nugget block structure in \eqref{eq:covariance-blocks}, and its signal block has the additional range separation in \eqref{eq:critical-covariance-block} when $d=4$.
Moreover, for every fixed $1\le p<\infty$,
\begin{equation}\label{eq:nugget-score-squares}
 \Delta_{n,t}=\frac1{2t_0^2\sqrt n}
             \sum_i(\varepsilon_{n,i}^2-t_0)+o_{L^p}(1).
\end{equation}
In particular, with $v_\varepsilon=t_0^2(\E_P\xi^4-1)$, $2t_0^2\Delta_{n,t}/\sqrt{v_\varepsilon}\Rightarrow N(0,1)$ uniformly. The deterministic assertions are uniform in the compact covariance family and in the allowed designs; the moment and distributional assertions are uniform also over $\mathcal P$.
\end{lemma}

\subsection{Feasible inference and parameter transformations}\label{sec:lemma-inference}

The covariance correction requires the fourth cumulant of the measurement errors, which are not observed directly. The fitted prediction residuals provide observable substitutes, from which we construct a root-$n$ consistent estimator of the fourth cumulant. 

\begin{lemma}[Residual estimation of the error fourth cumulant]
\label{lem:kurtosis}
Assume \Cref{ass:param,ass:sharp,ass:noise}, and let
$\widehat\theta_n$ be a global MLE. Define
\[
 \widehat A_n=\widehat t_n C_n(\widehat\theta_n)^{-1},\qquad
 \widehat r_n=\widehat A_nY_n,\qquad
 \widehat v_{n,i}=\widehat t_n(\widehat A_n)_{ii}.
\]
Thus $\widehat r_n$ is the observed response minus the fitted best linear predictor of the field, and $\widehat v_{n,i}$ is its model-based marginal variance. Set
\begin{equation}\label{eq:residual-cumulant}
 \widetilde\kappa_{4,n}
 =\frac{\displaystyle\sum_{i=1}^n
       \{\widehat r_{n,i}^{\,4}
          -6\widehat v_{n,i}\widehat r_{n,i}^{\,2}
          +3\widehat v_{n,i}^{\,2}\}}
       {\displaystyle\sum_{i,j=1}^n(\widehat A_n)_{ij}^{\,4}}.
\end{equation}
Then
\begin{equation}\label{eq:kurtosis-rate}
 \widetilde\kappa_{4,n}-\kappa_4(P,t_0)
       =O_{P,\mathcal K\times\mathcal P}(n^{-1/2}).
\end{equation}
For iid sites, this bound holds conditionally on the good designs of Lemma~\ref{lem:randomgood}, with common constants, and unconditionally. It also holds when smoothness is fixed at its true value. In low dimension no consistency of the fitted range is required.

If the errors are Gaussian, there is the sharper conclusion
\begin{equation}\label{eq:cumulant-gaussian-limit}
 \sqrt n\,\widetilde\kappa_{4,n}
       \Rightarrow N(0,24t_0^4),
\end{equation}
uniformly over the compact covariance family and allowed designs, and conditionally on good random designs and unconditionally as above. This is the limiting variance of the empirical fourth cumulant when the Gaussian errors themselves are observed. More precisely, let $\widetilde\kappa_{4,n}^{\,0}$ denote the estimator with all covariance
parameters held at the truth, put $A_0=t_0C_0^{-1}$, and write $d_n=\sum_{i,j}(A_0)_{ij}^4$. At Gaussian errors,
\begin{equation}\label{eq:cumulant-gaussian-oracle-variance}
 \Var_0(\widetilde\kappa_{4,n}^{\,0})
 =\frac{24t_0^4}{d_n}=\frac{24t_0^4}{n}\{1+o(1)\}.
\end{equation}
\end{lemma}

The local likelihood expansion needs its observed Hessian to be close to the deterministic information, uniformly on bounded $h$-sets. Once the fitted $h$ is tight, this same comparison permits evaluation of the curvature at the fit for Wald inference. In low dimension, control of the mixed range block also justifies using the identifiable block of the full inverse. Adding the estimated cumulant correction then
gives the sampling covariance needed for Wald inference.

\begin{lemma}[Local information and observed Hessian]\label{lem:studentization}
Assume \Cref{ass:param,ass:sharp,ass:noise}.
\begin{enumerate}[label=(\roman*),ref=(\roman*)]
\item\label{it:studentization-local}
\textbf{Local information and Hessian.}
For every fixed $H$, the expected
information satisfies, in high dimension,
\begin{equation}\label{eq:localIfisher}
 \sup_{\norm h\le H}
 \norm{E_n^\T\calI_n(\theta_0+E_nh)E_n-J_n}_{\op}\longrightarrow0.
\end{equation}
In low dimension the same assertion holds for the identifiable block, uniformly over $a\in[a_-,a_+]$ at $(\vartheta_0+E_nh,a)$. For the observed identifiable Hessian the analogous assertion is
\begin{equation}\label{eq:localIobserved}
 \sup_{\norm h\le H}
 \norm{-E_n^\T\nabla_{\mathrm{id}}^2\ell_n(\theta_0+E_nh)E_n-J_n}_{\op}
 \longrightarrow_{\PP}0,
\end{equation}
again with a uniform nuisance supremum in $d\le3$. The assertions are uniform over $\mathcal K\times\mathcal P$. The same results hold after restricting smoothness to its true value.

\item\label{it:studentization-mixed}
\textbf{Range information and mixed derivatives.}
For $d\le3$, write the full expected information in blocks for $(\vartheta,a)$ and put $\mathcal L_n(h,a)=\ell_n(\vartheta_0+E_nh,a)$. On $\Gamma_H=\{\norm h\le H\}\times[a_-,a_+]$, evaluated at $\theta=(\vartheta_0+E_nh,a)$,
\[
 c\le\calI_{aa,n}(\theta)\le C_H,\qquad
 \sup_{\Gamma_H}\norm{E_n^\T\calI_{\vartheta a,n}(\theta)}\longrightarrow0.
\]
The mixed observed derivatives satisfy
\[
 \sup_{\Gamma_H}
 \left\{\norm{\nabla_h\partial_a\mathcal L_n(h,a)}
             +\norm{\nabla_h\partial_a^2\mathcal L_n(h,a)}\right\}
 \longrightarrow_{\PP}0.
\]
These assertions are also uniform on compact sets of truths and hold after restricting smoothness.

\item\label{it:studentization-inverse}
\textbf{Substitution at the fit and the full inverse.}
At a tight global fit let $\widehat{\mathcal A}_n$ denote full expected curvature or the full negative Hessian in spectral coordinates. In low dimension let $\widehat H_n$ be either its identifiable block or its range Schur complement; in high dimension use the full matrix. These choices satisfy
\begin{equation*}
 E_n^\top\widehat H_nE_n=J_n+o_P(1).
\end{equation*}
For the observed Schur complement this uses Lemma~\ref{lem:range-curvature}. In particular the relevant inverses exist and their identifiable covariance blocks are positive definite with probability tending to one. Define
\begin{equation}\label{eq:feasible-sandwich}
 \widehat\Sigma_n^\vartheta=\widehat H_n^{-1}
                    +\frac{\widetilde\kappa_{4,n}}n e_te_t^\top.
\end{equation}
Then
\begin{equation}\label{eq:covariance-substitution}
 E_n^{-1}\widehat\Sigma_n^\vartheta E_n^{-\top}=V_n+o_P(1).
\end{equation}
The same assertions hold when the nuisance is any range coordinate $q$ permitted in Lemma~\ref{lem:range-curvature}, including endpoint fits. In particular this includes $a,\log a,l,\log l$. The proof of Corollary~\ref{cor:uncertainty} transfers the core covariance estimate to each admissible identifiable vector in the chosen full coordinates.
\end{enumerate}
\end{lemma}

Validity of the full inverse negative Hessian for Wald inference requires one additional check: a small mixed block has a negligible Schur correction only if the reciprocal of the nuisance curvature is controlled. The following lemma gives this
control at the joint maximizer, including an endpoint range fit.

\begin{lemma}[Observed range curvature at the joint MLE]
\label{lem:range-curvature}
Assume \Cref{ass:param,ass:sharp,ass:noise} and $d\le3$. For any measurable global maximizer $\widehat\theta_n$ and compact set of truths $\mathcal K$ as above, let $a=\alpha(q)$ be a fixed $C^4$ reparameterization of the enclosing range interval $I_*=[a_-,a_+]$, with $|\alpha'|$ bounded away from zero. All its derivatives through order four are bounded. This includes $q=a,\log a,l,\log l$. Hold $(\eta,\nu,t)$ fixed when
differentiating. Then
\[
 \lim_{\delta\downarrow0}\limsup_{n\to\infty}
 \sup_{\theta_0\in\mathcal K,\ P\in\mathcal P}
 P_{\theta_0,P}\!\left\{
       |\partial_q^2\ell_n(\widehat\theta_n)|\le\delta\right\}=0.
\]
Thus $|\partial_q^2\ell_n(\widehat\theta_n)|^{-1}=O_{\PP}(1)$, uniformly over $\mathcal K\times\mathcal P$, with the reciprocal set to infinity when the denominator is zero. The assertion allows the fitted range at either endpoint of its feasible interval $I(\widehat\vartheta_n)$ and does not require its consistency. It also holds for the known-smoothness MLE, with the identifiable coordinates restricted accordingly, and conditionally for almost every iid design.
\end{lemma}

The parameters converge at different rates, and the normalizer and coordinate Jacobians can mix these rates. Transferring the asymptotic normality of $\widehat h_n$ therefore requires control of the normalized Taylor remainders. The next lemma identifies the parameterizations for which the joint linear approximation holds and proves marginal linear approximations for every listed identifiable coordinate.

\begin{lemma}[Multirate parameter transformations]\label{lem:logtransform}
Under \Cref{ass:param,ass:sharp,ass:noise}, let $\zeta=g(\vartheta)$ be an identifiable vector in Table~\ref{tab:parameterizations}, other than the original-scale
variance--lengthscale vector when $4\le d\le8$. Then
\[
 E_n^{-1}Dg(\vartheta_0)^{-1}
      \{g(\widehat\vartheta_n)-g(\vartheta_0)\}
       -\widehat h_n=o_P(1).
\]
This remainder tends to zero in every fixed positive moment, and
\begin{equation}\label{eq:jacobianrelative}
 E_n^{-1}Dg(\vartheta_0)^{-1}
       Dg(\widehat\vartheta_n)E_n\longrightarrow_PI_{k_d}.
\end{equation}
The inverse of the matrix in \eqref{eq:jacobianrelative} also converges to $I_{k_d}$. The conclusions hold for the auxiliary log chart in \eqref{eq:logcharts} as well and are uniform over $\mathcal K\times\mathcal P$.

For every original-scale variance--lengthscale coordinate in $d\ge4$, the marginal expansion is valid:
\begin{equation}\label{eq:raw-marginal-remainder}
 \frac{g_j(\widehat\vartheta_n)-g_j(\vartheta_0)
       -Dg_j(\vartheta_0)E_n\widehat h_n}
      {\rho_{j,n}}=o_P(1),
\end{equation}
also in every fixed positive moment, uniformly.
\end{lemma}

\subsection{Random designs and supporting lemmas}\label{sec:lemma-supporting}

For random locations, the deterministic arguments can be applied conditionally once coverage and sampling hold with common constants. The next lemma supplies a single event on which both properties hold, uniformly over the parameter rectangle. Its exponentially small failure probability also controls exceptional designs when
the conditional tail and moment bounds are averaged over locations.

\begin{lemma}[Quantified random-design comparison]\label{lem:randomgood}
Under \Cref{ass:param} and the iid alternative in \Cref{ass:sharp}, for every fixed
compact parameter rectangle there are events $\mathcal G_n$ depending only on the locations and constants $c,C,\delta>0$ such that
\begin{equation}\label{eq:randomfailure}
 P_X(\mathcal G_n^c)\le Ce^{-cn^\delta}.
\end{equation}
On $\mathcal G_n$, a fixed polynomial-coverage bound and the Sobolev sampling inequality \eqref{eq:samplingfunctiontarget} hold with deterministic constants. The prediction bound \eqref{eq:predictionrisk} consequently holds uniformly with every
parameter in the rectangle taken as truth. The finite-$n$ information, preliminary localization, and local probability bounds proved from those two properties have the same constants for all such designs.
\end{lemma}

The remaining supporting results supply the analytic and probabilistic steps used above. The spectral covariance identity expresses covariances of field integrals and sampled field values as weighted Fourier inner products.

\begin{lemma}[Spectral covariance and Fourier responses]
\label{lem:spectral-covariance}
Under \Cref{ass:param}, let $g_1,g_2$ be compactly supported $L^2(\R^d;\mathbb C)$ functions, or finite linear combinations of such functions and point masses. With $Z_g$ and the Fourier convention of Appendix~\ref{sec:proofnotation},
\begin{equation}\label{eq:spectralcovariance}
 \Cov_\theta(Z_{g_1},Z_{g_2})
 =\int_{\R^d}f_\theta(\omega)\widehat g_1(\omega)
                         \overline{\widehat g_2(\omega)}\,d\omega.
\end{equation}
In particular, at the observation sites,
\begin{equation}\label{eq:sampledspectralcovariance}
 \Cov_\theta\bigl(F(x_{n,1:n})\bigr)
 =\int_{\R^d}f_\theta(\omega)
       e_\omega(x_{n,1:n})e_\omega(x_{n,1:n})^*\,d\omega.
\end{equation}
The normalized Fourier weights in \eqref{eq:fourierweights} satisfy
\begin{equation}\label{eq:fouriernormalizations}
 \begin{aligned}
 \int_{\R^d}f_0(\omega)|\ell_\lambda(\omega)|^2\,d\omega&=1,\\
 \E_0|\widehat U_{\lambda,n}|^2
 &=\int_{\R^d}f_0(\omega)|p_{\lambda,n}(\omega)|^2\,d\omega
       +t_0\sum_i|b_{\lambda,n,i}|^2.
 \end{aligned}
\end{equation}
\end{lemma}

The sampling argument rests on a deterministic fact: values at sites covering a cube control the $L^2$ norm of a smooth function there, up to an approximation error measured by its Sobolev norm. As the coverage becomes finer, this error decreases. The following local inequality is the approximation step used to prove the sampling
and prediction bounds in Lemma~\ref{lem:sampling}.

\begin{lemma}[Sobolev sampling from coverage on cubes]
\label{lem:cubesampling}
Let $Q_-,Q_+,\Omega$ be fixed cubes with $\overline{Q_-}\subset Q_+$ and $\overline{Q_+}\subset\Omega$, and fix $s>d/2$. For sites $x_1,\ldots,x_n\in\Omega$, put
\[
 h=\sup_{x\in Q_+}\min_{1\le i\le n}\norm{x-x_i}.
\]
There are constants $C,h_*>0$, depending only on $s,d$ and the cubes, such that, whenever $h\le h_*$ and $f\in H^s(\Omega)$,
\begin{equation}\label{eq:sobolevsampling}
 \norm f_{L^2(Q_-)}^2
 \le C\left\{h^{2s}\norm f_{H^s(\Omega)}^2
                 +h^d\sum_{i=1}^n|f(x_i)|^2\right\}.
\end{equation}
\end{lemma}

The next result gives general conditions for the Hessian approximation used in the likelihood Taylor expansion. Differentiating the Gaussian criterion separates its Hessian into deterministic trace terms and centered quadratic forms. Small Frobenius norms control the latter at each parameter value, and a Lipschitz bound makes that control uniform over the local parameter set. The Mat\'ern verification uses the bounds collected in Lemma~\ref{lem:derivatives}.

\begin{lemma}[Uniform Hessian from Lipschitz matrix bounds]
\label{lem:gaussian-hessian}
Let $\Gamma\subset\R^k$ be a fixed compact box containing the true index, where $k$ is the total dimension of $u$. Use $u=h$ without a nuisance and $u=(h,a)$ with a nuisance; in the latter case, $\Gamma$ is the local $h$-box times a compact interval $I$. Let $M(u)$ be the positive definite truth-whitened covariance, twice continuously differentiable on a neighborhood of $\Gamma$. Let $X$ be the whitened observation vector of Lemma~\ref{lem:gaussian} and define the covariance criterion, up to an additive constant, by
\[
 \mathcal L_n(u)=-\tfrac12\log\det M(u)-\tfrac12X^\T M(u)^{-1}X.
\]
Put $R=M^{-1}$ and $M_i=\partial_{h_i}M$. Let $B_i=M_i(0)$ without a nuisance and $B_i=M_i(0,a_0)$ with a nuisance, where $M(0)=I_n$ and $M(0,a_0)=I_n$, respectively. Define
\[
 A_{ij}=RM_{ij}R-RM_iRM_jR-RM_jRM_iR,\qquad
 J_{n,ij}=\tfrac12\tr(B_iB_j),\qquad
 \mathscr I_{n,ij}(u)=\tfrac12\tr(RM_iRM_j).
\]
Suppose $\norm R_\op+\norm{M-I}_\Frob+\max_i\norm{M_i}_\Frob\le C$ and $\norm{M_i-B_i}_\Frob+\norm{M_i}_\op=o(1)$, uniformly on the box. In addition, assume
\[
 \epsilon_n\coloneqq \sup_{u\in\Gamma}\max_{i,j}\|A_{ij}(u)\|_\Frob
       \longrightarrow0,\qquad
 \|A_{ij}(u)-A_{ij}(v)\|_\Frob\le L\|u-v\|,
 \quad u,v\in\Gamma,
\]
where $L$ is independent of $n$ and the indices. A sufficient condition for this Lipschitz bound is that the first partial derivatives of $A_{ij}$ are uniformly bounded in Frobenius norm. For every fixed $1\le p<\infty$,
\begin{equation}\label{eq:commonHessianconclusion}
 \sup_u\|\mathscr I_n(u)-J_n\|_\op=o(1),\qquad
 \sup_u\norm{\nabla_h^2\mathcal L_n(u)+J_n}_\op=o_{L^p}(1).
\end{equation}
All assertions are uniform over the noise class and over any family for which the displayed bounds and convergences are uniform.
\end{lemma}

Score stability has a similar general formulation. The difference between two baseline identifiable scores is the integral of a mixed identifiable--nuisance derivative. The matrix norm bounds below control its centered fluctuation, while the cross-trace bound
controls its mean. This yields a uniform comparison over the nuisance interval and is applied in Lemma~\ref{lem:nuisance}.

\begin{lemma}[Nuisance-score comparison from covariance derivatives]
\label{lem:score-comparison}
Fix $1\le p<\infty$ and a compact interval $I$, and let $h\in\R^k$, with $k$ fixed. Suppose $X_n\in\R^n$ satisfies $\E X_nX_n^\T=I_n$ and, for every real symmetric matrix $A$,
\[
 Q(A)=X_n^\T AX_n-\tr A,\qquad
 \|Q(A)\|_{L^p}\le C_p\|A\|_\Frob,
\]
with $C_p$ independent of $n$. Let $M_n(h,a)$ be a positive definite, twice continuously differentiable covariance matrix on a neighborhood of $\{0\}\times I$, with $M_n(0,a_0)=I_n$ for some $a_0\in I$. Define the Gaussian covariance criterion and its baseline score by
\[
 \mathcal L_n(h,a)=-\tfrac12\log\det M_n(h,a)
          -\tfrac12X_n^\T M_n(h,a)^{-1}X_n,\qquad
 s_n(a)=\nabla_h\mathcal L_n(0,a),\qquad \Delta_n=s_n(a_0).
\]
At $(0,a)$ abbreviate $M=M_n(0,a)$, $R=M^{-1}$, $M_i=\partial_{h_i}M_n(0,a)$,
$M_a=\partial_aM_n(0,a)$, and $M_{ia}=\partial_{h_i a}M_n(0,a)$. Suppose
\[
 \sup_{a\in I}\{\|R\|_\op+\|M-I_n\|_\Frob+\|M_a\|_\Frob\}\le C,
\]
and put
\[
 \epsilon_n=\sup_{a\in I}\max_{1\le i\le k}
       \{\|M_i\|_\op+\|M_{ia}\|_\Frob\},\qquad
 \rho_n=\sup_{a\in I}\max_{1\le i\le k}
       |\tr(RM_iRM_a)|.
\]
Then
\begin{equation}\label{eq:score-comparison-general}
 \left\|\sup_{a\in I}\|s_n(a)-\Delta_n\|\right\|_{L^p}
 \le C_{p,k,I}(\epsilon_n+\rho_n).
\end{equation}
In particular, if $\epsilon_n+\rho_n\to0$, the score difference is $o_{L^p}(1)$ uniformly over $I$. A sufficient condition for $\rho_n\to0$ is
\[
 \sup_{a\in I}\max_i\|M_i\|_\op\|M_a\|_1\longrightarrow0,
\]
where $\|\cdot\|_1$ is the trace norm. The conclusions are uniform over any family for which the displayed bounds are uniform. For Gaussian observations, and for the observations in Lemma~\ref{lem:gaussian}, the quadratic-form moment bound holds
for every fixed finite $p$. Thus the lemma applies to covariance models beyond the Mat\'ern family; Lemma~\ref{lem:nuisance} verifies its assumptions for the range nuisance here.
\end{lemma}

The following lemma is the general result used to prove Lemma~\ref{lem:range-curvature}.

\begin{lemma}[Nonvanishing curvature at a fitted scalar nuisance]
\label{lem:nuisance-curvature}
Let $I=[a_-,a_+]$ be a fixed nondegenerate compact interval and $p$ a fixed positive integer. Let $\mathcal L_n(h,a)$ be a random criterion, defined on $\mathcal H_n\times I$, where $\mathcal H_n\subset\R^p$ contains every fixed ball about zero eventually. At $h$, let the feasible nuisance interval be $I_n(h)=[\ell_n^a(h),u_n^a(h)]\subset I$, with deterministic endpoints and $\ell_n^a(h)<u_n^a(h)$. Assume, for each fixed $H$,
\[
 \rho_{n,H}\coloneqq \sup_{\|h\|\le H}
  \max\{|\ell_n^a(h)-\ell_n^a(0)|,
        |u_n^a(h)-u_n^a(0)|\}\longrightarrow0.
\]
The special case $I_n(h)=I$ is allowed. Assume $\mathcal L_n$ is twice continuously differentiable in $a$, with jointly continuous derivatives, and $G_n(a)\coloneqq \mathcal L_n(0,a)$ has $C^3$ paths. Let the measurable fitted indices satisfy
\[
 \widehat h_n=O_P(1),\qquad
 \widehat a_n\in\argmax_{a\in I_n(\widehat h_n)}
                                      \mathcal L_n(\widehat h_n,a).
\]
Fix a $C^2$ bijection $a=\alpha(q)$ from a compact interval $J$ onto $I$, with $\inf_J|\alpha'|>0$, and set
\[
 \widetilde{\mathcal L}_n(h,q)=\mathcal L_n(h,\alpha(q)),\quad
 \widetilde G_n(q)=G_n(\alpha(q)),\quad
 \mathcal J_n(h)=\alpha^{-1}(I_n(h)),\quad
 \widehat q_n=\alpha^{-1}(\widehat a_n).
\]
Suppose the following conditions hold.
\begin{enumerate}[label=(\roman*),ref=(\roman*)]
\item\label{it:nuisance-curvature-density}
\textbf{Bounded score densities.} For some $n_0$, every $G_n'(a)$
with $n\ge n_0$ and $a\in I$ has a density satisfying
\[
 \sup_{n\ge n_0}\sup_{a\in I,\,x\in\R}p_{G_n'(a)}(x)\le C.
\]
\item\label{it:nuisance-curvature-third}
\textbf{Tight derivative bounds.}
\[
 T_n\coloneqq \max_{1\le j\le3}\sup_{a\in I}|G_n^{(j)}(a)|=O_P(1).
\]
These bounds also control the modulus of continuity of $\widetilde G_n''$ when the feasible endpoints move.
\item\label{it:nuisance-curvature-endpoints}
\textbf{Endpoint anti-concentration.} For the deterministic endpoints of the baseline feasible interval $\mathcal J_n(0)$,
\[
 \lim_{\epsilon\downarrow0}\limsup_{n\to\infty}
 \max_{q_e\in\partial \mathcal J_n(0)}
 P\{|\widetilde G_n''(q_e)|\le\epsilon\}=0.
\]
Here
$\widetilde G_n''(q)=(\alpha'(q))^2G_n''(\alpha(q))
                   +\alpha''(q)G_n'(\alpha(q))$.
\item\label{it:nuisance-curvature-transfer}
\textbf{Transfer to bounded local sets.} For every fixed $H$,
\[
 d_{n,H}\coloneqq \sup_{\|h\|\le H,\ a\in I}\max_{j=1,2}
     |\partial_a^j\mathcal L_n(h,a)-G_n^{(j)}(a)|=o_P(1).
\]
\end{enumerate}
Then
\begin{equation}\label{eq:nuisance-curvature-general}
 \lim_{\delta\downarrow0}\limsup_{n\to\infty}
 P\{|\partial_q^2\widetilde{\mathcal L}_n(\widehat h_n,\widehat q_n)|
                                                    \le\delta\}=0.
\end{equation}
Equivalently,
$|\partial_q^2\widetilde{\mathcal L}_n(\widehat h_n,\widehat q_n)|^{-1}
 =O_P(1)$, with reciprocal infinity at zero. Taking $J=I$ and $\alpha(q)=q$ gives the result in the original nuisance coordinate. All conclusions are uniform over any family of experiments for which the density bound, tightness, and probability limits in the hypotheses are uniform.
\end{lemma}

After conditioning on the measurement errors, nuisance derivatives are quadratic polynomials in the Gaussian field coordinates. The following small-ball bound limits the probability that such a polynomial lies near a specified value. Its uniformity in the linear and constant terms is useful because these terms depend on the conditioned errors; the bound supplies an anti-concentration step in the nuisance-curvature argument.

\begin{lemma}[Small balls for Gaussian quadratic polynomials]
\label{lem:quadratic-smallball}
Let $G\sim N(0,I_N)$, $A=A^\top$, $b\in\R^N$, and $c\in\R$.
If $\tau=\|A\|_{\rm F}>0$, then
\[
 \sup_zP\{|G^\top AG+b^\top G+c-z|\le\epsilon\}
 \le C\sqrt{\epsilon/\tau}.
\]
The bound is independent of $N,b,c$. It may therefore be applied conditionally when the linear and constant coefficients are random and independent of $G$.
\end{lemma}

The localization argument measures separation of Gaussian covariance models through their affinity. The following formula expresses that separation using the eigenvalues of the relative covariance matrix.

\begin{lemma}[Gaussian affinity and covariance separation]
\label{lem:gaussian-affinity}
For probability laws $P,Q$, define their Hellinger affinity by
\[
 \Aff(P,Q)=\int\sqrt{\frac{dP}{d\mu}\frac{dQ}{d\mu}}\,d\mu,
 \qquad \mu=P+Q.
\]
Let $C_1,C_2$ be positive definite $n$-by-$n$ matrices, and let $\lambda_1,\ldots,\lambda_n$ be the eigenvalues of $C_1^{-1/2}C_2C_1^{-1/2}$. Then
\begin{equation}\label{eq:gaussian-affinity-formula}
 \begin{split}
 \Aff\{N(0,C_1),N(0,C_2)\}
 &=\frac{\det(C_1)^{1/4}\det(C_2)^{1/4}}
                  {\det\{(C_1+C_2)/2\}^{1/2}},\\
 -\log\Aff\{N(0,C_1),N(0,C_2)\}
 &=\frac12\sum_{k=1}^n
                    \log\cosh\!\left(\frac12\log\lambda_k\right).
 \end{split}
\end{equation}
There are absolute constants $0<c<C<\infty$ such that
\begin{equation}\label{eq:logcosh-comparison}
 c\min\{x^2,|x|\}\le\log\cosh x
                 \le C\min\{x^2,|x|\},\qquad x\in\R.
\end{equation}
Thus the negative log affinity is comparable to $\sum_k\min\{(\log\lambda_k)^2,|\log\lambda_k|\}$, with constants independent of the dimension and of the covariance matrices. If, in addition, $C_1,C_2\succeq t_-I$, then
\begin{equation}\label{eq:precision-affinity}
 -\log\Aff\{N(0,C_1),N(0,C_2)\}
 \ge\frac{t_-^2}{16}\|C_1^{-1}-C_2^{-1}\|_\Frob^2.
\end{equation}
\end{lemma}

The following result shows that applying the same randomized transformation to two probability laws cannot decrease affinity. We apply this result to obtain negative log affinity lower bounds for the Gaussian reference laws $G_{\theta,n}$ by passing to simpler statistics. Lemma~\ref{lem:affinitytails} converts these bounds into probability bounds under the actual error law.

\begin{lemma}[Affinity under statistics and added noise]
\label{lem:affinity-processing}
Let $P,Q$ be probability laws on the same space. Applying the same randomized transformation to both laws cannot decrease affinity. Precisely, let $\mathsf K(x,\cdot)$ be the conditional distribution of an output given input $x$, and let $P\mathsf K,Q\mathsf K$ denote the resulting output laws. Then
\begin{equation}\label{eq:affinity-processing}
 \Aff(P\mathsf K,Q\mathsf K)\ge\Aff(P,Q).
\end{equation}
In particular, if $T$ is any measurable statistic and $P_T,Q_T$ are its laws under $P,Q$, then
\[
 -\log\Aff(P,Q)\ge-\log\Aff(P_T,Q_T).
\]
If $C_1,C_2\succ0$ and $H\succeq0$, adding the same independent $N(0,H)$ noise gives
\begin{equation}\label{eq:affinity-common-noise}
 \Aff\{N(0,C_1+H),N(0,C_2+H)\}
       \ge\Aff\{N(0,C_1),N(0,C_2)\}.
\end{equation}
\end{lemma}

\subsection{Fractional moments and supremum bounds}
The next three lemmas supply the Ibragimov--Has'minskii arguments used in Lemma~\ref{lem:localization}. 

The first lemma shows how the Frobenius distance of the whitened covariance from the identity, together with its Frobenius increments, controls the moments and increments of a small power of the Gaussian likelihood ratio.

\begin{lemma}[Local fractional moments of the Gaussian likelihood ratio]
\label{lem:local-fractional}
Let $X$ be the whitened observation vector of Lemma~\ref{lem:gaussian}, so $\E XX^\T=I_n$. Fix $0<b\le B<\infty$. For each real symmetric $M$ with $bI_n\preceq M\preceq BI_n$, put
\[
 \mathcal L_X(M)=-\tfrac12\log\det M
                 -\tfrac12X^\T(M^{-1}-I_n)X.
\]
There are $s_*>0$ and $c>0$, depending only on $b,B$ and the common sub-Gaussian bound, such that
\[
 \E e^{s\mathcal L_X(M)}
       \le e^{-cs\|M-I_n\|_\Frob^2},\qquad 0\le s\le s_*.
\]
For every fixed $m\ge1$ and $0<s\le s_*/2$, set $\xi(M)=e^{s\mathcal L_X(M)/m}$. If $M,N$ both satisfy the same covariance bounds, then
\[
 \|\xi(M)-\xi(N)\|_{L^m}\le C_{m,s}\|M-N\|_\Frob.
\]
Consequently, suppose $u=h$ or $u=(h,a)$ indexes a continuous covariance family satisfying these bounds and, for $c_0>0$ and $C_0,L\ge0$,
\[
 \|M(u)-I_n\|_\Frob^2\ge c_0\|h\|^2-C_0,
 \qquad \|M(u)-M(v)\|_\Frob\le L\|u-v\|.
\]
Then the continuous process $\xi(u)=\xi(M(u))$ satisfies
\[
 \E|\xi(u)|^m\le A e^{-c_1\|h\|^2},\qquad
 \|\xi(u)-\xi(v)\|_{L^m}\le C_{m,s}L\|u-v\|,
\]
where $A,c_1>0$ are uniform in $n$ and the noise law whenever the displayed constants are uniform.
\end{lemma}

The next very general lemma shows how pointwise control of $L^m$ norms and the $L^m$ Lipschitz constant of a stochastic process on a compact set yields a bound on the $L^m$ norm of the supremum of the process.

\begin{lemma}[Supremum over a compact convex set]
\label{lem:compact-supremum}
Let $k\ge1$ and $m>k$. Let $\Gamma\subset\R^k$ be nonempty, compact, and convex, with diameter at most $\delta>0$. Suppose a real process $(\xi(u))_{u\in\Gamma}$ has continuous paths and satisfies
\[
 \sup_{u\in\Gamma}\|\xi(u)\|_{L^m}\le M,\qquad
 \|\xi(u)-\xi(v)\|_{L^m}\le L\|u-v\|.
\]
Then
\[
 \E\sup_{u\in\Gamma}|\xi(u)|^m
 \le C_{k,m}\{M^m+(L\delta)^kM^{m-k}\}.
\]
\end{lemma}

The final lemma gives a supremum bound whose constants do not depend on the diameter of the index set. This allows its application to expanding local parameter sets and yields the tail bound for $\widehat h_n$.

\begin{lemma}[A maximal inequality with a quadratic moment envelope]
\label{lem:quadratic-envelope}
Let $u=h\in\R^p$, or $u=(h,a)\in\R^p\times I$ with $I$ a fixed compact interval, and let $k$ be the total dimension. Suppose $\mathcal U$ is a nonempty compact convex set of such indices, $m>k$, and $(\xi(u))_{u\in\mathcal U}$ has continuous paths. Assume, for $A>0$, $c_0>0$, and $L\ge0$,
\[
 \E|\xi(u)|^m\le A e^{-c_0\|h\|^2},\qquad
 \|\xi(u)-\xi(v)\|_{L^m}\le L\|u-v\|.
\]
Then, with a supremum over the empty set interpreted as zero,
\[
 \E\sup_{\substack{u\in\mathcal U\\\|h\|>H}}|\xi(u)|^m
 \le C\{A+L^kA^{1-k/m}\}e^{-c_1H^2},\qquad H\ge0.
\]
Here $c_1>0$ and $C$ depend only on $k,m,c_0$ and the length of $I$ when present. In particular they do not depend on the diameter of $\mathcal U$ in the $h$ coordinates.
\end{lemma}

\section{Proofs of the main-body theorems and corollaries}
\label{app:proof-main}
We first prove Theorem~\ref{thm:normalized}, then explain its main constructions in the Gaussian errors, deterministic-design case. The remaining proofs follow the order of the results in Section~\ref{sec:main}. All auxiliary lemmas are stated in Appendix~\ref{app:main} and proved in Appendix~\ref{app:lemmas}.

\subsection{Proof of Theorem~\ref{thm:normalized}}\label{app:normalizedproof}
We first work at fixed observation locations satisfying the design bounds. All estimates are uniform over $\theta_0\in\mathcal K$ and $P\in\mathcal P$, with the compactness and interiority convention in Appendix~\ref{sec:proofnotation}. Steps 1--5 prove part (i) and the tail bound in part (iii); Step 6 proves the remaining moment assertions, Step 7 proves part (ii), and Step 8 gives the random-design assertions in part (iv).

\noindent\textbf{1. Tightness and curvature bounds.}
Lemma~\ref{lem:tails} gives
\[
 P_{\theta_0,P}(\|\widehat h_n\|>H)
 \le Ce^{-cH^2}+Ce^{-cn^\delta},\qquad H\ge0.
\]
This is \eqref{eq:tails}; in particular, $\widehat h_n=O_P(1)$ uniformly. We will use the same bound to prove uniform integrability after establishing convergence in distribution.

Use the definitions of $J_n,\Delta_n,K_n,V_n$ in \eqref{eq:normalizedscoreinfo}. Put $X=C_{n,0}^{-1/2}Y_n$ and recall the derivatives of the truth-whitened covariance $M$ from Lemma~\ref{lem:derivatives}:
\[
 B_j=\begin{cases}
       \partial_{h_j}M(0),&d\ge4,\\
       \partial_{h_j}M(0,a_0),&d\le3,
      \end{cases}
 \qquad 1\le j\le k_d.
\]
Thus $\Delta_{n,j}=\tfrac12\{X^\T B_jX-\tr B_j\}$ and $\E\Delta_n=0$. Lemma~\ref{lem:information} gives $J_n\asymp I$, and the quadratic-form variance bound in Lemma~\ref{lem:gaussian} gives $K_n\asymp J_n$. Consequently $V_n=J_n^{-1}K_nJ_n^{-1}\asymp I$, proving \eqref{eq:Jbounds}.

\noindent\textbf{2. Taylor expansion of the likelihood.}
For $d\le3$, put $I=[a_-,a_+]$ and
\[
 \mathcal L_n(h,a)=\ell_n(\vartheta_0+E_nh,a),\qquad
 s_n(a)=\nabla_h\mathcal L_n(0,a),\qquad s_n(a_0)=\Delta_n.
\]
The allowed pairs are $\{(h,a):a\in I_n(h)\}$ whenever $\vartheta_0+E_nh\in U$.
Every fixed $h$-ball lies in this projection for all sufficiently large $n$, uniformly, by \eqref{eq:feasible-endpoint-motion}. We evaluate the analytic criterion on the whole enclosing interval $I=I_*$ when deriving bounds, and impose $a\in I_n(h)$ only when profiling. This distinction permits Taylor expansion with $a$ fixed
even if a point on the segment is outside $\Theta$: the model and criterion are smooth on the enclosing rectangle $\mathcal R$. For $d\ge4$, use $\mathcal L_n(h)=\ell_n(\theta_0+E_nh)$. Since $\operatorname{dist}(\mathcal K,\Theta^c)>0$ and
$\|E_n\|_\op\to0$, every fixed $h$-ball is feasible eventually. In this case range is already a coordinate of $h$.

Fix $0<H<\infty$. In low dimension, Taylor's formula with integral remainder, applied along the segment from $0$ to $h$ with $a$ fixed, gives the exact identity
\begin{equation}\label{eq:likelihoodTaylorIntegral}
 \mathcal L_n(h,a)=\mathcal L_n(0,a)+s_n(a)^\T h
       +\int_0^1(1-s)h^\T\nabla_h^2\mathcal L_n(sh,a)h\,ds.
\end{equation}
The $L^p$ conclusion proved in Lemma~\ref{lem:studentization}\ref{it:studentization-local} gives, for every fixed $1\le p<\infty$,
\[
 \sup_{\substack{\|h\|\le H\\a\in I}}
       \|\nabla_h^2\mathcal L_n(h,a)+J_n\|_{\op}=o_{L^p}(1).
\]
Lemma~\ref{lem:nuisance} supplies the other input:
\[
 \sup_{a\in I}\|s_n(a)-\Delta_n\|=o_{L^p}(1).
\]
The Hessian is therefore close to $-J_n$, and the baseline score is close to the score at the true range. In high dimension the same Hessian bound holds with $a$ omitted and the score comparison is exact.

Put $Q_n(h)=\Delta_n^\T h-\tfrac12h^\T J_nh$. Subtracting this quadratic from \eqref{eq:likelihoodTaylorIntegral} gives
\[
 \begin{aligned}
 \epsilon_{n,H}
 &\coloneqq \sup_{\substack{\|h\|\le H\\a\in I}}
       |\mathcal L_n(h,a)-\mathcal L_n(0,a)-Q_n(h)|\\
 &\le H\sup_{a\in I}\|s_n(a)-\Delta_n\|
       +\frac{H^2}{2}
         \sup_{\substack{\|h\|\le H\\a\in I}}
                \|\nabla_h^2\mathcal L_n(h,a)+J_n\|_{\op}
   =o_{L^p}(1).
 \end{aligned}
\]
The factor $1/2$ is $\int_0^1(1-s)\,ds$. Thus, uniformly on the fixed ball and over the entire range interval,
\[
 \mathcal L_n(h,a)=\mathcal L_n(0,a)+Q_n(h)+o_{L^p}(1).
\]
Without a nuisance, omit $a$ and the first term in the bound for $\epsilon_{n,H}$.

\noindent\textbf{3. Profiling over the feasible range intervals.}
In low dimension define
\[
 \mathcal L_n^{\mathrm{prof}}(h)=\sup_{a\in I_n(h)}\mathcal L_n(h,a),
 \qquad G_n(a)=\mathcal L_n(0,a).
\]
The interval can now depend on $h$. Its endpoint displacement is at most $L_*H\|E_n\|_\op$ on the fixed ball, by \eqref{eq:feasible-endpoint-motion}. Project a point of $I_n(h)$ onto $I_n(0)$ by clipping it to the two endpoints. Its displacement
is no larger than that bound, and the same is true in the reverse direction. The mean value theorem therefore gives
\[
 \sup_{\|h\|\le H}
 \left|\sup_{a\in I_n(h)}G_n(a)
             -\sup_{a\in I_n(0)}G_n(a)\right|
 \le L_*H\|E_n\|_\op\sup_{a\in I_*}|G_n'(a)|
 =o_{L^p}(1),
\]
where the last step is \eqref{eq:nuisance-baseline-moments}. Thus moving the feasible nuisance interval changes the baseline profile value by a negligible amount.

Taking the supremum of the Taylor expansion over $I_n(h)$ gives
\[
 \begin{split}
 \sup_{\|h\|\le H}
  |\mathcal L_n^{\mathrm{prof}}(h)-\mathcal L_n^{\mathrm{prof}}(0)-Q_n(h)|& \le\epsilon_{n,H}
       +L_*H\|E_n\|_\op\sup_{a\in I_*}|G_n'(a)|\\
 &      \eqqcolon \bar\epsilon_{n,H}=o_{L^p}(1).
 \end{split}
\]
Equivalently, $\mathcal L_n^{\mathrm{prof}}(h)=\mathcal L_n^{\mathrm{prof}}(0)+Q_n(h)+o_{L^p}(1)$ uniformly on the fixed ball. In high dimension set $\mathcal L_n^{\mathrm{prof}}=\mathcal L_n$ and $\bar\epsilon_{n,H}=\epsilon_{n,H}$. Global maximization over $\Theta$ makes $\widehat h_n$ a global maximizer of this actual profile criterion over its feasible projection.

\noindent\textbf{4. Compare the two maximizers on a fixed ball.}
Since $J_n\succ0$, the unique maximizer of $Q_n$ over $\R^{k_d}$ is $h_n^*=J_n^{-1}\Delta_n$. It has mean zero and
\[
 \Var(h_n^*)=J_n^{-1}K_nJ_n^{-1}=V_n\preceq CI.
\]
Chebyshev's inequality gives $P(\|h_n^*\|>H)\le C/H^2$, so $h_n^*=O_P(1)$ uniformly. Together with the tail bound for $\widehat h_n$, this shows that the event
\[
 \mathcal E_{n,H}
 =\{\|\widehat h_n\|\le H,\ \|h_n^*\|\le H\}
\]
has complementary probability at most $Ce^{-cH^2}+Ce^{-cn^\delta}+C/H^2$. For fixed $H$ and sufficiently large $n$, the whole ball is allowed. On $\mathcal E_{n,H}$, therefore, $\widehat h_n$ maximizes the profile criterion on that ball, and $h_n^*$ maximizes $Q_n$ there.

Taylor expansion of $Q_n$ about $h_n^*$ is exact, because $Q_n$ is quadratic and $\nabla Q_n(h_n^*)=0$:
\[
 Q_n(h)-Q_n(h_n^*)
 =-\tfrac12(h-h_n^*)^\T J_n(h-h_n^*).
\]
On $\mathcal E_{n,H}$, profile maximization gives $\mathcal L_n^{\mathrm{prof}}(\widehat h_n)\ge \mathcal L_n^{\mathrm{prof}}(h_n^*)$. The uniform expansion on the ball then yields
\[
 -2\bar\epsilon_{n,H}
 \le Q_n(\widehat h_n)-Q_n(h_n^*)
 =-\tfrac12(\widehat h_n-h_n^*)^\T
                  J_n(\widehat h_n-h_n^*)\le0.
\]
The last inequality also follows directly from maximization of $Q_n$. If $J_n\succeq cI$, it follows that, for any $\eta>0$,
\[
 P(\|\widehat h_n-h_n^*\|>\eta)
 \le P(\mathcal E_{n,H}^c)
       +P(\bar\epsilon_{n,H}>c\eta^2/4).
\]
First let $n\to\infty$ with $H$ fixed, and then let $H\to\infty$. Both terms tend to zero in this order, uniformly. This proves
\[
 \widehat h_n=J_n^{-1}\Delta_n+o_P(1),
\]
the estimator expansion in \eqref{eq:MLEexpansion}.

\noindent\textbf{5. The score central limit theorem and standardization.}
For any deterministic unit vector $v$, put $q=K_n^{-1/2}v$ and
$B(q)=\sum_jq_jB_j$. Then
\[
 v^\T K_n^{-1/2}\Delta_n
 =\tfrac12\{X^\T B(q)X-\tr B(q)\},\qquad
 \Var(v^\T K_n^{-1/2}\Delta_n)=1.
\]
The bounds $J_n,K_n\asymp I$ give $\|B(q)\|_{\Frob}^2=2q^\T J_nq\le C$. Also $\|B(q)\|_{\op}\to0$ uniformly by the normalized derivative bounds in Lemma~\ref{lem:derivatives}. The quadratic-form central limit theorem in Lemma~\ref{lem:gaussian}, followed by Cram\'er--Wold, proves $K_n^{-1/2}\Delta_n\Rightarrow N_{k_d}(0,I)$. Applying the argument along arbitrary sequences of truths, noise laws, and unit directions gives the stated uniform convergence.

The estimator expansion implies
\[
 V_n^{-1/2}\widehat h_n
 =O_nK_n^{-1/2}\Delta_n+o_P(1),\qquad
 O_n=V_n^{-1/2}J_n^{-1}K_n^{1/2}.
\]
Here $O_nO_n^\T=I$, so $O_n$ is orthogonal. Every sequence of these deterministic matrices has a convergent subsequence, and the standard normal law is invariant under its limiting rotation. This compactness argument justifies the normal limit even though $O_n$ varies with $n$, and proves the final convergence in \eqref{eq:MLEexpansion}.

\noindent\textbf{6. Uniform integrability and moment convergence.}
Put $T_n=\|\widehat h_n\|$. Lemma~\ref{lem:tails} proves \eqref{eq:allmoments} by integrating the tail bound: for every $\lambda>0$ and $0<\alpha<\min\{2,2\delta\}$,
\[
 \sup_{n\ge n_0,\,\theta_0,P}\E e^{\lambda T_n^\alpha}<\infty.
\]
Every fixed polynomial moment is therefore uniformly bounded as well.

Let $F$ be continuous with $|F(x)|\le C_F e^{\lambda_F\|x\|^{\alpha_F}}$ for some fixed
$\lambda_F>0$ and $0<\alpha_F<\min\{2,2\delta\}$. Apply \eqref{eq:allmoments} with $2\lambda_F$ in place of $\lambda$:
\[
 \sup_{n\ge n_0,\,\theta_0,P}\E|F(\widehat h_n)|^2<\infty,
 \qquad
 \sup_{n\ge n_0,\,\theta_0,P}
 \E\!\left[|F(\widehat h_n)|
             \mathbf1_{\{|F(\widehat h_n)|>A\}}\right]
       \le\frac{C_F'}A\longrightarrow0.
\]
Thus $F(\widehat h_n)$ is uniformly integrable over the truths and noise laws as well as $n$. The same calculation applies to $F(V_n^{1/2}Z)$: $V_n\preceq CI$ and $\alpha_F<2$ imply $\sup_n\E_Ze^{2\lambda_F\|V_n^{1/2}Z\|^{\alpha_F}}<\infty$.

To combine this with convergence in distribution when $V_n$ need not converge, take any sequence of truths and noise laws. Every subsequence has a further subsequence on which $V_n\to V\succ0$. Along it, \eqref{eq:MLEexpansion} gives $\widehat h_n\Rightarrow V^{1/2}Z$. Continuity of $F$ and uniform integrability then give
\[
 \E F(\widehat h_n)\longrightarrow \E_ZF(V^{1/2}Z),
 \qquad
 \E_ZF(V_n^{1/2}Z)\longrightarrow \E_ZF(V^{1/2}Z).
\]
Their difference therefore tends to zero along every such subsequence. If the uniform conclusion failed, a sequence with its difference bounded away from zero would contradict this argument. This proves \eqref{eq:momentcomparison}. Taking $F(x)=x_i$ and $F(x)=x_ix_j$ gives the stated first- and second-moment comparisons. In particular,
$\Cov(\widehat h_n)=V_n+o(1)$. The covariance blocks needed for part (ii) are obtained in the next step.

Finally, $\|V_n^{-1/2}\widehat h_n\|\le C T_n$. The bound \eqref{eq:allmoments} holds for every fixed $\lambda$, so it also bounds $\E e^{\lambda\|V_n^{-1/2}\widehat h_n\|^\alpha}$ by replacing $\lambda$ with $\lambda C^\alpha$. The same uniform
integrability argument applied to the standardized convergence in \eqref{eq:MLEexpansion} proves \eqref{eq:normalizedmoments} for the stated function class.

For these fixed designs the estimator expansion remainder is also $o_{L^p}(1)$ for every fixed finite $p$. Indeed, Lemma~\ref{lem:gaussian} and $\|J_n^{-1}\|_\op\le C$ give $\sup \E\|h_n^*\|^{2p}\le C_p$. The tail bound supplies the same bound for $\widehat h_n$. If $D_n=\|\widehat h_n-h_n^*\|$, then
\[
 \E D_n^p\le\eta^p+(\E D_n^{2p})^{1/2}P(D_n>\eta)^{1/2}.
\]
For each $\eta>0$, convergence in probability makes the second term tend uniformly to zero. Letting $\eta\downarrow0$ proves the claim. For random locations this strengthened remainder assertion is conditional on designs with the common bounds.

\noindent\textbf{7. Covariance structure.}
Lemma~\ref{lem:noise-diagonals} proves the covariance correction \eqref{eq:kurtosis-correction}, the signal--nugget block structure \eqref{eq:covariance-blocks}, and, when $d=4$, the additional range separation in \eqref{eq:critical-covariance-block}.
These are the assertions in part (ii).

\noindent\textbf{8. Random locations.}
Lemma~\ref{lem:randomgood} supplies events $\mathcal G_n$, depending only on the locations, whose complements have summable probabilities. On these events all the preceding deterministic bounds have common constants for every truth and noise law. Almost every realized design therefore satisfies the same conditional proof eventually, including the global-localization step and the tail-to-moment argument in
Step~6. This proves parts (i)--(iii) conditionally. In particular, the expectations in \eqref{eq:momentcomparison} are conditional response expectations, and $\E_Z$ averages only over an independent standard normal vector, with $V_n$ fixed at the realized design. The local $L^p$ bounds used in the Taylor and score comparisons are conditional on these locations as well. For the joint-law assertions, Lemma~\ref{lem:gaussian} gives $cJ_n\preceq K_n\preceq CJ_n$ at every design, so $J_n$ and $K_n$ have the same nullspace. If $J_n$ is singular, define each of $J_n^{-1}\Delta_n$, $K_n^{-1/2}\Delta_n$, and $V_n^{-1/2}\widehat h_n$ to be the zero vector. On $\mathcal G_n$ the inverse matrices exist for all sufficiently large $n$. For every $\epsilon>0$, conditioning on the locations gives
\[
 \begin{aligned}
 &~~~~\sup_{\theta_0,P}
   P_{\theta_0,P}\!\left(\|\widehat h_n-J_n^{-1}\Delta_n\|>\epsilon\right)\\
 &\le P_X(\mathcal G_n^c)
 +\sup_{x\in\mathcal G_n}\sup_{\theta_0,P}
   P_{\theta_0,P}\!\left(\|\widehat h_n-J_n^{-1}\Delta_n\|>\epsilon
      \,\middle|\, x_{n,1:n}=x\right)\longrightarrow0 .
 \end{aligned}
\]
The first term tends to zero by \eqref{eq:randomfailure}; the conditional expansion on good designs bounds the second term uniformly. This proves the joint $o_P(1)$ remainder without a moment bound for $J_n^{-1}\Delta_n$ on the other designs. The random-design part of Lemma~\ref{lem:tails} proves the unconditional bound \eqref{eq:tails} by splitting over $\mathcal G_n$ and its complement. In particular, with $T_n=\|\widehat h_n\|$, the unconditional errors are uniformly tight.

Let $W_n$ denote either $K_n^{-1/2}\Delta_n$ or $V_n^{-1/2}\widehat h_n$, and put $\Phi=N_{k_d}(0,I_{k_d})$. Conditioning on the locations and using $\|f\|_\infty\le1$ in the definition of $d_{\mathrm{BL}}$ gives
\[
 \sup_{\theta_0,P}
 d_{\mathrm{BL}}\bigl(\mathcal L_{\theta_0,P}(W_n),\Phi\bigr)
 \le 2P_X(\mathcal G_n^c)
 +\sup_{x\in\mathcal G_n}\sup_{\theta_0,P}
 d_{\mathrm{BL}}\bigl(
    \mathcal L_{\theta_0,P}(W_n\mid x_{n,1:n}=x),\Phi\bigr).
\]
The first term tends to zero by \eqref{eq:randomfailure}; the second tends to zero by the uniform conditional result on good designs. This proves both joint-law normal limits.

For the moment assertions, the same random-design part of Lemma~\ref{lem:tails} integrates this unconditional tail bound and proves \eqref{eq:allmoments} under the joint law. Applying that bound with $2\lambda$ gives uniform integrability for the stated test functions, as in Step~6, and uniform bounds for all fixed polynomial
moments.

The same reasoning applies to the standardized errors: wherever the normalizer is defined,
\[
 V_n^{-1}=J_nK_n^{-1}J_n\preceq CJ_n\preceq CI,
 \qquad
 \|V_n^{-1/2}\widehat h_n\|\le CT_n,
\]
by the uniform quadratic-form variance bound and the universal curvature upper bound. Use the stated zero convention if the normalizer is undefined. The tail integral therefore also bounds the exponential moments of the standardized errors. Uniform integrability and their unconditional convergence in distribution prove
\eqref{eq:normalizedmoments} under the joint law. The matrices in the conditional comparison always use the realized locations. This proves all the joint-law assertions. \qed

\subsection{Proof explanation}\label{app:proof-outline}

We give an extended explanation of the proof of \Cref{thm:normalized} for Gaussian errors and deterministic sites, motivating all the key constructions of the proof.
The extensions to random sites and uniformly sub-Gaussian errors require technical adjustments to this argument. Throughout, let $\theta_0\in\mathcal K$ denote the true value and center the derivations around $\theta_0$. Conceptually, we carry out the proof assuming a fixed truth. The bounds are uniform in $\theta_0\in\mathcal K$, which gives the uniform conclusions. Write $P_{\theta,n}=N(0,C_n(\theta))$ and let $p_{\theta,n}$ denote its density; $P_0$ and $\E_0$ refer to the true law.

\noindent\textbf{Preliminary consistency.}
First, we prove the preliminary localization
\[
 \frac{|\widehat\eta_n-\eta_0|}{\log n}
       +|\widehat\nu_n-\nu_0|+|\widehat t_n-t_0|
       =O_P(n^{-\delta}),
\]
which is Theorem~\ref{thm:preliminary}. The method of proof is the affinity method, a version of which was also used in the RBF proof of \citet{qaqishli2026rbf}.
The argument starts with the following calculation:
\[
 \begin{aligned}
 P_0\{\ell_n(\theta)-\ell_n(\theta_0)\ge-1\}
 &=P_0\{e^{\{\ell_n(\theta)-\ell_n(\theta_0)\}/2}\ge e^{-1/2}\}\\
 &=P_0\!\left\{
      \sqrt{\frac{p_{\theta,n}(Y_n)}{p_{\theta_0,n}(Y_n)}}
                                             \ge e^{-1/2}\right\}\\
 &\le e^{1/2}\E_0
      \sqrt{\frac{p_{\theta,n}(Y_n)}{p_{\theta_0,n}(Y_n)}}\\
 &=e^{1/2}\int\sqrt{p_{\theta,n}(y)}\sqrt{p_{\theta_0,n}(y)}\,dy\\
 &=e^{1/2}\Aff(P_{\theta,n},P_{\theta_0,n})
   =e^{1/2-A_n(\theta)},
 \end{aligned}
\]
where $A_n(\theta)=-\log\Aff(P_{\theta,n},P_{\theta_0,n})$. The negative log affinity, also called the Bhattacharyya distance, therefore controls the probability that an alternative has log-likelihood at least $\ell_n(\theta_0)-1$.

We obtain the following key affinity bound. Put $\Delta_\eta=\eta-\eta_0$, $\Delta_\nu=\nu-\nu_0$, and $\Delta_t=t-t_0$. For a slowly growing frequency $B_n=n^\epsilon$, the bound is
\[
 |\Delta_\eta-2\Delta_\nu\log(2\pi B_n)|
       +|\Delta_\nu|+|\Delta_t|>Cn^{-\delta}
 \quad\Longrightarrow\quad
 A_n(\theta)\ge cn^{\delta_1}
\]
for some $\delta,\delta_1>0$. If we write
\[
 d_n(\theta)=|\Delta_\eta-2\Delta_\nu\log(2\pi B_n)|
                         +|\Delta_\nu|+|\Delta_t|,
\]
then, whenever $d_n(\theta)>Cn^{-\delta}$,
\[
 P_0\{\ell_n(\theta)-\ell_n(\theta_0)\ge-1\}
       \le Ce^{-cn^{\delta_1}}.
\]
It is therefore exponentially unlikely that such a parameter has log-likelihood even within one unit below that at the truth. A crude mesh and union bound argument then shows, after enlarging $C$, that
\[
 P_0\!\left\{
  \sup_{d_n(\theta)>Cn^{-\delta}}\ell_n(\theta)
                      \ge\ell_n(\theta_0)-\tfrac12\right\}
       \le Ce^{-cn^{\delta_1}}.
\]
To obtain this supremum bound, cover the parameter set by a sufficiently fine mesh with polynomially many points. Outside an exponentially unlikely event, a polynomial bound on the likelihood gradient makes the rounding error at most $1/4$. Thus the event in the preceding display implies that a mesh point, still separated from the truth after enlarging $C$, has log-likelihood at least $\ell_n(\theta_0)-1$. Applying the preceding pointwise bound to these mesh points and taking a union bound proves the claimed supremum bound: the polynomial number of points is absorbed by the exponential bound. A global maximizer has likelihood at least that of the truth, so
$d_n(\widehat\theta_n)=O_P(n^{-\delta})$. The preliminary localization follows on writing $\Delta_\eta=\{\Delta_\eta-2\Delta_\nu\log(2\pi B_n)\}
                 +2\Delta_\nu\log(2\pi B_n)$.

To prove the key affinity bound, we use the fact that discarding information increases affinity. For any statistic $T$, the affinity between its distributions is at least that between the full observation distributions. Consequently, for any collection $\mathcal T$ of statistics,
\[
 A_n(\theta)\ge\sup_{T\in\mathcal T}
       -\log\Aff\{\mathcal L_\theta(T),\mathcal L_0(T)\},
\]
where $\mathcal L_\theta(T)$ denotes the law of $T$ under $P_{\theta,n}$. We may therefore use whichever statistic distinguishes a given alternative and has an affinity that is easy to compute or bound. Let $Q$ have side length $L$ and let
$\chi\in C_c^\infty(Q)$, $0\le\chi\le1$, equal one on a cube $Q_0\Subset Q$ of side length $L/m$, for a fixed integer $m\ge3$. The statistics are constructed from estimates of the Fourier coefficients
\[
 X_\lambda=L^{-d/2}\int_Q\chi(x)F(x)e^{-2\pi i\lambda^\top x}\,dx,
 \qquad \lambda\in\mathscr L=(m/L)\mathbb Z^d.
\]
These form a subcollection of the Fourier series coefficients of $\chi F$ on $Q$ and also play a central role in the preprint of \citet{pillai2026}. For preliminary consistency, we use the nearest-site quadrature approximations
\[
 \widetilde X_{\lambda,n}
 =L^{-d/2}\int_Q\chi(x)Y_{n,i(x)}e^{-2\pi i\lambda^\top x}\,dx,
\]
where $i(x)$ is the index of a site closest to $x$. The heuristic, supported by Lemma~\ref{lem:windows}, is that the $X_\lambda$ are weakly correlated, with
\[
 \E_\theta|X_\lambda|^2
 =c_\chi f_\theta(\lambda)\{1+O(\langle\lambda\rangle^{-2})\},
 \qquad c_\chi=L^{-d}\int_Q\chi(x)^2\,dx.
\]
For $B>0$, set
\[
 \widetilde{\mathbf X}_{B,n}
 =\bigl(\sqrt{2}\,\Re\widetilde X_{\lambda,n},
        \sqrt{2}\,\Im\widetilde X_{\lambda,n}\bigr)
                                  _{\lambda\in\Lambda_B^+},
 \qquad
 \Lambda_B^+=\mathscr L\cap B(Be_1,\omega_*B/4),
\]
with fixed small relative bandwidth $\omega_*>0$. For slowly growing $B$, the relative $L^2$ recovery error tends to zero. More precisely, Lemma~\ref{lem:windows} gives, for any fixed $0<\zeta<\min\{\nu_-,1\}$,
\[
 \sup_{\theta\in\Theta}\sup_{\lambda\in\Lambda_B^+}
 \frac{\E_\theta|\widetilde X_{\lambda,n}-X_\lambda|^2}
      {\E_\theta|X_\lambda|^2}
 \le C(h_n^{2\zeta}+h_n^d)B^{2\nu_++d}=o(1),
\]
provided $B$ grows sufficiently slowly. Here $0<\nu_-<\nu_+<\infty$ are fixed smoothness bounds for an enclosing compact parameter rectangle. The square root of this ratio is the relative $L^2$ error. The same bound holds for every nonzero
real linear combination of the coefficient vector.

In the proof of Lemma~\ref{lem:preliminary}, the collection $\mathcal T$ consists of $\widetilde{\mathbf X}_{B_n,n}$, $\widetilde{\mathbf X}_{2B_n,n}$, and a vector of
$p_n=\lfloor n^\xi\rfloor$ disjoint pair differences $(Y_{n,i}-Y_{n,j})/\sqrt2$, at distances at most $Ch_n$. The exponents $\epsilon$ and $\xi$ are chosen sufficiently small. The Fourier vectors at $B_n$ and $2B_n$ distinguish $\eta$ and $\nu$.
The pair differences have joint covariance close to $tI_{p_n}$, so they are approximately iid $N(0,t)$ and distinguish the nugget variance. Choosing the frequencies and number of pairs to grow slowly keeps the approximation errors small enough for the affinity bound. This completes the proof of preliminary consistency.

Preliminary consistency implies, after decreasing $\delta$ if necessary, that every maximizer lies in
\[
 \mathcal N_n=\{\theta\in\Theta:
   |\Delta_\eta-2L_n\Delta_\nu|\le n^{-\delta}\log n,
   \ |\Delta_\nu|\le n^{-\delta},\ |\Delta_t|\le n^{-\delta}\},
 \qquad L_n=\log R_n,
\]
with probability tending to one; its failure probability is at most $Ce^{-cn^{\delta_1}}$. The rest of the proof therefore only needs to work within $\mathcal N_n$. For the remainder of this outline, assume $d\ge4$ to avoid the complication of treating the non-identifiable range as a nuisance parameter. The argument with this nuisance is similar, using additional lemmas that control its effect.

\noindent\textbf{Information from recoverable Fourier coefficients.}
The next step is to find the normalizing matrix $E_n$. After the change of coordinates $\theta=\theta_0+E_nh$, write $J_n$ for the Fisher information at $h=0$. Our aim is to prove $J_n\asymp I$. For example, in a regular iid model, $E_n=n^{-1/2}I$. With such a choice of $E_n$, one typically expects $\widehat h_n=O_P(1)$. In our setting, establishing $J_n\asymp I$ is a key part of proving this assertion. It is convenient to work with the whitened covariance matrix
\[
 M(h)=C_{n,0}^{-1/2}C_n(\theta_0+E_nh)C_{n,0}^{-1/2}.
\]
We have
\[
 [J_n]_{jk}=\tfrac12\tr(B_jB_k),\qquad B_j=\partial_{h_j}M(0),
 \qquad h^\top J_nh=\tfrac12\|B(h)\|_\Frob^2,
\]
where $B(h)=\sum_jh_jB_j$. The spectral calculations of Lemmas~\ref{lem:upper} and \ref{lem:derivatives} show that our choice of $E_n$ gives $\|B(h)\|_\Frob^2\le C\|h\|^2$. The more difficult part is the Fisher information lower bound in Lemma~\ref{lem:information}, which gives $\|B(h)\|_\Frob^2\ge c\|h\|^2$. These two calculations motivate the choice of $E_n$ and establish $J_n\asymp I$.

To prove the lower bound, we start from the Cram\'er--Rao information inequality: for any fixed real vector statistic $Z$ with $\Cov_0(Z)\preceq CI$,
\[
 h^\top J_nh\ge C^{-1}\|D_Zh\|^2,\qquad
 D_Zh=\left.\frac{d}{ds}\E_{\theta_0+sE_nh}Z\right|_{s=0}.
\]
The statistic, including any coefficients chosen using $\theta_0$, is held fixed while differentiating its expectation. The derivative of the mean is its covariance with the directional score, so the inequality follows from Cauchy--Schwarz.

We construct the statistics from the best linear predictions of the Fourier coefficients. Put
\[
 U_\lambda=\frac{X_\lambda}{\sqrt{v_0(\lambda)}},\qquad
 v_0(\lambda)=\E_0|X_\lambda|^2,\qquad
 \widehat U_{\lambda,n}=\E_0(U_\lambda\mid Y_n).
\]
These Gaussian conditional expectations are linear functions of the observations. Lemma~\ref{lem:sampling} gives, for
$\|\lambda\|\le\epsilon R_n$ and fixed small $\epsilon>0$,
\[
 \E_0|U_\lambda-\widehat U_{\lambda,n}|^2
       \le Cn^{-1}\langle\lambda\rangle^{2\nu_0+d}
       \le C\epsilon^{2\nu_0+d}
\]
for all sufficiently large $n$. The left side is exactly the relative mean squared recovery error for $X_\lambda$, since
\[
 \E_0|U_\lambda-\widehat U_{\lambda,n}|^2
 =\frac{\E_0|X_\lambda-\widehat X_{\lambda,n}|^2}{v_0(\lambda)}.
\]
Here $\widehat X_{\lambda,n}=\E_0(X_\lambda\mid Y_n)$. This explains the role of $R_n=n^{1/(2\nu_0+d)}$: taking $\epsilon$ small guarantees recovery with small relative error for $\|\lambda\|\le\epsilon R_n$. The error bound tends to zero for
$\|\lambda\|=o(R_n)$; at higher frequencies it no longer guarantees small relative error.

For comparison, the nearest-site quadrature bound from Lemma~\ref{lem:windows} gives the absolute mean squared error bound
\[
 \E_0|\widetilde X_{\lambda,n}-X_\lambda|^2
 \le C(h_n^{2\zeta}+h_n^d)
 \le C(n^{-2\zeta/d}+n^{-1}).
\]
The two terms come from replacing the field within each cell by its value at the nearest site and from the observation noise, respectively; the last inequality uses the sharp deterministic coverage assumption. Dividing by $v_0(\lambda)\asymp\langle\lambda\rangle^{-2\nu_0-d}$ therefore gives
\[
 \frac{\E_0|\widetilde X_{\lambda,n}-X_\lambda|^2}{v_0(\lambda)}
 \le C(n^{-2\zeta/d}+n^{-1})
                         \langle\lambda\rangle^{2\nu_0+d}.
\]
At $\|\lambda\|\asymp\epsilon R_n$, this upper bound is only
\[
 C\epsilon^{2\nu_0+d}(n^{1-2\zeta/d}+1),
\]
which grows with $n$ for fixed $\epsilon>0$: here $2\zeta/d<1$ because $d\ge4$ and $\zeta<1$. Thus the nearest-site upper bound does not establish small relative error at the required frequencies. This is a limitation of the available bound, not a lower bound on the quadrature error itself. The quadrature error bound suffices for the slowly growing frequencies used in preliminary consistency. For the frequencies used in the information lower bound, we need the sharper error bound for best linear prediction.

For $d>4$, choose a sufficiently small relative bandwidth $\omega_*\in(0,1/4)$ and define
\[
 t_j=2^j\epsilon/8,\qquad
 \Lambda_{n,j}=\left\{\lambda\in\mathscr L:\lambda_1>0,
       \left|\frac{\|\lambda\|}{t_jR_n}-1\right|\le\omega_*\right\},
 \qquad j=0,1,2.
\]
These are three narrow, disjoint bands at successively doubled radii below $\epsilon R_n$. With $N_{n,j}=\#\Lambda_{n,j}$, use the signal statistics
\[
 Z_{n,j}=N_{n,j}^{-1/2}\sum_{\lambda\in\Lambda_{n,j}}
       \{|\widehat U_{\lambda,n}|^2-\E_0|\widehat U_{\lambda,n}|^2\},
 \qquad j=0,1,2.
\]
The weak-correlation heuristic explains the bound $\Cov_0(Z_n)\preceq CI$, proved using
Lemmas~\ref{lem:windows} and \ref{lem:sampling}.

The choice of $Z_{n,j}$ is motivated by the score in an ideal model that treats the Fourier coefficients as independent, with variances given by the spectral density.
After absorbing the fixed window factor $c_\chi$ into the coefficients, this model is
\[
 X_\lambda^{\mathrm{id}}\stackrel{\mathrm{ind}}\sim
 N_{\mathbb C}(0,f_{\theta_0+E_nh}(\lambda)),\qquad
 U_\lambda^{\mathrm{id}}=\frac{X_\lambda^{\mathrm{id}}}
                                  {\sqrt{f_0(\lambda)}}.
\]
Here $N_{\mathbb C}(0,v)$ has independent real and imaginary parts, each with variance $v/2$. Differentiating its log-likelihood on band $j$ gives
\[
 \nabla_h\ell_{n,j}^{\mathrm{id}}(0)
 =\sum_{\lambda\in\Lambda_{n,j}}
    \left.\nabla_h\log f_{\theta_0+E_nh}(\lambda)\right|_{h=0}
       (|U_\lambda^{\mathrm{id}}|^2-1).
\]
The log spectral derivatives vary little across a sufficiently narrow band. Thus, up to these nearly constant derivative factors, the score is the band energy $\sum_{\lambda\in\Lambda_{n,j}}(|U_\lambda^{\mathrm{id}}|^2-1)$. This motivates $Z_{n,j}$: replace the actual normalized coefficients $U_\lambda$ by their predictions, center each squared magnitude at its expectation under the truth, and divide by $\sqrt{N_{n,j}}$ so that the variance is of order one.

To calculate the mean response $D_{Z_n}h$, first consider the latent energies. Write $v_\theta(\lambda)=\E_\theta|X_\lambda|^2$ and $\theta=\theta_0+E_nh$. Since $\E_\theta|U_\lambda|^2=v_\theta(\lambda)/v_0(\lambda)$,
\[
 \begin{aligned}
 \left.\nabla_h \E_\theta|U_\lambda|^2\right|_{h=0}
 &=\frac{\left.\nabla_h v_\theta(\lambda)\right|_{h=0}}
         {v_0(\lambda)}\\
 &=\left.\nabla_h\log v_\theta(\lambda)\right|_{h=0}
   \approx\left.\nabla_h\log f_\theta(\lambda)\right|_{h=0}.
 \end{aligned}
\]
The last approximation is made precise by the log-derivative bounds in Lemma~\ref{lem:windows}. Thus the derivatives of $\log f_\theta$ in $h$ coordinates describe how the band energies respond to parameter changes. For a signal direction $h=(h_u,h_v,h_a,0)$, put $m_0(\omega)=a_0+4\pi^2\|\omega\|^2$ and $\beta_0=\nu_0+d/2$.
Then the directional log spectral derivative at the truth is
\[
 g_h(\omega)=\left.\frac{d}{ds}
                 \log f_{\theta_0+sE_nh}(\omega)\right|_{s=0}
 =\frac{h_u+h_v\{2\log R_n-\log m_0(\omega)\}}{r_n}
       -\frac{\beta_0h_a}{b_nm_0(\omega)}.
\]
For $d>4$, $r_n/b_n=R_n^2$, so on a fixed annulus in $z$,
\[
 r_ng_h(R_nz)
 =h_u-h_v\log(4\pi^2\|z\|^2)
       -\frac{\beta_0h_a}{4\pi^2\|z\|^2}+o(1)\|h\|.
\]
The simplicity of the above formula is due to the choice of $E_n$, as we now explain. Since $\theta=\theta_0+E_nh$, the chain rule gives
\[
 \partial_{h_i}
 =\sum_j(E_n)_{ji}\partial_{\theta_j}.
\]
Holding the centered amplitude $u$ fixed gives
\[
 \left.\partial_v\right|_u
 =2L_n\partial_\eta+\partial_\nu,
 \qquad
 \partial_{h_v}
 =r_n^{-1}(2L_n\partial_\eta+\partial_\nu).
\]
This combination cancels the diverging $-2\log R_n$ term in the smoothness derivative and yields the displayed formula for $r_ng_h(R_nz)$. The three response shapes are constant, logarithmic, and reciprocal-square; the three bands distinguish their coefficients. We prove that $\|D_{Z_n}h\|\ge c\|h\|$ for signal directions $h$ with $h_t=0$, which gives the signal information lower bound.

At $d=4$, the normalized range cross-information with centered amplitude and smoothness tends to zero. Two bands suffice for the lower bound on the amplitude--smoothness information block. The additional statistic for the range information lower bound uses $\Lambda_n^a=\{\lambda\in\mathscr L:\lambda_1>0,
1\le\|\lambda\|\le\epsilon R_n\}$ and is
\[
 Z_{n,a}=\sum_{\lambda\in\Lambda_n^a}
 \frac{|\widehat U_{\lambda,n}|^2-\E_0|\widehat U_{\lambda,n}|^2}
      {b_nm_0(\lambda)}.
\]
In the same ideal independent model, the normalized range score on this frequency set is
\[
 \begin{aligned}
 \partial_{h_a}\ell_n^{\mathrm{id}}(0)
 &=\frac1{b_n}\sum_{\lambda\in\Lambda_n^a}
       \left.\partial_a\log f_\theta(\lambda)\right|_{\theta=\theta_0}
          (|U_\lambda^{\mathrm{id}}|^2-1)\\
 &=-\frac{\beta_0}{b_n}\sum_{\lambda\in\Lambda_n^a}
             \frac{|U_\lambda^{\mathrm{id}}|^2-1}{m_0(\lambda)}.
 \end{aligned}
\]
Thus, apart from the fixed factor $-\beta_0$, $Z_{n,a}$ is the predicted version of this range score, using the actual normalizations and centering each squared magnitude at its expectation under the truth. It combines frequencies from $1$ to $\epsilon R_n$ because the range information is spread over logarithmically many frequency scales. Finally, in every $d\ge4$, the nugget entry tends to $(2t_0^2)^{-1}$ and its cross-information with the signal coordinates converges to $0$. This completes the lower bound for the full $J_n$.

\noindent\textbf{Localization at the information scale.}
For the supremum argument, enlarge the preliminary set to $\mathcal N_n^\square$, defined using the enclosing rectangle $\mathcal R$ in Appendix~\ref{sec:proofnotation}. Its normalized image is compact and convex and contains every feasible candidate in $\mathcal N_n$. The covariance estimates hold on this enlarged set, so the argument does not require $\Theta$ to be convex. On $\mathcal N_n^\square$, we have $M(h)\asymp I$, which implies
\[
 A_n(\theta_0+E_nh)\asymp\|M(h)-I\|_\Frob^2.
\]
Expanding the covariance gives
\[
 M(h)-I=B(h)+\mathcal R_n(h),\qquad
 \|M(h)-I\|_\Frob\ge\|B(h)\|_\Frob-\|\mathcal R_n(h)\|_\Frob.
\]
The critical Fisher information lower bound gives $\|B(h)\|_\Frob\ge c\|h\|$, and the derivative estimates in Lemma~\ref{lem:derivatives} give $\|\mathcal R_n(h)\|_\Frob\le\epsilon_n\|h\|$, with $\epsilon_n\to0$. We therefore obtain
\[
 \|M(h)-I\|_\Frob^2\ge c\|h\|^2,
 \qquad A_n(\theta_0+E_nh)\ge c\|h\|^2.
\]

Just as in the RBF proof \citep{qaqishli2026rbf}, we use this separation to obtain a tail bound for $\widehat h_n$. Write $R(h)$ for the likelihood ratio $p_{\theta_0+E_nh,n}(Y_n)/p_{\theta_0,n}(Y_n)$. Take $q=5$, fix $0<s\le\min\{1/2,s_*/2\}$ with $s_*$ from Lemma~\ref{lem:local-fractional}, and put $\xi_n(h)=R(h)^{s/q}$. Since $\E_0 R(h)^{1/2}=e^{-A_n(\theta_0+E_nh)}$ and $2s\le1$,
Jensen's inequality gives
\[
 \|\xi_n(h)\|_{L^q(P_0)}^q=\E_0 R(h)^s
 \le\{\E_0 R(h)^{1/2}\}^{2s}
 =e^{-2s A_n(\theta_0+E_nh)}\le e^{-c\|h\|^2}.
\]
For this same power, Lemma~\ref{lem:local-fractional} and the Frobenius Lipschitz bound give
\[
 \|\xi_n(h)-\xi_n(h')\|_{L^q(P_0)}
       \le C\|M(h)-M(h')\|_\Frob\le C\|h-h'\|.
\]
Combining this pointwise moment bound with the uniform Lipschitz bound, the Ibragimov--Has'minskii argument in Lemma~\ref{lem:quadratic-envelope} yields, with
$\mathcal H_n=\{h:\theta_0+E_nh\in\mathcal N_n^\square\}$,
\[
 \left\|\sup_{\substack{h\in\mathcal H_n\\\|h\|>H}}
                 \xi_n(h)\right\|_{L^q(P_0)}\le Ce^{-cH^2}.
\]
Since $\xi_n(0)=1$ and a global maximizer has at least the likelihood of the truth,
\[
 \begin{aligned}
 P_0(\|\widehat h_n\|>H)
 &\le P_0(\widehat\theta_n\notin\mathcal N_n)
       +P_0\!\left\{\sup_{\substack{h\in\mathcal H_n\\\|h\|>H}}
                                      \xi_n(h)\ge1\right\}\\
 &\le Ce^{-cn^{\delta_1}}
       +\left\|\sup_{\substack{h\in\mathcal H_n\\\|h\|>H}}
                         \xi_n(h)\right\|_{L^q(P_0)}^q\\
 &\le Ce^{-cn^{\delta_1}}+Ce^{-cH^2}.
 \end{aligned}
\]
This is stronger than $\widehat h_n=O_P(1)$.

\noindent\textbf{Quadratic likelihood and normal approximation.}
The final step is a Taylor expansion in $h$ coordinates. Since $\widehat h_n=O_P(1)$, we may work on a fixed ball $B_H(0)$ that contains it with arbitrarily high probability when $H$ is sufficiently large. Set $\mathcal L_n(h)=\ell_n(\theta_0+E_nh)$. Taylor's formula with integral remainder gives
\[
 \mathcal L_n(h)=\mathcal L_n(0)+\Delta_n^\top h
       +\int_0^1(1-s)h^\top\nabla_h^2\mathcal L_n(sh)h\,ds.
\]
Lemmas~\ref{lem:derivatives} and \ref{lem:gaussian-hessian} show that
\[
 \sup_{\|h\|\le H}\|\nabla_h^2\mathcal L_n(h)+J_n\|_\op=o_P(1),
\]
so, uniformly on the ball,
\[
 \mathcal L_n(h)-\mathcal L_n(0)
 =\Delta_n^\top h-\tfrac12h^\top J_nh+o_P(1)
 =Q_n(h)+o_P(1).
\]
The quadratic $Q_n$ is maximized at $h_n^*=J_n^{-1}\Delta_n$. Since $\Cov_0(h_n^*)=J_n^{-1}\preceq CI$, this maximizer is also $O_P(1)$. Both maximizers therefore lie in a fixed, eventually feasible ball with high probability. Comparing their likelihood values gives
\[
 0\le Q_n(h_n^*)-Q_n(\widehat h_n)
   =\mathcal L_n(h_n^*)-\mathcal L_n(\widehat h_n)+o_P(1)\le o_P(1).
\]
Taylor expansion of the quadratic $Q_n$ about its maximizer gives the exact identity
\[
 Q_n(h_n^*)-Q_n(\widehat h_n)
 =\tfrac12(\widehat h_n-h_n^*)^\top J_n(\widehat h_n-h_n^*).
\]
Together with $J_n\succeq cI$, this shows that $\widehat h_n=J_n^{-1}\Delta_n+o_P(1)$.
Each standardized score contrast is a sum of independent centered squared Gaussian variables whose largest coefficient tends to zero. The Lyapunov central limit theorem and Cram\'er--Wold give
\[
 J_n^{-1/2}\Delta_n\Rightarrow N_4(0,I_4),\qquad
 J_n^{1/2}\widehat h_n\Rightarrow N_4(0,I_4).
\]
The tail bound gives uniform integrability and therefore also the stated moment convergences.

\subsection{Proof of Corollary~\ref{cor:rates}}\label{app:ratesproof}
We first work at fixed locations. Theorem~\ref{thm:normalized} gives, uniformly over the truth and noise law,
\[
 \|\widehat h_n\|_{L^p}\le C_p\quad(1\le p<\infty),\qquad
 \E(\widehat h_n\widehat h_n^\T)=V_n+o(1),\qquad V_n\asymp I.
\]
For a scalar coordinate $g$, write $b_{g,n}=Dg(\vartheta_0)E_n$ and $H_{g,n}(\vartheta)=E_n^\T\nabla_\vartheta^2g(\vartheta)E_n$. Taylor's formula gives
\begin{equation}\label{eq:coordinate-taylor-recipe}
 \begin{aligned}
 g(\widehat\vartheta_n)-g(\vartheta_0)
 &=b_{g,n}\widehat h_n+\mathcal R_{g,n},\\
 \mathcal R_{g,n}
 &=\int_0^1(1-s)\widehat h_n^\T
       H_{g,n}(\vartheta_0+sE_n\widehat h_n)
       \widehat h_n\,ds.
 \end{aligned}
\end{equation}
If $\sup_\vartheta\|H_{g,n}(\vartheta)\|_{\mathrm{op}} \le C\epsilon_{g,n}$ on the enclosing parameter rectangle, then
\[
 |\mathcal R_{g,n}|\le C\epsilon_{g,n}\|\widehat h_n\|^2,
 \qquad
 \|\mathcal R_{g,n}\|_{L^p}\le C_p\epsilon_{g,n},
 \quad 1\le p<\infty.
\]
The segment in \eqref{eq:coordinate-taylor-recipe} is contained in this rectangle even when the optimization set is not convex.

Write $B(\nu)=\partial_\nu\log A_d(\nu)$ and $B_0=B(\nu_0)$. Table~\ref{tab:coordinate-expansions} gives the derivative rows and remainder orders. All rows and matrices below use the order $(h_u,h_v,h_a,h_t)$. For $d\le3$, omit the third coordinate
and the formulas involving $a,l$, or $\sigma^2$.

\begin{table}[!htbp]
\centering\small
\renewcommand{\arraystretch}{1.35}
\caption{Derivative rows and Taylor remainder orders:
$g(\widehat\vartheta_n)-g(\vartheta_0)
=b_{g,n}\widehat h_n+O_{L^p}(\epsilon_{g,n})$, uniformly for
every fixed $1\le p<\infty$. Rows involving $a,l$, or
$\sigma^2$ apply when $d\ge4$.}
\label{tab:coordinate-expansions}
\begin{tabular}{@{}lcl@{}}\toprule
$g$ & $b_{g,n}=Dg(\vartheta_0)E_n$ & $\epsilon_{g,n}$\\\midrule
$\eta$ & $(r_n^{-1},\ 2L_nr_n^{-1},\ 0,\ 0)$ & $0$\\
$\nu$ & $(0,\ r_n^{-1},\ 0,\ 0)$ & $0$\\
$a$ & $(0,\ 0,\ b_n^{-1},\ 0)$ & $0$\\
$t=\tau^2$ & $(0,\ 0,\ 0,\ n^{-1/2})$ & $0$\\
$\log\nu$ & $(0,\ (\nu_0r_n)^{-1},\ 0,\ 0)$ & $r_n^{-2}$\\
$\log a$ & $(0,\ 0,\ (a_0b_n)^{-1},\ 0)$ & $b_n^{-2}$\\
$\log t$ & $(0,\ 0,\ 0,\ t_0^{-1}n^{-1/2})$ & $n^{-1}$\\
$\log c$ & $(r_n^{-1},\ (2L_n-B_0)r_n^{-1},\ 0,\ 0)$ & $r_n^{-2}$\\
$\log l$ & $(0,\ 0,\ -(2a_0b_n)^{-1},\ 0)$ & $b_n^{-2}$\\
$\log\sigma^2$
 & $(r_n^{-1},\ (2L_n-B_0-\log a_0)r_n^{-1},\ -\nu_0/(a_0b_n),\ 0)$
 & $b_n^{-2}$\\\midrule
$c$ & $c_0(r_n^{-1},\ (2L_n-B_0)r_n^{-1},\ 0,\ 0)$
 & $(1+L_n)^2r_n^{-2}$\\
$l$ & $(0,\ 0,\ -l_0/(2a_0b_n),\ 0)$ & $b_n^{-2}$\\
$\sigma^2$
 & $\sigma_0^2(r_n^{-1},\ (2L_n-B_0-\log a_0)r_n^{-1},\ -\nu_0/(a_0b_n),\ 0)$
 & $b_n^{-2}$\\\bottomrule
\end{tabular}
\end{table}
\FloatBarrier

The complete transformed Hessians, evaluated at $\vartheta=(\eta,\nu,a,t)$, are as follows. In these formulas, $r_n,b_n,L_n$, and $E_n$ remain fixed at the truth. Put
\[
 q_c=2L_n-B(\nu),\qquad q_\sigma=2L_n-B(\nu)-\log a.
\]
For the original coordinates and the logarithmic coordinates,
\[
 \begin{aligned}
 H_{\eta,n}&=H_{\nu,n}=H_{a,n}=H_{t,n}=0_{4\times4},\\
 H_{\log\nu,n}&=\operatorname{diag}
       \left(0,-\frac1{\nu^2r_n^2},0,0\right),\\
 H_{\log a,n}&=\operatorname{diag}
       \left(0,0,-\frac1{a^2b_n^2},0\right),\\
 H_{\log t,n}&=\operatorname{diag}
       \left(0,0,0,-\frac1{t^2n}\right),\\
 H_{\log c,n}&=\operatorname{diag}
       \left(0,-\frac{B'(\nu)}{r_n^2},0,0\right),\\
 H_{\log l,n}&=\operatorname{diag}
       \left(0,0,\frac1{2a^2b_n^2},0\right),\\
 H_{\log\sigma^2,n}
 &=\begin{pmatrix}
 0&0&0&0\\
 0&-B'(\nu)/r_n^2&-1/(ar_nb_n)&0\\
 0&-1/(ar_nb_n)&\nu/(a^2b_n^2)&0\\
 0&0&0&0
 \end{pmatrix}.
 \end{aligned}
\]
For the original-scale variance--lengthscale coordinates,
\[
 H_{c,n}=\frac{c}{r_n^2}
 \begin{pmatrix}
 1&q_c&0&0\\
 q_c&q_c^2-B'(\nu)&0&0\\
 0&0&0&0\\
 0&0&0&0
 \end{pmatrix},
 \qquad
 H_{l,n}=\operatorname{diag}
       \left(0,0,\frac{3l}{4a^2b_n^2},0\right),
\]
\[
 H_{\sigma^2,n}=\sigma^2
 \begin{pmatrix}
 r_n^{-2}&q_\sigma r_n^{-2}&-\nu/(ar_nb_n)&0\\
 q_\sigma r_n^{-2}&\{q_\sigma^2-B'(\nu)\}/r_n^2
       &-(1+\nu q_\sigma)/(ar_nb_n)&0\\
 -\nu/(ar_nb_n)&-(1+\nu q_\sigma)/(ar_nb_n)
       &\nu(\nu+1)/(a^2b_n^2)&0\\
 0&0&0&0
 \end{pmatrix}.
\]

The rates in Table~\ref{tab:rates} follow directly from these rows and Hessians, using $b_n(1+L_n)/r_n\to0$ when $d\ge4$. More precisely, for each listed rate $\rho_{g,n}$,
$\|b_{g,n}\|\asymp\rho_{g,n}$ and $\epsilon_{g,n}=o(\rho_{g,n})$. Hence
\[
 \|g(\widehat\vartheta_n)-g(\vartheta_0)\|_{L^p}
 =O(\rho_{g,n}),\qquad 1\le p<\infty,
\]
and the second-moment comparison gives
\begin{equation}\label{eq:coordinate-mse-expansion}
 \E|g(\widehat\vartheta_n)-g(\vartheta_0)|^2
 =b_{g,n}V_nb_{g,n}^\T+o(\rho_{g,n}^2)
 \asymp\rho_{g,n}^2.
\end{equation}
The centered amplitude error is exactly $r_n^{-1}\widehat h_{n,u}$, so it has rate $r_n^{-1}$ as well.

For random locations, the same proof applies conditionally on the good designs of Lemma~\ref{lem:randomgood}, with uniform constants, eventually almost surely. The unconditional moment bounds of Theorem~\ref{thm:normalized} give the unconditional upper orders.

To obtain the unconditional lower orders, let $\mathcal G_n$ be the good-design event and put $D_{g,n}=g(\widehat\vartheta_n)-g(\vartheta_0)$. The conditional version of \eqref{eq:coordinate-mse-expansion} gives, uniformly over the designs in $\mathcal G_n$,
\[
 \E\bigl(|D_{g,n}|^2\mid x_{n,1:n}\bigr)
 =b_{g,n}V_n(x_{n,1:n})b_{g,n}^\T+o(\rho_{g,n}^2).
\]
Here $V_n(x_{n,1:n})\succeq cI$ and $\|b_{g,n}\|\asymp\rho_{g,n}$ uniformly, so the leading term is at least $c_0\rho_{g,n}^2$. For all sufficiently large $n$, the absolute value of the uniform remainder is at most $c_0\rho_{g,n}^2/2$. Consequently, on $\mathcal G_n$,
\[
 \E\bigl(|D_{g,n}|^2\mid x_{n,1:n}\bigr)
 \ge\frac{c_0}{2}\rho_{g,n}^2.
\]
The rate $\rho_{g,n}$ depends on $n$ and the truth, but not on the realized locations. Nonnegativity and iterated expectation therefore give
\[
 \begin{aligned}
 \E|D_{g,n}|^2
 &=\E_X\!\left[\E\bigl(|D_{g,n}|^2\mid x_{n,1:n}\bigr)\right]\\
 &\ge\E_X\!\left[\mathbf1_{\mathcal G_n}
          \E\bigl(|D_{g,n}|^2\mid x_{n,1:n}\bigr)\right]\\
 &\ge\frac{c_0}{2}\rho_{g,n}^2P_X(\mathcal G_n)
 \ge\frac{c_0}{4}\rho_{g,n}^2
 \end{aligned}
\]
for all sufficiently large $n$, since $P_X(\mathcal G_n)\to1$. Together with the upper bound, this proves $\E|D_{g,n}|^2\asymp\rho_{g,n}^2$, uniformly.
\qed

\subsection{Proof of Theorem~\ref{thm:minimax}}\label{app:minimaxproof}
Corollary~\ref{cor:rates} already gives the pointwise mean squared error comparison for the MLE. Since each rate in Table~\ref{tab:rates} is nondecreasing in $\nu$, it also gives
\[
 \inf_{T_n}\sup_{\theta\in\mathcal K,\,P\in\mathcal P}
       \E_{\theta,P}|T_n-g(\theta)|^2
 \le \sup_{\theta\in\mathcal K,\,P\in\mathcal P}
       \E_{\theta,P}|g(\widehat\theta_n)-g(\theta)|^2
 \le C\rho_{g,n}(\bar\nu)^2.
\]
For the reverse inequality, we construct two parameter values whose $g$ values differ by the asserted rate, but whose observation laws are too close for a test to distinguish reliably. It is enough to use the standard Gaussian member of $\mathcal P$, since restricting the noise class can only decrease the supremum risk.

\noindent\textbf{1. Relating estimation error to two Gaussian laws.}
First fix the observation locations. For two parameter values $\theta_0,\theta_1$, write $G_j=N(0,C_n(\theta_j))$, and let $\E_j$ denote expectation under $G_j$. Put $\Delta=|g(\theta_1)-g(\theta_0)|$. Given any estimator $T_n$, choose the hypothesis whose $g$ value is closer to $T_n$; denote this test by $\phi_n\in\{0,1\}$. Whenever the test chooses the wrong hypothesis, the estimation error is at least $\Delta/2$. Hence
\[
 \begin{aligned}
 \max_{j=0,1}\E_j|T_n-g(\theta_j)|^2
 &\ge \frac12\sum_{j=0}^1\E_j|T_n-g(\theta_j)|^2\\
 &\ge \frac{\Delta^2}{8}
       \{G_0(\phi_n=1)+G_1(\phi_n=0)\}\\
 &\ge \frac{\Delta^2}{8}\{1-\TV(G_0,G_1)\}.
 \end{aligned}
\]
The last inequality follows from $G_0(\phi_n=1)+G_1(\phi_n=0) =1-\{G_1(\phi_n=1)-G_0(\phi_n=1)\}$ and the definition of total variation. In particular,
\begin{equation}\label{eq:minimax-two-point}
 \TV(G_0,G_1)\le\tfrac12
 \quad\Longrightarrow\quad
 \max_{j=0,1}\E_j|T_n-g(\theta_j)|^2\ge\frac{\Delta^2}{16}.
\end{equation}

To ensure the total variation bound, take $\vartheta_1=\vartheta_0+E_nh$ with $\|h\|\le1$, holding $a=a_0$ fixed when $d\le3$. All normalizations in this construction are evaluated at the baseline $\theta_0$. Write
\[
 M(h)=C_n(\theta_0)^{-1/2}C_n(\theta_1)C_n(\theta_0)^{-1/2},
\]
using $M(h,a_0)$ in the notation of Lemma~\ref{lem:derivatives} when $d\le3$. Parts \ref{it:covariance-basic} and \ref{it:covariance-local} of that lemma give $cI\preceq M(h)\preceq CI$ and
\[
 \|M(h)-I\|_\Frob
 =\left\|\int_0^1DM(uh)[h]\,du\right\|_\Frob
 \le C\|h\|.
\]
The Gaussian Kullback--Leibler formula and \eqref{eq:local-contrast-comparison} give
\[
 \KL(G_0,G_1)
 =\tfrac12\{\tr(M(h)^{-1})-n+\log\det M(h)\}
 \asymp\|M(h)-I\|_\Frob^2\le C\|h\|^2.
\]
By Pinsker's inequality, $\TV(G_0,G_1)\le\{\KL(G_0,G_1)/2\}^{1/2}$. Choose a fixed $\alpha_0\in(0,1]$ small enough that $\|h\|\le\alpha_0$ implies $\TV(G_0,G_1)\le1/2$.
These bounds are uniform over baseline parameters in the compact set and over all locations in the bounded observation region.

\noindent\textbf{2. Separating each parameter at its stated rate.}
The derivative rows and Hessians in the proof of Corollary~\ref{cor:rates} give
\[
 \|Dg(\vartheta_0)E_n\|\asymp\rho_{g,n}(\nu_0),\qquad
 \sup_\vartheta\|E_n^\T\nabla^2g(\vartheta)E_n\|_{\mathrm{op}}
       \le C\epsilon_{g,n},\qquad
 \epsilon_{g,n}=o\{\rho_{g,n}(\nu_0)\}.
\]
In low dimension, derivatives are taken in the identifiable coordinates, with range fixed. Choose any fixed $0<\alpha\le\alpha_0$ and set
\[
 h=\alpha\,
    \frac{(Dg(\vartheta_0)E_n)^\T}
         {\|Dg(\vartheta_0)E_n\|}.
\]
Then $\|h\|=\alpha$, and the first-order change in $g$ is $Dg(\vartheta_0)E_nh=\alpha\|Dg(\vartheta_0)E_n\|$. Applying the same Taylor formula as in \eqref{eq:coordinate-taylor-recipe}, now at this deterministic increment, gives
\[
 \begin{aligned}
 g(\theta_1)-g(\theta_0)
 &=\alpha\|Dg(\vartheta_0)E_n\|
       +O(\alpha^2\epsilon_{g,n}),\\
 |g(\theta_1)-g(\theta_0)|
 &\ge c\alpha\rho_{g,n}(\nu_0)
 \end{aligned}
\]
for all sufficiently large $n$. The derivative row thus determines the required direction for every coordinate in Table~\ref{tab:rates}, using one calculation.
Combining this separation with \eqref{eq:minimax-two-point} gives a lower bound $c\alpha^2\rho_{g,n}(\nu_0)^2$ for the maximum risk at these two parameter values, for every estimator.

\noindent\textbf{3. Placing the alternatives in the required sets.}
For the global minimax bound, keep all baseline coordinates other than smoothness at fixed interior values of the rectangle $\mathcal K$, and put
\[
 \nu_{0,n}=\bar\nu-\frac{2\alpha}{r_n(\bar\nu)},
 \qquad r_n(\nu)=n^{d/(4\nu+2d)}.
\]
Use Step~2 with this baseline and with $E_n$ evaluated there. The smoothness displacement satisfies
\[
 |\nu_1-\nu_{0,n}|
 =\frac{|h_v|}{r_n(\nu_{0,n})}
 \le\frac{\alpha}{r_n(\nu_{0,n})}
 \le\frac{\alpha}{r_n(\bar\nu)}.
\]
Thus neither hypothesis crosses the upper smoothness face of $\mathcal K$, and both stay above its lower face for large $n$. The other coordinates remain in $\mathcal K$ because their baseline values have fixed positive margins and $\|E_nh\|\le\alpha\|E_n\|_{\mathrm{op}}\to0$. Moreover,
\[
 (\bar\nu-\nu_{0,n})\log n
       =\frac{2\alpha\log n}{r_n(\bar\nu)}\longrightarrow0.
\]
The explicit rate formulas consequently give $\rho_{g,n}(\nu_{0,n})/\rho_{g,n}(\bar\nu)\to1$. Since both hypotheses belong to $\mathcal K$, Step~2 proves
\[
 \inf_{T_n}\sup_{\theta\in\mathcal K,\,P\in\mathcal P}
       \E_{\theta,P}|T_n-g(\theta)|^2
       \ge c\rho_{g,n}(\bar\nu)^2.
\]
Together with the upper bound at the start of the proof, this establishes the global minimax comparison and the matching worst-case risk of the MLE.

For the local assertion, fix an interior truth $\theta_0$ and $H>0$. Apply Steps~1 and~2 on a compact neighborhood of this truth, taking $\alpha=\min\{\alpha_0,H/2\}$.
Both $0$ and $h$ lie in the normalized ball of radius $H$. Since $\theta_0$ is interior and $E_nh\to0$, the two parameter values belong to $\Theta$ for all sufficiently large $n$; range is fixed at $a_0$ when $d\le3$. They therefore belong to
$\mathcal U_{n,H}(\theta_0)$. The same two-point bound gives
\[
 \inf_{T_n}\sup_{\theta\in\mathcal U_{n,H}(\theta_0),\,P\in\mathcal P}
       \E_{\theta,P}|T_n-g(\theta)|^2
       \ge c_H\rho_{g,n}(\nu_0)^2,
\]
where $c_H>0$ is independent of $n$ and of the estimator. Dividing by the squared rate and taking the lower limit proves \eqref{eq:localminimax}.

\noindent\textbf{4. Random locations.}
The covariance upper bounds used in Step~1 hold for arbitrary sites in the bounded region. The lower-bound argument therefore holds conditionally for every fixed realization of the locations, with the same constants. The alternative in Step~2 depends on the baseline and $n$, but not on the realized locations.

For the unconditional minimax lower bound, write $G_j^x$ for the conditional Gaussian response law and $G_j^{\mathrm{joint}}(dx,dy)=P_X(dx)G_j^x(dy)$ for the joint
law of locations and observations. The location distribution is the same under both parameters, so
\[
 \begin{aligned}
 \KL(G_0^{\mathrm{joint}},G_1^{\mathrm{joint}})
 =\E_X\!\left[\KL(G_0^X,G_1^X)\right]\le C\|h\|^2.
 \end{aligned}
\]
Pinsker's inequality and the testing argument in Step~1 now apply to any estimator based on both the locations and observations. The parameter separations are unchanged, so Steps~2 and~3 give the same global and local lower bounds under the joint law.
The conditional and unconditional MLE risk comparisons, including their uniformity over the noise class, are already supplied by Corollary~\ref{cor:rates}. This completes all assertions. \qed

\subsection{Proof of Corollary~\ref{cor:log}}\label{app:logproof}
We first work at fixed locations and use the forward map $\zeta=g(\vartheta)$. Table~\ref{tab:coordinate-expansions} gives the rows $Dg_j(\vartheta_0)E_n$ and their remainder orders. In particular,
\[
 [\Sigma_n^\zeta]_{jj}
 =Dg_j(\vartheta_0)E_nV_nE_n^{\T}Dg_j(\vartheta_0)^\T
 \asymp\|Dg_j(\vartheta_0)E_n\|^2\asymp\rho_{j,n}^2.
\]

\noindent\textbf{Comparison with Gaussian information.}
The chain rule, with the nuisance fixed when $d\le3$, gives
\[
 \begin{aligned}
 \mathcal I_n^\zeta
 =Dg(\vartheta_0)^{-\T}\mathcal I_n
                    Dg(\vartheta_0)^{-1},~
 Dg(\vartheta_0)E_nJ_n^{-1}E_n^{\T}Dg(\vartheta_0)^\T
 =(\mathcal I_n^\zeta)^{-1},~
 Dg(\vartheta_0)E_ne_t=n^{-1/2}q_\zeta(t_0)e_t.
 \end{aligned}
\]
These identities and \eqref{eq:kurtosis-correction} prove \eqref{eq:reporting-information} and \eqref{eq:reporting-information-correction}. For Gaussian errors,
$V_n=J_n^{-1}$ exactly, giving the stated exact equality. In low dimension, Lemma~\ref{lem:studentization} gives
\[
 E_n^{-1}\left\{[(\mathcal I_n^{(\vartheta,a)})^{-1}]_{\vartheta\vartheta}
                 -\mathcal I_n^{-1}\right\}E_n^{-\T}=o(1).
\]
The normalized mixed block vanishes and the nuisance information is bounded below, so the Schur correction vanishes. Changing $a$ to any listed nuisance coordinate does not change that Schur complement. Multiplication by $Dg(\vartheta_0)$ proves the full-inverse replacement.

\noindent\textbf{Standardized expansions.}
For every marginal coordinate, the Taylor table in the proof of Corollary~\ref{cor:rates} gives, for each fixed $1\le p<\infty$,
\[
 \frac{\widehat\zeta_{n,j}-\zeta_{0,j}}
      {\sqrt{[\Sigma_n^\zeta]_{jj}}}
 =\frac{Dg_j(\vartheta_0)E_nV_n^{1/2}}{\sqrt{[\Sigma_n^\zeta]_{jj}}}
       V_n^{-1/2}\widehat h_n+o_{L^p}(1).
\]
The row multiplying $V_n^{-1/2}\widehat h_n$ has norm one.

For a vector covered by the joint assertion, let $\mathcal R_n^\zeta$ be the forward Taylor remainder in Lemma~\ref{lem:logtransform}. Its Hessian table and \eqref{eq:common-chart-remainder-bound} give $\mathcal R_n^\zeta=o_{L^p}(1)$, and the exact relation is
\[
 T_n^\zeta\coloneqq E_n^{-1}Dg(\vartheta_0)^{-1}(\widehat\zeta_n-\zeta_0)
 =\widehat h_n+\mathcal R_n^\zeta.
\]
Consequently
\[
 \begin{aligned}
 (\Sigma_n^\zeta)^{-1/2}(\widehat\zeta_n-\zeta_0)
 &=O_n^\zeta V_n^{-1/2}\widehat h_n+o_{L^p}(1),\\
 O_n^\zeta&=(\Sigma_n^\zeta)^{-1/2}Dg(\vartheta_0)E_nV_n^{1/2},
 \qquad O_n^\zeta(O_n^\zeta)^\T=I.
 \end{aligned}
\]
Theorem~\ref{thm:normalized} now proves the marginal and joint normal limits. The coefficient rows and orthogonal matrices may vary with $n$: compactness of the unit sphere and the orthogonal group, followed by rotational invariance of the standard normal law, gives the conclusion along every sequence of truths and noise laws. This also proves uniformity.

\noindent\textbf{Tails and moments.}
The same derivative rows give the deterministic marginal bound
\[
 \frac{|\widehat\zeta_{n,j}-\zeta_{0,j}|}
      {\sqrt{[\Sigma_n^\zeta]_{jj}}}
 \le C\|\widehat h_n\|.
\]
For the joint bound, regard $h$ as a deterministic normalized parameter increment and use the map
\[
 T_n^\zeta(h)=E_n^{-1}Dg(\vartheta_0)^{-1}
                  \{g(\vartheta_0+E_nh)-g(\vartheta_0)\}.
\]
Since $T_n^\zeta(0)=0$ and $DT_n^\zeta(0)=I$, the same second-order Taylor formula applies to every coordinate:
\[
 T_{n,i}^\zeta(h)-h_i
 =\int_0^1(1-u)h^\T H_{i,n}^\zeta(uh)h\,du,
 \qquad 1\le i\le k_d,
\]
where the Hessians $H_{i,n}^\zeta$ are listed in Lemma~\ref{lem:logtransform}. The first-coordinate bounds are recorded in Table~\ref{tab:joint-tail-remainders}. In its
log-scale variance--lengthscale row, we use
\[
 |h_v\log(a/a_0)|\le C|h_v|\,|a-a_0|
                 =C\frac{|h_vh_a|}{b_n}.
\]

\begin{table}[!htbp]
\centering\small
\renewcommand{\arraystretch}{1.4}
\caption{Deterministic second-order bounds for $|T_{n,1}^\zeta(h)-h_u|$, up to a uniform constant, on the compact parameter set.}
\label{tab:joint-tail-remainders}
\begin{tabular}{@{}ll@{}}\toprule
Parameterization & Bound\\\midrule
Spectral & $0$\\
Spectral, log scale & $L_nh_v^2/r_n$\\
Microergodic & $(1+L_n)^2\|h\|^2/r_n$\\
Microergodic, log scale & $(1+L_n)h_v^2/r_n$\\
Variance--lengthscale, log scale
 & $|h_vh_a|/b_n+(1+L_n)h_v^2/r_n$\\
Variance--lengthscale, $d>8$
 & $(r_n/b_n^2+b_n^{-1})\|h\|^2$\\\bottomrule
\end{tabular}
\end{table}
\FloatBarrier

For every coordinate after the first, its Hessian is diagonal with at most one nonzero entry. The entries for the smoothness, range and nugget coordinates are bounded by $C/r_n$, $C/b_n$ and $C/\sqrt n$, respectively. The same Taylor formula therefore
gives
\[
 \sum_{i=2}^{k_d}|T_{n,i}^\zeta(h)-h_i|
 \le C\left(\frac{h_v^2}{r_n}+\frac{h_a^2}{b_n}
                         +\frac{h_t^2}{\sqrt n}\right),
\]
with the range term omitted when $d\le3$. A coordinate that is unchanged by the transformation has zero remainder. On the compact parameter set,
\[
 \frac{|h_v|}{r_n}=|\nu-\nu_0|\le C,\qquad
 \frac{|h_a|}{b_n}=|a-a_0|\le C,\qquad
 \frac{|h_t|}{\sqrt n}=|t-t_0|\le C.
\]
Consequently the three quadratic terms are bounded by $C|h_v|$, $C|h_a|$ and $C|h_t|$, and hence
\[
 \sum_{i=2}^{k_d}|T_{n,i}^\zeta(h)|\le C\|h\|.
\]
The mixed term $|h_vh_a|/b_n$ in the first-coordinate table is likewise at most $C|h_v|$. The first coordinate also has terms with different normalizations, such as
$r_nh_a^2/b_n^2$, which do not reduce to a uniform linear bound in this way. Retaining those quadratic terms yields, for every joint vector under consideration,
\[
 \begin{gathered}
 \|T_n^\zeta(h)\|\le C(\|h\|+a_n^\zeta\|h\|^2),\\
 a_n^\zeta=
 \begin{cases}
 r_n/b_n^2+b_n^{-1},&\text{original-scale variance--lengthscale vector, }d>8,\\
 (1+L_n)^2/r_n,&\text{all other joint vectors}.
 \end{cases}
 \end{gathered}
\]
In every case $a_n^\zeta\le Cn^{-\gamma}$ for some uniform $\gamma>0$. Since
$\|(\Sigma_n^\zeta)^{-1/2}(\widehat\zeta_n-\zeta_0)\| =\|V_n^{-1/2}T_n^\zeta(\widehat h_n)\|$ and $V_n^{-1/2}$ is bounded, on $\{\|\widehat h_n\|\le n^{\gamma/2}\}$ the joint standardized error is at most $C\|\widehat h_n\|$. Thus \eqref{eq:tails} gives
\[
 \begin{aligned}
 P\!\left(\|(\Sigma_n^\zeta)^{-1/2}
               (\widehat\zeta_n-\zeta_0)\|>H\right)
 &\le P(\|\widehat h_n\|>H/C)
       +P(\|\widehat h_n\|>n^{\gamma/2})\\
 &\le Ce^{-cH^2}+Ce^{-cn^\gamma}+Ce^{-cn^\delta}.
 \end{aligned}
\]
This is the asserted tail bound, with exponent $\min(\gamma,\delta)$. The marginal bound follows directly from \eqref{eq:tails}.

The deterministic bounds above and the moment bounds for $\widehat h_n$ give uniform bounds for every fixed moment of the standardized errors. Uniform integrability therefore upgrades their convergence in distribution to convergence of expectations
for every fixed continuous function of polynomial growth. In particular, \eqref{eq:coordinate-mse-expansion} gives
\[
 \frac{\E(\widehat\zeta_{n,j}-\zeta_{0,j})^2}
      {[\Sigma_n^\zeta]_{jj}}\longrightarrow1.
\]

\noindent\textbf{Random locations.}
All preceding comparisons hold conditionally on the good designs of Lemma~\ref{lem:randomgood}, with the same constants, eventually almost surely. Averaging conditional bounded-Lipschitz convergence and using $P_X(\mathcal G_n^c)\to0$ proves unconditional convergence in distribution. For unconditional tails and moments, consider any design on which $J_n$ is nonsingular. The bounds $K_n\succeq cJ_n$ and $J_n\preceq CI$ give
\[
 V_n^{-1}=J_nK_n^{-1}J_n\preceq CJ_n\preceq C'I.
\]
In particular, $V_n\succeq c'I$. For the marginal normalizer,
put $b_{j,n}=Dg_j(\vartheta_0)E_n$. Then
\[
 [\Sigma_n^\zeta]_{jj}
 =b_{j,n}V_nb_{j,n}^{\top}
 \ge c'\|b_{j,n}\|^2
 \asymp \rho_{j,n}^2.
\]
This lower bound is all that the marginal tail estimate uses. For the joint estimate, the preceding deterministic bound on $T_n^\zeta(\widehat h_n)$ applies since $\|V_n^{-1/2}\|_{\mathrm{op}}\le C$. Thus both estimates follow from the unconditional tail and moment bounds for $\widehat h_n$. On designs where $J_n$ is singular, set the marginal and joint standardized reporting errors to zero, as in Theorem~\ref{thm:normalized}. Uniform integrability then proves the unconditional standardized moment assertions as well. \qed

\subsection{Proof of Theorem~\ref{thm:preliminary}}
\label{app:preliminaryproof}
Work first at fixed sites satisfying Assumption~\ref{ass:coverage}. Lemma~\ref{lem:preliminary} places every global MLE in $\mathcal N_n$ with exceptional probability at most $Ce^{-cn^{\delta_1}}$, uniformly over $\mathcal K\times\mathcal P$.
When $\widehat\theta_n\in\mathcal N_n$, with $q_n=n^{-\delta}$,
\[
 |\widehat v|\le q_n,\qquad
 |\widehat t_n-t_0|\le q_n,\qquad
 |\widehat u|\le q_n\log n.
\]
Since $\widehat\eta_n-\eta_0=\widehat u+2L_n\widehat v$ and $L_n=(\log n)/(2\nu_0+d)$, it follows that
\[
 \frac{|\widehat\eta_n-\eta_0|}{\log n}
 +|\widehat\nu_n-\nu_0|+|\widehat t_n-t_0|
 \le\left(3+\frac{2}{2\nu_-+d}\right)n^{-\delta}.
\]
This proves \eqref{eq:preliminary-global-bound} and the stated uniform stochastic rates.

For the moment bounds, denote the left side of the preceding display by $D_n$. Compactness of $\Theta$ gives $D_n\le C$ for $n\ge2$, even on the exceptional event. Thus, for each fixed $1\le p<\infty$,
\[
 \begin{aligned}
 \E_{\theta_0,P}D_n^p
 &=\E_{\theta_0,P}[D_n^p\mathbf1_{\{\widehat\theta_n\in\mathcal N_n\}}]
   +\E_{\theta_0,P}[D_n^p\mathbf1_{\{\widehat\theta_n\notin\mathcal N_n\}}]\\
 &\le C_p n^{-p\delta}+C_p e^{-cn^{\delta_1}}
 \le C_p' n^{-p\delta}.
 \end{aligned}
\]
Each parameter error is bounded by $D_n$, with an additional factor $\log n$ for $\eta$, proving the $L^p$ assertions.

The identity $\log c=\eta-\log A_d(\nu)$ and bounded derivatives on the compact parameter intervals give
\[
 \begin{aligned}
 |\log\widehat c_n-\log c_0|
   &\le|\widehat\eta_n-\eta_0|+C|\widehat\nu_n-\nu_0|,\\
 |\log\widehat\nu_n-\log\nu_0|
   +|\log\widehat t_n-\log t_0|
   &\le C\{|\widehat\nu_n-\nu_0|+|\widehat t_n-t_0|\}.
 \end{aligned}
\]
The exponential map also has bounded derivative on the possible $\log c$ interval, so
$|\widehat c_n-c_0|\le C|\log\widehat c_n-\log c_0|$. These inequalities prove the transformed rates in probability and in every stated moment.

For iid sites, Lemma~\ref{lem:randomgood} provides events $\mathcal G_n$ on which one fixed polynomial-coverage bound holds with deterministic constants, and $P_X(\mathcal G_n^c)\le Ce^{-cn^\rho}$ for some $\rho>0$. The preceding proof therefore gives the same conditional bounds on every $\mathcal G_n$. Since the failure probabilities are
summable, almost every site sequence belongs to $\mathcal G_n$ eventually. This proves the conditional assertion. For the joint law, splitting over $\mathcal G_n$ gives
\[
 \begin{aligned}
 P_{\theta_0,P}(D_n>Cn^{-\delta})
 &\le \E_X\!\left[
   \mathbf1_{\mathcal G_n}
   P_{\theta_0,P}(D_n>Cn^{-\delta}\mid x_{n,1:n})\right]
   +P_X(\mathcal G_n^c)\\
 &\le Ce^{-cn^{\delta_1}}+Ce^{-cn^\rho}
 \le C'e^{-c'n^{\min(\delta_1,\rho)}}.
 \end{aligned}
\]
Likewise $\E D_n^p\le C_p n^{-p\delta}+C_pP_X(\mathcal G_n^c) \le C_p'n^{-p\delta}$. The same transformation inequalities finish the unconditional assertions. \qed

\subsection{Proof of Corollary~\ref{cor:uncertainty}}
\label{app:uncertaintyproof}
\noindent\textbf{Covariance estimation in spectral coordinates.}
Work first in the core spectral coordinates $\vartheta$. In low dimension use the nuisance $\rho$ already present in the selected row of Table~\ref{tab:parameterizations}; in high dimension omit it. Lemmas~\ref{lem:studentization} and~\ref{lem:range-curvature} apply to each of $\rho=a,\log a,l,\log l$. If $\widehat\Sigma_n^{\vartheta,\rho}$ is the identifiable block of the inverse full curvature, plus $\widetilde\kappa_{4,n}e_te_t^\top/n$,
those lemmas give
\begin{equation}\label{eq:relative-reporting-covariance}
 E_n^{-1}\widehat\Sigma_n^{\vartheta,\rho}E_n^{-\top}
       =V_n+o_P(1).
\end{equation}
This holds for full Gaussian Fisher information and full negative observed Hessian, including a boundary nuisance fit. In particular the corrected identifiable covariance is positive definite with probability tending to one. The notation suppresses $\rho$ in high dimension.

\noindent\textbf{Transformation into the fitted coordinates.}
Let $\zeta=g(\vartheta)$ be the selected forward map. In low dimension, $g$ is independent of the nuisance $\rho$. Gaussian Fisher information transforms exactly by congruence. In high dimension, the fit is interior with probability tending to one, so its score is zero and the observed Hessian also transforms exactly. In low dimension, a moving nuisance boundary can prevent an identifiable score from vanishing; we retain that term below.

In low dimension, let $\widehat H^{\vartheta,\rho}$ and $\widehat H^{\zeta,\rho}$ denote the full negative Hessians. View the same criterion in the two coordinate systems, so that $\ell_n(\vartheta,\rho)=\ell_n(g(\vartheta),\rho)$. Since $g$ does not depend on $\rho$, the first derivatives are
\[
 \partial_{\vartheta_i}\ell_n
 =\sum_j\partial_{\zeta_j}\ell_n\,\partial_{\vartheta_i}g_j,
 \qquad
 \nabla_\vartheta\ell_n=Dg(\vartheta)^\T\nabla_\zeta\ell_n.
\]
Differentiating the first identity once more gives
\[
 \begin{aligned}
 \partial_{\vartheta_k}\partial_{\vartheta_i}\ell_n
 &=\sum_{j,l}\partial_{\zeta_j\zeta_l}^2\ell_n\,
       \partial_{\vartheta_i}g_j\,\partial_{\vartheta_k}g_l
   +\sum_j\partial_{\zeta_j}\ell_n\,
       \partial_{\vartheta_i\vartheta_k}^2g_j,\\
 \partial_\rho\partial_{\vartheta_i}\ell_n
 &=\sum_j\partial_{\zeta_j\rho}^2\ell_n\,
       \partial_{\vartheta_i}g_j.
 \end{aligned}
\]
The derivatives of $\ell_n$ on the right are evaluated at $(g(\vartheta),\rho)$. Taking negatives and evaluating at the fit therefore yields
\[
 \begin{aligned}
 \widehat H^{\vartheta,\rho}_{\vartheta\vartheta}
 &=Dg(\widehat\vartheta_n)^\T
       \widehat H^{\zeta,\rho}_{\zeta\zeta}
       Dg(\widehat\vartheta_n)-\mathcal T_n,\\
 \mathcal T_n
 &=\sum_j\partial_{\zeta_j}\ell_n(\widehat\theta_n)
                   D^2g_j(\widehat\vartheta_n),\\
 \widehat H^{\vartheta,\rho}_{\vartheta\rho}
 &=Dg(\widehat\vartheta_n)^\T
       \widehat H^{\zeta,\rho}_{\zeta\rho},\qquad
 \widehat H^{\vartheta,\rho}_{\rho\rho}
 =\widehat H^{\zeta,\rho}_{\rho\rho}.
 \end{aligned}
\]
Substituting these identities into the Schur complement gives
\[
 \begin{aligned}
 &\widehat H^{\zeta,\rho}_{\zeta\zeta}
  -\widehat H^{\zeta,\rho}_{\zeta\rho}
       (\widehat H^{\zeta,\rho}_{\rho\rho})^{-1}
       \widehat H^{\zeta,\rho}_{\rho\zeta}=Dg(\widehat\vartheta_n)^{-\top}
       \Bigl\{\widehat H^{\vartheta,\rho}_{\vartheta\vartheta}
        -\widehat H^{\vartheta,\rho}_{\vartheta\rho}
          (\widehat H^{\vartheta,\rho}_{\rho\rho})^{-1}
          \widehat H^{\vartheta,\rho}_{\rho\vartheta}
        +\mathcal T_n\Bigr\}Dg(\widehat\vartheta_n)^{-1}.
 \end{aligned}
\]
The inverses in this identity exist with probability tending to one by Lemma~\ref{lem:range-curvature}.

Put $\widehat s_h=E_n^\T\nabla_\vartheta\ell_n(\widehat\theta_n)$. Use $\mathcal L_n(h,a)$ and $s_n(a)$ from Lemma~\ref{lem:nuisance}. The score representation
$\Delta_{n,j}=\tfrac12Q(B_j)$, the variance bound in Lemma~\ref{lem:gaussian}, and $J_n\preceq CI$ from Lemma~\ref{lem:information} give
\[
 \E\|\Delta_n\|^2
 \le C\sum_j\|B_j\|_{\mathrm F}^2
 =2C\operatorname{tr}J_n\le C'.
\]
Thus $\Delta_n=O_P(1)$, and Lemma~\ref{lem:nuisance} gives $\sup_{a\in I_*}\|s_n(a)\|=O_P(1)$. Lemma~\ref{lem:studentization}\ref{it:studentization-local},
specifically \eqref{eq:localIobserved}, gives, for every fixed $H$,
\[
 \sup_{\|h\|\le H,\ a\in I_*}
       \|\nabla_h^2\mathcal L_n(h,a)\|_{\mathrm{op}}=O_P(1).
\]
Since $\widehat h_n=O_P(1)$ by Lemma~\ref{lem:tails}, integration along the segment from $0$ to $\widehat h_n$ yields
\[
 \widehat s_h
 =s_n(\widehat a_n)
   +\int_0^1\nabla_h^2\mathcal L_n(s\widehat h_n,\widehat a_n)
          \widehat h_n\,ds
 =O_P(1).
\]
These bounds apply on the enclosing parameter rectangle, so the segment need not lie in the optimization set.

The score chain rule at the fit gives
\[
 \widehat s_h
 =E_n^\T Dg(\widehat\vartheta_n)^\T
       \nabla_\zeta\ell_n(\widehat\theta_n),
 \qquad
 \nabla_\zeta\ell_n(\widehat\theta_n)^\T
 =\widehat s_h^\T E_n^{-1}Dg(\widehat\vartheta_n)^{-1}.
\]
Here $D^2g(\vartheta)[x,y]$ is the column vector with entries $x^\T D^2g_j(\vartheta)y$. For unit vectors $u,v$, the definition of $\mathcal T_n$ and the preceding identity therefore give
\[
 \begin{aligned}
 u^\T E_n^\T\mathcal T_nE_nv
 &=\sum_j\partial_{\zeta_j}\ell_n(\widehat\theta_n)\,
       (E_nu)^\T D^2g_j(\widehat\vartheta_n)(E_nv)\\
 &=\nabla_\zeta\ell_n(\widehat\theta_n)^\T
       D^2g(\widehat\vartheta_n)[E_nu,E_nv]\\
 &=\widehat s_h^\T E_n^{-1}Dg(\widehat\vartheta_n)^{-1}
       D^2g(\widehat\vartheta_n)[E_nu,E_nv].
 \end{aligned}
\]
The forward Hessians in Lemma~\ref{lem:logtransform} and its Jacobian comparison imply, for the admissible low-dimensional parameterizations,
\[
 \sup_{\|u\|,\|v\|\le1}
 \left\|E_n^{-1}Dg(\widehat\vartheta_n)^{-1}
             D^2g(\widehat\vartheta_n)[E_nu,E_nv]\right\|
 =O_P\left\{\frac{(1+L_n)^2}{r_n}+\frac1{\sqrt n}\right\}
 =o_P(1).
\]
Indeed, replacing the fitted derivative inverse in this display by $Dg(\vartheta_0)^{-1}$ gives the Hessians displayed there; the two expressions differ by left multiplication by
\[
 \left[E_n^{-1}Dg(\vartheta_0)^{-1}
                 Dg(\widehat\vartheta_n)E_n\right]^{-1}
 =I+o_P(1).
\]
Consequently $\|E_n^\T\mathcal T_nE_n\|_{\mathrm{op}}=o_P(1)$. The normalized core Schur complement is $J_n+o_P(1)$, so inversion preserves this negligible error. For either curvature,
\[
 \begin{aligned}
 \widehat\Sigma_n^\zeta
 &=Dg(\widehat\vartheta_n)
       \widehat\Sigma_n^{\vartheta,\rho}
       Dg(\widehat\vartheta_n)^\T+Dg(\widehat\vartheta_n)E_n\,o_P(1)\,
       E_n^{\T}Dg(\widehat\vartheta_n)^\T,\\
 \Sigma_n^\zeta
 &=Dg(\vartheta_0)E_nV_nE_n^{\T}Dg(\vartheta_0)^\T.
 \end{aligned}
\]
For Gaussian Fisher information, or at an interior stationary fit, the additional congruence remainder is zero. The cumulant term transforms as stated because $Dg(\vartheta)e_t=q_\zeta(t)e_t$. The same argument applies to direct identifiable-block inversion. Thus the covariance formula and positivity remain valid at a nuisance boundary.

\noindent\textbf{Relative consistency and marginal Wald intervals.}
For a vector covered by \eqref{eq:joint-reporting-normal}, Lemma~\ref{lem:logtransform} gives directly
\[
 E_n^{-1}Dg(\vartheta_0)^{-1}
              Dg(\widehat\vartheta_n)E_n=I+o_P(1).
\]
Combining \eqref{eq:relative-reporting-covariance} with the congruence remainder above first gives
\[
 \widehat\Sigma_n^\zeta
 =Dg(\widehat\vartheta_n)E_n\{V_n+o_P(1)\}
        E_n^\T Dg(\widehat\vartheta_n)^\T.
\]
Consequently,
\[
 \begin{aligned}
 &~~~~E_n^{-1}Dg(\vartheta_0)^{-1}\widehat\Sigma_n^\zeta
        Dg(\vartheta_0)^{-\top}E_n^{-\top}\\
 &=\{E_n^{-1}Dg(\vartheta_0)^{-1}
                 Dg(\widehat\vartheta_n)E_n\}
        \{V_n+o_P(1)\}       \{E_n^{-1}Dg(\vartheta_0)^{-1}
                 Dg(\widehat\vartheta_n)E_n\}^\T\\
 &=\{I+o_P(1)\}\{V_n+o_P(1)\}\{I+o_P(1)\}^\T\\
&=V_n+o_P(1).
 \end{aligned}
\]
The same transformation of $\Sigma_n^\zeta$ is exactly $V_n$. Subtracting therefore gives
\[
 E_n^{-1}Dg(\vartheta_0)^{-1}
       (\widehat\Sigma_n^\zeta-\Sigma_n^\zeta)
       Dg(\vartheta_0)^{-\top}E_n^{-\top}=o_P(1).
\]
To pass to symmetric whitening, use $V_n\succeq cI$:
\[
 \begin{aligned}
 \Sigma_n^\zeta
 &=Dg(\vartheta_0)E_nV_nE_n^\T Dg(\vartheta_0)^\T\\
 &\succeq cDg(\vartheta_0)E_nE_n^\T Dg(\vartheta_0)^\T,
\end{aligned}
\]
Whitening this inequality and taking the largest eigenvalue gives
\[
 \|(\Sigma_n^\zeta)^{-1/2}Dg(\vartheta_0)E_n\|_{\mathrm{op}}^2
 \le c^{-1}.
\]
Factoring the covariance error through $Dg(\vartheta_0)E_n$ now proves
\begin{equation}\label{eq:relative-covariance-all-joint}
 \begin{aligned}
 &~~~~\left\|(\Sigma_n^\zeta)^{-1/2}
       (\widehat\Sigma_n^\zeta-\Sigma_n^\zeta)
       (\Sigma_n^\zeta)^{-1/2}\right\|_{\mathrm{op}}\\
 &\le
       \|(\Sigma_n^\zeta)^{-1/2}Dg(\vartheta_0)E_n\|_{\mathrm{op}}^2
       \left\|E_n^{-1}Dg(\vartheta_0)^{-1}
       (\widehat\Sigma_n^\zeta-\Sigma_n^\zeta)
       Dg(\vartheta_0)^{-\top}E_n^{-\top}\right\|_{\mathrm{op}}\\
       &=o_P(1).
 \end{aligned}
\end{equation}

For marginal variances, recall that the row Jacobian is $Dg_j(\vartheta)=\nabla g_j(\vartheta)^\T$ and that $D^2g_j$ is the derivative of the column gradient $\nabla g_j$. Transposing the row difference and integrating this gradient along the segment gives
\[
 \begin{aligned}
 \bigl(\{Dg_j(\widehat\vartheta_n)-Dg_j(\vartheta_0)\}E_n\bigr)^\T&=E_n^\T
       \{\nabla g_j(\widehat\vartheta_n)-\nabla g_j(\vartheta_0)\}\\
 &=\int_0^1 E_n^\T
       D^2g_j(\vartheta_0+sE_n\widehat h_n)E_n
       \widehat h_n\,ds.
 \end{aligned}
\]
The transformed Hessians $E_n^\T D^2g_j(\vartheta)E_n$ displayed in the proof of Corollary~\ref{cor:rates} are bounded by $C\epsilon_{g_j,n}$ throughout this segment.
Transposition does not change the Euclidean norm, so, for every listed coordinate,
\[
 \begin{aligned}
 \|\{Dg_j(\widehat\vartheta_n)-Dg_j(\vartheta_0)\}E_n\|
 &\le\int_0^1
       \|E_n^\T D^2g_j(\vartheta_0+sE_n\widehat h_n)E_n\|_{\mathrm{op}}
       \|\widehat h_n\|\,ds\\
 &\le C\epsilon_{g_j,n}\|\widehat h_n\|
       =o_P(\rho_{j,n}),
 \end{aligned}
\]
using $\widehat h_n=O_P(1)$ and $\epsilon_{g_j,n}=o(\rho_{j,n})$ from Table~\ref{tab:coordinate-expansions}. Therefore the covariance formula above gives
\[
 [\widehat\Sigma_n^\zeta]_{jj}
 =Dg_j(\vartheta_0)E_nV_nE_n^{\T}Dg_j(\vartheta_0)^\T
       +o_P(\rho_{j,n}^2)
 =[\Sigma_n^\zeta]_{jj}+o_P(\rho_{j,n}^2).
\]
Since $[\Sigma_n^\zeta]_{jj}\asymp\rho_{j,n}^2$, the ratio of the estimated to true linearized variance tends to one. This calculation also covers the original-scale variance--lengthscale vector when $4\le d\le8$. Corollary~\ref{cor:log} now gives uniform marginal Wald coverage for every admissible parameterization.

\noindent\textbf{Joint studentization.}
For joint studentization, abbreviate $\Sigma=\Sigma_n^\zeta$ and $\widehat\Sigma=\widehat\Sigma_n^\zeta$, and put $T_n=\widehat\Sigma^{-1/2}\Sigma^{1/2}$ on the event that $\widehat\Sigma$ is positive definite. By
\eqref{eq:relative-covariance-all-joint},
\[
 T_n^\top T_n
 =\{\Sigma^{-1/2}\widehat\Sigma\Sigma^{-1/2}\}^{-1}
 =I+o_P(1).
\]
Thus the smallest and largest singular values of $T_n$ both converge to one in probability. The eigenvalues of $T_n$ itself are positive real numbers, since $T_n$ is similar to the positive definite matrix $\widehat\Sigma^{-1/4}\Sigma^{1/2}\widehat\Sigma^{-1/4}$. Each eigenvalue lies between the smallest and largest singular values, so $\tr(T_n)=k_d+o_P(1)$. Consequently,
\[
 \|T_n-I\|_{\mathrm F}^2
 =\tr(T_n^\top T_n)-2\tr(T_n)+k_d=o_P(1).
\]
The joint conclusion of Corollary~\ref{cor:log} proves \eqref{eq:studentizednormal} and the chi-squared limit.

\noindent\textbf{Transformed regions and random designs.}
A joint region $\mathcal R_n$ in the core coordinates can be mapped to the desired coordinates by the same forward map $g$. Since $g$ is one-to-one, it preserves the coverage event exactly:
\[
 \{\zeta_0\in g(\mathcal R_n)\}
       =\{\vartheta_0\in\mathcal R_n\}.
\]
This proves coverage of transformed joint regions, including those in the original variance--lengthscale coordinates, and of exponentiated marginal intervals. All estimates hold conditionally on good designs with common constants; bounded convergence gives the unconditional coverage statements. \qed

\subsection{Proof of Corollary~\ref{cor:known}}
\label{app:knownproof}
Restrict the likelihood to $\Theta\cap\{\nu=\nu_0\}$. In normalized coordinates this sets the smoothness \emph{increment} $h_v$ to zero. The truth belongs to this slice, so its global maximizer has likelihood at least that of the truth. The mesh bound in the proof of Lemma~\ref{lem:preliminary} rules out all distant candidates with this
property, including those on the slice. The restricted maximizer is therefore in $\mathcal N_n$ with the same preliminary probability bound.

Let $I$ index the free identifiable coordinates of $h$: $(u,t)$ for $d\le3$ and $(u,a,t)$ for $d\ge4$. Differentiating the restricted criterion selects the score $\Delta_{n,I}$, the expected negative Hessian $J_{n,II}$, and the score covariance $K_{n,II}$. In particular, $J_n\asymp I_{k_d}$ gives $J_{n,II}\asymp I_{|I|}$. The tail and Taylor arguments in the proof of Theorem~\ref{thm:normalized} now apply on the slice $h_v=0$. When $d\le3$, fix $H<\infty$. For all sufficiently large $n$, every $(h_u,h_t)$ with $\|(h_u,h_t)\|\le H$ gives an identifiable parameter
\[
 \vartheta=(\eta_0+h_u/r_n,\ \nu_0,\ t_0+h_t/\sqrt n)\in U.
\]
At this point the restricted likelihood is profiled over $I(\vartheta)$. Its endpoints may depend on $(h_u,h_t)$, and \eqref{eq:feasible-endpoint-motion} controls their displacement uniformly on the ball. Hence the normalized identifiable estimation error has expansion $J_{n,II}^{-1}\Delta_{n,I}+o_P(1)$, with covariance benchmark
\[
 V_n^{\mathrm{res}}=J_{n,II}^{-1}K_{n,II}J_{n,II}^{-1}.
\]
The restriction is made before inversion: in general this covariance is not $(V_n)_{II}$. The same tail bound gives the moment assertions. For the minimax lower bounds, assume $N(0,1)\in\mathcal P$ and apply the two-point argument in the proof of
Theorem~\ref{thm:minimax} within the slice $\nu=\nu_0$. For each reported scalar coordinate $g$, let
\[
 b_{g,n,I}=\{Dg(\vartheta_0)E_n\}_I
\]
be its derivative row restricted to the free identifiable coordinates. Choose the alternative with $h_v=0$ and
\[
 h_I=\alpha\,\frac{b_{g,n,I}^{\top}}{\|b_{g,n,I}\|},
\]
holding $a=a_0$ when $d\le3$. As in Step~1 of that proof, the Gaussian Kullback--Leibler divergence is at most $C\alpha^2$. Choose $\alpha>0$ small enough that the total variation distance is at most $1/2$. By the Taylor expansion in Step~2, the two
values of $g$ differ by $\alpha\|b_{g,n,I}\|\{1+o(1)\}$. Hence \eqref{eq:minimax-two-point} gives a squared-error lower bound of order $\|b_{g,n,I}\|^2$.

With the smoothness column removed, the derivative rows in Table~\ref{tab:coordinate-expansions} have norm of order $r_n^{-1}$ for $\eta,c,\log c$; $b_n^{-1}$ for $a,\log a,l,\log l,\sigma^2,\log\sigma^2$ when $d\ge4$; and $n^{-1/2}$ for $t,\log t$. For the global bound, choose the baseline in the relative interior of the rectangle of free
parameters, so both alternatives remain in the rectangle for all sufficiently large $n$. For the local bound, take the given interior truth as the baseline and also choose $\alpha<H$. The alternatives do not depend on the locations, so Step~4 of
the same proof gives the joint-law lower bounds for random designs. The transformations are calculated by setting $h_v=0$ in the increments used in the proof of Corollary~\ref{cor:rates}, before forming the derivative rows and remainder bounds in Table~\ref{tab:coordinate-expansions}. In particular, write $\widehat h_n$ for the normalized error of this restricted fit,
\[
 \widehat\eta_n-\eta_0
 =\log\widehat c_n-\log c_0=r_n^{-1}\widehat h_{n,u},\qquad
 \widehat c_n-c_0
 =c_0r_n^{-1}\widehat h_{n,u}+O_{L^p}(r_n^{-2})
\]
for every fixed $1\le p<\infty$. Hence all three amplitude parameters have rate $r_n^{-1}$, by the same moment calculation. The nugget cumulant calculation and the feasible covariance argument are unchanged. This proves the corollary for the global
fit with known smoothness.
\qed

\section{Proofs of the lemmas}\label{app:lemmas}
The proofs follow the lemma numbering in Appendix~\ref{app:main}. For the sub-Gaussian extension, we invoke established concentration, normal approximation, and anti-concentration results, checking their hypotheses and uniform constants below. The covariance calculations specific to this model are given explicitly. The supporting spectral, sampling, Hessian, score-comparison, and nuisance-curvature arguments are collected in Lemmas~\ref{lem:spectral-covariance}--\ref{lem:nuisance-curvature}; their proofs may be read when first invoked. The fractional-moment and supremum tools used for localization are stated in Lemmas~\ref{lem:local-fractional}--\ref{lem:quadratic-envelope} and proved at the end of this appendix. For the feasible-inference argument, parts \ref{it:studentization-local}--\ref{it:studentization-mixed} of Lemma~\ref{lem:studentization} are proved independently of Lemma~\ref{lem:range-curvature}. The latter uses those local bounds, and then supplies the range-curvature input for part~\ref{it:studentization-inverse} of Lemma~\ref{lem:studentization}.

\subsection{Proof of Lemma~\ref{lem:gaussian}}
We first express each quadratic form in independent coordinates. This makes its variance explicit and also supplies the coordinates used in the concentration and central limit arguments.

\noindent\textbf{Independent coordinates and the variance calculation.}
Put $N=2n$ and
\[
 Z=(G^\T,\xi^\T)^\T\in\R^N,\qquad
 L=C_0^{-1/2}(S_0^{1/2},\sqrt{t_0}I_n).
\]
The coordinates of $Z$ are independent, centered, and have variance one. Since $C_0=S_0+t_0I_n$,
\[
 LL^\T=C_0^{-1/2}(S_0+t_0I_n)C_0^{-1/2}=I_n,
 \qquad X=LZ.
\]
For a real symmetric $A$, let $H=L^\T AL$. Cyclicity of the trace gives
\[
 \tr H=\tr(ALL^\T)=\tr A,\qquad
 \|H\|_\Frob^2
 =\tr(L^\T ALL^\T AL)=\tr(A^2)=\|A\|_\Frob^2.
\]
Also $\|L\|_\op=1$, so $\|H\|_\op\le\|A\|_\op$. The identity $A=LHL^\T$ gives the reverse inequality. Thus
\[
 Q(A)=Z^\T HZ-\tr H,\qquad
 \|H\|_\Frob=\|A\|_\Frob,\qquad
 \|H\|_\op=\|A\|_\op.
\]

Separate the diagonal and off-diagonal terms:
\[
 Z^\T HZ-\tr H
   =\sum_iH_{ii}(Z_i^2-1)+2\sum_{i<j}H_{ij}Z_iZ_j.
\]
Every summand is centered, and distinct summands in this display are uncorrelated. For example, two products $Z_iZ_j$ and $Z_kZ_l$ have zero covariance unless the unordered pairs coincide: otherwise some coordinate occurs exactly once, and independence and centering make the expectation zero. A diagonal term and an off-diagonal term are
also uncorrelated. Even when they share an index,
\[
 \E\{(Z_i^2-1)Z_iZ_j\}
 =\E\{(Z_i^2-1)Z_i\}\,\E Z_j=0,\qquad i\ne j.
\]
Consequently, with $\underline v=\min\{2,v_*\}$,
\[
 \begin{aligned}
 \Var(Z^\T HZ)
 &=4\sum_{i<j}H_{ij}^2+\sum_i(\E Z_i^4-1)H_{ii}^2\\
 &=2\sum_{i\ne j}H_{ij}^2+\sum_i(\E Z_i^4-1)H_{ii}^2\\
 &=2\|H\|_\Frob^2+\sum_i(\E Z_i^4-3)H_{ii}^2\\
 &\ge\underline v\left\{\sum_{i\ne j}H_{ij}^2+\sum_iH_{ii}^2\right\}
   =\underline v\|H\|_\Frob^2.
 \end{aligned}
\]
For the last inequality, use the second line: its off-diagonal coefficient is $2$, and each diagonal coefficient is at least $\underline v$. The Gaussian coordinates have $\E Z_i^4-1=2$, and the noise coordinates have $\E Z_i^4-1\ge v_*$. The uniform fourth-moment bound likewise gives $\Var(Z^\T HZ)\le C\|H\|_\Frob^2$. The norm identity above proves the asserted variance bounds for $Q(A)$.

For two matrices $A,B$, write $H_A=L^\T AL$ and $H_B=L^\T BL$. The same pairwise covariance calculation gives, explicitly,
\[
 \begin{aligned}
 \Cov\{Q(A),Q(B)\}
 &=2\sum_{i\ne j}(H_A)_{ij}(H_B)_{ij}
       +\sum_i(\E Z_i^4-1)(H_A)_{ii}(H_B)_{ii}\\
 &=2\tr(H_AH_B)
       +\sum_i(\E Z_i^4-3)(H_A)_{ii}(H_B)_{ii}.
 \end{aligned}
\]
Here $\tr(H_AH_B)=\tr(AB)$ because $LL^\T=I_n$. The first $n$ coordinates of $Z$ are Gaussian, so their terms in the last sum vanish. The bottom-right block of $H_A$ is
$t_0C_0^{-1/2}AC_0^{-1/2}$, and similarly for $H_B$. The remaining sum is therefore
\[
 t_0^2(\E_P\xi^4-3)
 \sum_{i=1}^n(C_0^{-1/2}AC_0^{-1/2})_{ii}
              (C_0^{-1/2}BC_0^{-1/2})_{ii}.
\]
Since its leading factor is $\kappa_4$, this proves \eqref{eq:quadratic-cumulant}.

\noindent\textbf{Concentration and moments.}
Let $K_\psi$ be a common upper bound on the sub-Gaussian norms of the coordinates of $Z$, supplied by Assumption~\ref{ass:noise}. The Hanson--Wright inequality of
\citet[Theorem~1.1]{rudelson2013hansonwright}, applied to $H$, states that for every threshold $s\ge0$,
\[
 P\{|Q(A)|>s\}
 \le2\exp\left[-c_{\rm HW}\min\left\{
       \frac{s^2}{K_\psi^4\|A\|_\Frob^2},
       \frac{s}{K_\psi^2\|A\|_\op}\right\}\right],
\]
where we used $\|H\|_\Frob=\|A\|_\Frob$ and $\|H\|_\op=\|A\|_\op$. The constant $c_{\rm HW}>0$ is absolute. If $A=0$ there is nothing to prove. Otherwise put
\[
 a=\frac{K_\psi^2}{\sqrt{c_{\rm HW}}}\|A\|_\Frob,
 \qquad
 b=\frac{K_\psi^2}{c_{\rm HW}}\|A\|_\op,
 \qquad s(z)=a\sqrt z+bz.
\]
For $z\ge0$, the two terms in the minimum satisfy
\[
 \frac{c_{\rm HW}s(z)^2}{K_\psi^4\|A\|_\Frob^2}
       =\frac{s(z)^2}{a^2}\ge z,
 \qquad
 \frac{c_{\rm HW}s(z)}{K_\psi^2\|A\|_\op}
       =\frac{s(z)}b\ge z.
\]
Thus $P\{|Q(A)|>a\sqrt z+bz\}\le2e^{-z}$. Enlarging one common constant in front of the two matrix norms gives \eqref{eq:quadconc}. The square-root and linear terms correspond to the Gaussian and exponential parts of the tail, respectively.

Here is the tail integration giving the moment bound. For $x>0$, choose
\[
 z_x=\min\left\{\frac{x^2}{4a^2},\frac{x}{2b}\right\}.
\]
Then $a\sqrt{z_x}\le x/2$ and $bz_x\le x/2$, so
\[
 \begin{aligned}
 P\{|Q(A)|>x\}
 &\le2e^{-z_x}\\
 &\le2\left\{e^{-x^2/(4a^2)}+e^{-x/(2b)}\right\}.
 \end{aligned}
\]
Consequently, for $p\ge1$,
\[
 \begin{aligned}
 \E|Q(A)|^p
 &=p\int_0^\infty x^{p-1}P\{|Q(A)|>x\}\,dx\\
 &\le2p\int_0^\infty x^{p-1}
             \{e^{-x^2/(4a^2)}+e^{-x/(2b)}\}\,dx\\
 &=2(2a)^p\Gamma(1+p/2)+2(2b)^p\Gamma(1+p).
 \end{aligned}
\]
In the last line, substitute $y=x^2/(4a^2)$ in the first integral and $y=x/(2b)$ in the second. Taking the $p$th root and using $\Gamma(1+p/2)^{1/p}\le C\sqrt p$ and
$\Gamma(1+p)^{1/p}\le Cp$ gives
\[
 \|Q(A)\|_{L^p}
 \le C\{a\sqrt p+bp\}
 \le C\{\sqrt p\|A\|_\Frob+p\|A\|_\op\}.
\]
\noindent\textbf{The exponential-moment bound.}
Retain the independent-coordinate representation $H=L^\T AL$ and write
\[
 D_H=\sum_iH_{ii}(Z_i^2-1),\qquad
 O_H=\sum_{i\ne j}H_{ij}Z_iZ_j,\qquad Q(A)=D_H+O_H.
\]
The centered variables $Z_i^2-1$ have moments $\E|Z_i^2-1|^r\le r!(C K_\psi^2)^r$ for integers $r\ge2$. Expanding their exponential and using their zero means gives
\[
 \begin{aligned}
 \E e^{tH_{ii}(Z_i^2-1)}
 &\le 1+\sum_{r=2}^\infty(CK_\psi^2|tH_{ii}|)^r\\
 &\le 1+C K_\psi^4t^2H_{ii}^2
 \le e^{C K_\psi^4t^2H_{ii}^2},
 \qquad CK_\psi^2|tH_{ii}|\le\tfrac12.
 \end{aligned}
\]
Independence therefore gives $\log \E e^{tD_H}\le C K_\psi^4t^2\sum_iH_{ii}^2$
when $K_\psi^2|t|\max_i|H_{ii}|\le c$. For the off-diagonal sum, the decoupling calculation in \citet[proof of Theorem~1.1, Steps~2--4]{rudelson2013hansonwright} gives
\[
 \log \E e^{tO_H}\le C K_\psi^4t^2\|H\|_\Frob^2,
 \qquad K_\psi^2|t|\|H\|_\op\le c.
\]
Their displayed moment-generating-function bound is for $K_\psi=1$ and $t>0$; rescaling the coordinates and replacing $H$ by $-H$ give the stated form. The two sums need not be independent: Cauchy--Schwarz gives
\[
 \E e^{tQ(A)}\le (\E e^{2tD_H})^{1/2}(\E e^{2tO_H})^{1/2}
          \le e^{C K_\psi^4t^2\|H\|_\Frob^2}
\]
after decreasing $c$. Since $\max_i|H_{ii}|\le\|H\|_\op$ and the two matrix norms of $H$ equal those of $A$, this proves the claimed bound. Only the uniform upper sub-Gaussian bound is used.

\noindent\textbf{A uniform central limit theorem.}
We use \citet[Theorem~4.1]{privault2022berryesseen}, which allows independent coordinates with different distributions. Write
\[
 W_i=H_{ii}(Z_i^2-1),\qquad
 W_{ij}=W_{ji}=2H_{ij}Z_iZ_j\quad(i\ne j).
\]
Then $Q(A)=\sum_iW_i+\sum_{i<j}W_{ij}$. These are the order-one and order-two terms of its Hoeffding decomposition: $\E W_i=0$ and $\E(W_{ij}\mid Z_i)=\E(W_{ij}\mid Z_j)=0$.
Put $\sigma^2=\Var Q(A)$. Apply their equation~(4.1) with maximum order two to the normalized terms $W_i/\sigma$ and $W_{ij}/\sigma$. With $d_{\rm K}$ denoting Kolmogorov distance, that equation becomes
\[
 d_{\rm K}\!\left(\frac{Q(A)}\sigma,N(0,1)\right)
 \le\frac C{\sigma^2}
       (\mathcal R_1+\mathcal R_2+\mathcal R_3
                       +\mathcal R_4+\mathcal R_5)^{1/2},
\]
where the five unnormalized terms are
\[
 \begin{aligned}
 \mathcal R_1&=\sum_i \E W_i^4,\\
 \mathcal R_2&=\sum_{i<j}\E W_{ij}^4,\\
 \mathcal R_3&=\sum_i \E\left[
       \left\{\sum_{j\ne i}\E(W_{ij}^2\mid Z_i)\right\}^2\right],\\
 \mathcal R_4&=\sum_{i\ne j} \E\left[
       \left\{\sum_{k\notin\{i,j\}}
                  \E(W_{ik}W_{jk}\mid Z_i,Z_j)\right\}^2\right],\\
 \mathcal R_5&=\sum_i \E\left[
       \left\{\sum_{j\ne i}\E(W_jW_{ij}\mid Z_i)\right\}^2\right].
 \end{aligned}
\]
The first sum in their (4.1) gives $\mathcal R_1,\mathcal R_2, \mathcal R_3$ from the order pairs $(1,0),(2,0),(2,1)$. Its second and third sums give $\mathcal R_4$ and $\mathcal R_5$, respectively. The factor $\sigma^{-2}$ comes from taking the square
root of expressions homogeneous of degree four in the normalized terms. We next evaluate these five expressions for our quadratic form.

Let
\[
 K=H-\diag(H_{11},\ldots,H_{NN}),\qquad
 c_i=H_{ii}\E Z_i^3,\qquad r_i=\sum_jH_{ij}^2.
\]
The uniform eighth-moment bound controls $\E(Z_i^2-1)^4$ and all fourth moments. Independence therefore gives
\[
 \begin{aligned}
 \mathcal R_1
 &=\sum_iH_{ii}^4\E(Z_i^2-1)^4
       \le C\sum_iH_{ii}^4,\\
 \mathcal R_2
 &=16\sum_{i<j}H_{ij}^4(\E Z_i^4)(\E Z_j^4)
       \le C\sum_{i<j}H_{ij}^4.
 \end{aligned}
\]
For the other three expressions, the relevant conditional expectations are, for distinct indices as indicated,
\[
 \begin{aligned}
 \E(W_{ij}^2\mid Z_i)&=4H_{ij}^2Z_i^2,\\
 \E(W_{ik}W_{jk}\mid Z_i,Z_j)&=4H_{ik}H_{jk}Z_iZ_j,
                   \qquad k\notin\{i,j\},\\
 \E(W_jW_{ij}\mid Z_i)
   &=2H_{ij}H_{jj}\E\{(Z_j^2-1)Z_j\}\,Z_i
     =2H_{ij}c_jZ_i.
 \end{aligned}
\]
Here $\E Z_k^2=1$ in the second line, and $\E Z_j=0$ in the third. It follows that
\[
 \begin{aligned}
 \mathcal R_3
 &=16\sum_i(\E Z_i^4)\left(\sum_{j\ne i}H_{ij}^2\right)^2
       \le C\sum_i r_i^2,\\
 \mathcal R_4
 &=16\sum_{i\ne j}
             \left(\sum_{k\notin\{i,j\}}H_{ik}H_{jk}\right)^2
   =16\sum_{i\ne j}(K^2)_{ij}^2
       \le16\|K^2\|_\Frob^2,\\
 \mathcal R_5
 &=4\sum_i\left(\sum_{j\ne i}H_{ij}c_j\right)^2
   =4\|Kc\|^2.
 \end{aligned}
\]
We used $\E(Z_i^2Z_j^2)=1$ for $i\ne j$ in $\mathcal R_4$. Also $H_{ii}^4\le r_i^2$ and
$\sum_jH_{ij}^4\le(\sum_jH_{ij}^2)^2=r_i^2$. Thus $\mathcal R_1+\mathcal R_2+\mathcal R_3\le C\sum_i r_i^2$, and the cited theorem yields
\[
 d_{\rm K}\!\left(\frac{Z^\T HZ-\tr H}{\sigma},N(0,1)\right)
 \le\frac C{\sigma^2}
 \left\{\sum_i\left(\sum_jH_{ij}^2\right)^2
              +\|K^2\|_\Frob^2+\|Kc\|^2\right\}^{1/2}.
\]
The constant depends only on the common moment bound, so it is uniform over the noise class.

For completeness, let $e_i$ be the $i$th coordinate vector. Since $H$ is symmetric,
\[
 r_i=\sum_jH_{ij}^2=\|He_i\|^2\le\|H\|_\op^2,
 \qquad \sum_i r_i=\sum_{i,j}H_{ij}^2=\|H\|_\Frob^2.
\]
Consequently,
\[
 \sum_i\left(\sum_jH_{ij}^2\right)^2
 =\sum_i r_i^2
 \le(\max_i r_i)\sum_i r_i
 \le\|H\|_\op^2\|H\|_\Frob^2.
\]
Moreover, $|H_{ii}|=|e_i^\T H e_i|\le\|H\|_\op$, so $\|K\|_\op\le2\|H\|_\op$ and $\|K\|_\Frob\le\|H\|_\Frob$. The remaining terms satisfy
\[
 \begin{aligned}
 \|K^2\|_\Frob^2
   &\le\|K\|_\op^2\|K\|_\Frob^2
    \le4\|H\|_\op^2\|H\|_\Frob^2,\\
 \|Kc\|^2
   &\le\|K\|_\op^2\sum_iH_{ii}^2(\E Z_i^3)^2
    \le C\|H\|_\op^2\|H\|_\Frob^2.
 \end{aligned}
\]
We therefore obtain
\[
 d_{\rm K}\!\left(\frac{Q(A)}\sigma,N(0,1)\right)
 \le C\sigma^{-2}\|H\|_\op\|H\|_\Frob.
\]
Under the CLT hypotheses, $\|H\|_\Frob=\|A\|_\Frob\le C$ and $c\le\sigma^2\le C\|H\|_\Frob^2$. These inequalities imply both $\sigma^2\asymp1$ and $\|H\|_\Frob\asymp1$. Hence
\[
 \sigma^{-2}\|H\|_\op\|H\|_\Frob
       \asymp\|H\|_\op=\|A\|_\op\longrightarrow0
\]
uniformly. Intuitively, the vanishing operator norm prevents any one direction from carrying a substantial part of the fluctuation. The cited theorem quantifies this through Stein's method, including the diagonal terms and possible nonzero third moments.

For the vector assertion, write $W_n=(Q(A_{n,1}),\ldots,Q(A_{n,k}))^\T$ and
$\mathcal K_n=\Cov(W_n)\asymp I$, with $k$ fixed. For a unit vector $v$, put $q=\mathcal K_n^{-1/2}v$. Then
\[
 v^\T\mathcal K_n^{-1/2}W_n
 =Q\!\left(\sum_jq_jA_{n,j}\right),\qquad
 \Var(v^\T\mathcal K_n^{-1/2}W_n)=1.
\]
Since $\|q\|\le C$, the matrix in this display has bounded Frobenius norm and vanishing operator norm, uniformly in $v$. The scalar bound and Cram\'er--Wold give the stated vector limit. Applying this argument along arbitrary sequences of truths, noise laws, and matrices proves uniformity.

For the error-only form,
\[
 Q_\varepsilon(A)=t_0\{\xi^\T A\xi-\tr A\}.
\]
Apply the same upper-bound arguments directly to the independent vector $\xi$ and the matrix $t_0A$. Since $t_0$ lies in a fixed positive compact interval, the constants remain uniform. These upper bounds do not require a positive variance lower bound.
\qed

\subsection{Proof of Lemma~\ref{lem:windows}}
We construct Fourier coefficients of the cut-off field and retain frequencies in a narrow band. The nearest-cell rule recovers their linear combinations under every candidate model. Fourier localization then controls their variances, while decay of the
covariances controls sums of their squared magnitudes.

\noindent\textbf{Recovering integrals.}
Uniformly in $\theta\in\Theta$,
$f_\theta(\omega)\le C\langle\omega\rangle^{-2\nu_--d}$. For $0<\zeta<\min(\nu_-,1)$, the spectral formula therefore gives
\begin{equation}\label{eq:variogram}
 \E_\theta|F(x)-F(y)|^2
 =\int f_\theta(\omega)|e^{2\pi i\omega^\T(x-y)}-1|^2\,d\omega
 \le C\norm{x-y}^{2\zeta}.
\end{equation}
Here we used $|e^{is}-1|\le C\min\{|s|,1\}\le C|s|^\zeta$ and the finite $2\zeta$ spectral moment. The integral $Z_g$ is well defined in mean square for $g\in L^2(Q)$; a continuous version of $F$ also exists by Gaussian moment bounds and Kolmogorov's criterion.

Let $\pi_n(x)$ be the site assigned to $x$ by the nearest-site rule. Then
\[
 \widehat Z_g-Z_g
 =\int_Qg(x)\{F(\pi_n(x))-F(x)\}\,dx
   +\sum_iw_{n,i}(g)\varepsilon_{n,i}.
\]
The two terms are independent. The first has standard deviation at most $Ch_n^\zeta\norm g_{L^1}\le Ch_n^\zeta\norm g_{L^2}$ by \eqref{eq:variogram}. Since each cell has volume at most $Ch_n^d$, the second has variance
\[
 t\sum_i\left|\int_{V_{n,i}}g(x)\,dx\right|^2
 \le Ch_n^d\sum_i\int_{V_{n,i}}|g(x)|^2\,dx
 =Ch_n^d\norm g_{L^2}^2.
\]
This proves \eqref{eq:intrecovery}. Applying it to $g=\sum_kb_kg_k$ proves the covariance bound, because $\norm g_{L^2}^2=\sum_k|b_k|^2$ for an orthonormal basis. In particular, the estimate holds for every linear combination with the same constant, which is why no dimension factor is lost.

\noindent\textbf{Constructing stable bands of Fourier coefficients.}
Use the cubes $Q_0\Subset Q$, cutoff, and single lattice from the statement. The exponentials with frequencies in $\mathscr L=(m/L)\Z^d$ are orthogonal both on $Q_0$, whose side is $L/m$, and on $Q$, whose side is $L$. Since $\chi=1$ on $Q_0$ and is supported in $Q$, these two orthogonality calculations give
\begin{equation}\label{eq:windownorm}
 m^{-d}\sum_\lambda|b_\lambda|^2
 \le\left\|L^{-d/2}\chi
           \sum_\lambda b_\lambda e^{2\pi i\lambda^\T x}\right\|_{L^2}^2
 \le\norm\chi_{L^\infty}^2\sum_\lambda|b_\lambda|^2.
\end{equation}
The factor $m^{-d}$ is $L^{-d}|Q_0|$. This explains the spacing: the lower bound uses a whole period cube on which the cutoff equals one. Keeping all frequencies of $L^{-1}\Z^d$ would not give this lower bound, since their trigonometric polynomials can concentrate where $\chi$ vanishes. Lattice-point counting gives $\#\Lambda_R^+\asymp R^d$, and \eqref{eq:windownorm} makes the coefficient map injective. The real and imaginary parts of $X_\lambda$ are precisely the integrals against the corresponding cutoff cosines and sines.

\noindent\textbf{Controlling mass outside the band.}
The Fourier transform is
\[
 \widehat g(\omega)=L^{-d/2}\sum_{\lambda\in\Lambda_R}
 b_\lambda\widehat\chi(\omega-\lambda).
\]
Since $\chi\in C_c^\infty(\R^d)$, its Fourier transform $\widehat\chi$ is a Schwartz function. In particular, for every $N\ge0$ there is $C_N$ such that $|\widehat\chi(\eta)|\le C_N\langle\eta\rangle^{-N}$. Moreover, the frequency clusters lie a distance of order $R$ inside $\mathcal A_R$. Consequently, for $R\ge1$,
\[
 \langle\omega-\lambda\rangle\ge c(R+\langle\omega\rangle),
 \qquad\omega\notin\mathcal A_R,\quad\lambda\in\Lambda_R.
\]
Cauchy--Schwarz in the finite sum, followed by \eqref{eq:windownorm} and $\#\Lambda_R\asymp R^d$, gives
\begin{align*}
 |\widehat g(\omega)|^2
 &\le C\left(\sum_\lambda|b_\lambda|^2\right)
          \left(\sum_\lambda|\widehat\chi(\omega-\lambda)|^2\right)\\
 &\le C_NR^d(R+\langle\omega\rangle)^{-2N}\norm g_{L^2}^2,
 \qquad\omega\notin\mathcal A_R.
\end{align*}
Let $B\ge0$ and $2N>2B+d$. Since $\langle\omega\rangle\le R+\langle\omega\rangle$, polar coordinates $r=\norm\omega$ and then $r=Rt$ give
\begin{align*}
 &~~~~\int_{\R^d\setminus\mathcal A_R}
       \langle\omega\rangle^{2B}|\widehat g(\omega)|^2\,d\omega\\
 &\le C_NR^d\norm g_{L^2}^2
       \int_{\R^d}\langle\omega\rangle^{2B}
                       (R+\langle\omega\rangle)^{-2N}\,d\omega\\
 &\le C_{N,B}R^d\norm g_{L^2}^2
       \int_0^\infty r^{d-1}(R+r)^{2B-2N}\,dr\\
 &=C_{N,B}R^{2d+2B-2N}\norm g_{L^2}^2
       \int_0^\infty t^{d-1}(1+t)^{2B-2N}\,dt.
\end{align*}
The surface area of the unit sphere is absorbed into $C_{N,B}$. The last integral is finite because $2N>2B+d$. For a prescribed $A\ge0$, choose $2N>2d+2B+A$; the bound is
then at most $C_{A,B}R^{-A}\norm g_{L^2}^2$. This proves \eqref{eq:windowleak} and defines precisely the Fourier leakage used below.

\noindent\textbf{Sobolev norms and Mat\'ern variances.}
Write the fixed interval of Sobolev orders as $[s_-,s_+]$. On $\mathcal A_R$, $cR\le\langle\omega\rangle\le CR$, so $\langle\omega\rangle^{2s}\asymp R^{2s}$ uniformly in this interval, including negative $s$. Plancherel and the leakage
bound with $A=1$, $B=0$ give, for all sufficiently large $R$,
\[
 \int_{\mathcal A_R}|\widehat g|^2
 =\norm g_{L^2}^2-\int_{\R^d\setminus\mathcal A_R}|\widehat g|^2
 \ge(1-CR^{-1})\norm g_{L^2}^2
 \ge\tfrac12\norm g_{L^2}^2.
\]
The same integral is at most $\norm g_{L^2}^2$. Hence the contribution from the band satisfies
\[
 cR^{2s}\norm g_{L^2}^2
 \le\int_{\mathcal A_R}\langle\omega\rangle^{2s}
                         |\widehat g(\omega)|^2\,d\omega
 \le CR^{2s}\norm g_{L^2}^2.
\]
To control the complement uniformly, take $B_*\coloneqq \max\{s_+,0\}$ and $A_*>\max\{0,-2s_-\}$. Since $\langle\omega\rangle\ge1$ and $s\le B_*$, the weighted leakage estimate gives
\begin{align*}
 \int_{\R^d\setminus\mathcal A_R}
       \langle\omega\rangle^{2s}|\widehat g(\omega)|^2\,d\omega
 &\le\int_{\R^d\setminus\mathcal A_R}
       \langle\omega\rangle^{2B_*}|\widehat g(\omega)|^2\,d\omega\\
 &\le CR^{-A_*}\norm g_{L^2}^2\\
 &\le CR^{-A_*-2s_-}R^{2s}\norm g_{L^2}^2.
\end{align*}
The factor $R^{-A_*-2s_-}$ tends to zero. Thus the complement is smaller than the band contribution uniformly in $s$, even when $s<0$ and $R^{2s}$ tends to zero. Adding the two contributions proves \eqref{eq:windowsobolev}.

The spectral covariance identity in Lemma~\ref{lem:spectral-covariance}, applied under parameter $\theta$, gives
\[
 \Var_\theta(Z_g)=\int f_\theta(\omega)|\widehat g(\omega)|^2\,d\omega.
\]
Uniformly on the parameter rectangle, $f_\theta(\omega)\asymp\langle\omega\rangle^{-2\beta}$. The Sobolev estimate with $s=-\beta$ proves \eqref{eq:windowvariance}.

For the variance ratio, write the spectral integral as $I_\theta+O_\theta$, its parts inside and outside $\mathcal A_R$. Taking the leakage exponent sufficiently large gives
\[
 0\le O_\theta\le C_AR^{-A}\Var_\theta(Z_g)
\]
uniformly in $\theta$. On the band, $I_\theta/I_{\theta'}$ is a weighted average of
$f_\theta/f_{\theta'}$, with nonnegative weights $f_{\theta'}|\widehat g|^2$. Bounding that average by its infimum and supremum, and then restoring the two tails, proves
\eqref{eq:windowratio}. For a basis $g_1,\ldots,g_N$, the covariance quotient in direction $b$ is the same variance ratio for $g=\sum_kb_kg_k$. The generalized-eigenvalue assertion follows by the min--max characterization of these quotients.

\noindent\textbf{Covariance and energy bounds for the individual coefficients.}
This is the localized Fourier calculation used in \citet[Sections 5--8]{pillai2026}; we give the details needed here in general dimension and with our Fourier convention.
Apply the spectral covariance identity in Lemma~\ref{lem:spectral-covariance}, with $f_0$ replaced by $f_\theta$, to the test functions $L^{-d/2}\chi(x)e^{-2\pi i\lambda^\T x}$. Using the evenness of $f_\theta$ gives
\begin{align*}
 \Sigma_\theta(\lambda,\mu)
 &\coloneqq \E_\theta[X_\lambda\overline{X_\mu}]\\
 &=L^{-d}\int f_\theta(\omega)
       \widehat\chi(\lambda-\omega)
       \overline{\widehat\chi(\mu-\omega)}\,d\omega,\\
 \Pi_\theta(\lambda,\mu)
 &\coloneqq \E_\theta[X_\lambda X_\mu]\\
 &=L^{-d}\int f_\theta(\omega)
       \widehat\chi(\lambda-\omega)
       \widehat\chi(\mu+\omega)\,d\omega.
\end{align*}
Here $\Pi_\theta(\lambda,\mu)=\Sigma_\theta(\lambda,-\mu)$, because $X_{-\mu}=\overline{X_\mu}$; the displayed formula also uses $\overline{\widehat\chi(\xi)}=\widehat\chi(-\xi)$ for real $\chi$. Setting $\mu=\lambda$ in the first identity proves the variance formula in \eqref{eq:fouriervariance}.

\noindent\textbf{A second-order expansion of the individual variances.}
The convolution in \eqref{eq:fouriervariance} smooths the spectrum over a frequency neighborhood of fixed width. At large frequencies the Mat\'ern spectrum varies slowly over that neighborhood. Moreover, the smoothing kernel
\[
 k_\chi(y)=L^{-d}|\widehat\chi(y)|^2
\]
is even, because $\chi$ is real. Consequently its first moments vanish, and the leading relative error has order $\langle\lambda\rangle^{-2}$ rather than
$\langle\lambda\rangle^{-1}$.

\begin{samepage}
Here are the details, including the more precise second-order term. For every fixed multi-index $\alpha$, differentiating $f_\theta(\xi)=e^\eta m_a(\xi)^{-\beta}$ gives
\[
 |\partial_\xi^\alpha f_\theta(\xi)|
 \le C_\alpha f_\theta(\xi)\langle\xi\rangle^{-|\alpha|},
 \qquad m_a(\xi)=a+4\pi^2\norm\xi^2.
\]
\end{samepage}
Indeed, each differentiation either lowers the power of $m_a$ and introduces a linear factor in $\xi$, or differentiates one of the existing polynomial factors; using
$m_a(\xi)\asymp\langle\xi\rangle^2$ gives the displayed bound. On $\norm y\le\norm\lambda/2$, all points of the segment from $\lambda-y$ to $\lambda+y$ have brackets comparable to $\langle\lambda\rangle$. Taylor expansion through degree three,
averaged at $y$ and $-y$, therefore gives
\[
 \frac{f_\theta(\lambda-y)+f_\theta(\lambda+y)}2
 =f_\theta(\lambda)+\frac12y^\T\nabla_\lambda^2
                                      f_\theta(\lambda)y
   +O\{f_\theta(\lambda)\langle\lambda\rangle^{-4}\norm y^4\}.
\]
The odd terms cancel. Integrate against $k_\chi(y)$. Its fourth moment is finite. We next bound the complementary integral before dividing by $f_\theta(\lambda)$. Write $\beta_+=\nu_++d/2$. The bracket inequality
\[
 \langle\lambda\rangle\le\sqrt2\langle\lambda-y\rangle
                                      \langle y\rangle
\]
implies $\langle\lambda-y\rangle^{-2\beta}
\le C\langle\lambda\rangle^{-2\beta}\langle y\rangle^{2\beta}$, uniformly for $\beta\le\beta_+$. For any prescribed $N\ge0$, Schwartz decay gives
\[
 0\le k_\chi(y)\le C_N
             \langle y\rangle^{-2\beta_+-N-d-1}.
\]
Moreover, $\norm y>\norm\lambda/2$ implies $\langle y\rangle\ge\langle\lambda\rangle/2$. Starting with the unnormalized tail, these inequalities give
\begin{align*}
 &~~~~\int_{\norm y>\norm\lambda/2}
                     k_\chi(y)f_\theta(\lambda-y)\,dy\\*
 &\asymp\int_{\norm y>\norm\lambda/2}
                  k_\chi(y)\langle\lambda-y\rangle^{-2\beta}\,dy\\*
 &\le C\langle\lambda\rangle^{-2\beta}
       \int_{\norm y>\norm\lambda/2}
                          k_\chi(y)\langle y\rangle^{2\beta}\,dy\\*
 &\le C_N\langle\lambda\rangle^{-2\beta}
       \int_{\norm y>\norm\lambda/2}\langle y\rangle^{-N-d-1}\,dy\\*
 &\le C_N\langle\lambda\rangle^{-2\beta-N}
                          \int_{\R^d}\langle y\rangle^{-d-1}\,dy\\*
 &\le C_N f_\theta(\lambda)\langle\lambda\rangle^{-N}.
\end{align*}
The integral in the penultimate line is finite. The last step uses $f_\theta(\lambda)\asymp\langle\lambda\rangle^{-2\beta}$. Dividing by $f_\theta(\lambda)>0$ gives the required relative tail bound
\[
 \int_{\norm y>\norm\lambda/2}
 k_\chi(y)\frac{f_\theta(\lambda-y)}{f_\theta(\lambda)}\,dy
 \le C_N\langle\lambda\rangle^{-N}.
\]
The same extraction of $\langle\lambda\rangle^{-N}$, after increasing the Schwartz exponent, gives for $j=0,2$
\[
 \int_{\norm y>\norm\lambda/2}\norm y^j k_\chi(y)\,dy
 \le C_N\langle\lambda\rangle^{-N}.
\]
\begin{samepage}
These control the omitted tails of the zeroth and second moments. Thus, with the constants $c_\chi$ and $M_\chi$ in the statement,
\[
 v_\theta(\lambda)
 =c_\chi f_\theta(\lambda)
  +\frac12\sum_{i,j}(M_\chi)_{ij}
                          \partial_{\lambda_i\lambda_j}f_\theta(\lambda)
  +O\{f_\theta(\lambda)\langle\lambda\rangle^{-4}\}.
\]
\end{samepage}
In particular, $v_\theta/(c_\chi f_\theta)=1+O(\langle\lambda\rangle^{-2})$.

\noindent\textbf{Logarithmic parameter derivatives.}
Start with the identity
\[
 \partial_j\log v_\theta(\lambda)-\partial_j\log f_\theta(\lambda)
 =\frac{\partial_j\{v_\theta(\lambda)/f_\theta(\lambda)\}}
              {v_\theta(\lambda)/f_\theta(\lambda)},
 \qquad j\in\{\nu,a\}.
\]
The second-order expansion bounds the denominator below by $c_\chi/2$ for large $\norm\lambda$. On bounded frequency sets, continuity and positivity of the variance integral give a uniform positive lower bound as well. It therefore suffices to bound the derivatives of the ratio in the numerator. This motivates writing
\[
 R_\theta(\lambda,y)
 =\frac{f_\theta(\lambda-y)}{f_\theta(\lambda)}
 =\left\{\frac{m_a(\lambda)}{m_a(\lambda-y)}\right\}^{\beta}.
\]
Then
\[
 \frac{v_\theta(\lambda)}{f_\theta(\lambda)}
       =\int k_\chi(y)R_\theta(\lambda,y)\,dy.
\]
Since $\partial_\nu\beta=1$, direct differentiation gives
\[
 \partial_\nu R_\theta
 =R_\theta\log\frac{m_a(\lambda)}{m_a(\lambda-y)},\qquad
 \partial_a R_\theta
 =\beta R_\theta
        \{m_a(\lambda)^{-1}-m_a(\lambda-y)^{-1}\}.
\]
Both expressions vanish at $y=0$. On $\norm y\le\norm\lambda/2$, $R_\theta$ is bounded,
$\nabla_yR_\theta=O(\langle\lambda\rangle^{-1})$, and $\nabla_y^2R_\theta=O(\langle\lambda\rangle^{-2})$. The other factors have the derivative orders
\[
 \begin{array}{c|ccc}
  &\text{function}&\text{first }y\text{ derivative}
                  &\text{second }y\text{ derivative}\\ \hline
  \log\{m_a(\lambda)/m_a(\lambda-y)\}
       &O(1)&O(\langle\lambda\rangle^{-1})
             &O(\langle\lambda\rangle^{-2})\\
  m_a(\lambda)^{-1}-m_a(\lambda-y)^{-1}
       &O(\langle\lambda\rangle^{-2})
       &O(\langle\lambda\rangle^{-3})
       &O(\langle\lambda\rangle^{-4}).
 \end{array}
\]
For example, $\nabla\log m_a(\xi)=8\pi^2\xi/m_a(\xi)$ and $\nabla m_a(\xi)^{-1}=-8\pi^2\xi/m_a(\xi)^2$; differentiating once more gives the last column. The product rule therefore bounds the Hessians in $y$ of $\partial_\nu R_\theta$
and $\partial_aR_\theta$ by $C\langle\lambda\rangle^{-2}$ and $C\langle\lambda\rangle^{-4}$, respectively. To write out the integration of the Taylor remainder, put
\[
 D_\lambda=\{y:\norm y\le\norm\lambda/2\},\qquad
 G_j(y)=\partial_jR_\theta(\lambda,y),\qquad
 s_\nu=2,\quad s_a=4,\quad j\in\{\nu,a\}.
\]
We have $G_j(0)=0$ and
$\sup_{y\in D_\lambda}\norm{\nabla_y^2G_j(y)}_\op
\le C\langle\lambda\rangle^{-s_j}$.
The first derivatives at zero are
\[
 \nabla_yG_\nu(0)=\frac{8\pi^2\lambda}{m_a(\lambda)},\qquad
 \nabla_yG_a(0)=-\frac{8\pi^2\beta\lambda}{m_a(\lambda)^2}.
\]
Taylor's formula with integral remainder, for $y\in D_\lambda$, is therefore
\[
 G_j(y)=\nabla_yG_j(0)^\T y
       +\int_0^1(1-t)y^\T\nabla_y^2G_j(ty)y\,dt.
\]
Both $D_\lambda$ and $k_\chi$ are invariant under $y\mapsto-y$, whereas the linear term changes sign. In particular,
\[
 \int_{D_\lambda}y k_\chi(y)\,dy=0,\qquad
 \int_{D_\lambda}k_\chi(y)\nabla_yG_j(0)^\T y\,dy=0.
\]
Only the remainder remains after integration. Since $ty\in D_\lambda$ for $0\le t\le1$, its contribution satisfies
\begin{align*}
 \left|\int_{D_\lambda}k_\chi(y)G_j(y)\,dy\right|
 &=\left|\int_{D_\lambda}k_\chi(y)
        \int_0^1(1-t)y^\T\nabla_y^2G_j(ty)y\,dt\,dy\right|\\*
 &\le C\langle\lambda\rangle^{-s_j}
        \int_{D_\lambda}\norm y^2k_\chi(y)\,dy
        \int_0^1(1-t)\,dt\\*
 &\le\frac C2\langle\lambda\rangle^{-s_j}
                         \int_{\R^d}\norm y^2k_\chi(y)\,dy
 \le C\langle\lambda\rangle^{-s_j}.
\end{align*}
The second moment is finite because $k_\chi$ is Schwartz.

For completeness, the derivative tails have arbitrarily fast decay as well. The bracket inequality, used in both directions, gives
\[
 \frac{m_a(\lambda)}{m_a(\lambda-y)}
 \le C\frac{\langle\lambda\rangle^2}{\langle\lambda-y\rangle^2}
 \le C\langle y\rangle^2,\qquad
 \frac{m_a(\lambda-y)}{m_a(\lambda)}\le C\langle y\rangle^2.
\]
Consequently
\[
 R_\theta(\lambda,y)\le C\langle y\rangle^{2\beta_+},\qquad
 \left|\log\frac{m_a(\lambda)}{m_a(\lambda-y)}\right|
       \le\log C+2\log\langle y\rangle.
\]
Also $m_a(\xi)\ge a_-$ for every $\xi$, so both reciprocals are at most $a_-^{-1}$. The triangle inequality gives $|m_a(\lambda)^{-1}-m_a(\lambda-y)^{-1}|\le2/a_-$. The formulas for the two integrands are
\[
 \begin{aligned}
 G_\nu(y)&=\partial_\nu R_\theta(\lambda,y)
       =R_\theta(\lambda,y)
                    \log\frac{m_a(\lambda)}{m_a(\lambda-y)},\\
 G_a(y)&=\partial_a R_\theta(\lambda,y)
       =\beta R_\theta(\lambda,y)
                    \{m_a(\lambda)^{-1}-m_a(\lambda-y)^{-1}\}.
 \end{aligned}
\]
Substituting the preceding bounds in these formulas gives
\[
 |G_\nu(y)|+|G_a(y)|
 \le C\langle y\rangle^{2\beta_+}(1+\log\langle y\rangle).
\]
This supplies an integrable bound after multiplication by $k_\chi$ and justifies differentiating under the integral. More specifically, Schwartz decay makes $k_\chi(y)|G_j(y)|$ at most $C_N\langle y\rangle^{-N-d-1}$ for every fixed $N$.
The same bracket extraction used for the variance tail therefore gives
\[
 \int_{D_\lambda^c}k_\chi(y)|G_j(y)|\,dy
 \le C_N\int_{\norm y>\norm\lambda/2}
                          \langle y\rangle^{-N-d-1}\,dy
 \le C_N\langle\lambda\rangle^{-N}.
\]
Combining the ball and its complement, and taking $N\ge4$, yields
\[
 \partial_\nu\{v_\theta/f_\theta\}
       =O(\langle\lambda\rangle^{-2}),\qquad
 \partial_a\{v_\theta/f_\theta\}
       =O(\langle\lambda\rangle^{-4}).
\]
Substitute these bounds in the identity at the start of this paragraph and use
\[
 \partial_\nu\log f_\theta=-\log m_a,\qquad
 \partial_a\log f_\theta=-\beta/m_a.
\]
The positive lower bound on $v_\theta/f_\theta$ then proves the stated logarithmic-derivative estimates. Bounded $\lambda$ are covered by enlarging the constants. All estimates are uniform on $\Theta$.

\noindent\textbf{Covariances and energy bounds.}

For the off-diagonal bound, use the polynomial comparison for $R_\theta$ just proved, once from $\lambda$ and once from $\mu$, and take their geometric mean:
\[
 f_\theta(\omega)
 \le C\sqrt{f_\theta(\lambda)f_\theta(\mu)}
       \langle\lambda-\omega\rangle^{\beta_+}
       \langle\mu-\omega\rangle^{\beta_+}.
\]
Put $a(\eta)=\langle\eta\rangle^{\beta_+}|\widehat\chi(\eta)|$. After division by the two standard deviations, the bound for $\Sigma_\theta$ is at most $C\int a(\eta)a(\eta+\mu-\lambda)\,d\eta$. For $\delta=\mu-\lambda$, the bracket inequality
$\langle\delta\rangle\le\sqrt2\langle\eta\rangle \langle\eta+\delta\rangle$ and Cauchy--Schwarz imply, for every $M>0$,
\begin{align*}
 \langle\delta\rangle^M\int a(\eta)a(\eta+\delta)\,d\eta
 &\le C_M\int
       \{\langle\eta\rangle^M a(\eta)\}
       \{\langle\eta+\delta\rangle^M a(\eta+\delta)\}\,d\eta\\*
 &\le C_M\norm{\langle\cdot\rangle^M a}_{L^2}^2<\infty.
\end{align*}
The last norm is finite because $\widehat\chi$ is Schwartz. This proves decay by $\langle\mu-\lambda\rangle^{-M}$. For $\Pi_\theta$, use $X_\mu=\overline{X_{-\mu}}$:
\[
 \Pi_\theta(\lambda,\mu)
 =\E_\theta[X_\lambda\overline{X_{-\mu}}]
 =\Sigma_\theta(\lambda,-\mu),\qquad
 v_\theta(-\mu)=v_\theta(\mu).
\]
The bound just proved, with $-\mu$ in place of $\mu$, therefore gives
\[
 \frac{|\Pi_\theta(\lambda,\mu)|}
      {\sqrt{v_\theta(\lambda)v_\theta(\mu)}}
 =\frac{|\Sigma_\theta(\lambda,-\mu)|}
       {\sqrt{v_\theta(\lambda)v_\theta(-\mu)}}
 \le C_M\langle\lambda+\mu\rangle^{-M}.
\]
This proves \eqref{eq:fouriercovdecay}.

Fix a finite $\Lambda\subset\mathscr L$ and put $U=(U_\lambda)_{\lambda\in\Lambda}$, where $U_\lambda=X_\lambda/\sqrt{v_\theta(\lambda)}$. Write
\[
 \Sigma=\E_\theta[UU^*],\qquad
 \Pi=\E_\theta[UU^\T]
\]
for its complex covariance and relation matrices. These are the normalized versions of $\Sigma_\theta$ and $\Pi_\theta$ above. For $M>d$, translations and reflections of the lattice give
\[
 \sup_{\lambda\in\mathscr L}\sum_{\mu\in\mathscr L}
       \{\langle\lambda-\mu\rangle^{-M}
                  +\langle\lambda+\mu\rangle^{-M}\}
 =2\sum_{\gamma\in\mathscr L}\langle\gamma\rangle^{-M}<\infty.
\]
Thus the absolute row and column sums of $\Sigma$ and $\Pi$ are bounded independently of $\Lambda$. In particular, $\norm\Sigma_\op+\norm\Pi_\op\le C$.

The Gaussian fourth-moment formula gives the real symmetric, entrywise nonnegative matrix
\[
 A_{\lambda\mu}
 \coloneqq \Cov_\theta(|U_\lambda|^2,|U_\mu|^2)
 =|\Sigma_{\lambda\mu}|^2+|\Pi_{\lambda\mu}|^2.
\]
\begin{samepage}
The squared decay bounds similarly imply $\sup_{\lambda\in\Lambda}\sum_{\mu\in\Lambda}A_{\lambda\mu}\le C$. For real weights, expand the variance and use $2|w_\lambda w_\mu|\le w_\lambda^2+w_\mu^2$:
\begin{align*}
 \Var_\theta\left\{\sum_{\lambda\in\Lambda}
                   w_\lambda(|U_\lambda|^2-1)\right\}
 &=\sum_{\lambda,\mu\in\Lambda}
                    w_\lambda w_\mu A_{\lambda\mu}\\*
 &\le\frac12\sum_{\lambda,\mu\in\Lambda}
                    (w_\lambda^2+w_\mu^2)A_{\lambda\mu}\\*
 &=\sum_{\lambda\in\Lambda}w_\lambda^2
                         \sum_{\mu\in\Lambda}A_{\lambda\mu}
 \le C\sum_{\lambda\in\Lambda}w_\lambda^2.
\end{align*}
\end{samepage}
The penultimate equality uses $A_{\lambda\mu}=A_{\mu\lambda}$. This proves \eqref{eq:fourierenergyvariance}, including signed weights.

Finally, write $U=U_{\mathrm R}+iU_{\mathrm I}$. Expanding $\Sigma=\E_\theta[UU^*]$ and $\Pi=\E_\theta[UU^\T]$ gives
\[
 \Cov_\theta\begin{pmatrix}U_{\mathrm R}\\U_{\mathrm I}\end{pmatrix}
 =\frac12
 \begin{pmatrix}
  \Re(\Sigma+\Pi)&\Im(\Pi-\Sigma)\\
  \Im(\Pi+\Sigma)&\Re(\Sigma-\Pi)
 \end{pmatrix},
\]
where real and imaginary parts are taken entrywise. For example, the sum of absolute values in a row of the upper two blocks is at most
\[
 \frac12\sum_{\mu\in\Lambda}
       \{|\Sigma_{\lambda\mu}+\Pi_{\lambda\mu}|
                    +|\Pi_{\lambda\mu}-\Sigma_{\lambda\mu}|\}
 \le\sum_{\mu\in\Lambda}
       (|\Sigma_{\lambda\mu}|+|\Pi_{\lambda\mu}|)\le C.
\]
The lower block rows obey the same bound. The real covariance matrix is symmetric and positive semidefinite, so its operator norm is at most $C$ as well. This allows arbitrary dependence between real and imaginary parts and includes the relation
$X_{-\lambda}=\overline{X_\lambda}$.

\noindent\textbf{Envelope for the normalized Fourier weights.}
The variance expansion \eqref{eq:fourierlocalvariance}, together with positivity and continuity on bounded frequency sets, gives $v_0(\lambda)\asymp f_0(\lambda)
\asymp\langle\lambda\rangle^{-2\beta_0}$. By evenness of $f_0$ and the spectral-ratio bound,
\[
 \frac{f_0(\omega)}{f_0(\lambda)}
 =\frac{f_0(\omega)}{f_0(-\lambda)}
 \le C\langle\omega+\lambda\rangle^{2\beta_0}.
\]
It follows from \eqref{eq:fourierweights} that
\[
 \begin{aligned}
 \sum_{\lambda\in\Lambda}
       f_0(\omega)|\ell_\lambda(\omega)|^2
 &=L^{-d}\sum_{\lambda\in\Lambda}
       \frac{f_0(\omega)}{v_0(\lambda)}
                         |\widehat\chi(\omega+\lambda)|^2\\
 &\le C\sum_{\lambda\in\mathscr L}
       \langle\omega+\lambda\rangle^{2\beta_0}
                         |\widehat\chi(\omega+\lambda)|^2
 \le C.
 \end{aligned}
\]
For the last step, Schwartz decay bounds the summand by $C_M\langle\omega+\lambda\rangle^{-M}$ for any fixed $M>d$. The sum of these brackets is uniformly bounded in $\omega$: translation by a lattice vector reduces $\omega$ to a fixed fundamental cell, where the tail is bounded by the convergent lattice sum $C\sum_{\lambda\in\mathscr L} \langle\lambda\rangle^{-M}$.

If $\norm\lambda\le B$, we also have $v_0(\lambda)^{-1}\le C\langle B\rangle^{2\beta_0}$. Using this bound in the first line instead gives
\[
 \begin{aligned}
 \sum_{\lambda\in\Lambda}
       f_0(\omega)|\ell_\lambda(\omega)|^2
 &\le C\langle B\rangle^{2\beta_0}f_0(\omega)
                \sum_{\lambda\in\mathscr L}
                            |\widehat\chi(\omega+\lambda)|^2\\
 &\le C\langle B\rangle^{2\beta_0}f_0(\omega).
 \end{aligned}
\]
Taking the smaller of the two bounds proves part~\ref{it:windows-envelope}. For $B=\epsilon R_n\ge1$, we have $\langle B\rangle^{2\beta_0}\le CR_n^{2\beta_0}=Cn$;
the bounded positive $t_0$ then gives its stated consequence.

\noindent\textbf{Random coverage.}
Partition $Q$ into cubes of side comparable to $s_n=(A\log n/n)^{1/d}$, where $A$ is a fixed constant to be chosen. There are at most $Cs_n^{-d}$ cells, and the density
lower bound makes the probability that a particular cell is empty at most $e^{-cns_n^d}$. Hence
\[
 P_X\{\text{some cell is empty}\}
 \le C\frac{n}{A\log n}\,n^{-cA}.
\]
Choose $A$ large enough to make this summable. Borel--Cantelli then gives a site in every cell eventually almost surely, and $h_n\le Cs_n$. All constants above can be chosen uniformly on the fixed compact parameter sets in the statement. \qed

\subsection{Proof of Lemma~\ref{lem:preliminary}}
The variance of a localized Fourier coefficient near radius $B$ is approximately a fixed cutoff constant times $e^\eta(4\pi^2B^2)^{-\nu-d/2}$. Thus one band probes
$\eta-2\nu\log(2\pi B)$, and comparing bands at $B$ and $2B$ separates amplitude from smoothness. We compare covariance matrices of the full correlated coefficient vectors through affinity in the auxiliary Gaussian experiment $G_{\theta,n}$. These are deterministic covariance comparisons. Lemma~\ref{lem:affinitytails} then controls the Gaussian criterion under the actual response law.

Here the nearest-cell coefficients $\widetilde X_{\lambda,n}$ are convenient because their recovery error is uniform over all $\theta\in\Theta$. We may let $B$ grow slowly; no sharp information rate is needed yet. The best linear predictions used later in Lemma~\ref{lem:information} instead optimize recovery at the fixed truth, at frequencies approaching the resolution $R_n$.

\noindent\textbf{The observed log variance ratio.}
Let $B\to\infty$ and form the real vector
\[
 \mathbf X_B=(\sqrt2\Re X_\lambda,\sqrt2\Im X_\lambda)
                                      _{\lambda\in\Lambda_B^+}.
\]
Its observable version $\widetilde{\mathbf X}_{B,n}$ uses $\widetilde X_{\lambda,n}$ from \eqref{eq:fourierquadrature}. Write $\Gamma_{B,n}(\theta)$ and $\widetilde\Gamma_{B,n}(\theta)$ for their covariance matrices. Their size is $N_B=2\#\Lambda_B^+\asymp B^d$. A real combination $b^\T\mathbf X_B$ equals $Z_{g_b}$ for a cutoff trigonometric polynomial $g_b\in\mathcal V_B$; \eqref{eq:windownorm} gives
$\norm{g_b}_{L^2}^2\asymp\norm b^2$. Our object of interest is
\[
 \mathcal Q_{B,n}(\theta;b)
 =\log\frac{b^\T\widetilde\Gamma_{B,n}(\theta)b}
                {b^\T\widetilde\Gamma_{B,n}(\theta_0)b},
 \qquad b\ne0.
\]
We will approximate it uniformly in $b$ and $\theta$ by
\[
 U_B(\theta)=\eta-\eta_0-2v\log(2\pi B),\qquad v=\nu-\nu_0.
\]
Uniformity in $b$ will control all generalized eigenvalues, which is what the Gaussian affinity calculation needs.

First pass from recovered to latent Fourier coefficients. Write $\widetilde{\mathbf X}_{B,n}=\mathbf X_B+E$. With $e=\gamma\min(2\zeta,d)>0$ and $\kappa_+=2\nu_++d$,
Lemma~\ref{lem:windows} gives
\begin{equation}\label{eq:bandcovsizes}
 \Var_\theta(b^\T E)\le Cn^{-e}\norm b^2,
 \qquad b^\T\Gamma_{B,n}(\theta)b
 \ge cB^{-2\nu-d}\norm b^2\ge cB^{-\kappa_+}\norm b^2.
\end{equation}
Expanding the variance of $b^\T(\mathbf X_B+E)$ gives
\[
 b^\T\widetilde\Gamma_{B,n}(\theta)b
 =b^\T\Gamma_{B,n}(\theta)b
   +2\Cov_\theta(b^\T\mathbf X_B,b^\T E)
   +\Var_\theta(b^\T E).
\]
For $b\ne0$, \eqref{eq:bandcovsizes} implies
\[
 \frac{\Var_\theta(b^\T E)}{b^\T\Gamma_{B,n}(\theta)b}
 \le\frac{Cn^{-e}\norm b^2}{cB^{-\kappa_+}\norm b^2}
 \le Cn^{-e}B^{\kappa_+}.
\]
Cauchy--Schwarz bounds the covariance term by $\{b^\T\Gamma_{B,n}(\theta)b\}^{1/2}
\{\Var_\theta(b^\T E)\}^{1/2}$. Dividing the variance expansion by $b^\T\Gamma_{B,n}(\theta)b$ therefore gives
\begin{equation*}
 \begin{gathered}
 \begin{aligned}
 \left|\frac{b^\T\widetilde\Gamma_{B,n}(\theta)b}
                  {b^\T\Gamma_{B,n}(\theta)b}-1\right|
 &\le2\left(\frac{\Var_\theta(b^\T E)}
                    {b^\T\Gamma_{B,n}(\theta)b}\right)^{1/2}
       +\frac{\Var_\theta(b^\T E)}{b^\T\Gamma_{B,n}(\theta)b}\\
 &\le C\rho_n(B),
 \end{aligned}\\
 \rho_n(B)=n^{-e/2}B^{\kappa_+/2}+n^{-e}B^{\kappa_+}.
 \end{gathered}
\end{equation*}
Provided $\rho_n(B)=o(1)$, taking logarithms under $\theta$ and $\theta_0$ gives the first approximation:
\[
 \mathcal Q_{B,n}(\theta;b)
 =\log\frac{\Var_\theta(Z_{g_b})}{\Var_0(Z_{g_b})}
   +O(\rho_n(B)).
\]

The variance-ratio bound in Lemma~\ref{lem:windows}, with leakage exponent $2$, places the middle term between the infimum and supremum of $\log(f_\theta/f_0)$ on $\mathcal A_B$, up to $O(B^{-2})$. Finally, factor the leading power from the spectrum:
\begin{align*}
 \log\frac{f_\theta(\omega)}{f_0(\omega)}
 ={}&\eta-\eta_0-2v\log(2\pi\norm\omega)-\beta\log\left(1+\frac a{4\pi^2\norm\omega^2}\right)
 +\beta_0\log\left(1+\frac {a_0}{4\pi^2\norm\omega^2}\right).
\end{align*}
On $\mathcal A_B$ we have $\log(\norm\omega/B)=O(\omega_*)$; the two range terms are
$O(B^{-2})$. Thus
\begin{equation*}
 \log(f_\theta/f_0)=U_B+O(\omega_*|v|+B^{-2})
 \quad\text{on }\mathcal A_B.
\end{equation*}
Combining the three approximations gives the useful final bound
\begin{equation}\label{eq:bandeigen}
 \sup_{b\ne0}|\mathcal Q_{B,n}(\theta;b)-U_B(\theta)|
 \le C\omega_*|v|+C\{B^{-2}+\rho_n(B)\}.
\end{equation}
The constant in the bandwidth term is uniform for sufficiently small bandwidths; the remaining constants may depend on the fixed bandwidth.

\noindent\textbf{Two bands give separation.}
Replacing $B$ by $2B$ proves the same estimate for the second statistic, with
$U_{2B}=U_B-2v\log2$ and errors of the same orders. The linear map $(U_B,v)\mapsto(U_B,U_{2B})$ is invertible, so
\[
 \max\{|U_B|,|U_{2B}|\}\ge c_*s,
 \qquad s=(U_B^2+v^2)^{1/2},
\]
where $c_*>0$ is fixed. Select $B_*\in\{B,2B\}$ for which $|U_{B_*}|\ge c_*s$, and put $\sigma=\operatorname{sgn}(U_{B_*})$. Since $|v|\le s$, we can first fix the bandwidth so small that $C\omega_*\le c_*/2$. The approximation \eqref{eq:bandeigen}, applied at $B_*$, then gives, for every $b\ne0$,
\[
\begin{aligned}
 \sigma\mathcal Q_{B_*,n}(\theta;b)
 &\ge |U_{B_*}|-
       |\mathcal Q_{B_*,n}(\theta;b)-U_{B_*}|\\
 &\ge c_*s-C\omega_*|v|-C\{B^{-2}+\rho_n(B)\}\\
 &\ge \tfrac12c_*s-C\{B^{-2}+\rho_n(B)\}
 \ge \tfrac14c_*s,
\end{aligned}
\]
where the last inequality holds whenever $s\ge C\{B^{-2}+\rho_n(B)\}$ with a sufficiently large fixed $C$. Thus every log covariance quotient on the selected band has the sign of $U_{B_*}$ and absolute value at least $cs$. Each generalized eigenvalue is such a covariance quotient at a generalized eigenvector, so its logarithm obeys the same bound.

The selected statistic has $N_{B_*}\asymp B^d$ coordinates. By Lemma~\ref{lem:gaussian-affinity}, each of its log generalized eigenvalues, whose magnitude is at least $cs$, contributes at least $c'\min\{s^2,s\}$ to its negative log affinity. Lemma~\ref{lem:affinity-processing} transfers this lower bound to the full observations. Consequently
\begin{equation}\label{eq:bandseparation}
 \begin{gathered}
 -\log\Aff(G_{\theta,n},G_{0,n})
 \ge cB^d\min\{s^2,s\},
 \text{ if }
 s=(U_B^2+v^2)^{1/2}\ge C\{B^{-2}+\rho_n(B)\}.
 \end{gathered}
\end{equation}
For each candidate we use just one band; independence between bands is unnecessary.

\noindent\textbf{Separating the nugget.}
The latent field changes little over a short distance, whereas the independent measurement errors do not. To exploit this, place a lattice of spacing $Ah_n$ in a smaller cube, with $A$ fixed and large. Select an observation within $h_n$ of each lattice point and pair neighboring lattice points without reuse. The selected sites are distinct, and this provides at least $ch_n^{-d}$ disjoint pairs at distance at most $Ch_n$.

Keep $p$ pairs and put $D_j=(Y_{i_j}-Y_{i'_j})/\sqrt2$. Their covariance is $tI_p+H_{\theta,n}$, where $H_{\theta,n}\succeq0$ is the covariance of the latent increments. By \eqref{eq:variogram},
\[
 \norm{H_{\theta,n}}_\op\le\tr H_{\theta,n}
 \le Cp h_n^{2\zeta}.
\]
For a unit vector $b\in\R^p$, write $d_\theta(b)=b^\T H_{\theta,n}b$. Then
$0\le d_\theta(b)\le Cp h_n^{2\zeta}$, and the corresponding covariance quotient is
\[
 q_\theta(b)=\frac{t+d_\theta(b)}{t_0+d_0(b)},\qquad
 \left|q_\theta(b)-\frac t{t_0}\right|
 =\frac{|t_0d_\theta(b)-t d_0(b)|}
         {t_0\{t_0+d_0(b)\}}
 \le Cp h_n^{2\zeta}.
\]
Here we used the fixed positive lower and upper bounds on the nugget. If $\sigma_t=\operatorname{sgn}(t-t_0)$, this yields
\[
 \sigma_t\{q_\theta(b)-1\}
 \ge \frac{|t-t_0|}{t_0}-Cp h_n^{2\zeta}
 \ge c|t-t_0|,
 \qquad |t-t_0|\ge Cp h_n^{2\zeta},
\]
where the constant in the threshold is sufficiently large. When $p h_n^{2\zeta}=o(1)$, as our choices below will ensure, all the quotients $q_\theta(b)$ lie in a fixed compact subinterval of $(0,\infty)$. On that interval the same scalar comparison gives
\[
 \log\cosh(\tfrac12\log q)\asymp(\log q)^2\asymp(q-1)^2.
\]
Applying Lemmas~\ref{lem:gaussian-affinity} and \ref{lem:affinity-processing} to these $p$ generalized eigenvalues gives
\[
 -\log\Aff(G_{\theta,n},G_{0,n})\ge cp|t-t_0|^2,
 \qquad |t-t_0|\ge Cp h_n^{2\zeta}.
\]

\begin{samepage}
We now want to choose the scales and a threshold $n^{-\delta}$ so that, uniformly over candidate parameters $\theta$,
\[
 \begin{gathered}
 |U_{B_n}(\theta)|+|\nu-\nu_0|+|t-t_0|>Cn^{-\delta}\\
 \Longrightarrow\quad
 -\log\Aff(G_{\theta,n},G_{0,n})\ge cn^{\delta_1}
 \quad\text{for some }\delta_1>0.
 \end{gathered}
\]
\end{samepage}
The likelihood--affinity bound \eqref{eq:likelihoodaffinity} will then make it exponentially unlikely that any such fixed candidate has likelihood close to the true likelihood. A parameter mesh will turn this into a statement about every
global maximizer.

\noindent\textbf{Choosing scales and covering the likelihood.}
Set $B=B_n=n^\epsilon$ and $p=p_n=\lfloor n^\xi\rfloor$, with $\epsilon,\xi>0$. The coverage bound $h_n\le Cn^{-\gamma}$ gives
\[
 \begin{aligned}
 \rho_n(B_n)=n^{-(e-\epsilon\kappa_+)/2}
                +n^{-(e-\epsilon\kappa_+)},\qquad
 p_nh_n^d\le Cn^{-(\gamma d-\xi)},\qquad
 p_nh_n^{2\zeta}\le Cn^{-(2\gamma\zeta-\xi)}.
 \end{aligned}
\]
Before choosing $\delta$, we need $\rho_n(B_n)\to0$ to justify the recovered log ratios, $p_n=o(h_n^{-d})$ to have enough disjoint pairs, and $p_nh_n^{2\zeta}\to0$ to make their latent covariance perturbation small. These three requirements hold on
choosing
\[
 0<\epsilon<e/\kappa_+,\qquad
 0<\xi<\min\{\gamma d,2\gamma\zeta\}.
\]

Next, the approximation errors must be smaller than the proposed separation threshold, while the affinity exponents at that threshold must grow. Thus we require
\[
 \begin{aligned}
 B_n^{-2}+\rho_n(B_n)&=o(n^{-\delta}),&
 p_nh_n^{2\zeta}&=o(n^{-\delta}),\\
 B_n^dn^{-2\delta}&\longrightarrow\infty,&
 p_nn^{-2\delta}&\longrightarrow\infty.
 \end{aligned}
\]
All four conditions follow by choosing $\delta>0$ so that
\[
 \delta<\min\left\{2\epsilon,
 \frac{e-\epsilon\kappa_+}{2},\frac{\epsilon d}{2},
 \frac\xi2,2\gamma\zeta-\xi\right\}.
\]
Every entry of this minimum is positive by the choices of $\epsilon$ and $\xi$.

To verify the desired implication, note that $|U_{B_n}|+|v|\le\sqrt2s$. Thus a candidate outside the stated threshold has either $s\ge cn^{-\delta}$ or
$|t-t_0|\ge cn^{-\delta}$. The error conditions ensure that the appropriate separation bound applies for all sufficiently large $n$. Since $n^{-\delta}\to0$, the first case gives $\min\{s^2,s\}\ge cn^{-2\delta}$, so \eqref{eq:bandseparation} yields
$-\log\Aff\ge cn^{\epsilon d-2\delta}$; in the second, the pair statistic gives $-\log\Aff\ge cn^{\xi-2\delta}$. Consequently the desired bound holds with
\[
 \delta_1=\min\{\epsilon d-2\delta,\xi-2\delta,1\}>0.
\]

It remains to rule out all distant candidates whose likelihood is at least that of the truth. Write
\[
 d_n(\theta)=|U_{B_n}(\theta)|+|\nu-\nu_0|+|t-t_0|,
\]
and fix $A>0$ large enough that the preceding affinity bound applies whenever $d_n(\theta)>An^{-\delta}$. We want to bound the probability of the event
\[
 \mathcal E_n=
 \{\exists\theta\in\Theta:\ d_n(\theta)>2An^{-\delta},\quad
                       \ell_n(\theta)\ge\ell_n(\theta_0)\}.
\]
We will show that any such candidate gives a distant point on a deterministic finite
mesh whose log-likelihood is at least $\ell_n(\theta_0)-1$, except on an exponentially unlikely event.

\emph{Controlling the likelihood under rounding.}
Under the actual law, put $X=C_{n,0}^{-1/2}Y_n$. Apply the sub-Gaussian quadratic-form bound \eqref{eq:quadconc} in Lemma~\ref{lem:gaussian} to $X^\T I_nX$ with $z=n$. Since $\norm{I_n}_\Frob=\sqrt n$ and $\norm{I_n}_\op=1$, it gives
\[
 P_0\{\norm X^2>Cn\}\le2e^{-n}.
\]
The marginal variances are uniformly bounded on the compact parameter set, so
$\norm{C_{n,0}}_\op\le\tr C_{n,0}\le Cn$ and hence $\norm{Y_n}^2\le Cn\norm X^2$. Consequently, for a sufficiently large fixed $C_Y$,
\[
 \mathcal Y_n=\{\norm{Y_n}^2\le C_Yn^2\},\qquad
 P_0(\mathcal Y_n^c)\le2e^{-n}.
\]

The covariance entries and their first parameter derivatives are uniformly bounded: for the signal part this follows by integrating the absolute derivatives of \eqref{eq:spectrum}. Thus, with $C_{n,j}=\partial_{\theta_j}C_n$, $\sup_{\mathcal R}\norm{C_{n,j}}_\op\le Cn$, while $\sup_{\mathcal R}\norm{C_n^{-1}}_\op\le t_-^{-1}$.
Differentiating \eqref{eq:likelihood} gives
\[
 \partial_{\theta_j}\ell_n
 =-\tfrac12\tr(C_n^{-1}C_{n,j})
   +\tfrac12Y_n^\T C_n^{-1}C_{n,j}C_n^{-1}Y_n.
\]
The trace term satisfies
\[
 \frac12|\tr(C_n^{-1}C_{n,j})|
 \le\frac n2\|C_n^{-1}\|_\op\|C_{n,j}\|_\op\le Cn^2.
\]
The quadratic term is bounded by $Cn\norm{Y_n}^2$. Hence, on $\mathcal Y_n$,
\[
 \sup_{\theta\in\mathcal R}\norm{\nabla\ell_n(\theta)}
 \le C\{n^2+n\norm{Y_n}^2\}\le Cn^3.
\]

Apply the preceding affinity separation on the enclosing family $\mathcal R$ from Appendix~\ref{sec:proofnotation}; its proof uses only the compact parameter bounds. Let $\mathcal M_n\subset\mathcal R$ be a rectangular mesh of spacing at most $n^{-K}$, including the endpoints in each coordinate, and let $\pi_n(\theta)$ be a nearest
mesh point. This point need not belong to $\Theta$: it is used only to bound the probability of a likelihood event, and the criterion is defined on all of $\mathcal R$. The segment joining it to $\theta$ lies in $\mathcal R$, so the gradient bound applies. Since the enclosing rectangle has four coordinates, $\#\mathcal M_n\le Cn^{4K}$ and $\norm{\pi_n(\theta)-\theta}\le Cn^{-K}$. The gradient bound and the formula for $U_{B_n}$ yield
\[
 \begin{aligned}
 \sup_{\theta\in\mathcal R}
 |\ell_n(\pi_n(\theta))-\ell_n(\theta)|
 &\le Cn^{3-K} &&\text{on }\mathcal Y_n,\\
 \sup_{\theta\in\mathcal R}
 |d_n(\pi_n(\theta))-d_n(\theta)|
 &\le Cn^{-K}\log n.
 \end{aligned}
\]
The factor $\log n$ comes from the coefficient $2\log(2\pi B_n)$ of $\nu-\nu_0$ in $U_{B_n}$. Thus any fixed $K>\max\{3,\delta\}$ makes the two errors $o(1)$ and $o(n^{-\delta})$, respectively. For a concrete choice, take $K=4$: our earlier choices give
$\delta<\xi/2<\gamma d/2\le1/2$. In particular, $\#\mathcal M_n\le Cn^{16}$.

\emph{From a distant candidate to a distant mesh point.}
On $\mathcal E_n\cap\mathcal Y_n$, take any $\theta$ witnessing $\mathcal E_n$ and set $\theta^*=\pi_n(\theta)$. For all sufficiently large $n$, the rounding bounds give
\[
 \begin{aligned}
 d_n(\theta^*)&\ge d_n(\theta)-o(n^{-\delta})>An^{-\delta},\\
 \ell_n(\theta^*)&\ge\ell_n(\theta)-o(1)\ge\ell_n(\theta_0)-1.
 \end{aligned}
\]
Thus the event of interest is contained in a finite union:
\[
 \mathcal E_n\subseteq\mathcal Y_n^c\ \cup
 \bigcup_{\substack{\theta\in\mathcal M_n\\
                     d_n(\theta)>An^{-\delta}}}
 \{\ell_n(\theta)-\ell_n(\theta_0)\ge-1\}.
\]
For each mesh point in this union, \eqref{eq:likelihoodaffinity} with $z=1$ and the affinity bound above give
\[
 P_0\{\ell_n(\theta)-\ell_n(\theta_0)\ge-1\}
 \le Ce^{1/2}\Aff(G_{\theta,n},G_{0,n})^c
 \le Ce^{-cn^{\delta_1}}.
\]
Taking the union bound therefore proves
\[
 P_0(\mathcal E_n)
 \le2e^{-n}+Cn^{16}e^{-cn^{\delta_1}}
 \le Ce^{-c'n^{\delta_1}}.
\]
Here $\delta_1\le1$ handles the first term, and the polynomial factor is absorbed by decreasing the exponential constant: $16\log n\le(c/2)n^{\delta_1}$ for all sufficiently large $n$. Every global MLE has likelihood at least that of the truth, so this also bounds the probability that some MLE has $d_n(\widehat\theta_n)>2An^{-\delta}$. Since $\eta-\eta_0=U_{B_n}+2v\log(2\pi B_n)$, this proves \eqref{eq:prelim}. Also $u=U_{B_n}+2\{\log(2\pi B_n)-L_n\}v$. For any fixed $0<\delta'<\delta$, the bounds
$|u|\le Cn^{-\delta}\log n$ and $|v|+|t-t_0|\le Cn^{-\delta}$ imply \eqref{eq:localset}
with $q_n=n^{-\delta'}$, since $Cn^{-\delta}\le n^{-\delta'}$ eventually. This gives \eqref{eq:prelimfailure} with the same exceptional event. Relabel $\delta'$ as $\delta$; \eqref{eq:prelim} continues to hold with this smaller exponent. All scale choices use only the fixed compact and coverage constants. \qed

\subsection{Proof of Lemma~\ref{lem:upper}}
The main step is a bound on the contribution of one frequency to the whitened covariance. Once that bound is available, all the upper estimates reduce to scalar integrals. This separates two effects: boundedness of the observation region limits how
much a frequency can contribute, and the nugget suppresses frequencies where the signal spectrum is small. No sampling lower bound is used in this proof.

\noindent\textbf{A bound for one frequency.}
Use the true covariance $C_{n,0}=S_{n,0}+t_0I_n$ and write
\[
 a_\omega=(e^{2\pi i\omega^\T x_{n,i}})_{i=1}^n,\qquad
 w_n(\omega)
   =f_0(\omega)a_\omega^*C_{n,0}^{-1}a_\omega.
\]
This nonnegative function measures the contribution of frequency $\omega$ after whitening. The spectral covariance formula is $S_{n,0}=\int f_0(\omega)a_\omega a_\omega^*\,d\omega$. Consequently,
\[
 \tr(S_{n,0}C_{n,0}^{-1})
       =\int_{\R^d}w_n(\omega)\,d\omega.
\]
We claim the pointwise estimate
\begin{equation}\label{eq:frequencyweight}
 w_n(\omega)
       \le C\min\{1,(n/t_0)f_0(\omega)\},
       \qquad \omega\in\R^d.
\end{equation}
The second bound is immediate: $C_{n,0}^{-1}\preceq t_0^{-1}I_n$ and $\norm{a_\omega}^2=n$ give $w_n(\omega)\le(n/t_0)f_0(\omega)$. The first bound is the only analytic ingredient to prove.

Because $D$ is bounded, choose $\chi_D\in C_c^\infty(\R^d)$ equal to one on a neighborhood of $\overline D$. Fourier inversion, evaluated at each site, gives
\[
 e^{2\pi i\omega^\T x}
 =\int_{\R^d}\widehat{\chi}_D(\xi-\omega)
                 e^{2\pi i\xi^\T x}\,d\xi,\qquad x\in D,
 \qquad
 a_\omega=\int_{\R^d}\widehat{\chi}_D(\xi-\omega)a_\xi\,d\xi.
\]
Thus any linear combination of the sampled wave at $\omega$ can be expressed using nearby frequencies, with rapidly decreasing weights. For $c\in\mathbb C^n$, weighted
Cauchy--Schwarz makes this statement quantitative:
\begin{align*}
 |c^*a_\omega|^2
 &\le
 \left(\int\frac{|\widehat{\chi}_D(\xi-\omega)|^2}{f_0(\xi)}
                       \,d\xi\right)
 \left(\int f_0(\xi)|c^*a_\xi|^2\,d\xi\right)=J(\omega)c^*S_{n,0}c,\\
 J(\omega)&\coloneqq \int
          \frac{|\widehat{\chi}_D(\xi-\omega)|^2}{f_0(\xi)}\,d\xi.
\end{align*}
Since $f_0(\xi)^{-1}\asymp\langle\xi\rangle^{2\beta_0}$, the Japanese-bracket inequality
$\langle\xi\rangle\le\sqrt2\langle\omega\rangle
\langle\xi-\omega\rangle$ gives
\begin{align*}
 J(\omega)
 &\asymp\int_{\R^d}\langle\xi\rangle^{2\beta_0}
                   |\widehat{\chi}_D(\xi-\omega)|^2\,d\xi\\
 &\le C\langle\omega\rangle^{2\beta_0}
       \int_{\R^d}\langle\xi-\omega\rangle^{2\beta_0}
                   |\widehat{\chi}_D(\xi-\omega)|^2\,d\xi\\
 &=C\langle\omega\rangle^{2\beta_0}
       \int_{\R^d}\langle\zeta\rangle^{2\beta_0}
                   |\widehat{\chi}_D(\zeta)|^2\,d\zeta\\
 &\le C\langle\omega\rangle^{2\beta_0}
 \asymp f_0(\omega)^{-1}.
\end{align*}
Here we substituted $\zeta=\xi-\omega$. The last integral is bounded uniformly for $\beta_0$ in its fixed compact interval, because $\widehat{\chi}_D$ is a Schwartz function.

Take the square root of the weighted Cauchy--Schwarz bound, set $c=C_{n,0}^{-1}a_\omega$, and use $S_{n,0}\preceq C_{n,0}$:
\begin{align*}
 a_\omega^*C_{n,0}^{-1}a_\omega
 &=(C_{n,0}^{-1}a_\omega)^*a_\omega\\
 &\le Cf_0(\omega)^{-1/2}
       \left(a_\omega^*C_{n,0}^{-1}S_{n,0}
                         C_{n,0}^{-1}a_\omega\right)^{1/2}\\
 &\le Cf_0(\omega)^{-1/2}
       \left(a_\omega^*C_{n,0}^{-1}a_\omega\right)^{1/2}.
\end{align*}
Dividing by the positive square root of $a_\omega^*C_{n,0}^{-1}a_\omega$ and squaring gives
\[
 a_\omega^*C_{n,0}^{-1}a_\omega\le Cf_0(\omega)^{-1}.
\]
Multiplication by $f_0(\omega)$ proves the first bound in \eqref{eq:frequencyweight}.
This argument only needs a smooth cutoff around the bounded set $D$; it imposes no boundary regularity on $D$.

\noindent\textbf{Effective dimension.}
Put $z=n/t_0>0$, recall $\kappa_0=2\nu_0+d=2\beta_0$, and set $R=z^{1/\kappa_0}$. Since
$f_0(\omega)\asymp(1+\norm\omega^2)^{-\beta_0}$, \eqref{eq:frequencyweight}, polar coordinates $\omega=r\vartheta$, and a split at $r=R$ give
\begin{align*}
 \tr(S_{n,0}C_{n,0}^{-1})
 &\le C\int_{\R^d}
       \min\{1,z(1+\norm\omega^2)^{-\beta_0}\}\,d\omega\\
 &=C|S^{d-1}|\int_0^\infty
       \min\{1,z(1+r^2)^{-\beta_0}\}r^{d-1}\,dr\\
 &\le C\left\{\int_0^R r^{d-1}\,dr
       +z\int_R^\infty(1+r^2)^{-\beta_0}r^{d-1}\,dr\right\}\\
 &\le C\left\{\frac{R^d}{d}
       +z\int_R^\infty r^{-2\nu_0-1}\,dr\right\}\\
 &=C\left\{\frac{R^d}{d}
       +\frac{zR^{-2\nu_0}}{2\nu_0}\right\}
 \le C z^{d/\kappa_0}.
\end{align*}
Since $t_0$ stays in a fixed positive compact interval and $M_n=R_n^d=n^{d/\kappa_0}$, this gives $\tr P_n=\tr(S_{n,0}C_{n,0}^{-1})\le CM_n$, proving \eqref{eq:effective}.

\noindent\textbf{Converting scalar integrals into matrix bounds.}
For a real multiplier $h$, the integral representation is
\[
 \mathcal B_n\bigl(h\bigr)
 =\int_{\R^d}h(\omega)f_0(\omega)
        (C_{n,0}^{-1/2}a_\omega)
        (C_{n,0}^{-1/2}a_\omega)^*\,d\omega.
\]
For $h\ge0$, taking the trace gives $\tr(\mathcal B_n\bigl(h\bigr))=\int hw_n$. For a signed $h$, the triangle inequality for the trace norm of this rank-one integral gives
$\norm{\mathcal B_n\bigl(h\bigr)}_1\le\int|h|w_n$. There is also a useful squared Frobenius estimate:
\begin{equation}\label{eq:weightedmatrixbounds}
 \norm{\mathcal B_n\bigl(h\bigr)}_1\le\int|h|w_n,\qquad
 \norm{\mathcal B_n\bigl(h\bigr)}_\Frob^2\le\int h^2w_n.
\end{equation}
For the second inequality, set $z_\omega=f_0(\omega)^{1/2}C_{n,0}^{-1/2}a_\omega$ within this calculation. Then $\mathcal B_n(h)=\int h z_\omega z_\omega^*$. For $u,v\in\mathbb C^n$, Cauchy--Schwarz gives
\begin{align*}
 |u^*\mathcal B_n(h)v|^2
 &\le\left(\int h^2|u^*z_\omega|^2\right)
       \left(\int|z_\omega^*v|^2\right)\\
 &=(u^*\mathcal B_n(h^2)u)(v^*P_nv)
 \le(u^*\mathcal B_n(h^2)u)\norm v^2.
\end{align*}
Taking the supremum over unit $v$, using that $\mathcal B_n(h)$ is Hermitian, and then taking the trace yields
\[
 \mathcal B_n(h)^2\preceq\mathcal B_n(h^2),\qquad
 \norm{\mathcal B_n(h)}_\Frob^2
 \le\tr\mathcal B_n(h^2)=\int h^2w_n.
\]
One can first truncate $h$ to make it bounded and then pass to the limit whenever the displayed integrals are finite.

\noindent\textbf{Range derivatives: where the dimension thresholds arise.}
For any fixed $p>0$, split the integral $I_{p,n}=\int m_0^{-p}w_n$ at radii $1$ and $R_n$. The bounds $w_n\le C$ below $R_n$ and $w_n\le Cnf_0$ above $R_n$ yield
\begin{align*}
 I_{p,n}
 &\le C\left\{1+\int_1^{R_n}r^{d-2p-1}\,dr
                   +n\int_{R_n}^\infty
                              r^{-2\nu_0-2p-1}\,dr\right\}\\
 &\le C\left\{1+\int_1^{R_n}r^{d-2p-1}\,dr
                   +R_n^{d-2p}\right\}.
\end{align*}
The last equality in order uses $n=R_n^{2\nu_0+d}$. The middle integral is bounded if $d<2p$, equals $\log R_n$ if $d=2p$, and is of order $R_n^{d-2p}$ if $d>2p$. Therefore
\begin{equation}\label{eq:scalarpowerbound}
 I_{p,n}\le C
 \begin{cases}
 1,&d<2p,\\
 \log n,&d=2p,\\
 R_n^{d-2p},&d>2p.
 \end{cases}
\end{equation}
The trace identity with $p=k$ proves \eqref{eq:rangetrace}. The squared Frobenius inequality in \eqref{eq:weightedmatrixbounds} with $h=m_0^{-k}$ uses $p=2k$ and proves
\eqref{eq:rangeupper}. The critical dimension is $4k$ for the Frobenius norm because the multiplier is squared.

\noindent\textbf{Centered logarithms: why there is no extra $\log n$ in the Frobenius bound.}
Recall $g_n=2\log R_n-\log m_0$. Its value can be of order $\log n$ at low frequencies, but most contributing frequencies have size comparable to $R_n$, where centering makes it bounded. To quantify this, substitute $\omega=R_nz$:
\[
 g_n(R_nz)=-\log(a_0R_n^{-2}+4\pi^2\norm z^2),\qquad
 w_n(R_nz)\le C\min\{1,\norm z^{-2\beta_0}\}.
\]
The second bound follows from \eqref{eq:frequencyweight} and $n=R_n^{2\beta_0}$. For $r=\norm z>0$,
\[
 |g_n(R_nz)|\le C(1+|\log r|).
\]
For $r\ge1$, this follows from $4\pi^2r^2\le a_0R_n^{-2}+4\pi^2r^2 \le(a_++4\pi^2)r^2$. For $0<r\le1$, the same expression is at least $4\pi^2r^2$ and at most $a_++4\pi^2$; these bounds control respectively its negative and positive
logarithms.

Thus polar coordinates give, for every fixed integer $j\ge0$,
\begin{align*}
 \int g_n^{2j}w_n
 &\le C R_n^d\left\{
       \int_0^1 r^{d-1}(1+|\log r|)^{2j}\,dr
       +\int_1^\infty r^{-2\nu_0-1}(1+\log r)^{2j}\,dr
                       \right\}\\
 &\le C_j M_n.
\end{align*}
The substitutions $s=-\log r$ on $(0,1)$ and $s=\log r$ on $(1,\infty)$ give, respectively,
\[
 \int_0^\infty e^{-ds}(1+s)^{2j}\,ds,\qquad
 \int_0^\infty e^{-2\nu_0s}(1+s)^{2j}\,ds.
\]
For any $c>0$, expanding the polynomial and integrating by parts $\ell$ times yields
\[
 \int_0^\infty e^{-cs}(1+s)^{2j}\,ds
 =\sum_{\ell=0}^{2j}\binom{2j}{\ell}
                         \frac{\ell!}{c^{\ell+1}}.
\]
Both expressions are finite; the second is uniformly bounded because $\nu_0\ge\nu_->0$. This proves \eqref{eq:logmoments}, and \eqref{eq:weightedmatrixbounds}
with $h=g_n^j$ proves the Frobenius part of \eqref{eq:logbounds}.

For the operator norm, use a larger cutoff $\norm\omega=n^A$, where $A$ is a fixed sufficiently large constant. Below this cutoff, $|g_n|^j\le C_A(\log n)^j$. Order preservation and $\mathcal B_n(1)\preceq I_n$ bound the operator norm of this part by $C_A(\log n)^j$. The trace norm of the complementary term is at most
\[
 Cn\int_{\norm\omega>n^A}f_0|g_n|^j\,d\omega
 \le C_A n^{1-2A\nu_-}(\log n)^j=o(1).
\]
This tail estimate follows from the same substitution and exponential integral just calculated. It proves the operator part of \eqref{eq:logbounds}.

\noindent\textbf{Mixed derivatives and uniformity.}
We first compute the multipliers, then bound their contributions below and above a frequency cutoff. Recall $f_\theta=e^\eta m_a^{-\beta}$ with $\beta=\nu+d/2$.
For the rising factorial $(\beta)_k=\beta(\beta+1)\cdots(\beta+k-1)$, with $(\beta)_0=1$, we have
\[
 \partial_a^k f_\theta=(-1)^k(\beta)_k m_a^{-k}f_\theta,
 \qquad D_{\nu,n}f_\theta=g_{n,a}f_\theta.
\]
Since $D_{\nu,n}\beta=1$ while $D_{\nu,n}m_a=D_{\nu,n}g_{n,a}=0$, the product rule gives
\[
 D_{\nu,n}^{\,j}\partial_a^k f_\theta
 =(-1)^k m_a^{-k}f_\theta
   \sum_{\ell=0}^{\min(j,k)}
       \binom{j}{\ell}
       \frac{d^\ell(\beta)_k}{d\beta^\ell}\,
       g_{n,a}^{\,j-\ell}.
\]
For example, $D_{\nu,n}\partial_a f_\theta =-m_a^{-1}f_\theta(1+\beta g_{n,a})$. All coefficients in the displayed polynomial are bounded on the compact parameter set for fixed $j,k$. Hence
\[
 \bigl|D_{\nu,n}^{\,j}\partial_a^k f_\theta\bigr|
 \le C_{j,k}f_\theta m_a^{-k}(1+|g_{n,a}|)^j.
\]
Also $\partial_\eta f_\theta=f_\theta$, so extra $\eta$ derivatives do not change these expressions. Differentiation under the spectral covariance integral is valid
because these multipliers are integrable uniformly on the compact parameter set; the tail estimate below makes this explicit.

Fix a constant $A>0$ and write $E_{n,A}=\{\omega:\norm\omega\le n^A\}$. Uniformly over the range interval, $m_a\asymp m_0$, and
\[
 \log\frac{f_\theta}{f_0}
 =(\eta-\eta_0)-(\nu-\nu_0)\log m_0
       -\beta\log(m_a/m_0).
\]
On $\mathcal N_n$, $|\eta-\eta_0|\le |u|+2L_n|v|\le Cq_n\log n=o(1)$. On $E_{n,A}$ we also have $|\log m_0|\le C_A\log n$, so $|\nu-\nu_0|\,|\log m_0|\le C_Aq_n\log n=o(1)$.
The last term in the spectral ratio is uniformly bounded because $m_a\asymp m_0$. It follows that, on $E_{n,A}$,
\[
 c_A\le \rho_\theta\coloneqq f_\theta/f_0\le C_A,\qquad
 |g_{n,a}|\le C_A\log n.
\]
For the real relative derivative multiplier
\[
 h_{j,k,\theta}
   =\frac{D_{\nu,n}^{\,j}\partial_a^k f_\theta}{f_0},
\]
we therefore obtain
\[
 |h_{j,k,\theta}|\mathbf1_{E_{n,A}}
 \le C_{j,k,A}(\log n)^j m_0^{-k}.
\]
The spectral covariance formula gives $B_{j,k,n}(\theta)=\mathcal B_n\bigl(h_{j,k,\theta}\bigr)$. Let $B^{\le}=\mathcal B_n\bigl(h_{j,k,\theta}\mathbf1_{E_{n,A}}\bigr)$, and suppose first that $k\ge1$. Applying the squared Frobenius estimate in \eqref{eq:weightedmatrixbounds} gives
\[
 \norm{B^{\le}}_\Frob^2
 \le\int_{E_{n,A}}|h_{j,k,\theta}|^2w_n
 \le C_{j,k,A}(\log n)^{2j}\int m_0^{-2k}w_n
 \le C_{j,k,A}L_{k,n}^2(\log n)^{2j}.
\]
The last step is \eqref{eq:scalarpowerbound} with $p=2k$. For the operator norm, use the pointwise multiplier bound together with $m_0^{-k}\le a_-^{-k}$ and $\mathcal B_n(1)\preceq I_n$. Order preservation gives
\[
 \norm{B^{\le}}_\Frob\le C_{j,k,A}L_{k,n}(\log n)^j,
 \qquad
 \norm{B^{\le}}_\op\le C_{j,k,A}(\log n)^j.
\]

For the complementary term $B^{>}=\mathcal B_n\bigl(h_{j,k,\theta}\mathbf1_{E_{n,A}^c}\bigr)$, use the trace-norm estimate in \eqref{eq:weightedmatrixbounds} and $w_n(\omega)\le(n/t_0)f_0(\omega)$:
\begin{align*}
 \norm{B^{>}}_1
 &\le\int_{\norm\omega>n^A}|h_{j,k,\theta}(\omega)|
                                      w_n(\omega)\,d\omega\\
 &\le \frac n{t_0}\int_{\norm\omega>n^A}
       \bigl|D_{\nu,n}^{\,j}\partial_a^k f_\theta(\omega)\bigr|
       \,d\omega\\
 &\le Cn\int_{n^A}^\infty
       r^{-2(\nu_-+k)-1}(\log n+\log r)^j\,dr\\
 &\le C_{j,k,A}n^{1-2A(\nu_-+k)}(\log n)^j.
\end{align*}
For the radial integral bound, use $f_\theta\le Cr^{-2\nu_--d}$, $m_a^{-k}\le Cr^{-2k}$, and $|g_{n,a}|\le C(\log n+\log r)$ for $r\ge n^A$. For the last inequality substitute $r=n^As$; the remaining integral is bounded by $C_{j,A}\int_1^\infty s^{-2(\nu_-+k)-1}(1+\log s)^j\,ds$, which is finite by the preceding calculation.

Given any fixed $N>0$, choose $2A\nu_->N+1$. The displayed tail is then $O(n^{-N})$, uniformly in $\theta\in\Theta$ and in any fixed finite collection of $j,k$. Since trace norm bounds both Frobenius and operator norms, adding this tail proves the claimed mixed derivative bounds. Taking the multiplier $m_0^{-k}g_n^j$ directly in the same cutoff argument proves the stated mixed multiplier estimate.

Finally, the sharper bound for pure centered derivatives uses the logarithmic moments, rather than the pointwise bound $|g_{n,a}|\le C\log n$. Here $k=0$ and $h_{j,0,\theta}=\rho_\theta g_{n,a}^j$. Since
\[
 |g_{n,a}-g_n|=|\log(m_0/m_a)|\le Cm_0^{-1}\le C,
\]
the boundedness of $\rho_\theta$ on $E_{n,A}$ and the matrix Cauchy--Schwarz inequality proved above yield
\begin{align*}
 \norm{\mathcal B_n\bigl(\rho_\theta g_{n,a}^j
                   \mathbf1_{E_{n,A}}\bigr)}_\Frob^2
 &\le\tr\{\mathcal B_n\bigl(\rho_\theta^2|g_{n,a}|^{2j}
                   \mathbf1_{E_{n,A}}\bigr)\}\\
 &\le C_j\tr\{\mathcal B_n\bigl(1+g_n^{2j}\bigr)\}
 \le C_jM_n.
\end{align*}
The centered-logarithm argument established this moment bound for every fixed $j$. The same tail estimate, now with $k=0$, therefore gives $\norm{B_{j,0,n}(\theta)}_\Frob\le C_j\sqrt{M_n}=C_jr_n$. The truncated operator norm is at most $C_j(\log n)^j$, and its tail is negligible as before. Every comparison used the fixed compact parameter bounds, the bounded spatial set, and $q_n\log n\to0$, which proves the asserted uniformity. \qed

\subsection{Proof of Lemma~\ref{lem:sampling}}
We want to predict the Fourier coefficients at frequencies where low-order quadrature need not be accurate enough. Linear least squares gives the covariance-optimal linear predictions. A variational formula then turns prediction error into the Sobolev sampling inequality. For iid sites, averaging many disjoint covering subsets proves
that inequality without a logarithmic loss.

\noindent\textbf{1. Express prediction error directly.}
First let $g\in L^2(Q)$ be real, extend it by zero when taking Fourier transforms, and put
\[
 Z_g=\int_Qg(x)F(x)\,dx,\qquad
 c_g=\Cov_0(Y_n,Z_g),\qquad v_g=\Var_0(Z_g).
\]
Here $c_g$ is a column vector. Minimizing the prediction mean squared error over linear predictors gives
\begin{equation*}
 \widehat Z_g^{\,0}=c_g^\T C_{n,0}^{-1}Y_n,\qquad
 \E_0|Z_g-\widehat Z_g^{\,0}|^2
       =v_g-c_g^\T C_{n,0}^{-1}c_g.
\end{equation*}
For the complex Fourier coefficients, use the weights in \eqref{eq:fourierweights}. The representations $U_\lambda=Z_{\check\ell_\lambda}$ and $\widehat U_{\lambda,n}=Z_{\check p_{\lambda,n}}  +b_{\lambda,n}\varepsilon_n$, together with independence of the signal and observation noise, give the more explicit error formula
\begin{equation}\label{eq:fourierresidual}
 \begin{split}
 \E_0|U_\lambda-\widehat U_{\lambda,n}|^2
 &=\int f_0(\omega)
       |\ell_\lambda(\omega)-p_{\lambda,n}(\omega)|^2\,d\omega
        +t_0\sum_i|b_{\lambda,n,i}|^2\\
 &=1-\E_0|\widehat U_{\lambda,n}|^2.
 \end{split}
\end{equation}
The second equality follows because the best linear predictor is uncorrelated with its residual. Thus prediction error measures both the error in the signal's Fourier weight and the noise introduced by the prediction weights.

\noindent\textbf{2. A variational formula reduces the problem to sampling.}
For real functions $u,v\in H^{\beta_0}(\R^d)$, define the inner product
\[
 \langle u,v\rangle_{n,0}
 =\int_{\R^d}\frac{\widehat u(\omega)
                          \overline{\widehat v(\omega)}}{f_0(\omega)}
                          \,d\omega
       +t_0^{-1}\sum_i u(x_{n,i})v(x_{n,i}).
\]
The first term induces a norm equivalent to the usual $H^{\beta_0}(\R^d)$ norm, and point evaluations are bounded because $\beta_0>d/2$. We claim
\begin{equation}\label{eq:predictionvariational}
 \E_0|Z_g-\widehat Z_g^{\,0}|^2
 =\sup_{u\ne0}
       \frac{|\int_Qg(x)u(x)\,dx|^2}{\langle u,u\rangle_{n,0}},
 \qquad u\in H^{\beta_0}(\R^d)\text{ real}.
\end{equation}

\begin{samepage}
To derive the formula, we first represent latent integrals $Z_g$ as vectors in a negative Sobolev space. All Hilbert spaces in the following argument are real. Give $\mathcal H_-\coloneqq H^{-\beta_0}(\R^d)$ the inner product
\[
 \langle p,q\rangle_-
 =\int_{\R^d}f_0(\omega)\widehat p(\omega)
                    \overline{\widehat q(\omega)}\,d\omega,
 \qquad \norm p_-^2=\langle p,p\rangle_-.
\]
\end{samepage}
Since $f_0(\omega)\asymp\langle\omega\rangle^{-2\beta_0}$, this norm is equivalent to the usual $H^{-\beta_0}$ norm. For real smooth compactly supported $p,q$, put $Z_p\coloneqq \int_{\R^d}p(x)F(x)\,dx$. The spectral covariance identity in Lemma~\ref{lem:spectral-covariance} gives
\[
 \E_0[Z_pZ_q]=\langle p,q\rangle_-.
\]
By density, $p\mapsto Z_p$ therefore extends to a linear isometry from $\mathcal H_-$ onto a closed subspace of $L^2(P_{\theta_0})$. This agrees with the preceding integral definition for $g\in L^2(Q)$.

\begin{samepage}
Dirac masses belong to $\mathcal H_-$ because
\[
 \widehat{\delta_x}(\omega)=e^{-2\pi i\omega^\T x},\qquad
 \norm{\delta_x}_-^2=\int_{\R^d}f_0(\omega)\,d\omega=K_0(0).
\]
\end{samepage}
Approximation by smooth mollifiers and mean-square continuity of $F$ give
\[
 Z_{\delta_x}=F(x)\qquad\text{in }L^2(P_{\theta_0}).
\]

Define the bounded operator
\[
 D=\begin{pmatrix}
       \delta_{x_{n,1}}&\cdots&\delta_{x_{n,n}}
    \end{pmatrix}:\R^n\longrightarrow\mathcal H_-,\qquad
 Db=\sum_i b_i\delta_{x_{n,i}}.
\]
Taking adjoints gives
\[
 D^*=\begin{pmatrix}
       \delta_{x_{n,1}}^*\\ \vdots\\ \delta_{x_{n,n}}^*
      \end{pmatrix},\qquad
 (D^*u)_i=\delta_{x_{n,i}}^*u
         =\langle u,\delta_{x_{n,i}}\rangle_-.
\]
Write $\varepsilon_n$ for the observation noise, independence gives
\[
 \begin{aligned}
 Z_g-b^\T Y_n&=Z_{g-Db}-b^\T\varepsilon_n,\\
 \E_0|Z_g-b^\T Y_n|^2
 &=\norm{g-Db}_-^2+t_0\norm b^2.
 \end{aligned}
\]
Orthogonal projection makes $\widehat Z_g^{\,0}$ the best linear predictor, so
\[
 \E_0|Z_g-\widehat Z_g^{\,0}|^2
 =\min_{b\in\R^n}\{\norm{g-Db}_-^2+t_0\norm b^2\}.
\]

In the product Hilbert space $\mathcal E=\mathcal H_-\oplus\R^n$, with $\norm{(u,a)}_{\mathcal E}^2=\norm u_-^2+\norm a^2$, put
\[
 z=(g,0),\qquad
 \mathcal M=\{(Db,\sqrt{t_0}\,b):b\in\R^n\}.
\]
The subspace $\mathcal M$ is finite-dimensional and closed. If $P_{\mathcal M^\perp}$ denotes orthogonal projection onto its orthogonal complement, then
\[
 \min_{b\in\R^n}\{\norm{g-Db}_-^2+t_0\norm b^2\}
 =\operatorname{dist}(z,\mathcal M)^2
 =\norm{P_{\mathcal M^\perp}z}_{\mathcal E}^2.
\]
For $w\in\mathcal M^\perp$, we have $\langle z,w\rangle_{\mathcal E} =\langle P_{\mathcal M^\perp}z,w\rangle_{\mathcal E}$. Cauchy--Schwarz, with equality at
$w=P_{\mathcal M^\perp}z$ when this vector is nonzero, gives
\[
 \norm{P_{\mathcal M^\perp}z}_{\mathcal E}^2
 =\sup_{\substack{w\in\mathcal M^\perp\\w\ne0}}
       \frac{|\langle z,w\rangle_{\mathcal E}|^2}
            {\norm w_{\mathcal E}^2}.
\]
To describe this orthogonal complement, write $w=(u,a)$. Since
\[
 \langle (u,a),(Db,\sqrt{t_0}\,b)\rangle_{\mathcal E}
 =(D^*u+\sqrt{t_0}\,a)^\T b,
\]
orthogonality for every $b$ is equivalent to $a=-t_0^{-1/2}D^*u$. Therefore
\[
 \mathcal M^\perp
 =\{(u,-t_0^{-1/2}D^*u):u\in\mathcal H_-\},
\]
and substitution into the preceding supremum yields
\[
 \E_0|Z_g-\widehat Z_g^{\,0}|^2
 =\sup_{\substack{u\in\mathcal H_-\\u\ne0}}
       \frac{|\langle g,u\rangle_-|^2}
            {\norm u_-^2+t_0^{-1}\norm{D^*u}^2}.
\]

Finally, change variables using the covariance multiplier:
\[
 v=\mathcal R_0u,\qquad
 \widehat v(\omega)=f_0(\omega)\widehat u(\omega).
\]
This is a bijection from $\mathcal H_-$ onto $H^{\beta_0}(\R^d)$, because
\[
 \int_{\R^d}\frac{|\widehat v|^2}{f_0}
 =\int_{\R^d}f_0|\widehat u|^2=\norm u_-^2,
\]
and the inverse is given by $\widehat u=\widehat v/f_0$. Plancherel and Fourier inversion give
\[
 \begin{aligned}
 \langle g,u\rangle_-
 &=\int_{\R^d}\widehat g\,\overline{f_0\widehat u}
   =\int_Qg(x)v(x)\,dx,\\
 (D^*u)_i
 &=\langle u,\delta_{x_{n,i}}\rangle_-
   =\int_{\R^d}f_0(\omega)\widehat u(\omega)
                  e^{2\pi i\omega^\T x_{n,i}}\,d\omega
   =v(x_{n,i}).
 \end{aligned}
\]
The last Fourier integral is absolutely convergent, since $\int|f_0\widehat u|\le\norm u_-(\int f_0)^{1/2}$. Thus the denominator becomes
\[
 \norm u_-^2+t_0^{-1}\norm{D^*u}^2
 =\int_{\R^d}\frac{|\widehat v|^2}{f_0}
        +t_0^{-1}\sum_i|v(x_{n,i})|^2
 =\langle v,v\rangle_{n,0}.
\]
Since $v$ ranges over all nonzero real functions in $H^{\beta_0}(\R^d)$, this proves
\eqref{eq:predictionvariational}.

Put $Q_-=Q$, the working cube already fixed with the cutoff and Fourier lattice. Choose $Q_+$ with $Q_-\Subset Q_+$ inside the cube on which the sampling assumption holds, and fix $s=\nu_-+d/2>d/2$. We now prove the sampling inequality and then apply \eqref{eq:predictionvariational} separately for each type of design.

\noindent\textbf{3. Deterministic sites: sampling and the prediction bound.}
The assumed fill distance on $Q_+$ is at most $Cn^{-1/d}$. Lemma~\ref{lem:cubesampling}, at the fixed order $s$, gives, for every $f\in H^s(\Omega)$,
\[
 \begin{aligned}
 \norm f_{L^2(Q_-)}^2
 &\le C\left\{n^{-2s/d}\norm f_{H^s(\Omega)}^2
               +n^{-1}\sum_i|f(x_{n,i})|^2\right\}\\
 &\le\frac Cn\left\{\norm f_{H^s(\Omega)}^2
                     +\sum_i|f(x_{n,i})|^2\right\},
 \end{aligned}
\]
where the last inequality uses $2s/d>1$. This is \eqref{eq:samplingfunctiontarget} with $Q_-$ in place of $Q$.

To apply it in the variational formula, let $u\in H^{\beta_0}(\R^d)$ be real. Its whole-space norm bounds the restriction norm, and $s\le\beta_0$ gives
\[
 \norm u_{H^s(\Omega)}^2
 \le\norm u_{H^s(\R^d)}^2
 =\int\langle\omega\rangle^{2s}|\widehat u(\omega)|^2\,d\omega
 \le C\int\frac{|\widehat u|^2}{f_0}.
\]
Since $t_0$ is bounded above and below, the sampling inequality therefore implies, for real $g\in L^2(Q_-)$,
\[
 \begin{aligned}
 \left|\int_{Q_-}gu\right|^2
 &\le\norm g_{L^2}^2\norm u_{L^2(Q_-)}^2\\
 &\le\frac Cn\norm g_{L^2}^2
       \left\{\norm u_{H^s(\Omega)}^2
                       +\sum_i|u(x_{n,i})|^2\right\}\\
 &\le\frac Cn\norm g_{L^2}^2\langle u,u\rangle_{n,0}.
 \end{aligned}
\]
Divide by $\langle u,u\rangle_{n,0}$ and take the supremum over nonzero real $u\in H^{\beta_0}(\R^d)$. Formula \eqref{eq:predictionvariational} yields
\[
 \E_0|Z_g-\widehat Z_g^{\,0}|^2
 \le\frac Cn\norm g_{L^2}^2,
\]
which proves \eqref{eq:predictionrisk} for this design.

\noindent\textbf{4. Iid sites: average sampling inequalities, then bound prediction.}
We want the approximation term in Lemma~\ref{lem:cubesampling} to have size $n^{-1}$. This suggests using a mesh with $\ell_n^{2s}\asymp n^{-1}$. Divide each side of $Q_+$ into $\lceil n^{1/(2s)}\rceil$ equal intervals. The resulting cells have side length
\[
 \ell_n\asymp n^{-1/(2s)},\qquad
 \#\{\text{cells}\}\le C\ell_n^{-d}.
\]
Let $p_*>0$ be the density lower bound on $Q_+$. For each cell $E$, its count $N_E=\#\{i:X_i\in E\}$ is binomial with mean $np_E\ge p_*n\ell_n^d$. Put
\[
 k_n=\left\lfloor\frac{p_*}{2}n\ell_n^d\right\rfloor.
\]
Because $s>d/2$, $n\ell_n^d\asymp n^{1-d/(2s)}\to\infty$, so $k_n\asymp n\ell_n^d$ for large $n$. The binomial bound $P\{N_E<np_E/2\}\le e^{-np_E/8}$ and a union bound give,
for the event $\mathcal C_n=\{N_E\ge k_n\text{ for every cell }E\}$, \begin{equation}\label{eq:cellcountprob}
 P_X(\mathcal C_n^c)
 \le C\ell_n^{-d}e^{-cn\ell_n^d}
 \le C'e^{-c'n^{\,1-d/(2s)}}.
\end{equation}
These probabilities are summable.

On $\mathcal C_n$, choose $k_n$ distinct observation indices from each cell. For $k=1,\ldots,k_n$, let $I_{n,k}$ contain the $k$th chosen index from every cell. The sets $I_{n,k}$ are disjoint, and each of them covers $Q_+$ with fill distance at most $\sqrt d\,\ell_n$. Conditional on the locations, Lemma~\ref{lem:cubesampling} applies to every such set:
\[
 \norm f_{L^2(Q_-)}^2
 \le C\left\{\ell_n^{2s}\norm f_{H^s(\Omega)}^2
       +\ell_n^d\sum_{i\in I_{n,k}}|f(X_i)|^2\right\}.
\]
Average these $k_n$ inequalities. Each selected observation appears only once in the double sum, so
\begin{align}
 \norm f_{L^2(Q_-)}^2
 &\le C\ell_n^{2s}\norm f_{H^s(\Omega)}^2
       +\frac{C\ell_n^d}{k_n}
          \sum_{k=1}^{k_n}\sum_{i\in I_{n,k}}|f(X_i)|^2
       \notag\\
 &\le\frac Cn\left\{\norm f_{H^s(\Omega)}^2
                    +\sum_{i=1}^n|f(X_i)|^2\right\}.
       \label{eq:averagedsampling}
\end{align}
The mesh choice gives $\ell_n^{2s}\asymp n^{-1}$, and the averaging gives $\ell_n^d/k_n\asymp n^{-1}$. This proves \eqref{eq:samplingfunctiontarget} on $\mathcal C_n$, simultaneously for every $f\in H^s(\Omega)$. Borel--Cantelli now proves
\eqref{eq:samplingfunctiontarget} eventually almost surely.

On $\mathcal C_n$, \eqref{eq:samplingfunctiontarget} and the same whole-space norm comparison used in the deterministic case give, for every real $u\in H^{\beta_0}(\R^d)$,
\[
 \left|\int_{Q_-}gu\right|^2
 \le\frac Cn\norm g_{L^2}^2\langle u,u\rangle_{n,0}.
\]
Taking the supremum in \eqref{eq:predictionvariational} proves \eqref{eq:predictionrisk} conditionally on every design in $\mathcal C_n$, and hence eventually almost surely.

The fixed choice $s=\nu_-+d/2$ makes the design event independent of the true parameter. Moreover, for every $\theta$ in the compact rectangle and $u\in H^{\nu+d/2}(\R^d)$,
\[
 \norm u_{H^s(\Omega)}^2
 \le\norm u_{H^s(\R^d)}^2
 \le C\int\frac{|\widehat u|^2}{f_\theta}
\]
with a common constant. The same calculation, with each such parameter as truth, proves the prediction bound uniformly over the rectangle in both design cases. This uniformity will be used in Lemma~\ref{lem:randomgood}.

\noindent\textbf{5. Normalize the Fourier prediction error.}
Recall that $Q_-=Q$ throughout. Steps 3 and 4 proved \eqref{eq:predictionrisk} for real $g$ in the two design cases. Applying this bound to the real and imaginary parts gives it for complex $g$, since their squared $L^2$ norms add to $\norm g_{L^2}^2$.
For $g=L^{-d/2}\chi e^{-2\pi i\lambda^\T x}$ this norm is a fixed constant. Dividing its prediction bound by $v_0(\lambda)\asymp(1+\norm\lambda^2)^{-\beta_0}$ proves
\eqref{eq:fourierpredictionrisk}. In particular, at $\norm\lambda\le\epsilon R_n$, the relative error is at most $C\epsilon^{\kappa_0}+o(1)$. Equation~\eqref{eq:fourierresidual} then shows that the prediction explains at least half the latent variance when $\epsilon$ is sufficiently small and $n$ sufficiently large.

\noindent\textbf{6. Bound fluctuations of the predicted Fourier energies.}
For a finite frequency set, let $V$ be the real vector formed by the real and imaginary parts of $U_\lambda$, and put $\widehat V=\Cov_0(V,Y_n)C_{n,0}^{-1}Y_n$.
Lemma~\ref{lem:windows} gives $\Cov_0(V)\preceq CI$. Linear prediction and its residual are uncorrelated, so
\[
 \Cov_0(\widehat V)\preceq\Cov_0(V)\preceq CI.
\]
Let $D$ be diagonal with each $w_\lambda$ repeated for the real and imaginary coordinates. Then $\sum_\lambda w_\lambda|\widehat U_{\lambda,n}|^2
=\widehat V^\T D\widehat V$ and $\norm D_\Frob^2=2\sum_\lambda w_\lambda^2$.
Write $A=\Cov_0(\widehat V)$, the quadratic-form bound in Lemma~\ref{lem:gaussian}, applied to the linear map from $Y_n$ to $\widehat V$, gives
\[
 \Var_0(\widehat V^\T D\widehat V)
 \le C\norm{A^{1/2}DA^{1/2}}_\Frob^2
 \le C\norm A_\op^2\norm D_\Frob^2
 \le C\sum_\lambda w_\lambda^2.
\]
This proves \eqref{eq:predictedenergyvariance}, including signed weights.

The same covariance bound also controls the Fourier responses of all predictions at once. Let $A_\Lambda$ be the matrix whose rows are $b_{\lambda,n}$, $\lambda\in\Lambda$, so that the complex vector of predictions is $A_\Lambda Y_n$.
The complex covariance bound in Lemma~\ref{lem:windows} and orthogonality of the prediction residual give
\[
 A_\Lambda C_{n,0}A_\Lambda^*
 =\Cov_0((\widehat U_{\lambda,n})_{\lambda\in\Lambda})
 \preceq\Cov_0((U_\lambda)_{\lambda\in\Lambda})\preceq CI.
\]
Thus $\norm{A_\Lambda C_{n,0}^{1/2}}_\op^2\le C$ and, for every vector $a\in\mathbb C^n$,
\[
 \norm{A_\Lambda a}^2
 =\norm{A_\Lambda C_{n,0}^{1/2}C_{n,0}^{-1/2}a}^2
 \le Ca^*C_{n,0}^{-1}a.
\]
Take $a=a_{-\omega}$, so that $(A_\Lambda a_{-\omega})_\lambda=p_{\lambda,n}(\omega)$.
Since $C_{n,0}$ is real, $a_{-\omega}=\overline{a_\omega}$ gives $a_{-\omega}^*C_{n,0}^{-1}a_{-\omega}  =a_\omega^*C_{n,0}^{-1}a_\omega$. Multiplying by $f_0(\omega)$ proves \eqref{eq:predictionenvelope}.

\qed

\subsection{Proof of Lemma~\ref{lem:information}}
Throughout this proof the probability law is the auxiliary Gaussian experiment $G_{\theta,n}$. Its Fisher information is exactly the deterministic Gaussian curvature $\calI_n$ needed in the Gaussian criterion analysis. This proof makes no Fisher-information assertion about the actual non-Gaussian response law.

The guiding approximation is that the localized Fourier coefficients behave like independent Gaussian variables with variances $c_\chi f_0(\lambda)$. A band containing $N$ such coefficients estimates its relative power with error of order $N^{-1/2}$.
We therefore average normalized squared coefficients in two or three bands. Their responses distinguish the constant, logarithmic, and reciprocal-square shapes of the spectral derivatives. The covariance bounds in Lemmas~\ref{lem:windows} and \ref{lem:sampling} give the variance control that this heuristic suggests, and the prediction bound lets us replace latent coefficients by observable statistics. At $d=4$, range information is spread over logarithmically many frequency scales, so we add one energy that uses all those scales.

Condition on a deterministic design satisfying the assumptions, or on an iid design for which Lemma~\ref{lem:sampling} holds eventually. Use its fixed interior cube and cutoff throughout. All prediction weights, frequency sets, and variance normalizations
are fixed at the truth when differentiating.

\noindent\textbf{1. Information is bounded below by the response of a statistic.}
For the normalized local coordinates, write
\[
 B_u=\frac{P_n}{r_n},\qquad
 B_v=\frac{\mathcal B_n(g_n)}{r_n},\qquad
 B_a=-\frac{\beta_0}{b_n}\mathcal B_n(m_0^{-1}),\qquad
 B_t=\frac{I-P_n}{t_0\sqrt n}.
\]
The range coordinate is included only when $d\ge4$. For a normalized parameter direction $k$, put $B(k)=\sum_i k_iB_i$ and $\xi_n=C_{n,0}^{-1/2}Y_n\sim N(0,I_n)$.
The directional score at the truth is
\[
 \dot\ell_k=\frac12\{\xi_n^\T B(k)\xi_n-\tr B(k)\},
 \qquad
 \E_0\dot\ell_k^2=k^\T J_nk
      =\frac12\norm{B(k)}_\Frob^2.
\]
For a fixed square-integrable statistic $Z=Z(Y_n)$, differentiation of its mean gives
$\dot\mu_Z(k)=\E_0[(Z-\E_0Z)\dot\ell_k]$. Cauchy--Schwarz gives
\[
 |\dot\mu_Z(k)|^2
 \le\Var_0(Z)\E_0\dot\ell_k^2
 =\Var_0(Z)\,k^\T J_nk.
\]
Consequently, for $\Var_0(Z)>0$,
\begin{equation}\label{eq:statisticinformation}
 k^\T J_nk\ge\frac{|\dot\mu_Z(k)|^2}{\Var_0(Z)}.
\end{equation}
Thus a statistic with bounded variance and a nonvanishing mean response supplies a lower information bound. For a real vector $Z$ with $\Cov_0(Z)\preceq CI$ and mean derivative $Dk$, apply this inequality to $b^\T Z$, where $b=Dk/\norm{Dk}$, to obtain
\begin{equation}\label{eq:vectorstatisticinformation}
 k^\T J_nk\ge C^{-1}\norm{Dk}^2.
\end{equation}
If $Dk=0$ this bound is immediate.

The upper bound is already available: Lemma~\ref{lem:upper} gives $\norm{B_i}_\Frob\le C$ for each signal coordinate, and $0\preceq P_n\preceq I$ gives $\norm{B_t}_\Frob\le t_0^{-1}$.
Consequently
\[
 k^\T J_nk\le\frac12\norm k^2\sum_i\norm{B_i}_\Frob^2
 \le C\norm k^2.
\]

\noindent\textbf{2. Prediction preserves the response of a weighted Fourier energy.}
For now $h$ is a signal direction and the nugget is fixed. Choose $0<\epsilon\le1$, to be made small below, and let
\[
 \Lambda_n=\{\lambda\in\mathscr L:\lambda_1>0,
                                      \ \norm\lambda\le\epsilon R_n\}.
\]
The half-space avoids conjugate duplication. For real weights $w=(w_\lambda)_{\lambda\in\Lambda_n}$ define the latent and observable energies
\begin{equation*}
 \begin{split}
 Z_n^{\rm lat}(w)&=\sum_\lambda w_\lambda(|U_\lambda|^2-1),\\
 Z_n(w)&=\sum_\lambda w_\lambda
       \{|\widehat U_{\lambda,n}|^2-\E_0|\widehat U_{\lambda,n}|^2\}.
 \end{split}
\end{equation*}
Their variances are at most $C\norm w_{\ell^2}^2$, by Lemmas~\ref{lem:windows} and \ref{lem:sampling}. The relative variance derivative of one latent coefficient is
\[
 q_\lambda(h)
 =\left.\frac{d}{ds}\log v_{\theta_0+sE_nh}(\lambda)\right|_{s=0}
 =\int g_h f_0|\ell_\lambda|^2.
\]
Since $\partial_\eta\log v_\theta=1$, the logarithmic-derivative bounds \eqref{eq:fourierlocalderivatives} give
\begin{equation}\label{eq:latenttangentapproximation}
 |q_\lambda(h)-g_h(\lambda)|
 \le C\left\{
       \frac{|h_v|}{r_n\langle\lambda\rangle^2}
       +\frac{|h_a|}{b_n\langle\lambda\rangle^4}\right\}.
\end{equation}
Here and below the range term is omitted when it is absent from $h$. This uses the normalized directional derivative derived in Section~\ref{sec:proofnotation}, applied to $\log f_\theta$ and $\log v_\theta$, respectively, at the truth.

The mean derivative of $Z_n^{\rm lat}(w)$ is $\sum_\lambda w_\lambda q_\lambda(h)$. For its predicted version, the spectral covariance identity gives
\[
 \E_\theta|\widehat U_{\lambda,n}|^2
 =\int f_\theta|p_{\lambda,n}|^2
                  +t_0\sum_i|b_{\lambda,n,i}|^2.
\]
The nugget term has zero derivative in a signal direction. Hence the difference between the two mean derivatives is
\[
 e_n(w,h)=\sum_\lambda w_\lambda
       \int g_h f_0\{|p_{\lambda,n}|^2-|\ell_\lambda|^2\}.
\]
We claim the useful bound
\begin{equation}\label{eq:traceerror}
 |e_n(w,h)|\le C\epsilon^{\kappa_0/2}
                    \norm w_{\ell^2}\norm h,
 \qquad \kappa_0=2\nu_0+d.
\end{equation}

To prove it, set
\[
 \vartheta_n(\omega)=\min\{1,(n/t_0)f_0(\omega)\},\qquad
 W_n(h)=\int g_h(\omega)^2\vartheta_n(\omega)\,d\omega.
\]
The polar-coordinate estimates from Lemma~\ref{lem:upper} give
\[
 \int\vartheta_n\le CR_n^d,\qquad
 \int g_n^2\vartheta_n\le CR_n^d,\qquad
 \int m_0^{-2}\vartheta_n\le Cb_n^2\quad(d\ge4).
\]
These are the spectral integrals estimated in \eqref{eq:effective}, \eqref{eq:logmoments}, and \eqref{eq:rangetrace}, respectively, with $w_n$ replaced by its upper envelope. Dividing by $r_n^2$ and $b_n^2$ yields
\begin{equation*}
 W_n(h)\le C\norm h^2.
\end{equation*}

The prediction envelope \eqref{eq:predictionenvelope} and \eqref{eq:frequencyweight} give $\sum_\lambda f_0|p_{\lambda,n}|^2\le C\vartheta_n$. The latent-weight envelope in
Lemma~\ref{lem:windows}\ref{it:windows-envelope}, with $B=\epsilon R_n$, gives the same bound for $\ell_\lambda$. Together,
\begin{equation*}
 \sum_{\lambda\in\Lambda_n}
       f_0\{|p_{\lambda,n}|^2+|\ell_\lambda|^2\}
 \le C\vartheta_n.
\end{equation*}
Lemma~\ref{lem:sampling} and \eqref{eq:fourierresidual} also give
\[
 \int f_0|p_{\lambda,n}-\ell_\lambda|^2
 \le Cn^{-1}\langle\lambda\rangle^{\kappa_0}
 \le C\epsilon^{\kappa_0},\qquad \lambda\in\Lambda_n,
\]
for all sufficiently large $n$, with $\epsilon$ fixed. Finally use $|p|^2-|\ell|^2=(p-\ell)\overline p+ \ell(\overline p-\overline\ell)$ and Cauchy--Schwarz over the
frequencies and the coefficients:
\[
\begin{aligned}
 |e_n(w,h)|
 &\le\left\{\sum_\lambda w_\lambda^2
              \int f_0|p_{\lambda,n}-\ell_\lambda|^2\right\}^{1/2}\left[
       \left\{\sum_\lambda\int g_h^2f_0|p_{\lambda,n}|^2\right\}^{1/2}
      +\left\{\sum_\lambda\int g_h^2f_0|\ell_\lambda|^2\right\}^{1/2}
                  \right]\\
 &\le C\epsilon^{\kappa_0/2}\norm w_{\ell^2}\sqrt{W_n(h)}
 \le C\epsilon^{\kappa_0/2}\norm w_{\ell^2}\norm h.
\end{aligned}
\]
This proves \eqref{eq:traceerror} for arbitrary signed weights and signal directions.

\noindent\textbf{3. Two or three bands distinguish the spectral derivative shapes.}
First take the amplitude and smoothness coordinates; include range also when $d>4$. Thus let $k=2$ for $d\le4$ and $k=3$ for $d>4$, with $h_a=0$ for the present argument when $d=4$. Choose a small relative bandwidth $\delta\in(0,1/4)$, and put
\[
 \rho=\epsilon/8,\qquad t_j=2^j\rho,\qquad
 \Lambda_{n,j}=\left\{\lambda\in\mathscr L:\lambda_1>0,
          \left|\frac{\norm\lambda}{t_jR_n}-1\right|\le\delta\right\},
 \quad 0\le j<k.
\]
The bands are disjoint and contained in $\Lambda_n$. Lattice cells have fixed side $m/L$; cells meeting a band boundary have total volume $O((t_jR_n)^{d-1})$, whereas the band has volume comparable to $(t_jR_n)^d$. Thus their counts satisfy, for fixed $\delta$,
\[
 N_{n,j}\coloneqq \#\Lambda_{n,j}\asymp(\epsilon R_n)^d,\qquad
 \frac{\sqrt{N_{n,j}}}{r_n}\asymp\epsilon^{d/2}.
\]
Use the $k$ observable band energies
\[
 Z_{n,j}=N_{n,j}^{-1/2}
       \sum_{\lambda\in\Lambda_{n,j}}
       \{|\widehat U_{\lambda,n}|^2-\E_0|\widehat U_{\lambda,n}|^2\}.
\]
For any $c\in\R^k$, $c^\T Z_n$ is the energy with weight $w_\lambda=c_j/\sqrt{N_{n,j}}$ on band $j$. Its squared weight norm is $\sum_jc_j^2$, so $\Cov_0(Z_n)\preceq CI_k$.

We next compute the response of the latent band energies. Write $D_n^{\rm lat}h$ for their mean derivative vector and $D_nh$ for the observable one. Its latent entries are
\[
 (D_n^{\rm lat}h)_j
   =\frac1{\sqrt{N_{n,j}}}\sum_{\lambda\in\Lambda_{n,j}}q_\lambda(h).
\]
The response shapes can be seen directly by writing $\omega=R_nz$:
\[
 r_ng_h(R_nz)
 =h_u-h_v\log\{a_0/R_n^2+4\pi^2\norm z^2\}
   -\mathbf1_{\{d>4\}}
       \frac{\beta_0h_a}{a_0/R_n^2+4\pi^2\norm z^2}.
\]
Here $r_n/(b_nR_n^2)=1$ for $d>4$. To see how to choose the coefficients of the response matrix, take $\norm\lambda=t_jR_nq$, where $t_j=2^j\rho$ and $|q-1|\le\delta$. Put
\[
 d_{n,j}(q)=\frac{a_0}{4\pi^2t_j^2R_n^2q^2},\qquad
 m_0(\lambda)=4\pi^2t_j^2R_n^2q^2\{1+d_{n,j}(q)\}.
\]
Substitution in the logarithmic and range terms gives
\[
 \begin{aligned}
 2\log R_n-\log m_0(\lambda)
 &=-\log(4\pi^2\rho^2)-2j\log2-2\log q
                                      -\log\{1+d_{n,j}(q)\},\\
 \frac{r_n}{b_nm_0(\lambda)}
 &=\frac{4^{-j}q^{-2}}{4\pi^2\rho^2\{1+d_{n,j}(q)\}}
                                      \qquad(d>4).
 \end{aligned}
\]
Consequently,
\[
 \begin{aligned}
 r_ng_h(\lambda)
 &=\{h_u-h_v\log(4\pi^2\rho^2)\}
       +j\{-2(\log2)h_v\}-2h_v\log q\\
 &\quad+4^{-j}\left\{-\frac{\beta_0h_a}{4\pi^2\rho^2}\right\}
                    q^{-2}\{1+d_{n,j}(q)\}^{-1}
       -h_v\log\{1+d_{n,j}(q)\},
 \end{aligned}
\]
where the range term is omitted when $k=2$. The coefficients of the three shapes $1,j,4^{-j}$ are therefore naturally
\[
 \alpha_1=h_u-h_v\log(4\pi^2\rho^2),\qquad
 \alpha_2=-2(\log2)h_v,\qquad
 \alpha_3=-\frac{\beta_0h_a}{4\pi^2\rho^2}\quad(d>4).
\]
Since $q\in[3/4,5/4]$, we have $0\le d_{n,j}(q)\le C(\rho R_n)^{-2}$.
The inequalities $\log(1+x)\le x$ and $|(1+x)^{-1}-1|\le x$ for $x\ge0$ then imply
\[
 \begin{aligned}
 r_ng_h(\lambda)
 &=\alpha_1+j\alpha_2-2h_v\log q
                  +4^{-j}\alpha_3q^{-2}+r_{n,\lambda}(h),\\
 |r_{n,\lambda}(h)|
 &\le\frac{C}{(\rho R_n)^2}
               (|h_v|+|\alpha_3|)
 \le\frac{C}{(\rho R_n)^2}\norm\alpha.
 \end{aligned}
\]
Again omit $\alpha_3$ when $k=2$. At the band centers $q=1$, the limiting responses are $\alpha_1+j\alpha_2+4^{-j}\alpha_3$. Within a band, $|\log q|+|q^{-2}-1|\le C\delta$, so
\[
 \left|r_ng_h(\lambda)
           -(\alpha_1+j\alpha_2+4^{-j}\alpha_3)\right|
 \le C\{\delta+(\rho R_n)^{-2}\}\norm\alpha.
\]

We must also absorb the smoothing error in \eqref{eq:latenttangentapproximation}. On these bands, $\langle\lambda\rangle\ge c\rho R_n$, and multiplication of that inequality by $r_n$ gives
\[
 \begin{aligned}
 r_n|q_\lambda(h)-g_h(\lambda)|
 &\le C\left\{\frac{|h_v|}{\langle\lambda\rangle^2}
       +\frac{r_n|h_a|}{b_n\langle\lambda\rangle^4}\right\}\\
 &\le C\left\{\frac{|h_v|}{\rho^2R_n^2}
                         +\frac{|h_a|}{\rho^4R_n^2}\right\}\\
 &\le\frac{C}{(\rho R_n)^2}
                     (|h_v|+\rho^{-2}|h_a|)
 \le\frac{C}{(\rho R_n)^2}\norm\alpha.
 \end{aligned}
\]
For the range term we used $r_n/b_n=R_n^2$ when $d>4$ and $h_a=-(4\pi^2\rho^2/\beta_0)\alpha_3$; otherwise that term is absent. The two errors therefore combine to give
\[
 r_nq_\lambda(h)=\alpha_1+j\alpha_2+4^{-j}\alpha_3
                     +\varepsilon_{n,\lambda}(h),\qquad
 |\varepsilon_{n,\lambda}(h)|
 \le C\{\delta+(\rho R_n)^{-2}\}\norm\alpha.
\]

The two- and three-band evaluation matrices are
\[
 V_2=\begin{pmatrix}1&0\\1&1\end{pmatrix},\qquad
 V_3=\begin{pmatrix}1&0&1\\1&1&1/4\\1&2&1/16\end{pmatrix},
 \qquad \det V_2=1,\quad\det V_3=9/16.
\]
Let $s_* =\min\{\sigma_{\min}(V_2),\sigma_{\min}(V_3)\}>0$. After averaging within each band, the preceding estimate becomes
\[
 \operatorname{diag}\left(\frac{r_n}{\sqrt{N_{n,j}}}\right)
               D_n^{\rm lat}h
 =V_k\alpha+e_n^{\rm band},\qquad
 \norm{e_n^{\rm band}}\le C\{\delta+(\rho R_n)^{-2}\}\norm\alpha.
\]
Choose $\delta$ so that $C\delta\le s_*/4$. For each fixed $\epsilon>0$, take $n$ sufficiently large that $C(\rho R_n)^{-2}\le s_*/4$. This absorbs both the variation
within a band and the smoothing error from \eqref{eq:latenttangentapproximation}, since
\[
 \norm{V_k\alpha+e_n^{\rm band}}
 \ge\norm{V_k\alpha}-\norm{e_n^{\rm band}}
 \ge s_*\norm\alpha-\frac{s_*}{2}\norm\alpha
 =\frac{s_*}{2}\norm\alpha.
\]
Using $\sqrt{N_{n,j}}/r_n\asymp\epsilon^{d/2}$ gives
\[
 \norm{D_n^{\rm lat}h}\ge c\epsilon^{d/2}\norm\alpha
       \ge\frac{c\epsilon^{d/2}}{1+|\log\epsilon|}\norm h.
\]
For the last inequality, solve for $h_v$ from $\alpha_2$, then for $h_u$ from $\alpha_1$, and, when present, for $h_a$ from $\alpha_3$. These formulas give
$\norm h\le C(1+|\log\epsilon|)\norm\alpha$.

To transfer this bound to the predicted energies, write $w^{(j)}_\lambda=N_{n,j}^{-1/2}\mathbf1_{\Lambda_{n,j}}(\lambda)$. For the combined statistic $c^\T Z_n$, the disjoint bands give
\[
 w_\lambda(c)=\sum_jc_jw^{(j)}_\lambda,\qquad
 \norm{w(c)}_{\ell^2}^2
 =\sum_j\sum_{\lambda\in\Lambda_{n,j}}\frac{c_j^2}{N_{n,j}}
 =\sum_jc_j^2=\norm c^2.
\]
By the definition of the response error,
\[
 e_n(w(c),h)=c^\T(D_nh-D_n^{\rm lat}h).
\]
If the vector on the right is nonzero, choose $c=(D_nh-D_n^{\rm lat}h)/\norm{D_nh-D_n^{\rm lat}h}$. Then $\norm{w(c)}_{\ell^2}=1$ and $e_n(w(c),h)=\norm{D_nh-D_n^{\rm lat}h}$. Applying \eqref{eq:traceerror} gives
\[
 \norm{D_nh-D_n^{\rm lat}h}
 =\sup_{\norm c=1}|e_n(w(c),h)|
 \le C\epsilon^{\kappa_0/2}\norm h.
\]
The same bound is immediate if the difference vector is zero. Relative to the preceding lower bound, the error is at most $C\epsilon^{\nu_0}(1+|\log\epsilon|)$. Choose $\epsilon$ so small that this is less than $1/2$, uniformly since $\nu_0\ge\nu_->0$. Equations \eqref{eq:vectorstatisticinformation} and $\Cov_0(Z_n)\preceq CI$ then give
\[
 h^\T J_{n,\mathrm{sig}}h
 \ge C^{-1}\norm{D_nh}^2\ge c\norm h^2.
\]
This proves the full signal lower bound for $d\ne4$ and the amplitude--smoothness principal-block bound for $d=4$.

\noindent\textbf{4. At $d=4$, collect the range information across scales.}
In dimension four the squared range shape has radial integral $\int_1^{R_n}r^3r^{-4}\,dr=\log R_n$. Every doubling of the frequency radius contributes the same order, so a single band near $R_n$ would lose the factor $\log R_n$. Use instead
\[
 \Lambda_n^a=\{\lambda\in\Lambda_n:\norm\lambda\ge1\},\qquad
 w_\lambda^a=\frac{\mathbf1_{\Lambda_n^a}(\lambda)}{b_nm_0(\lambda)},
 \qquad Z_{n,a}=Z_n(w^a),\qquad b_n=\sqrt{\log n}.
\]
For every sufficiently large dyadic radius $2^j$, the half-shell $2^j\le\norm\lambda<2^{j+1}$ has order $2^{4j}$ lattice points, and $m_0(\lambda)^{-2}\asymp2^{-4j}$ on it. Each complete shell therefore contributes an amount bounded above and below by positive constants. There are order $\log(\epsilon R_n)$ such shells, and the finitely many low-frequency points and the last partial shell contribute at most a constant. Consequently
\[
 \sum_{\lambda\in\Lambda_n^a}m_0(\lambda)^{-2}
       \asymp\log R_n,\qquad
 \norm{w^a}_{\ell^2}^2
       \asymp\frac{\log R_n}{\log n}
       =\frac1{\kappa_0}\asymp1.
\]
Here $\epsilon$ is fixed before $n\to\infty$; the comparison constants can be chosen independently of small $\epsilon$ once $n$ is sufficiently large. Thus $\Var_0(Z_{n,a})\le C$.

In the unit normalized range direction $e_a$, the latent mean derivative is, by \eqref{eq:latenttangentapproximation},
\[
\begin{aligned}
 \sum_\lambda w_\lambda^a q_\lambda(e_a)
 &=-\frac{\beta_0}{b_n^2}
               \sum_{\lambda\in\Lambda_n^a}m_0(\lambda)^{-2}
       +O\left\{\frac1{b_n^2}
               \sum_{\lambda\in\Lambda_n^a}
                  m_0(\lambda)^{-1}\langle\lambda\rangle^{-4}\right\}\\
 &=-\frac{\beta_0}{b_n^2}
               \sum_{\lambda\in\Lambda_n^a}m_0(\lambda)^{-2}
       +O((\log n)^{-1}).
\end{aligned}
\]
The last sum in braces is bounded, since its summand is at most $C\langle\lambda\rangle^{-6}$ and $d=4$. The magnitude of this response is bounded below by a positive constant. The prediction error in its response is at most
$C\epsilon^{\kappa_0/2}\norm{w^a}_{\ell^2}$ by \eqref{eq:traceerror}. Decrease the fixed $\epsilon$, if necessary, to keep this error below half the lower bound. The information inequality \eqref{eq:statisticinformation} then gives $J_{n,aa}\ge c$.

\smallskip
\noindent\emph{Combining the range and the two-band bounds.}
It remains to control the range cross-information; separate lower bounds for the two principal blocks would not by themselves give a lower bound for their joint matrix. We show directly that $\tr(B_uB_a)$ and $\tr(B_vB_a)$ tend to zero. Put $Q_n=\mathcal B_n(m_0^{-1})$, so $B_a=-\beta_0Q_n/\sqrt{\log n}$, and split the range multiplier at
$\rho_n=R_n/\log n$, which tends to infinity:
\[
 Q_n=Q_n^{\le}+Q_n^{>},\qquad
 Q_n^{\le}=\mathcal B_n\bigl(m_0^{-1}
       \mathbf1_{\{\norm\omega\le\rho_n\}}\bigr).
\]
Both parts are positive. Below $\rho_n$ we use the trace norm and the small operator norms of $B_u$ and $B_v$. Above $\rho_n$ only $O(\log\log n)$ of the squared range norm
remains, compared with its full order $\log n$.

To verify these statements, use the frequency weight $w_n$ from \eqref{eq:frequencyweightdefinition}. Equations \eqref{eq:frequencyweight} and
\eqref{eq:weightedmatrixbounds}, followed by polar coordinates in dimension four, give
\begin{align*}
 \norm{Q_n^{\le}}_1
 &=\int_{\norm\omega\le\rho_n}m_0^{-1}w_n(\omega)\,d\omega
 \le C\left(1+\int_1^{\rho_n}r\,dr\right)
 \le C\rho_n^2,\\
 \norm{Q_n^{>}}_\Frob^2
 &\le\int_{\norm\omega>\rho_n}m_0^{-2}w_n(\omega)\,d\omega\\*
 &\le C\left\{\int_{\rho_n}^{R_n}\frac{dr}{r}
       +n\int_{R_n}^\infty r^{-2\nu_0-5}\,dr\right\}\\*
 &=C\left\{\log(R_n/\rho_n)
       +\frac{nR_n^{-2\nu_0-4}}{2\nu_0+4}\right\}
 \le C(1+\log\log n).
\end{align*}
The two terms in the last line use $R_n/\rho_n=\log n$ and $n=R_n^{2\nu_0+4}$, respectively.

For $B=B_u$ or $B_v$, Lemma~\ref{lem:upper} gives $\norm B_\Frob\le C$ and $\norm B_\op\le C\log n/r_n$. Since $r_n=R_n^2$ at $d=4$, trace-norm duality on the first
part and Frobenius Cauchy--Schwarz on the second give
\begin{align*}
 \left|\tr\left(B\frac{Q_n}{\sqrt{\log n}}\right)\right|
 &\le\frac{\norm B_\op\norm{Q_n^{\le}}_1
            +\norm B_\Frob\norm{Q_n^{>}}_\Frob}
                {\sqrt{\log n}}\\
 &\le C\frac{(\log n)\rho_n^2}{R_n^2\sqrt{\log n}}
       +C\sqrt{\frac{1+\log\log n}{\log n}}\\
 &=C(\log n)^{-3/2}
       +C\sqrt{\frac{1+\log\log n}{\log n}}
 \longrightarrow0.
\end{align*}
Multiplication by $-\beta_0/2$ gives the corresponding cross-information with $B_a$.

The two-band bound gives $J_{n,(u,v)}\succeq cI_2$, and the weighted energy gives $J_{n,aa}\ge c$. Since the cross block tends to zero, their joint matrix is bounded below by $c'I_3$ for all sufficiently large $n$. This proves the remaining signal lower
bound and the stated additional block decomposition at $d=4$.

Together with Step 3, this proves in every dimension that
\[
 \norm{\mathcal B_n(g_h)}_\Frob^2
 =2h^\T J_{n,\mathrm{sig}}h\ge c\norm h^2,
\]
which is \eqref{eq:signalcoercivity}.

\noindent\textbf{5. Handle the nugget and finish the information bound.}
Since $\tr P_n=O(M_n)=o(n)$, the term $I-P_n$ retains nearly all $n$ observation directions. The formula $B_t=(I-P_n)/(t_0\sqrt n)$ makes this precise:
\[
 J_{tt}=\frac1{2t_0^2n}
       \{n-2\tr P_n+\tr P_n^2\}
       =\frac1{2t_0^2}+O(M_n/n)
       \longrightarrow\frac1{2t_0^2}.
\]
We used $0\preceq P_n\preceq I$, so $\tr P_n^2\le\tr P_n\le CM_n$, and $M_n/n=R_n^{-2\nu_0}\to0$.

For $B_j=B_u,B_v$, or $B_a$ when present, trace-norm duality gives
\[
 |J_{jt}|\le\frac12\norm{B_j}_1\norm{B_t}_\op
       \le Cn^{-1/2}\norm{B_j}_1.
\]
For amplitude, $\norm{B_u}_1=\tr(P_n)/r_n\le Cr_n$. For centered smoothness, the positive moment bounds and Cauchy--Schwarz with measure $w_n(\omega)\,d\omega$ give
\[
 \norm{B_v}_1\le\frac1{r_n}\int|g_n|w_n
 \le\frac1{r_n}
       \left(\int g_n^2w_n\right)^{1/2}
       \left(\int w_n\right)^{1/2}
 \le\frac{CM_n}{r_n}=Cr_n.
\]
Consequently,
\[
 |J_{ut}|+|J_{vt}|\le C\frac{r_n}{\sqrt n}
       =CR_n^{-\nu_0}\longrightarrow0.
\]
For $d\ge4$, the positive trace bound for $m_0^{-1}$ similarly gives
\[
 |J_{at}|\le C\frac{R_n^{d-2}}{b_n\sqrt n}
 =\begin{cases}
 O\!\left(R_n^{-\nu_0}/\sqrt{\log n}\right),&d=4,\\
 O(R_n^{-\nu_0}),&d>4,
 \end{cases}
 \qquad\longrightarrow0.
\]

Let $J_{n,\mathrm{sig}}$ be the information block for amplitude and smoothness, including range when $d\ge4$, and let $j_t$ collect its nugget cross-entries.
We have proved $J_{n,\mathrm{sig}}\succeq c_sI$, $J_{tt}\to(2t_0^2)^{-1}$, and $\norm{j_t}\to0$. Thus
\[
 J_n=
 \begin{pmatrix}J_{n,\mathrm{sig}}&0\\0&J_{tt}\end{pmatrix}
 +\begin{pmatrix}0&j_t\\j_t^\T&0\end{pmatrix}
 \succeq\{\min(c_s,J_{tt})-\norm{j_t}\}I
 \succeq cI
\]
for all sufficiently large $n$; the off-diagonal matrix has operator norm $\norm{j_t}$. This proves the bound for $J_n$ in \eqref{eq:Jbounds}. The cross-term limits just established give the stated orthogonality assertions and the two block decompositions in the lemma. \qed

\subsection{Proof of Lemma~\ref{lem:derivatives}}
All segment estimates below are taken on the enclosing set $\mathcal N_n^\square$ and then restricted to the feasible parameters in $\mathcal N_n$.

\noindent\textbf{Covariance comparison on the preliminary set.}
Fix a sufficiently large $A$. The frequency-cutoff estimates in the proof of Lemma~\ref{lem:upper} give $f_\theta/f_0\asymp1$ on $\norm\omega\le n^A$, uniformly on
$\mathcal N_n$, and a sampled spectral tail of operator norm $O(n^{1-2A\nu_-})=o(1)$. To express the comparison after whitening, put $\chi_n=\mathbf1_{\{\norm\omega\le n^A\}}$ and
\[
 E_{\theta,n}
   =\mathcal B_n\bigl((1-\chi_n)\frac{f_\theta}{f_0}\bigr),
 \qquad
 E_{0,n}=\mathcal B_n\bigl(1-\chi_n\bigr).
\]
These matrices are positive semidefinite, and the tail bound together with $C_{n,0}\succeq t_0I_n$ gives $\sup_{\theta\in\mathcal N_n}\norm{E_{\theta,n}}_\op=o(1)$. Since $f_\theta/f_0\asymp1$ on the support of $\chi_n$ and $t\asymp t_0$, uniformly on $\mathcal N_n$,
\begin{align*}
 M(\theta)
 &=\mathcal B_n\bigl(\frac{f_\theta}{f_0}\bigr)+tC_{n,0}^{-1}\\
 &=\mathcal B_n\bigl(\chi_n\frac{f_\theta}{f_0}\bigr)
       +tC_{n,0}^{-1}+E_{\theta,n}\\
 &\asymp \mathcal B_n\bigl(\chi_n\bigr)+t_0C_{n,0}^{-1}+E_{\theta,n}\\
 &=I_n-E_{0,n}+E_{\theta,n}
 =I_n+o_{\op}(1)\asymp I_n.
\end{align*}
Here $\asymp$ denotes two-sided Loewner bounds, and the last line uses $\mathcal B_n\bigl(1\bigr)+t_0C_{n,0}^{-1}=I_n$. Thus $cI_n\preceq M(\theta)\preceq CI_n$
uniformly on $\mathcal N_n$.

\noindent\textbf{Derivatives in the centered coordinates.}
The exact spectral ratio in $(u,v,a,t)$ is
\begin{equation}\label{eq:localratioexact}
 \frac{f_\theta}{f_0}
 =\exp\{u+vg_n(\omega)\}
   \left(\frac{m_0(\omega)}{a+4\pi^2\norm\omega^2}\right)^{\beta_0+v}.
\end{equation}
Use the derivative multipliers $g_{n,a}$ and $m_a^{-1}$ from the mixed-derivative calculation in Lemma~\ref{lem:upper}.

Lemma~\ref{lem:upper} now gives the pure $u,v$ Frobenius bounds $Cr_n$ through order three, and the stated first-derivative operator bounds. To see the mixed terms explicitly, with $m_a=a+4\pi^2\norm\omega^2$,
\[
 \partial_{ua}f_\theta=-\frac{\beta_0+v}{m_a}f_\theta,
 \quad
 \partial_{va}f_\theta
 =-\frac{1+(\beta_0+v)g_{n,a}}{m_a}f_\theta,
 \quad
 \partial_{aa}f_\theta
 =\frac{(\beta_0+v)(\beta_0+v+1)}{m_a^2}f_\theta.
\]
The corresponding matrices $M_{ua},M_{va},M_{aa}$ have Frobenius norms at most $Cb_n$,
$Cb_n\log n$, and $CL_{2,n}$ for $d\ge4$. Higher mixed derivatives have the same structure: additional range derivatives add inverse powers of $m_a$, and additional
$v$ derivatives add centered logarithms. Nugget differentiation is simpler: $M_t=C_{n,0}^{-1}$ and all second derivatives involving $t$ are zero. For $d\le3$, the same estimates in Lemma~\ref{lem:upper} use $L_{1,n}=1$ and give $\norm{M_a}_\Frob+\norm{M_{ua}}_\Frob\le C$ and $\norm{M_{va}}_\Frob\le C\log n$, uniformly on $\mathcal N_n$.

\noindent\textbf{The high-dimensional Taylor remainder.}
Suppose $d\ge4$ and write $\delta a=a-a_0$. In this paragraph, $M(u,v,a,t)$ denotes the map in the centered, unnormalized coordinates; $M(h)$ denotes its expression in normalized coordinates. Write $M_{u,0},M_{v,0},M_{a,0}$ for derivatives at $(0,0,a_0,t_0)$. Then
\[
 B(h)=uM_{u,0}+vM_{v,0}+\delta a\,M_{a,0}
                         +(t-t_0)C_{n,0}^{-1}.
\]
To bound the remainder on the entire preliminary neighborhood, expand in $(u,v)$ at the current $a$, compare those first derivatives with their values at $a_0$, and expand the baseline in $a$. Taylor's formula gives the exact decomposition
\[
 M(h)-I_n-B(h)=\mathcal R_n(h)=R_{uv}+R_{\mathrm{mix}}+R_{aa},
\]
where
\begin{align*}
 R_{uv}
 &=\int_0^1(1-s)
   [u^2M_{uu}+2uvM_{uv}+v^2M_{vv}](su,sv,a,t_0)\,ds,\\
 R_{\mathrm{mix}}
 &=\delta a\int_0^1
   [uM_{ua}+vM_{va}](0,0,a_0+s\delta a,t_0)\,ds,\\
 R_{aa}
 &=(\delta a)^2\int_0^1(1-s)
                  M_{aa}(0,0,a_0+s\delta a,t_0)\,ds.
\end{align*}
There is no nugget remainder because $M$ is affine in $t$. All interpolation points remain in the enclosing preliminary set $\mathcal N_n^\square$, where the same derivative bounds hold; they need not belong to $\Theta$. The bounds in part (i) therefore give
\begin{align*}
 \norm{\mathcal R_n(h)}_\Frob
 \le C\{&r_n(|u|+|v|)^2
 +b_n|\delta a|(|u|+|v|\log n)
 +L_{2,n}|\delta a|^2\}.
\end{align*}

For $h=(r_nu,r_nv,b_n\delta a,\sqrt n(t-t_0))^\T$, these three terms are bounded respectively by $C(|u|+|v|)\norm h$, $C(|u|+|v|\log n)\norm h$, and $C(L_{2,n}/b_n)\norm h$, using boundedness of $\delta a$ for the last one. Since $|u|\le q_n\log n$ and $|v|\le q_n$ on $\mathcal N_n$,
\[
 \norm{\mathcal R_n(h)}_\Frob
 \le C\left(q_n\log n+\frac{L_{2,n}}{b_n}\right)\norm h
 \eqqcolon \epsilon_n\norm h,\qquad \epsilon_n\longrightarrow0.
\]
Here we used
\begin{equation}\label{eq:L2ratios}
 L_{2,n}/b_n\to0,\qquad L_{2,n}/b_n^2\to0,
 \qquad d\ge4.
\end{equation}
For example, $L_{2,n}$ is bounded below dimension eight, of order $\sqrt{\log n}$ at eight, and of order $R_n^{(d-8)/2}$ above eight; comparison with $b_n$ gives both
claims.

For the full covariance map, the normalized first-derivative bounds give $\norm{DM(h)[k]}_\Frob\le C\norm k$ throughout $\mathcal H_n$. Integration along segments proves \eqref{eq:fullcovlip}. This proves parts (ii)--(iii).

\noindent\textbf{Range perturbations in low dimension.}
Let $d\le3$. We obtain the needed covariance and tangent differences by integrating the range derivatives just bounded. In normalized coordinates, at the fixed identifiable truth,
\[
 D_n(a)=M(0,a)-M(0,a_0)=\int_{a_0}^a M_a(0,b)\,db.
\]
Since the range interval has fixed length and $\sup_b\norm{M_a(0,b)}_\Frob\le C$, this gives $\sup_a\norm{D_n(a)}_\Frob\le C$. For the identifiable tangents, the same identity gives
\[
 \partial_{h_i}M(0,a)-B_i
       =\int_{a_0}^a\partial_a\partial_{h_i}M(0,b)\,db.
\]
The normalization is $\partial_{h_1}=r_n^{-1}\partial_u$, $\partial_{h_2}=r_n^{-1}\partial_v$, and $\partial_{h_3}=n^{-1/2}\partial_t$. Thus the preceding bounds on $M_{ua}$ and $M_{va}$ prove \eqref{eq:rangetangentbounds}. The third difference is zero because $M_t=C_{n,0}^{-1}$ is independent of range.

\noindent\textbf{Low-dimensional increments in normalized coordinates.}
Let $h=(r_nu,r_nv,\sqrt n(t-t_0))^\T$. At each fixed $a$, Taylor's formula gives
\[
 M(h,a)=M(0,a)+D_hM(0,a)[h]+T_n(h,a),
 \quad
 T_n(h,a)=\int_0^1(1-s)D_h^2M(sh,a)[h,h]\,ds.
\]
Since $M(0,a)=I_n+D_n(a)$ and $B(h)=D_hM(0,a_0)[h]$,
\eqref{eq:Mlowexpansion} holds with
\[
 \mathcal R_n(h,a)
 =T_n(h,a)+\{D_hM(0,a)-D_hM(0,a_0)\}[h].
\]
The pure $(u,v)$ second-derivative bounds control the first term, and the tangent-difference bounds control the second:
\begin{align*}
 \norm{T_n(h,a)}_\Frob
       &\le Cr_n(|u|+|v|)^2
        \le Cq_n\log n\,\norm h,\\
 \norm{\{D_hM(0,a)-D_hM(0,a_0)\}[h]}_\Frob
       &\le C\frac{\log n}{r_n}\norm h.
\end{align*}
Thus \eqref{eq:Mlowglobalbounds} holds with $\epsilon_n=C(q_n\log n+\log n/r_n)\to0$; the bound on $D_n(a)$ was proved in the preceding paragraph.

Finally, the normalized identifiable derivatives satisfy $\norm{D_hM(h,a)[k]}_\Frob\le C\norm k$, and $\norm{\partial_aM(h,a)}_\Frob\le C$ by part (i). Integrating along the segment between $(h,a)$ and $(h',a')$ proves \eqref{eq:Mlowlip}. This completes part (iv).

For later use in parts (v)--(vii), first record one common low-dimensional mixed-derivative calculation. On $\Gamma_H$, write $\rho_\theta=f_\theta/f_0$ and choose $0<\zeta<4-d$. Equation~\eqref{eq:localratioexact} gives $\rho_\theta\le C_H\langle\omega\rangle^{2|v|}$, and eventually $4|v|\le\zeta$ uniformly. A normalized amplitude or smoothness derivative supplies a factor $r_n^{-1}$, and a centered smoothness derivative adds at most $C(\log n+\log\langle\omega\rangle)$. Let $\gamma$ be a multi-index in these two normalized coordinates, with $1\le|\gamma|\le2$, and let $j\ge1$ with $|\gamma|+j\le4$. The Frobenius bound \eqref{eq:weightedmatrixbounds} gives
\[
 \begin{aligned}
 \|\partial_h^\gamma\partial_a^jM\|_\Frob^2
 &\le C_Hr_n^{-2|\gamma|}
   \int_{\R^d}\langle\omega\rangle^{-4j+\zeta}
       (\log n+\log\langle\omega\rangle)^{2|\gamma|}\,d\omega\\
 &\le C_Hr_n^{-2|\gamma|}(\log n)^{2|\gamma|}.
 \end{aligned}
\]
The integral is finite because $4j-\zeta>d$; its fixed logarithmic moments are finite as well. Derivatives involving nugget vanish because $M$ is affine in $t$ and its nugget term is independent of the signal parameters. Taking square roots proves \eqref{eq:localnuisancederivatives}.

\noindent\textbf{Normalized derivatives for the likelihood Hessian.}
We prove part (v). In high dimension, differentiation in $h_1,h_2,h_3,h_4$ contributes the factors $r_n^{-1},r_n^{-1},b_n^{-1},n^{-1/2}$, respectively. Part (i) therefore bounds the first derivatives in Frobenius norm and, using $r_n\le\sqrt n$, gives
\[
 \max_i\norm{M_i}_\op
 \le C\left(\frac{\log n}{r_n}+\frac1{b_n}\right).
\]
For second derivatives it gives
\[
 \begin{aligned}
 \max_{i,j\le2}\norm{M_{ij}}_\Frob&\le C/r_n,\\
 \norm{M_{13}}_\Frob&\le C/r_n,\qquad
 \norm{M_{23}}_\Frob\le C\log n/r_n,\\
 \norm{M_{33}}_\Frob&\le CL_{2,n}/b_n^2,
 \qquad M_{i4}=0.
 \end{aligned}
\]
These all vanish. Lemma~\ref{lem:upper} gives the third derivatives as well; for $i,j,k\le2$,
\[
 \begin{aligned}
 \norm{M_{ijk}}_\Frob&\le C/r_n^2,\\
 \norm{M_{ij3}}_\Frob&\le C(\log n)^2/r_n^2,\\
 \norm{M_{i33}}_\Frob&\le CL_{2,n}\log n/(r_nb_n^2),\\
 \norm{M_{333}}_\Frob&\le CL_{3,n}/b_n^3.
 \end{aligned}
\]
Third derivatives involving $h_4$ vanish. By \eqref{eq:rangeupper}, $L_{3,n}\le CL_{2,n}\le Cb_n$. Since $b_n\ge1$ and $(\log n)/r_n\le1$ eventually, every second- and third-derivative bound displayed above is at most
\[
 C\left(\frac{\log n}{r_n}+\frac{L_{2,n}}{b_n^2}\right).
\]
This tends to zero by \eqref{eq:L2ratios} and gives the stated high-dimensional rate.

In low dimension the normalized coordinates are $(r_nu,r_nv,\sqrt n(t-t_0))$. Pure signal derivatives of orders two and three have Frobenius bounds $C/r_n$ and $C/r_n^2$,
respectively, uniformly in $a$; derivatives involving nugget vanish. Thus their sum is at most $C/r_n$, and part (i) gives $\max_i\norm{M_i}_\op\le C(\log n)/r_n$. The nuisance bounds follow from the common mixed-derivative calculation above and $\norm{M_a}_\Frob\le C$ from part (i).

To compare the tangents with those at the truth, integrate their derivatives along a local segment. In high dimension,
\[
 \norm{M_i(h)-B_i}_\Frob
 =\left\|\int_0^1\sum_jh_jM_{ij}(sh)\,ds\right\|_\Frob
 \le C_H\left(\frac{\log n}{r_n}+\frac{L_{2,n}}{b_n^2}\right)
\]
uniformly on $\norm h\le H$. The right side is at most $C_H\{(\log n)/r_n+b_n^{-1}\}$ because $L_{2,n}\le Cb_n$. In low dimension, integration at fixed $a$ gives $C_H/r_n$;
the tangent differences in part (iv) then give $\norm{M_i(h,a)-B_i}_\Frob\le C_H(\log n)/r_n$. Finally, parts (ii) and (iv) give $\norm{M(h)-I_n}_\Frob\le C_H$ or
\[
 \norm{M(h,a)-I_n}_\Frob
 \le\norm{M(h,a)-M(0,a)}_\Frob+\norm{D_n(a)}_\Frob\le C_H.
\]
This proves all of \eqref{eq:studentderiv}.

\noindent\textbf{Matrices in the identifiable Hessian.}
We prove part~\ref{it:covariance-hessian}. Work on $\Gamma_H$ and write $u=h$ or $u=(h,a)$, with $k$ its total dimension.

Parts (i) and (v) give $\|R\|_\op\le C$, $\max_i\|M_i\|_\Frob\le C$, $\max_i\|M_i\|_\op=o(1)$, and $\max_{i,j}\|M_{ij}\|_\Frob=o(1)$. The definition of $A_{ij}$ therefore yields
\[
 \|A_{ij}\|_\Frob
 \le C\{\|M_{ij}\|_\Frob
          +\|M_i\|_\op\|M_j\|_\Frob
          +\|M_j\|_\op\|M_i\|_\Frob\}=o(1)
\]
uniformly on $\Gamma_H$.

For the increment bound, write $M_\ell=\partial_{u_\ell}M$,
$M_{i\ell}=\partial_{u_\ell}M_i$, and $M_{ij\ell}=\partial_{u_\ell}M_{ij}$; thus $u_\ell$ may be a normalized coordinate or the unnormalized range nuisance. Equation \eqref{eq:studentderiv} and the nuisance bounds in part~\ref{it:covariance-nuisance} give, uniformly,
\[
 \max_\ell\|M_\ell\|_\op\le C,\qquad
 \max_{i,\ell}\|M_{i\ell}\|_\Frob
       +\max_{i,j,\ell}\|M_{ij\ell}\|_\Frob=o(1).
\]
Use $\partial_{u_\ell}R=-RM_\ell R$ and differentiate the three products defining $A_{ij}$. Applying $\|UVW\|_\Frob\le\|U\|_\op\|V\|_\Frob\|W\|_\op$ to each
resulting term gives
\[
 \begin{aligned}
 \|\partial_{u_\ell}A_{ij}\|_\Frob\le C\bigl[&
       \|M_{ij\ell}\|_\Frob
       +\|M_\ell\|_\op\|M_{ij}\|_\Frob+\|M_{i\ell}\|_\Frob\|M_j\|_\op
       +\|M_i\|_\op\|M_{j\ell}\|_\Frob\\
 &+\|M_\ell\|_\op\{
       \|M_i\|_\Frob\|M_j\|_\op
       +\|M_j\|_\Frob\|M_i\|_\op\}\bigr]=o(1).
 \end{aligned}
\]
Each term contains a factor tending to zero; all the other factors are bounded. Integrating along the segment from $v$ to $u$ gives
\[
 A_{ij}(u)-A_{ij}(v)
 =\int_0^1\sum_{\ell=1}^k(u_\ell-v_\ell)
       \partial_{u_\ell}A_{ij}(v+s(u-v))\,ds.
\]
Taking Frobenius norms and using $\sum_\ell|u_\ell-v_\ell|\le\sqrt{k}\|u-v\|$ yields
\[
 \|A_{ij}(u)-A_{ij}(v)\|_\Frob
 \le\sqrt{k}\max_\ell\sup_{w\in\Gamma_H}
       \|\partial_{u_\ell}A_{ij}(w)\|_\Frob\,\|u-v\|.
\]
The coefficient even tends to zero, and in particular is bounded by a fixed $L$. Each small factor above is bounded by the rate in part~\ref{it:covariance-hessian}: $C_H(\log n)/r_n$ in low dimension and $C_H\{(\log n)/r_n+b_n^{-1}\}$ in high dimension. The same inequalities therefore prove the quantitative bounds in that
part, using covariance derivatives only through third order.

\noindent\textbf{Nuisance derivatives, traces, and mixed Hessian matrices.}
We prove part~\ref{it:covariance-nuisance}. Let $d\le3$ and work on $\Gamma_H$. Part (i), the normalized first-derivative bounds in (v), and the common mixed-derivative calculation above give
\begin{equation*}
 \begin{aligned}
 \sup_{\Gamma_H}\{\norm R_\op+\norm{M-I_n}_\Frob
                         +\norm{M_a}_\Frob\}&\le C_H,\\
 \epsilon_n\coloneqq \sup_{\Gamma_H}\max_i
       \{\norm{M_i}_\op+\norm{M_{ia}}_\Frob\}
       &\le C_H\frac{\log n}{r_n}=o(1).
 \end{aligned}
\end{equation*}
For the trace norms, put $\theta=(\vartheta_0+E_nh,a)$, $m_a=a+4\pi^2\norm\omega^2$, and $\beta=\nu+d/2$. Then
\begin{equation*}
 \rho_\theta=\frac{f_\theta}{f_0},\qquad
 -M_a=\beta\mathcal B_n(\rho_\theta/m_a),\qquad
 \partial_a^jM=(-1)^j(\beta)_j
                 \mathcal B_n(\rho_\theta m_a^{-j}),
\end{equation*}
where $(\beta)_j=\beta(\beta+1)\cdots(\beta+j-1)$. The displayed multipliers have a fixed sign, so their trace norms are their absolute traces.

On $\norm\omega\le n^A$, the covariance comparison argument above gives $\rho_\theta\le C_H$, and $m_a\asymp m_0$. The contribution of the complement to $\|M_a\|_1$ is at most
\[
 \frac n{t_0}\int_{\norm\omega>n^A}
                 |\partial_af_\theta(\omega)|\,d\omega
 \le Cn^{1-2A(\nu_-+1)}=o(1)
\]
for sufficiently large fixed $A$, by the tail estimate in Lemma~\ref{lem:upper}. Thus \eqref{eq:rangetrace} with $k=1$ gives
\begin{equation*}
 \sup_{\Gamma_H}\norm{M_a}_1
 =\sup_{\Gamma_H}(-\tr M_a)
 \le C_HT_{d,n}
 =C_H\begin{cases}
 1,&d=1,\\
 \log n,&d=2,\\
 R_n,&d=3.
 \end{cases}
\end{equation*}
For higher range derivatives, a direct integrable bound suffices. Use the bound $\rho_\theta\le C_H\langle\omega\rangle^{2|v|}$ and $0<\zeta<4-d$ from the common mixed-derivative calculation above; eventually $4|v|\le\zeta$ uniformly. The trace formula and $w_n\le C$ give, for $2\le j\le4$,
\[
 \|\partial_a^jM\|_1
 =(\beta)_j\int\rho_\theta m_a^{-j}w_n
 \le C_H\int_{\R^d}\langle\omega\rangle^{-2j+\zeta}\,d\omega
 \le C_H.
\]
In particular this proves the trace-norm bound for $M_{aa}$.

Trace-norm duality now gives the cross-information bound directly on the whole local set:
\begin{equation*}
 \begin{aligned}
 \rho_n\coloneqq \sup_{\Gamma_H}\max_i|\tr(RM_iRM_a)|
 &\le C_H\frac{\log n}{r_n}\sup_{\Gamma_H}\norm{M_a}_1\\
 &\le C_H\begin{cases}
 (\log n)R_n^{-1/2},&d=1,\\
 (\log n)^2R_n^{-1},&d=2,\\
 (\log n)R_n^{-1/2},&d=3,
 \end{cases}
 \qquad\longrightarrow0.
 \end{aligned}
\end{equation*}
Write $\mathscr I_{ia}=\tfrac12\tr(RM_iRM_a)$. Differentiating $R=M^{-1}$ in $a$, then applying Frobenius Cauchy--Schwarz and trace-norm duality, yields
\[
 \begin{aligned}
 |\partial_a\mathscr I_{ia}|
 &\le C_H\{\norm{M_{ia}}_\Frob\norm{M_a}_\Frob
       +\norm{M_i}_\op(\norm{M_{aa}}_1+\norm{M_a}_\Frob^2)\}\\
 &\le C_H\frac{\log n}{r_n}.
 \end{aligned}
\]

The range information is
\[
 \calI_{aa,n}(\theta)
 =\tfrac12\tr(RM_aRM_a)
 =\tfrac12\norm{R^{1/2}M_aR^{1/2}}_\Frob^2\le C_H.
\]
For a lower bound, let $s(\theta)=K_\theta(0)$, so that $C_{11}=s+t$ and $\partial_aC_{11}=-\nu s/a$ by \eqref{eq:inverse}. The unit vector $w=C_n(\theta)^{1/2}e_1/\sqrt{C_{11}}$ gives
\[
 \begin{aligned}
 \calI_{aa,n}(\theta)
 &=\tfrac12\norm{C_n(\theta)^{-1/2}C_a
                         C_n(\theta)^{-1/2}}_\Frob^2\\
 &\ge\tfrac12\left(
 w^\T C_n(\theta)^{-1/2}C_aC_n(\theta)^{-1/2}w\right)^2
 =\tfrac12\left(\frac{\nu s}{a(s+t)}\right)^2\ge c>0.
 \end{aligned}
\]
The last constant is uniform on the parameter rectangle. Thus even one marginal observation gives a positive amount of range information.

It remains to bound $A_{ia}$ and $\partial_aA_{ia}$ and their first partial derivatives. Start with
\[
 A_{ia}=RM_{ia}R-RM_iRM_aR-RM_aRM_iR.
\]
In a differentiated term from $RM_{ia}R$, the central covariance derivative still has at least one normalized $h$ derivative and one range derivative. Its Frobenius norm is at most $C_H(\log n)/r_n$ by the preceding bounds. In a term from either of the other two products, there are two possibilities. If a derivative hits $M_i$, it produces such a mixed derivative, or a second normalized derivative with Frobenius norm $O_H(r_n^{-1})$. Otherwise $\|M_i\|_\op\le C_H(\log n)/r_n$ multiplies a range
derivative with bounded Frobenius norm. The derivatives of $R$ that occur here are products of bounded inverse and covariance factors, so they have bounded operator norm. Applying $\|UVW\|_\Frob\le\|U\|_\op\|V\|_\Frob\|W\|_\op$ term by term proves
\[
 \sup_{\Gamma_H}\max_i\left\{
 \norm{A_{ia}}_\Frob+\norm{\partial_aA_{ia}}_\Frob
 +\max_\ell\norm{\partial_{u_\ell}A_{ia}}_\Frob
 +\max_\ell\norm{\partial_{u_\ell}\partial_aA_{ia}}_\Frob
 \right\}\le C_H\frac{\log n}{r_n}.
\]
The largest covariance derivative used here has order four: two identifiable and two range derivatives, or one identifiable and three range derivatives. Integration along segments, as in part~\ref{it:covariance-hessian}, gives the stated Lipschitz bounds.

\noindent\textbf{Identifiable tangents and covariance separation.}
Assume \Cref{ass:sharp}. Lemma~\ref{lem:information} and $\|B(h)\|_\Frob^2=2h^\T J_nh$ give both tangent bounds in part~\ref{it:covariance-separation}. Work on the entire preliminary set $\mathcal H_n$ for the remaining assertions. In high dimension,
\[
 M(h)-I_n=B(h)+\mathcal R_n(h),\qquad
 \|B(h)\|_\Frob^2=2h^\T J_nh\ge c\|h\|^2,
 \qquad \|\mathcal R_n(h)\|_\Frob\le\epsilon_n\|h\|,
\]
where $\epsilon_n\to0$ uniformly, by Lemma~\ref{lem:information} and parts (iii)--(iv).
The triangle inequality therefore gives $\|M(h)-I_n\|_\Frob\ge c'\|h\|$ for all sufficiently large $n$. In low dimension the additional term is $D_n(a)=M(0,a)-I_n$, with $\sup_a\|D_n(a)\|_\Frob\le C$. The same calculation gives $\|M(h,a)-I_n\|_\Frob\ge(c'\|h\|-C)_+$. Squaring and using $(x-y)_+^2\ge x^2/2-y^2$ yields
\[
 \|M(u)-I_n\|_\Frob^2\ge c\|h\|^2-C.
\]

The same decompositions, the tangent upper bound, and the triangle inequality give $\|M(h)-I_n\|_\Frob\le C\|h\|$ in high dimension and $\|M(h,a)-I_n\|_\Frob\le C(1+\|h\|)$ in low dimension. This proves all the separation bounds.

\noindent\textbf{Gaussian signal bounds for nuisance derivatives.}
We prove part~\ref{it:covariance-signal}. Assume \Cref{ass:sharp} and $d\le3$, fix $h=0$, and write $C_0=C_{n,0}$, $P_0=C_0^{-1/2}S_{n,0}C_0^{-1/2}$, and $R=M(0,a)^{-1}$.
Recall that $0\preceq P_0\preceq I$ and $M\asymp I$ uniformly. We first establish one weighted Cauchy--Schwarz inequality. It will show both that applying $R$ preserves the directions carrying signal and that $A_q$ cannot have substantial Frobenius norm concentrated entirely in directions where $P_0$ is small.

\smallskip
\noindent\emph{A weighted Cauchy--Schwarz inequality.}
For any bounded real even multiplier $g$, the Cauchy--Schwarz calculation proving \eqref{eq:weightedmatrixbounds}, before its $P_0$ factor is bounded by the identity, gives
\[
 |y^\T\mathcal B_n(g)x|^2
 \le (y^\T P_0y)\{x^\T\mathcal B_n(g^2)x\},
 \qquad x,y\in\R^n.
\]
If $y\in\ker P_0$, this inequality shows that $y^\T\mathcal B_n(g)x=0$ for every $x$. Since $\mathcal B_n(g)$ is symmetric, it therefore vanishes on $\ker P_0$ and has range contained in $\operatorname{ran}P_0$.

This applies to $M-I,M_a,M_q$, and $M_{qq}$: each is $\mathcal B_n(g)$ for a bounded real even multiplier. Consequently $M$ and $R$ act as the identity on $\ker P_0$, and $B_a$ and $A_q$ vanish there. In calculations involving inverse powers of $P_0$, we therefore work on $\operatorname{ran}P_0$ and invert $P_0$ only on that subspace; we do not assume these inverses are uniformly bounded. Restricting the covariance derivatives, $B_a$, and $A_q$ to this subspace preserves their Frobenius norms
and nonzero eigenvalues. The other formulas may be read on $\R^n$. This also covers designs with repeated observation sites.

In the preceding Cauchy--Schwarz inequality, take $y=P_0^{-1}\mathcal B_n(g)x$. Then
\[
 \begin{aligned}
 \{x^\T\mathcal B_n(g)P_0^{-1}\mathcal B_n(g)x\}^2
 &\le \{x^\T\mathcal B_n(g)P_0^{-1}\mathcal B_n(g)x\}
           \{x^\T\mathcal B_n(g^2)x\}.
 \end{aligned}
\]
Canceling the common nonnegative factor gives the first bound below; when that factor is zero the bound is immediate. Taking traces then gives the second:
\begin{equation}\label{eq:range-weighted-factorization}
 \begin{aligned}
 \mathcal B_n(g)P_0^{-1}\mathcal B_n(g)
     &\preceq\mathcal B_n(g^2),\\
 \|P_0^{-1/2}\mathcal B_n(g)\|_\Frob^2
     &=\tr\{\mathcal B_n(g)P_0^{-1}\mathcal B_n(g)\}\\
     &\le\tr\mathcal B_n(g^2)
       =\int_{\R^d}g(\omega)^2w_n(\omega)\,d\omega.
 \end{aligned}
\end{equation}

\smallskip
\noindent\emph{The inverse covariance preserves signal directions.}
At $h=0$, the spectral ratio $f_{\eta_0,\nu_0,a}/f_0$ is uniformly bounded above and below. Thus the multiplier in
\[
 M-I=\mathcal B_n\left(\frac{f_{\eta_0,\nu_0,a}}{f_0}-1\right)
\]
is uniformly bounded. Since $P_0^{-1}\succeq I$ on $\operatorname{ran}P_0$, \eqref{eq:range-weighted-factorization}
implies
\[
 (M-I)^2\preceq(M-I)P_0^{-1}(M-I)\preceq CP_0.
\]
In particular, $\|(I-M)x\|\le C\|P_0^{1/2}x\|$. The identity $R=I+R(I-M)$ now gives
\[
 \begin{aligned}
 \|P_0^{1/2}Rx\|
 &\le\|P_0^{1/2}x\|
       +\|P_0^{1/2}R\|_\op\|(I-M)x\|\\
 &\le C\|P_0^{1/2}x\|,
 \end{aligned}
\]
because $\|P_0^{1/2}\|_\op\le1$ and $\|R\|_\op\le C$. Squaring proves $RP_0R\preceq CP_0$. Also $0\preceq-M_a\preceq CP_0$, because $-\partial_af_{\eta_0,\nu_0,a}/f_0$ is nonnegative and uniformly bounded. Congruencing this last order bound by $R$ yields
\begin{equation}\label{eq:range-signal-domination}
 RP_0R\preceq CP_0,\qquad
 0\preceq B_a\coloneqq -RM_aR\preceq CP_0.
\end{equation}

\smallskip
\noindent\emph{Two Gaussian directions have uniformly positive variance.}
Choose two fixed distinct points $x^{(1)},x^{(2)}$ in the coverage cube. Coverage supplies indices $i_{n,1},i_{n,2}$ such that $x_{n,i_{n,j}}\to x^{(j)}$; the indices are distinct for all sufficiently large $n$. If $e_1,\ldots,e_n$ are the standard basis
vectors of $\R^n$, define the two-by-$n$ matrix
\[
 E=\begin{pmatrix}e_{i_{n,1}}^\T\\ e_{i_{n,2}}^\T\end{pmatrix}
   =\begin{pmatrix}e_{i_{n,1}}&e_{i_{n,2}}\end{pmatrix}^{\!\T},
 \qquad
 Ey=\begin{pmatrix}y_{i_{n,1}}\\y_{i_{n,2}}\end{pmatrix}.
\]
The spectral density $-\partial_af_{\eta_0,\nu_0,a}$ is strictly positive. Its covariance at two distinct points is positive definite: for $b\in\R^2\setminus\{0\}$, its quadratic form is
\[
 \int_{\R^d}(-\partial_af_{\eta_0,\nu_0,a})(\omega)
       \left|b_1e^{2\pi i\omega^\T x^{(1)}}
                  +b_2e^{2\pi i\omega^\T x^{(2)}}\right|^2d\omega>0.
\]
The two distinct exponential functions are linearly independent. Continuity in the sites and parameters, and compactness of their parameter ranges, therefore give
\[
 E(-C_a)E^\T\succeq cI_2,\qquad
 t_-I_2\preceq EC_0E^\T\preceq CI_2,
\]
where $C_a$ denotes the range derivative of the observation covariance at the current value $a$. We can apply min--max directly on the two-dimensional subspace $\operatorname{ran}(C_0^{1/2}E^\T)$. For $b\ne0$,
\[
 \frac{(C_0^{1/2}E^\T b)^\T(-M_a)(C_0^{1/2}E^\T b)}
      {\|C_0^{1/2}E^\T b\|^2}
 =\frac{b^\T E(-C_a)E^\T b}{b^\T EC_0E^\T b}\ge c.
\]
Hence $\lambda_2(-M_a)\ge c$. Since $R\succeq cI$, the nonzero eigenvalues of $B_a=R(-M_a)R$ equal those of
\[
 (-M_a)^{1/2}R^2(-M_a)^{1/2}\succeq c^2(-M_a).
\]
Thus $\lambda_2(B_a)\ge c$ as well. By \eqref{eq:range-signal-domination},
\[
 B_a^{1/2}P_0B_a^{1/2}\succeq C^{-1}B_a^2.
\]
Its nonzero eigenvalues equal those of $P_0^{1/2}B_aP_0^{1/2}$. The latter matrix therefore has two eigenvalues bounded below by a uniform positive constant.

\smallskip
\noindent\emph{A nonnegligible nuisance Hessian retains Gaussian signal.}
Let $a=\alpha(q)$ be the coordinate change in the statement and write $m_a(\omega)=a+4\pi^2\|\omega\|^2$. At the current parameter $\theta=(\eta_0,\nu_0,\alpha(q),t_0)$, the Mat\'ern derivative formulas are
\[
 \partial_af_\theta=-\frac{\beta_0}{m_a}f_\theta,
 \qquad
 \partial_a^2f_\theta
    =\frac{\beta_0(\beta_0+1)}{m_a^2}f_\theta.
\]
The chain rule gives
\[
 \partial_qf_\theta=\alpha'\partial_af_\theta,
 \qquad
 \partial_q^2f_\theta
    =(\alpha')^2\partial_a^2f_\theta+\alpha''\partial_af_\theta.
\]
The bounded coordinate derivatives and spectral ratio consequently imply
\[
 \left|\frac{\partial_qf_\theta(\omega)}{f_0(\omega)}\right|
 +\left|\frac{\partial_q^2f_\theta(\omega)}{f_0(\omega)}\right|
 \le C\langle\omega\rangle^{-2}.
\]
Since $\partial_q^jM=\mathcal B_n\{(\partial_q^jf_\theta)/f_0\}$ for $j=1,2$, \eqref{eq:range-weighted-factorization} gives
\[
 \begin{aligned}
 \|P_0^{-1/2}M_q\|_\Frob^2
 +\|P_0^{-1/2}M_{qq}\|_\Frob^2
 &\le C\int_{\R^d}\langle\omega\rangle^{-4}w_n(\omega)\,d\omega\\
 &\le C\int_{\R^d}\langle\omega\rangle^{-4}\,d\omega\le C.
 \end{aligned}
\]
Here we used $w_n\le C$ from \eqref{eq:frequencyweight} and $d\le3$. Also, \eqref{eq:range-signal-domination} implies
\[
 \|P_0^{-1/2}RP_0^{1/2}\|_\op^2
 =\|P_0^{-1/2}RP_0RP_0^{-1/2}\|_\op\le C.
\]
Insert $P_0^{1/2}P_0^{-1/2}$ between the leading $R$ and the covariance derivative in each term of $A_q=RM_{qq}R-2RM_qRM_qR$. The preceding bounds, together with $\|R\|_\op+\|M_q\|_\op\le C$, yield
\[
 \begin{aligned}
 \|P_0^{-1/2}A_q\|_\Frob
 &\le C\left\{\|P_0^{-1/2}M_{qq}\|_\Frob
      +\|P_0^{-1/2}M_q\|_\Frob\|M_q\|_\op\right\}\\
 &\le C.
 \end{aligned}
\]
This weighted bound prevents $A_q$ from concentrating its Frobenius norm in directions with arbitrarily little Gaussian signal.

It remains to apply Frobenius Cauchy--Schwarz twice. With $\langle U,V\rangle_\Frob=\tr(U^\T V)$,
\[
 \begin{aligned}
 \|A_q\|_\Frob^2
 &=\langle P_0^{-1/2}A_q,P_0^{1/2}A_q\rangle_\Frob\\
 &\le\|P_0^{-1/2}A_q\|_\Frob\|P_0^{1/2}A_q\|_\Frob
 \le C\|P_0^{1/2}A_q\|_\Frob.
 \end{aligned}
\]
For the second application, symmetry of $A_q$ and cyclicity of trace give
\[
 \begin{aligned}
 &\langle P_0^{-1/2}A_qP_0^{1/2},
                  P_0^{1/2}A_qP_0^{1/2}\rangle_\Frob\\
 &\qquad=\tr(P_0^{1/2}A_q^2P_0^{1/2})
       =\tr(A_qP_0A_q)=\|P_0^{1/2}A_q\|_\Frob^2.
 \end{aligned}
\]
Since $\|P_0^{1/2}\|_\op\le1$, the first matrix in this inner product has Frobenius norm at most $C$. Therefore
\[
 \|P_0^{1/2}A_q\|_\Frob^2
 \le C\|P_0^{1/2}A_qP_0^{1/2}\|_\Frob,
 \qquad
 \|A_q\|_\Frob^4
 \le C^3\|P_0^{1/2}A_qP_0^{1/2}\|_\Frob.
\]
In particular,
\begin{equation}\label{eq:range-gaussian-part}
 \|A_q\|_\Frob\ge\eta>0
 \quad\Longrightarrow\quad
 \|P_0^{1/2}A_qP_0^{1/2}\|_\Frob
       \ge C^{-3}\eta^4\eqqcolon c_\eta>0.
\end{equation}
All bounds are uniform in the locations and parameters under the stated assumptions. This proves part~\ref{it:covariance-signal}.

\qed

\subsection{Proof of Lemma~\ref{lem:localization}}
\noindent\textbf{The deterministic likelihood contrast.}
Write $\ell_n(h,a)=\ell_n(\vartheta_0+E_nh,a)$ in low dimension and $\ell_n(h)=\ell_n(\theta_0+E_nh)$ in high dimension. For this proof let $\mathcal H_n$ be the normalized image of the enlarged preliminary set $\mathcal N_n^\square$, rather than only the feasible candidates. It is compact and convex: it is the affine
inverse image of $\mathcal R$ intersected with the linear inequalities in \eqref{eq:localset}. It contains the true index, $(0,a_0)$ or $0$, and every feasible candidate in $\mathcal N_n$. All covariance bounds in Lemma~\ref{lem:derivatives} hold on this enclosing set. A supremum over it therefore bounds the event generated by a fit over any compact $\Theta\subset\mathcal R$, without assuming $\Theta$ is convex.

Put $X=C_{n,0}^{-1/2}Y_n$ and abbreviate the likelihood increment by
\[
 \mathcal L_n(u)=\ell_n(u)-\ell_n(u_0)
       =-\tfrac12\log\det M(u)
        -\tfrac12X^\T\{M(u)^{-1}-I_n\}X,
\]
where $u=(h,a)$, $u_0=(0,a_0)$ in low dimension and $u=h$, $u_0=0$ otherwise. Since $\E XX^\T=I_n$, its negative expectation is
\[
 \mathcal D_n(u)\coloneqq -\E\mathcal L_n(u)
       =\tfrac12\{\tr M(u)^{-1}-n+\log\det M(u)\}
       \asymp\|M(u)-I_n\|_\Frob^2.
\]
The last comparison is \eqref{eq:local-contrast-comparison}, using the covariance bounds of Lemma~\ref{lem:derivatives}.

Lemma~\ref{lem:derivatives}\ref{it:covariance-separation} gives $\|M(u)-I_n\|_\Frob^2\ge c\|h\|^2-C$ throughout $\mathcal H_n$. Consequently
\begin{equation*}
 \mathcal D_n(u)\ge c\|h\|^2-C.
\end{equation*}
Thus a bounded nuisance perturbation cannot remove the quadratic separation as $\|h\|$ grows.

\noindent\textbf{Fractional moments and the supremum.}
We apply the Ibragimov--Has'minskii argument to a small power of the Gaussian likelihood ratio. Let $k$ be the total dimension of $u$ and fix $m>k$. Choose the fixed $s>0$ supplied by Lemma~\ref{lem:local-fractional} and set
\[
 \xi_n(u)=e^{s\mathcal L_n(u)/m},\qquad u\in\mathcal H_n.
\]
Besides the covariance comparison and separation just verified, Lemma~\ref{lem:derivatives} gives
\[
 \|M(u)-M(v)\|_\Frob\le C\|u-v\|,
\]
by \eqref{eq:fullcovlip} or \eqref{eq:Mlowlip}. The latter uses the equivalence of $\|h-h'\|+|a-a'|$ with the Euclidean distance on the fixed-dimensional joint index set. Lemma~\ref{lem:local-fractional} therefore gives, with uniform constants,
\[
 \E|\xi_n(u)|^m\le A e^{-c_0\|h\|^2},\qquad
 \|\xi_n(u)-\xi_n(v)\|_{L^m}\le C\|u-v\|.
\]
The maximal inequality in Lemma~\ref{lem:quadratic-envelope} now yields
\[
 \E\sup_{\substack{u\in\mathcal H_n\\\|h\|>H}}
           |\xi_n(u)|^m\le Ce^{-cH^2}.
\]
If $\mathcal L_n(u)\ge0$, then $|\xi_n(u)|^m=e^{s\mathcal L_n(u)}\ge1$. Markov's inequality consequently gives
\[
 \begin{aligned}
 P\{\exists u\in\mathcal H_n:\|h\|>H,\ \mathcal L_n(u)\ge0\}
 &\le P\!\left\{\sup_{\substack{u\in\mathcal H_n\\\|h\|>H}}
                          |\xi_n(u)|^m\ge1\right\}\le Ce^{-cH^2}.
 \end{aligned}
\]
Every global maximizer has likelihood at least that of the truth. If it lies in $\mathcal N_n$ and its normalized error exceeds $H$, it therefore belongs to the event in this display. Adding the exceptional probability from Lemma~\ref{lem:preliminary} proves \eqref{eq:tails}, including maximization over the varying feasible
nuisance intervals.

\noindent\textbf{From the tail bound to moments.}
Put $T_n=\|\widehat h_n\|$. On the compact parameter set, $T_n\le C_\Theta\sqrt n$: $r_n^2\log^2n=o(n)$ uniformly, $b_n^2\le Cr_n^2$, and the nugget contribution is at most $Cn$. Fix $\lambda>0$ and $0<\alpha<\min\{2,2\delta\}$. Tonelli's theorem and \eqref{eq:tails} give directly
\[
 \begin{aligned}
 \E e^{\lambda T_n^\alpha}
 &=1+\int_0^{C_\Theta\sqrt n}
      \lambda\alpha H^{\alpha-1}e^{\lambda H^\alpha}
                     P(T_n>H)\,dH\\
 &\le1+C\lambda\alpha\int_0^\infty
      H^{\alpha-1}e^{\lambda H^\alpha-cH^2}\,dH
       +Ce^{-cn^\delta}
          \bigl(e^{\lambda C_\Theta^\alpha n^{\alpha/2}}-1\bigr)\\
 &\le C_{\lambda,\alpha}
       +C e^{\lambda C_\Theta^\alpha n^{\alpha/2}-cn^\delta}
 \le C'_{\lambda,\alpha}.
 \end{aligned}
\]
The integral is finite because $\alpha>0$ gives integrability at zero and $\alpha<2$ makes $-cH^2$ dominate $\lambda H^\alpha$ at infinity. The exceptional term tends to zero since $\alpha/2<\delta$. All constants are uniform, proving \eqref{eq:allmoments}.

Fixing smoothness at its true value only restricts the set of candidates whose criterion is at least its value at the truth. The preliminary localization and the supremum bound therefore remain valid for the known-smoothness fit.

\noindent\textbf{Random designs.}
Let $\mathscr X_n=\sigma(x_{n,1:n})$ and take the location event $\mathcal G_n\in\mathscr X_n$ from Lemma~\ref{lem:randomgood}. Decrease $\delta>0$ if necessary so that it is no larger than both the preceding tail exponent and the design-failure exponent in \eqref{eq:randomfailure}. Thus $P_X(\mathcal G_n^c)\le Ce^{-cn^\delta}$.
On $\mathcal G_n$, the preceding arguments apply conditionally on $\mathscr X_n$ with the same deterministic constants for every such design:
\[
 \begin{aligned}
 P(T_n>H\mid\mathscr X_n)
     &\le Ce^{-cH^2}+Ce^{-cn^\delta},\qquad H\ge0,\\
 \E(e^{\lambda T_n^\alpha}\mid\mathscr X_n)
     &\le C_{\lambda,\alpha},\qquad
       \lambda>0,\quad 0<\alpha<\min\{2,2\delta\}.
 \end{aligned}
\]
Here unconditional probabilities and expectations include the locations. Splitting over $\mathcal G_n$ and its complement gives
\[
 \begin{aligned}
 P(T_n>H)
 &\le \E\!\left[\mathbf1_{\mathcal G_n}
                   P(T_n>H\mid\mathscr X_n)\right]
                  +P_X(\mathcal G_n^c)\\
 &\le Ce^{-cH^2}+Ce^{-cn^\delta}+P_X(\mathcal G_n^c)\\
 &\le C'e^{-c'H^2}+C'e^{-c'n^\delta}.
 \end{aligned}
\]
For the exponential moment, the deterministic bound $T_n\le C_\Theta\sqrt n$ also holds on $\mathcal G_n^c$. Hence
\[
 \begin{aligned}
 \E e^{\lambda T_n^\alpha}
 &=\E\!\left[\mathbf1_{\mathcal G_n}
             \E(e^{\lambda T_n^\alpha}\mid\mathscr X_n)\right]
       +\E\!\left[\mathbf1_{\mathcal G_n^c}e^{\lambda T_n^\alpha}\right]\\
 &\le C_{\lambda,\alpha}
       +e^{\lambda C_\Theta^\alpha n^{\alpha/2}}P_X(\mathcal G_n^c)\\
 &\le C_{\lambda,\alpha}
       +C e^{\lambda C_\Theta^\alpha n^{\alpha/2}-cn^\delta}
 \le C'_{\lambda,\alpha}.
 \end{aligned}
\]
Again the exceptional term tends to zero because $\alpha/2<\delta$. These are the unconditional versions of \eqref{eq:tails} and \eqref{eq:allmoments}, with constants uniform over $\mathcal K\times\mathcal P$. \qed

\subsection{Proof of Lemma~\ref{lem:nuisance}}
Apply Lemma~\ref{lem:score-comparison} with $I=[a_-,a_+]$ and $X_n=C_{n,0}^{-1/2}Y_n$, so $\E X_nX_n^\T=I_n$. Its quadratic-form moment assumption follows from Lemma~\ref{lem:gaussian}. The covariance family is smooth, positive definite, and satisfies $M(0,a_0)=I_n$. Lemma~\ref{lem:derivatives}, parts \ref{it:covariance-local} and \ref{it:covariance-nuisance}, gives the required bounds on $R$, $M-I_n$, and $M_a$, and, in the notation of Lemma~\ref{lem:score-comparison},
\[
 \epsilon_n\le C\frac{\log n}{r_n},\qquad
 \rho_n\le C\frac{\log n}{r_n}T_{d,n}\longrightarrow0.
\]
These estimates hold on the whole local set in Lemma~\ref{lem:derivatives}, so in particular at $(0,a)$ for every $a\in I$. Lemma~\ref{lem:score-comparison} therefore yields, for each fixed
$1\le p<\infty$,
\[
 \left\|\sup_{a\in I}\|s_n(a)-\Delta_n\|\right\|_{L^p}
       \le C_p(\epsilon_n+\rho_n)=o(1),
\]
which is \eqref{eq:nuisancescorecomparison}.

For the baseline derivative bounds, abbreviate $G_n(a)=\mathcal L_n(0,a)$, $R=M^{-1}$, and $Q(A)=X_n^\top AX_n-\tr A$. Differentiation gives
\[
 \begin{split}
 G_n'(a)
 &=\tfrac12Q(RM_aR)
       +\tfrac12\tr\{(I_n-M)RM_aR\},\\
 G_n''(a)
 &=\tfrac12Q(A_{aa})-\tfrac12\tr(RM_aRM_a)
       +\tfrac12\tr\{(I_n-M)A_{aa}\},\\
 A_{aa}&=RM_{aa}R-2RM_aRM_aR.
 \end{split}
\]
Lemma~\ref{lem:derivatives} bounds $M-I_n$ and $M_a$ in Frobenius norm, $M_{aa}$ in trace norm, and $R$ in operator norm. Hence $RM_aR$ and $A_{aa}$ have bounded Frobenius norms, and all deterministic terms in these displays are bounded. The quadratic-form moment inequality in Lemma~\ref{lem:gaussian} proves $\sup_a\|G_n^{(j)}(a)\|_{L^p}\le C_p$ for $j=1,2$.

For the third and fourth derivatives, every matrix product in the quadratic or trace terms contains either a covariance derivative of order at least two or at least two factors $M_a$. The former has bounded trace norm by Lemma~\ref{lem:derivatives}; for the latter use $\|AB\|_1\le\|A\|_\Frob\|B\|_\Frob$, while the remaining factors
have bounded operator norms. Thus the corresponding trace and Frobenius norms are bounded, and the same moment inequality gives
\[
 \max_{1\le j\le4}\sup_{a\in I_*}
                 \|G_n^{(j)}(a)\|_{L^p}\le C_p.
\]
For $j=1,2,3$, the fundamental theorem of calculus yields
\[
 \sup_{a\in I_*}|G_n^{(j)}(a)|
 \le |G_n^{(j)}(a_-)|+\int_{I_*}|G_n^{(j+1)}(a)|\,da.
\]
Minkowski's inequality now proves \eqref{eq:nuisance-baseline-moments}. Fixing smoothness removes its score coordinate and preserves these bounds. \qed

\subsection{Proof of Lemma~\ref{lem:affinitytails}}
The quadratic-form exponential-moment bound is supplied by Lemma~\ref{lem:gaussian}. The deterministic affinity identities and monotonicity are collected in Lemmas~\ref{lem:gaussian-affinity} and \ref{lem:affinity-processing}. We combine them below by first conditioning on the errors and integrating the Gaussian field.

\noindent\textbf{Integrating out the Gaussian field.}
Put $D=C^{-1}-C_0^{-1}$ and $s=1/2$. Write $F=S_0^{1/2}G$, where $G\sim N(0,I_n)$ is independent of $\varepsilon$, and define
\[
 A_s=I+sS_0^{1/2}DS_0^{1/2},\qquad
 h_\varepsilon=sS_0^{1/2}D\varepsilon.
\]
The matrix $A_s$ is positive definite. Indeed, because $C_0=S_0+t_0I$,
\[
 0\preceq S_0^{1/2}C_0^{-1}S_0^{1/2}
       =S_0(S_0+t_0I)^{-1}\preceq I,
\]
and therefore
\[
 A_s=I+sS_0^{1/2}C^{-1}S_0^{1/2}
          -sS_0^{1/2}C_0^{-1}S_0^{1/2}
       \succeq(1-s)I.
\]
Conditionally on $\varepsilon$, expand the quadratic form in the Gaussian integral:
\begin{align*}
 \E_F[e^{sL_C}\mid\varepsilon]
 &=\left(\frac{\det C_0}{\det C}\right)^{s/2}(2\pi)^{-n/2}
   \int_{\R^n}\exp\left\{-\frac12g^\T g
       -\frac s2(S_0^{1/2}g+\varepsilon)^\T
                    D(S_0^{1/2}g+\varepsilon)\right\}\,dg\\*
 &=\left(\frac{\det C_0}{\det C}\right)^{s/2}(2\pi)^{-n/2}
     e^{-s\varepsilon^\T D\varepsilon/2}
   \int_{\R^n}\exp\left\{-\frac12g^\T A_sg
                            -h_\varepsilon^\T g\right\}\,dg.
\end{align*}
The cross term is $-s\varepsilon^\T DS_0^{1/2}g =-h_\varepsilon^\T g$. Completing the square means using the identity
\[
 g^\T A_sg+2h_\varepsilon^\T g
 =(g+A_s^{-1}h_\varepsilon)^\T A_s
       (g+A_s^{-1}h_\varepsilon)
       -h_\varepsilon^\T A_s^{-1}h_\varepsilon.
\]
With $z=A_s^{1/2}(g+A_s^{-1}h_\varepsilon)$, the Jacobian is $dg=\det(A_s)^{-1/2}\,dz$. Hence
\begin{align*}
 (2\pi)^{-n/2}\int_{\R^n}
     e^{-g^\T A_sg/2-h_\varepsilon^\T g}\,dg&=e^{h_\varepsilon^\T A_s^{-1}h_\varepsilon/2}
       \det(A_s)^{-1/2}(2\pi)^{-n/2}
                      \int_{\R^n}e^{-\norm z^2/2}\,dz\\*
 &=\det(A_s)^{-1/2}
             e^{h_\varepsilon^\T A_s^{-1}h_\varepsilon/2}.
\end{align*}
Finally,
\[
 h_\varepsilon^\T A_s^{-1}h_\varepsilon
 =s^2\varepsilon^\T DS_0^{1/2}A_s^{-1}S_0^{1/2}D\varepsilon.
\]
Thus, on defining
\[
 B=D-sDS_0^{1/2}A_s^{-1}S_0^{1/2}D,
\]
the conditional expectation is
\[
 \begin{aligned}
 \E_F[e^{sL_C}\mid\varepsilon]
 &=\left(\frac{\det C_0}{\det C}\right)^{s/2}
   \det(A_s)^{-1/2}
   \exp\left\{-\frac s2\varepsilon^\T D\varepsilon
          +\frac{s^2}{2}\varepsilon^\T
               DS_0^{1/2}A_s^{-1}S_0^{1/2}D\varepsilon\right\}\\
 &=\left(\frac{\det C_0}{\det C}\right)^{s/2}
   \det(I+sS_0^{1/2}DS_0^{1/2})^{-1/2}
                  e^{-s\varepsilon^\T B\varepsilon/2}.
 \end{aligned}
\]

\noindent\textbf{Determining the normalization by Gaussian errors.}
Denote the deterministic prefactor just obtained by
\[
 K_s=\left(\frac{\det C_0}{\det C}\right)^{s/2}
                         \det(A_s)^{-1/2}>0.
\]
Introduce hypothetical errors $\varepsilon^{\mathrm G}\sim N(0,t_0I)$, independent of $F$. Then $Y^{\mathrm G}=F+\varepsilon^{\mathrm G}\sim N(0,C_0)$. Let $p_C$ and $p_0$ be the densities of $N(0,C)$ and $N(0,C_0)$, respectively. The contrast evaluated at a vector $y$ is $L_C(y)=\log\{p_C(y)/p_0(y)\}$. Since $s=1/2$,
\begin{align*}
 \E_{F,\varepsilon^{\mathrm G}}
            e^{sL_C(F+\varepsilon^{\mathrm G})}
 &=\int_{\R^n}p_0(y)
                    \left\{\frac{p_C(y)}{p_0(y)}\right\}^{1/2}\,dy\\*
 &=\int_{\R^n}\sqrt{p_C(y)p_0(y)}\,dy
   =\Aff\{N(0,C),N(0,C_0)\}=e^{-A_C}.
\end{align*}
Averaging the conditional identity over these Gaussian errors therefore gives
\[
 e^{-A_C}=K_sJ_s,\qquad
 J_s=\E_{\varepsilon^{\mathrm G}}
       e^{-s(\varepsilon^{\mathrm G})^\T B\varepsilon^{\mathrm G}/2}.
\]
In particular $J_s$ is finite. Diagonalize the symmetric matrix $B=U\operatorname{diag}(b_1,\ldots,b_n)U^\T$ and write $\varepsilon^{\mathrm G}=\sqrt{t_0}\,UZ$ with independent standard Gaussian coordinates $Z_j$. Then
\[
 J_s=\prod_{j=1}^n(2\pi)^{-1/2}
       \int_{\R}\exp\{-\tfrac12(1+st_0b_j)z^2\}\,dz.
\]
Finiteness forces every $1+st_0b_j$ to be positive: otherwise the corresponding one-dimensional integral is infinite. Thus $I+st_0B\succ0$, and evaluating the Gaussian integrals gives
\[
 J_s=\prod_{j=1}^n(1+st_0b_j)^{-1/2}
      =\det(I+st_0B)^{-1/2},\qquad
 K_s=e^{-A_C}\det(I+st_0B)^{1/2}.
\]
Substituting this expression for $K_s$ in the conditional identity gives an exact formula before any inequality is applied. To bound its determinant, use $\log(1+x)\le x$ for $x>-1$:
\[
 \log\det(I+st_0B)=\sum_j\log(1+st_0b_j)
                  \le st_0\sum_jb_j=st_0\tr B.
\]
Consequently, with $s=1/2$,
\begin{equation}\label{eq:conditional-affinity}
 \begin{aligned}
 \E_F[e^{L_C/2}\mid\varepsilon]
 &=\exp\left\{-A_C+\frac12\log\det(I+\tfrac{t_0}{2}B)
                             -\frac14\varepsilon^\T B\varepsilon\right\}\\
 &\le\exp\left\{-A_C+\frac{t_0}{4}\tr B
                             -\frac14\varepsilon^\T B\varepsilon\right\}
   =\exp\{-A_C-\tfrac14Q_\varepsilon(B)\},
 \end{aligned}
\end{equation}
where $Q_\varepsilon(B)=\varepsilon^\T B\varepsilon-t_0\tr B$. This is centered under the actual error law, since
\[
 \E[\varepsilon^\T B\varepsilon]
 =\tr\{B\,\E[\varepsilon\varepsilon^\T]\}=t_0\tr B.
\]

\noindent\textbf{Controlling the error quadratic form.}
Only two norms of $B$ are needed. Set $W=I-sS_0C_0^{-1}$ and $H=sW^{-1}S_0$. Resolvent algebra gives
\[
 B=W^{-1}\{(C+H)^{-1}-(C_0+H)^{-1}\}W^{-1}.
\]
For verification, $I+sS_0D=W(I+HC^{-1})$, $C_0+H=W^{-1}C_0$, and multiplication gives this formula from $B=D(I+sS_0D)^{-1}$. To bound its operator norm, note that $S_0$ commutes with $C_0=S_0+t_0I$ and
\[
 0\preceq S_0C_0^{-1}\preceq I,
 \qquad (1-s)I\preceq W\preceq I.
\]
Consequently $\|W^{-1}\|_\op\le(1-s)^{-1}=2$. The matrices $W^{-1}$ and $S_0$ also commute, so $H=sW^{-1}S_0\succeq0$. Since $C,C_0\succeq t_-I$,
\[
 C+H\succeq t_-I,\qquad C_0+H\succeq t_-I,
\]
and each inverse has operator norm at most $t_-^{-1}$. The displayed representation of $B$ therefore gives
\[
 \begin{aligned}
 \|B\|_\op
 &\le\|W^{-1}\|_\op^2
       \bigl\{\|(C+H)^{-1}\|_\op
                    +\|(C_0+H)^{-1}\|_\op\bigr\}\\
 &\le 4(t_-^{-1}+t_-^{-1})=8t_-^{-1}.
 \end{aligned}
\]
This bound is uniform in $n$ and the candidate covariance $C$.

Apply the precision bound in Lemma~\ref{lem:gaussian-affinity} to $C+H$ and $C_0+H$. They are both at least $t_-I$. Lemma~\ref{lem:affinity-processing} shows that adding the common Gaussian noise with covariance $H$ increases affinity, so
\[
\begin{aligned}
 \|B\|_\Frob^2
 &\le\|W^{-1}\|_\op^4
        \|(C+H)^{-1}-(C_0+H)^{-1}\|_\Frob^2\\
 &\le C\bigl[-\log\Aff\{N(0,C+H),N(0,C_0+H)\}\bigr]
 \le CA_C.
\end{aligned}
\]
\noindent\textbf{A fractional moment and Markov's inequality.}
The matrices above retain $s=1/2$. For a second power $0<q\le1/2$, conditional Jensen and \eqref{eq:conditional-affinity} give
\[
 \begin{aligned}
 \E_F[e^{qL_C}\mid\varepsilon]
 &=\E_F[(e^{L_C/2})^{2q}\mid\varepsilon]\\
 &\le\bigl(\E_F[e^{L_C/2}\mid\varepsilon]\bigr)^{2q}
 \le e^{-2qA_C-(q/2)Q_\varepsilon(B)}.
 \end{aligned}
\]
By Lemma~\ref{lem:gaussian}, the centered error quadratic form satisfies
\[
 \log \E_\varepsilon e^{uQ_\varepsilon(B)}
       \le Cu^2\|B\|_\Frob^2,
 \qquad |u|\|B\|_\op\le c.
\]
Since $\|B\|_\op\le8t_-^{-1}$, this applies with $u=-q/2$ for every sufficiently small fixed $q>0$, uniformly in $n$ and $C$. Averaging the conditional bound over the errors and using $\|B\|_\Frob^2\le CA_C$, we obtain
\[
 \begin{aligned}
 \E e^{qL_C}
 &\le e^{-2qA_C}\E_\varepsilon e^{-(q/2)Q_\varepsilon(B)}\\
 &\le\exp\{-2qA_C+Cq^2\|B\|_\Frob^2\}
 \le\exp\{(-2q+C_1q^2)A_C\}.
 \end{aligned}
\]
Choose one such $q$ with $q\le1/2$ and $C_1q\le1$. Then $\E e^{qL_C}\le e^{-qA_C}$, and Markov's inequality gives
\[
 P(L_C\ge-z)\le e^{qz}\E e^{qL_C}
             \le e^{qz-qA_C}\le e^{z/2-qA_C},\qquad z\ge0.
\]
This proves \eqref{eq:likelihoodaffinity}, with $c=q$ and $C_*=1$.
\qed

\subsection{Proof of Lemma~\ref{lem:noise-diagonals}}
\noindent\textbf{Signal diagonals.}
Write $F_n=F(x_{n,1:n})$, $C_0=S_0+t_0I$, and $P_n=S_0C_0^{-1}$. The best linear predictor of $F_n$ from $Y_n$ is
\[
 \widehat F_n=S_0C_0^{-1}Y_n.
\]
Indeed, $\Cov(F_n,Y_n)=S_0$ and $\Cov(Y_n)=C_0$, so the residual $F_n-\widehat F_n$ is uncorrelated with every linear function of $Y_n$. Its covariance is
\[
 \begin{aligned}
 \Cov(F_n-\widehat F_n)
 &=S_0-S_0C_0^{-1}S_0\\
 &=S_0C_0^{-1}(C_0-S_0)=t_0P_n.
 \end{aligned}
\]
Thus the minimum mean squared error for linear prediction at site $i$ is $t_0(P_n)_{ii}$. These identities use only second moments.

Fix $0<\zeta<\min\{\nu_-,1\}$. The uniform increment bound \eqref{eq:variogram} gives
$\E|F(x)-F(y)|^2\le C\|x-y\|^{2\zeta}$. Partition an enclosing box into $K_\rho<\infty$ cells of diameter at most $\rho$, assigning boundary points by a fixed rule. For a cell $Q$, let $I_Q=\{j:x_{n,j}\in Q\}$ and $m_Q=\#I_Q$. If $i\in I_Q$, consider the linear predictor $\widetilde F_i=m_Q^{-1}\sum_{j\in I_Q}Y_{n,j}$. Independence of field and errors, followed by Cauchy--Schwarz for the finite average, gives
\[
 \begin{aligned}
 \E|F(x_{n,i})-\widetilde F_i|^2
 &=\E\left|\frac1{m_Q}\sum_{j\in I_Q}
                    \{F(x_{n,i})-F(x_{n,j})\}\right|^2
                       +\frac{t_0}{m_Q}\\
 &\le\frac1{m_Q}\sum_{j\in I_Q}
              \E|F(x_{n,i})-F(x_{n,j})|^2+\frac{t_0}{m_Q}\\
 &\le C\rho^{2\zeta}+\frac{t_+}{m_Q}.
 \end{aligned}
\]
Optimality of $\widehat F_n$ therefore implies
\[
 t_0(P_n)_{ii}\le C\rho^{2\zeta}+t_+/m_Q.
\]
For any fixed $\delta>0$, choose $\rho$ so that $C\rho^{2\zeta}\le t_-\delta/2$, and choose an integer $m_\delta\ge2t_+/(t_-\delta)$. If $m_Q\ge m_\delta$, then
$t_0(P_n)_{ii}\le t_-\delta$, hence $(P_n)_{ii}\le\delta$. The cells with $m_Q<m_\delta$ contain at most
\[
 N_\delta\coloneqq K_\rho(m_\delta-1)
\]
observation indices in total. Consequently all but at most $N_\delta$ sites have $(P_n)_{ii}\le\delta$. The choices of $\rho,m_\delta,N_\delta$ are independent of $n$ and of the true parameter in the compact set. Repeated sites are counted by their observation indices; their errors are still independent, so the same argument applies.

For a real signal multiplier $g$, put
\[
 D_g=\int g(\omega)f_0(\omega)a_\omega a_\omega^*\,d\omega,
 \qquad
 A_g=C_0^{-1}D_gC_0^{-1}
     =C_0^{-1/2}\mathcal B_n(g)C_0^{-1/2}.
\]
With $c_i=C_0^{-1}e_i$, $(A_g)_{ii}=\int gf_0|a_\omega^*c_i|^2$. Cauchy--Schwarz in this integral gives
\[
 |(A_g)_{ii}|^2\le q_ir_i,\qquad
 q_i=(C_0^{-1}S_0C_0^{-1})_{ii},\qquad
 r_i=(C_0^{-1}D_{g^2}C_0^{-1})_{ii}\ge0.
\]
Since $C_0^{-1}=t_0^{-1}(I-P_n)$ and $0\preceq P_n\preceq I$,
\[
 q_i=t_0^{-1}(P_n-P_n^2)_{ii}\le t_0^{-1}(P_n)_{ii}.
\]
Also $\mathcal B_n(g^2)\succeq0$ and $C_0^{-1}\preceq t_-^{-1}I$, so
\[
 \begin{aligned}
 \sum_i r_i
 &=\tr\{C_0^{-1/2}\mathcal B_n(g^2)C_0^{-1/2}\}\\
 &=\tr\{C_0^{-1}\mathcal B_n(g^2)\}
 \le t_-^{-1}\tr\mathcal B_n(g^2)
 =t_-^{-1}\int g^2w_n\le C.
 \end{aligned}
\]
For the normalized multipliers in \eqref{eq:signaltangent}, the last bound follows from Lemma~\ref{lem:upper}: use \eqref{eq:effective} for $g=1/r_n$, \eqref{eq:logmoments} with $j=1$ for $g=g_n/r_n$, and \eqref{eq:rangetrace} with $k=2$ for $g=-\beta_0/(b_nm_0)$ when $d\ge4$, together with \eqref{eq:rates}--\eqref{eq:arate}. Furthermore, Lemma~\ref{lem:derivatives}(v), specifically
\eqref{eq:studentderiv} at the truth, gives
\[
 \|A_g\|_\op
 \le t_-^{-1}\|\mathcal B_n(g)\|_\op=o(1).
\]
The corresponding Frobenius estimate from \eqref{eq:logbounds} and \eqref{eq:rangeupper} is only $\|A_g\|_\Frob\le t_-^{-1}\|\mathcal B_n(g)\|_\Frob=O(1)$. We instead split the diagonal entries according to $(P_n)_{ii}$:
\[
 \begin{aligned}
 \sum_i|(A_g)_{ii}|^2
 &\le\sum_{i:(P_n)_{ii}\le\delta}q_ir_i
        +\sum_{i:(P_n)_{ii}>\delta}|(A_g)_{ii}|^2\\
 &\le\frac{\delta}{t_0}\sum_i r_i
        +N_\delta\|A_g\|_\op^2
 \le C\delta+N_\delta\|A_g\|_\op^2.
 \end{aligned}
\]
First let $n\to\infty$ with $\delta$ fixed, and then let $\delta\downarrow0$. This proves \eqref{eq:diagonal-negligible} uniformly.

\noindent\textbf{The nugget diagonal and covariance correction.}
Since $C_0^{-1}=t_0^{-1}(I-P_n)$ and $\tr P_n=O(R_n^d)=o(n)$,
\[
 \norm{C_0^{-2}-t_0^{-2}I}_\Frob^2\le C\tr P_n=o(n).
\]
This proves the nugget diagonal limit. Cauchy--Schwarz makes all signal--nugget cumulant cross terms vanish. Equation \eqref{eq:quadratic-cumulant} proves the first part of \eqref{eq:score-covariance-correction}. The Gaussian curvature has $J_{tt}\to(2t_0^2)^{-1}$ and signal--nugget cross entries tending to zero by Lemma~\ref{lem:information}. Since $J_n\asymp I$, its inverse satisfies
\[
 J_n^{-1}
 =\begin{pmatrix}J_{n,\mathrm{sig}}^{-1}&0\\0&2t_0^2\end{pmatrix}+o(1),
 \qquad J_n^{-1}e_t=2t_0^2e_t+o(1).
\]
Multiplication of the score covariance formula on both sides by $J_n^{-1}$ therefore gives
\[
 \begin{aligned}
 V_n
 &=J_n^{-1}
    +\frac{\kappa_4}{4t_0^4}(J_n^{-1}e_t)(J_n^{-1}e_t)^\T+o(1)\\
 &=J_n^{-1}+\kappa_4e_te_t^\T+o(1)\\
 &=\begin{pmatrix}
       J_{n,\mathrm{sig}}^{-1}&0\\0&t_0^2(\E_P\xi^4-1)
    \end{pmatrix}+o(1).
 \end{aligned}
\]
At $d=4$, Lemma~\ref{lem:information} additionally gives
\[
 J_{n,\mathrm{sig}}^{-1}
 =\begin{pmatrix}J_{n,(u,v)}^{-1}&0\\0&J_{n,aa}^{-1}\end{pmatrix}+o(1).
\]
Inversion is continuous uniformly here because all the principal blocks have eigenvalues bounded above and away from zero. These calculations prove the covariance correction and both asserted block structures, without requiring a limit for the design-dependent signal information.

\noindent\textbf{The nugget score as a sum of noise squares.}
Write $Y=F_n+\varepsilon_n$ and $D=C_0^{-2}-t_0^{-2}I$. The exact score formula gives
\[
 \begin{aligned}
 2\sqrt n\,\Delta_{n,t}
 &=Y^\T C_0^{-2}Y-\tr C_0^{-1}\\
 &=t_0^{-2}\sum_i(\varepsilon_{n,i}^2-t_0)
     +\{\varepsilon_n^\T D\varepsilon_n-t_0\tr D\}+\{F_n^\T C_0^{-2}F_n-\tr(S_0C_0^{-2})\}
     +2F_n^\T C_0^{-2}\varepsilon_n.
 \end{aligned}
\]
The expectations cancel because $\tr(S_0C_0^{-2})+t_0\tr C_0^{-2}=\tr C_0^{-1}$.
The centered noise remainder has $L^p$ norm at most $C_p\|D\|_\Frob=o(\sqrt n)$ for every fixed finite $p$, by Lemma~\ref{lem:gaussian} and the bound on $D$ above.

For the field and cross terms write $F_n=S_0^{1/2}G$ and $\varepsilon_n=\sqrt{t_0}\xi$, with $G$ and $\xi$ independent. If $\lambda_i$ are
the eigenvalues of $S_0$, commutation of $S_0$ and $C_0$ gives
\[
 \begin{aligned}
 \|S_0^{1/2}C_0^{-2}S_0^{1/2}\|_\Frob^2
 &=\sum_i\frac{\lambda_i^2}{(\lambda_i+t_0)^4}
   \le C\sum_i\frac{\lambda_i}{\lambda_i+t_0}
   =C\tr P_n,\\
 t_0\|S_0^{1/2}C_0^{-2}\|_\Frob^2
 &=t_0\sum_i\frac{\lambda_i}{(\lambda_i+t_0)^4}
   \le C\tr P_n.
 \end{aligned}
\]
The first bound controls the centered Gaussian quadratic form. For the cross term use the symmetric off-diagonal block matrix on the independent vector $(G^\T,\xi^\T)^\T$:
\[
 2F_n^\T C_0^{-2}\varepsilon_n
 =\begin{pmatrix}G\\\xi\end{pmatrix}^{\!\T}
   \begin{pmatrix}
    0&\sqrt{t_0}S_0^{1/2}C_0^{-2}\\
    \sqrt{t_0}C_0^{-2}S_0^{1/2}&0
   \end{pmatrix}
   \begin{pmatrix}G\\\xi\end{pmatrix}.
\]
Its squared Frobenius norm is bounded by $C\tr P_n$ by the second estimate. The quadratic-form moment bound thus gives $L^p$ norms at most $C_p\sqrt{\tr P_n}=o(\sqrt n)$ for both terms. Division by $2\sqrt n$ proves \eqref{eq:nugget-score-squares}
with an $o_{L^p}(1)$ remainder, uniformly in the noise class.

The distributional assertion follows directly from the quadratic-form central limit theorem in Lemma~\ref{lem:gaussian}, applied to $B_t=C_0^{-1}/\sqrt n$: we have
$\|B_t\|_\Frob\le t_-^{-1}$ and $\|B_t\|_\op\le t_-^{-1}/\sqrt n$, and the covariance calculation above gives
\[
 \Var(\Delta_{n,t})=\frac{v_\varepsilon}{4t_0^4}+o(1),
 \qquad v_\varepsilon=t_0^2(\E_P\xi^4-1)\asymp1.
\]
Thus $2t_0^2\Delta_{n,t}/\sqrt{v_\varepsilon}\Rightarrow N(0,1)$ uniformly. \qed

\subsection{Proof of Lemma~\ref{lem:kurtosis}}
We first hold the covariance parameters at the truth. The correction in \eqref{eq:residual-cumulant} then gives an exactly unbiased estimator. We bound its variance and subsequently control the effect of substituting the global fit. Steps 5--7 establish the additional Gaussian variance and normal limit. All constants
below are uniform over $\mathcal K\times\mathcal P$ and the allowed designs.

\noindent\textbf{1. The correction at the true covariance.}
Suppress the subscript $n$ and put
\[
 C_0=S_0+t_0I_n,\qquad A_0=t_0C_0^{-1}=I_n-P_n,\qquad
 r_0=A_0Y_n,\qquad v_{0,i}=t_0(A_0)_{ii}.
\]
Since $S_0C_0^{-1}Y_n$ is the best linear predictor of the field at the sites, $r_0$ is its residual vector. Its covariance is
\[
 \Cov_0(r_0)=A_0C_0A_0=t_0A_0,
 \qquad \E_0r_{0,i}^2=v_{0,i}.
\]
Write $Y_n=F_n+\varepsilon_n$. The Gaussian vector $A_0F_n$ has zero fourth cumulant and is independent of $A_0\varepsilon_n$. Independence of the errors therefore gives
\[
 \E_0r_{0,i}^4-3v_{0,i}^2
 =\operatorname{cum}_4(r_{0,i})
 =\kappa_4(P,t_0)\sum_j(A_0)_{ij}^4.
\]
For a vector $r$, a vector $v$ of marginal variances, and a matrix $A$, write
\[
 T(r,v)=\sum_i(r_i^4-6v_ir_i^2+3v_i^2),\qquad
 \mathfrak d(A)=\sum_{i,j}A_{ij}^4.
\]
It follows that
\begin{equation}\label{eq:oracle-cumulant-mean}
 \E_0T(r_0,v_0)=\kappa_4(P,t_0)\mathfrak d(A_0).
\end{equation}
The numerator uses the fourth Hermite polynomial:
\[
 r_i^4-6v_ir_i^2+3v_i^2
 =v_i^2\operatorname{He}_4(r_i/\sqrt{v_i}),\qquad
 \operatorname{He}_4(x)=x^4-6x^2+3.
\]

The denominator has order $n$. Indeed, for any symmetric matrix $0\preceq A\preceq I_n$, all its entries have absolute value at most one, and
\[
 \mathfrak d(A)\le\sum_{i,j}A_{ij}^2=\tr A^2\le n,
 \qquad
 \mathfrak d(A)\ge\sum_i A_{ii}^4
          \ge n-4\tr(I_n-A).
\]
The last inequality uses $x^4\ge1-4(1-x)$ for $0\le x\le1$. By Lemma~\ref{lem:upper},
\begin{equation}\label{eq:oracle-cumulant-denominator}
 n-4CM_n\le\mathfrak d(A_0)\le n,\qquad
 M_n=n^{d/(2\nu_0+d)}=o(n).
\end{equation}

\noindent\textbf{2. A variance bound for the numerator.}
Represent the residual as $r_0=HZ$, where
\[
 Z=\begin{pmatrix}G\\\xi\end{pmatrix},\qquad
 H=A_0\begin{pmatrix}S_0^{1/2}&\sqrt{t_0}I_n\end{pmatrix}.
\]
The $2n$ coordinates of $Z$ are independent and centered, have variance one, and obey a common sub-Gaussian bound. Also
\[
 HH^\T=t_0A_0,\qquad \|H\|_\op\le C,\qquad
 \|H\|_\Frob^2\le Cn.
\]
In particular, all fixed moments of each residual coordinate are uniformly bounded. Apply \citet[Theorem~1.4]{adamczak2015concentration} to the degree-four polynomial $\Phi(z)=T(Hz,v_0)$. Every partition norm in that theorem is at most the Hilbert--Schmidt norm, so its $p=2$ case, followed by squaring and Cauchy--Schwarz over the four derivative orders, gives
\[
 \Var_0\{\Phi(Z)\}
 \le C\sum_{k=1}^4
       \|\E_0D^k\Phi(Z)\|_{\mathrm{HS}}^2.
\]
Here the Hilbert--Schmidt norm of a tensor is the square root of the sum of its squared entries. The constant depends only on the degree and the common sub-Gaussian bound.

Let $h_i^\T$ denote row $i$ of $H$. Differentiating $\Phi$ gives
\[
 \begin{aligned}
 \E_0D\Phi(Z)&=4H^\T\E_0(r_0^{\circ3}),\\
 \E_0D^2\Phi(Z)
 &=12\sum_i(\E_0r_{0,i}^2-v_{0,i})h_i^{\otimes2}=0,\\
 \E_0D^3\Phi(Z)&=24\sum_i\E_0r_{0,i}\,h_i^{\otimes3}=0,\\
 \E_0D^4\Phi(Z)&=24\sum_i h_i^{\otimes4},
 \end{aligned}
\]
where $r_0^{\circ3}$ denotes coordinatewise cubing. Bounded sixth moments give
\[
 \|\E_0D\Phi(Z)\|^2
 \le16\|H\|_\op^2\sum_i\E_0|r_{0,i}|^6\le Cn.
\]
For the fourth derivative, tensor inner products and $\max_{i,j}|(HH^\T)_{ij}|\le\|HH^\T\|_\op\le C$ give
\[
 \begin{aligned}
 \|\E_0D^4\Phi(Z)\|_{\mathrm{HS}}^2
 &=24^2\sum_{i,j}\langle h_i,h_j\rangle^4
 \le C\|HH^\T\|_\Frob^2\\
 &\le C\|H\|_\op^2\|H\|_\Frob^2\le Cn.
 \end{aligned}
\]
The cited theorem therefore yields
\[
 \Var_0\{T(r_0,v_0)\}\le Cn.
\]
Together with \eqref{eq:oracle-cumulant-mean}--\eqref{eq:oracle-cumulant-denominator},
this gives the known-covariance bound
\begin{equation}\label{eq:oracle-cumulant-rate}
 \left\|\frac{T(r_0,v_0)}{\mathfrak d(A_0)}-\kappa_4(P,t_0)
       \right\|_{L^2}\le Cn^{-1/2}.
\end{equation}

\noindent\textbf{3. Substituting the fitted covariance.}
Lemma~\ref{lem:tails} gives $\widehat h_n=O_P(1)$ and $\widehat t_n-t_0=O_P(n^{-1/2})$, uniformly. We show that
\begin{equation}\label{eq:residual-plugin-bounds}
 \|\widehat A_n-A_0\|_\Frob=O_P(1),\qquad
 \|\widehat r_n-r_0\|_2=O_P(1),\qquad
 \|\widehat v_n-v_0\|_2=O_P(1).
\end{equation}
Fix a local radius $H_*>0$. Use $u=h$ in high dimension and $u=(h,a)$ in low dimension, on the compact convex set $\Gamma_{H_*}$ of Lemma~\ref{lem:derivatives}. Put
\[
 A(u)=t(u)C_n(u)^{-1},\qquad R(u)=M(u)^{-1},\qquad
 D(u)=\{A(u)-A_0\}C_0^{1/2}
     =C_0^{-1/2}\{t(u)R(u)-t_0I_n\}.
\]
Lemma~\ref{lem:derivatives} gives, on this set,
\[
 \|R(u)\|_\op\le C,\qquad
 \|M(u)-I_n\|_\Frob+
       \max_\ell\|\partial_{u_\ell}M(u)\|_\Frob\le C_{H_*}.
\]
The derivatives here include the unnormalized nuisance $a$ when $d\le3$. Also $\sqrt n\,|t(u)-t_0|\le H_*$; the only nonzero derivative of $t$ is $\partial_{h_t}t=n^{-1/2}$. Using $R-I_n=-R(M-I_n)$ and
\[
 \partial_{u_\ell}D
 =C_0^{-1/2}\{(\partial_{u_\ell}t)R
                 -tR(\partial_{u_\ell}M)R\},
\]
we obtain
\begin{equation}\label{eq:residual-plugin-matrix}
 \sup_{u\in\Gamma_{H_*}}
 \left\{\|D(u)\|_\Frob+
               \max_\ell\|\partial_{u_\ell}D(u)\|_\Frob\right\}
 \le C_{H_*}.
\end{equation}
For example, the term containing $\partial_{h_t}t$ has Frobenius norm at most $Cn^{-1/2}\|R\|_\Frob\le C$.

To evaluate these bounds at the random fit, write $X=C_0^{-1/2}Y_n$. For any deterministic matrix $B$ and $p\ge2$, the sub-Gaussian moment bound for a linear combination gives
\[
 \left\|\|BX\|_2\right\|_{L^p}^2
 \le\sum_i\|(BX)_i\|_{L^p}^2
 \le C_p\sum_i\|B_{i\cdot}\|_2^2
 =C_p\|B\|_\Frob^2.
\]
Indeed, $X=C_0^{-1/2}(S_0^{1/2},\sqrt{t_0}I_n)Z$, and the matrix multiplying $Z$ has orthonormal rows. Apply this estimate and \eqref{eq:residual-plugin-matrix} to the continuous process $u\mapsto\|D(u)X\|_2$. Its pointwise $L^p$ norms and its $L^p$ Lipschitz constant are bounded by $C_{H_*,p}$. Lemma~\ref{lem:compact-supremum}, with $p$ larger than the fixed index dimension and as large as needed for the moments below, yields
\[
 \sup_{u\in\Gamma_{H_*}}\|D(u)X\|_2=O_P(1).
\]
On $\{\|\widehat h_n\|\le H_*\}$ this bounds
$\|\widehat r_n-r_0\|_2$; tightness then lets $H_*$ increase. The first assertion in \eqref{eq:residual-plugin-bounds} follows in the same way from $A(u)-A_0=D(u)C_0^{-1/2}$. Finally, $0\preceq\widehat A_n\preceq I_n$ and
\[
 \|\widehat v_n-v_0\|_2
 \le\sqrt n\,|\widehat t_n-t_0|
          +t_0\|\widehat A_n-A_0\|_\Frob=O_P(1),
\]
which proves the last assertion.

\noindent\textbf{4. The resulting rate.}
Let $e=\widehat r_n-r_0$. The bounded sixth moments of $r_{0,i}$ and \eqref{eq:residual-plugin-bounds} imply
\[
 \begin{aligned}
 \sum_i|\widehat r_{n,i}^4-r_{0,i}^4|
 &\le C\|e\|_2
        \left\{\sum_i(|r_{0,i}|+|e_i|)^6\right\}^{1/2}
       =O_P(\sqrt n),\\
 \sum_i|\widehat r_{n,i}^2-r_{0,i}^2|
 &\le\|e\|_2(2\|r_0\|_2+\|e\|_2)=O_P(\sqrt n).
 \end{aligned}
\]
Here $\sum_i|e_i|^6\le\|e\|_2^6=O_P(1)$ and $\sum_i|r_{0,i}|^6=O_P(n)$. Since $v_{0,i}$ and $\widehat v_{n,i}$ are bounded by a fixed constant, we also have
\[
 \begin{aligned}
 \left|\sum_i\widehat v_{n,i}\widehat r_{n,i}^2
                -\sum_iv_{0,i}r_{0,i}^2\right|
 &\le C\sum_i|\widehat r_{n,i}^2-r_{0,i}^2|
      +\|\widehat v_n-v_0\|_2\left(\sum_i r_{0,i}^4\right)^{1/2}
       =O_P(\sqrt n),\\
 \left|\sum_i\widehat v_{n,i}^2-\sum_iv_{0,i}^2\right|
 &\le\|\widehat v_n-v_0\|_2
            (\|\widehat v_n\|_2+\|v_0\|_2)=O_P(\sqrt n).
 \end{aligned}
\]
Thus $T(\widehat r_n,\widehat v_n)-T(r_0,v_0)=O_P(\sqrt n)$. For the denominator, every matrix on the segment between $A_0$ and $\widehat A_n$ is a positive contraction. Differentiating $\mathfrak d$ along this segment and applying Cauchy--Schwarz gives
\[
 |\mathfrak d(\widehat A_n)-\mathfrak d(A_0)|
 \le4\sqrt n\,\|\widehat A_n-A_0\|_\Frob=O_P(\sqrt n),
\]
because $\sum_{i,j}|A_{ij}|^6\le\|A\|_\Frob^2\le n$ there. In particular, $\mathfrak d(\widehat A_n)/n\to_P1$ uniformly. Combining these bounds with \eqref{eq:oracle-cumulant-rate} yields
\[
 \begin{split}
 T(\widehat r_n,\widehat v_n)
       -\kappa_4\mathfrak d(\widehat A_n)
 &=T(r_0,v_0)-\kappa_4\mathfrak d(A_0)+\{T(\widehat r_n,\widehat v_n)-T(r_0,v_0)\}
       -\kappa_4\{\mathfrak d(\widehat A_n)-\mathfrak d(A_0)\}\\
 &=O_P(\sqrt n).
 \end{split}
\]
Division by $\mathfrak d(\widehat A_n)$ proves \eqref{eq:kurtosis-rate}.

For random sites, every bound above holds conditionally on the good designs $\mathcal G_n$ of Lemma~\ref{lem:randomgood}, with the same constants. More explicitly, for each fixed $H_*$ the moment bounds and the tail bound for $\widehat h_n$ give, for $z\ge1$,
\[
 \sup_{\theta_0,P}
 P_{\theta_0,P}\{\sqrt n\,|\widetilde\kappa_{4,n}-\kappa_4|>z
                         \mid x_{1:n}\}
 \le Ce^{-cH_*^2}+C_{H_*}z^{-2}+o(1)
 \quad\text{on }\mathcal G_n.
\]
The $o(1)$ includes the vanishing probability that the denominator is smaller than $n/2$. Averaging adds at most $P_X(\mathcal G_n^c)=o(1)$. First choose $H_*$ large and then $z$ large to obtain the unconditional uniform $O_P(n^{-1/2})$ assertion. With smoothness fixed, the same covariance bounds and the tightness in Lemma~\ref{lem:tails} apply on the restricted local sets, so the proof is unchanged.

\noindent\textbf{5. Exact variance under Gaussian errors.}
For the remainder of the proof assume Gaussian errors, so that $\kappa_4=0$, and put $d_n=\mathfrak d(A_0)$ and $K=t_0A_0$. For centered jointly Gaussian variables $U,V$ with variances $v,w$ and covariance $k$, their moment-generating function gives
\[
 \E\{e^{sU-vs^2/2}e^{zV-wz^2/2}\}=e^{ksz}.
\]
Comparing coefficients of $s^4z^4$ gives
\[
 \E\{(U^4-6vU^2+3v^2)(V^4-6wV^2+3w^2)\}=24k^4.
\]
Applying this identity to each pair of residuals yields
\[
 \Var_0\{T(r_0,v_0)\}
 =24\sum_{i,j}K_{ij}^4=24t_0^4d_n.
\]
Division by $d_n^2$ proves \eqref{eq:cumulant-gaussian-oracle-variance}.

\noindent\textbf{6. The oracle Gaussian limit.}
Put
\[
 W_n=\sum_i(\varepsilon_i^4-6t_0\varepsilon_i^2+3t_0^2).
\]
Under the Gaussian law,
\[
 \varepsilon_n\mid Y_n
 \sim N(r_0,t_0I_n-t_0A_0).
\]
In particular, its $i$th conditional variance is $t_0-v_{0,i}$. Using the second and fourth moments of this conditional normal distribution gives
\[
 \begin{aligned}
 \E_0(\varepsilon_i^4-6t_0\varepsilon_i^2+3t_0^2\mid Y_n)
 &=r_{0,i}^4+6r_{0,i}^2(t_0-v_{0,i})
       +3(t_0-v_{0,i})^2
       -6t_0\{r_{0,i}^2+t_0-v_{0,i}\}+3t_0^2\\
 &=r_{0,i}^4-6v_{0,i}r_{0,i}^2+3v_{0,i}^2.
 \end{aligned}
\]
Consequently
\begin{equation}\label{eq:cumulant-conditional-expectation}
 \E_0(W_n\mid Y_n)=T(r_0,v_0).
\end{equation}
Independence of the Gaussian errors gives $\Var_0(W_n)=24t_0^4n$. By \eqref{eq:cumulant-conditional-expectation}, conditional expectation is an orthogonal projection, so
\[
 \E_0|W_n-T(r_0,v_0)|^2
 =\Var_0(W_n)-\Var_0\{T(r_0,v_0)\}
 =24t_0^4(n-d_n)=o(n).
\]
The Lyapunov central limit theorem for the iid summands of $W_n$ therefore gives
\[
 \frac{T(r_0,v_0)}{\sqrt n}
 =\frac{W_n}{\sqrt n}+o_{L^2}(1)
 \Rightarrow N(0,24t_0^4).
\]
The convergence is uniform over the compact parameter family: after division by $t_0^2$, the summand distribution is fixed, and $(n-d_n)/n\to0$ uniformly.

\noindent\textbf{7. Fitting the covariance does not change the Gaussian limit.}
We sharpen the bound in Step 4 to
\begin{equation}\label{eq:cumulant-gaussian-plugin}
 T(\widehat r_n,\widehat v_n)-T(r_0,v_0)=O_P(1)
 \quad\text{under Gaussian errors}.
\end{equation}
Use the compact local sets $\Gamma_{H_*}$ from Step 3, and now write $X=C_0^{-1/2}Y_n\sim N(0,I_n)$ and
\[
 L(u)=A(u)C_0^{1/2},\qquad L_0=A_0C_0^{1/2},\qquad
 w(u)=\diag\{L(u)L(u)^\T\},\qquad
 v(u)=t(u)\diag A(u).
\]
Here $w(u)$ is the true marginal variance vector of $L(u)X$, whereas $v(u)$ is the candidate variance vector. Equation~\eqref{eq:residual-plugin-matrix} and the same bounds for $v(u)$ give, uniformly on $\Gamma_{H_*}$,
\[
 \|L(u)\|_\op+\|L(u)-L_0\|_\Frob+\|v(u)-v_0\|_2
 +\max_\ell\{\|\partial_{u_\ell}L(u)\|_\Frob
                    +\|\partial_{u_\ell}v(u)\|_2\}\le C_{H_*}.
\]
The entries of $v(u)$ are uniformly bounded. Also
\[
 \|w(u)-v(u)\|_2
 \le(\|L(u)\|_\op+\|L_0\|_\op)\|L(u)-L_0\|_\Frob
                  +\|v(u)-v_0\|_2\le C_{H_*}.
\]
These vectors are Lipschitz in $u$ with bounded constants. For a deterministic candidate, Gaussian fourth moments give the exact cancellation
\[
 \E_0 T(L(u)X,v(u))=3\|w(u)-v(u)\|_2^2.
\]
Thus this expectation and its Lipschitz constant are bounded on $\Gamma_{H_*}$, rather than of order $\sqrt n$.

The centered fluctuations have the same bounded orders. If $L,L'$ have bounded operator norms and $v,v'$ have bounded entries, the Gaussian Poincar\'e inequality gives
\begin{equation}\label{eq:cumulant-gaussian-polynomial}
 \Var\{T(LX,v)-T(L'X,v')\}
 \le C\{\|L-L'\|_\Frob^2+\|v-v'\|_2^2\}.
\end{equation}
To check this bound, put $r=LX$, $r'=L'X$, and $b=r^{\circ3}-3v\circ r$, where $\circ$ denotes entrywise multiplication. The gradient of the difference is
\[
 4(L-L')^\T b
 +4L'^\T\{r^{\circ3}-r'^{\circ3}
             -3v\circ(r-r')-3(v-v')\circ r'\}.
\]
Write $K=LL^\T$ and $w=\diag K$, Gaussian moments give
\[
 \Cov(r^{\circ3})
 =9\diag(w)K\diag(w)+6K^{\circ3},
 \qquad \|\Cov(b)\|_\op\le C.
\]
For example, $K^{\circ3}$ is a compression of $K^{\otimes3}$ to the orthonormal vectors $e_i\otimes e_i\otimes e_i$, so its operator norm is at most $\|K\|_\op^3$.
It follows that the expected squared norm of the first gradient term is at most $C\|L-L'\|_\Frob^2$. For the other term, Gaussian moments and Cauchy--Schwarz give
\[
 \E|r_i^3-r_i'^3|^2
 \le\{\E|r_i-r_i'|^4\}^{1/2}
      \{\E|r_i^2+r_ir_i'+r_i'^2|^4\}^{1/2}
 \le C\E|r_i-r_i'|^2.
\]
Summing this bound and bounding the two remaining entrywise products proves \eqref{eq:cumulant-gaussian-polynomial}. The Gaussian Poincar\'e inequality used here is the $p=2$ case of the Gaussian Sobolev bound discussed in \citet[equation~(3)]{adamczak2015concentration}.

For a centered Gaussian polynomial of degree at most four, its $L^p$ norm, for each fixed $p\ge2$, is at most $C_p$ times its $L^2$ norm; see \citet[Theorem~1.3]{adamczak2015concentration}. Apply this fact to \eqref{eq:cumulant-gaussian-polynomial}. The process
\[
 u\longmapsto T(L(u)X,v(u))-T(L_0X,v_0)
\]
has bounded pointwise $L^p$ norms and a bounded $L^p$ Lipschitz constant, including its expectation as calculated above. Lemma~\ref{lem:compact-supremum}, followed by the tightness of $\widehat h_n$, proves \eqref{eq:cumulant-gaussian-plugin}. Since $\mathfrak d(\widehat A_n)/n\to_P1$, Step 6 now yields
\[
 \sqrt n\,\widetilde\kappa_{4,n}
 =\frac{n}{\mathfrak d(\widehat A_n)}
       \left\{\frac{T(r_0,v_0)}{\sqrt n}+o_P(1)\right\}
 \Rightarrow N(0,24t_0^4).
\]
This proves \eqref{eq:cumulant-gaussian-limit}. The Gaussian limit is conditional on every good random design with the same uniform estimates; averaging adds at most
$P_X(\mathcal G_n^c)=o(1)$ to the distributional comparison. The proof with fixed smoothness uses the same restricted local sets as in Step 4. \qed

\subsection{Proof of Lemma~\ref{lem:studentization}}
\label{app:studentizationproof}
\noindent\textbf{Part (i): local information and Hessian.}
Fix $H$ and keep the local map $E_n$ evaluated at the truth. Use $M(h)$ in high dimension and $M(h,a)$ in low dimension, with the nuisance held fixed when taking $h$ derivatives. Write $\mathcal L_n(h)=\ell_n(\theta_0+E_nh)$ or $\mathcal L_n(h,a)=\ell_n(\vartheta_0+E_nh,a)$, respectively. Lemma~\ref{lem:derivatives}, parts
\ref{it:covariance-local} and \ref{it:covariance-hessian}, verifies the hypotheses of Lemma~\ref{lem:gaussian-hessian} on a fixed box containing the local ball, and over the range interval when present.

To identify its information matrix with the one in the statement, put $M_i=\partial_{h_i}M$ and $R=M^{-1}$. Cyclicity of the trace and the affine local parameter map give
\[
 \mathscr I_{ij}(h,a)=\tfrac12\tr(RM_iRM_j)
 =\left[E_n^\T\mathcal I_{n,\mathrm{id}}(\vartheta_0+E_nh,a)E_n\right]_{ij}
\]
in low dimension; omit $a$ and use the full parameter in high dimension. Thus the deterministic conclusion of Lemma~\ref{lem:gaussian-hessian} is \eqref{eq:localIfisher}. For the observed Hessian, the same affine change of coordinates gives
\[
 \nabla_h^2 \mathcal L_n(h,a)
 =E_n^\T\nabla_{\mathrm{id}}^2
       \ell_n(\vartheta_0+E_nh,a)E_n.
\]
The stochastic conclusion of that lemma is therefore exactly \eqref{eq:localIobserved}. This proves part~\ref{it:studentization-local}.

\noindent\textbf{Part (ii): range information and mixed derivatives.}
Let $d\le3$, $I=[a_-,a_+]$, and $\Gamma_H=\{\norm h\le H\}\times I$. All matrices in this paragraph are evaluated at $\theta=(\vartheta_0+E_nh,a)$, and $R=M^{-1}$. The information entries are
\[
 \begin{aligned}
 \calI_{aa,n}(\theta)
   &=\tfrac12\tr(RM_aRM_a)
     =\tfrac12\|R^{1/2}M_aR^{1/2}\|_\Frob^2,\\
 [E_n^\T\calI_{\vartheta a,n}(\theta)]_i
   &=\mathscr I_{ia}(h,a)\coloneqq \tfrac12\tr(RM_iRM_a).
 \end{aligned}
\]
Lemma~\ref{lem:derivatives}\ref{it:covariance-nuisance} bounds the first quantity above and away from zero, and makes the second uniformly $o(1)$ on $\Gamma_H$. This proves the expected-information assertions.

For the observed derivatives, put $X=C_{n,0}^{-1/2}Y_n$ and $Q(A)=X^\T AX-\tr A$. The exact mixed Hessian identity gives
\[
 \begin{aligned}
 A_{ia}&=RM_{ia}R-RM_iRM_aR-RM_aRM_iR,\\
 \mathcal L_{n,h_i a}-\E_0\mathcal L_{n,h_i a}&=\tfrac12Q(A_{ia}),\\
 \E_0\mathcal L_{n,h_i a}
 &=-\mathscr I_{ia}+\tfrac12\tr\{(I_n-M)A_{ia}\}.
 \end{aligned}
\]
The trace and norm bounds in Lemma~\ref{lem:derivatives}\ref{it:covariance-nuisance} imply
\[
 |\E_0\mathcal L_{n,h_i a}|
 \le |\mathscr I_{ia}|
       +\tfrac12\|I_n-M\|_\Frob\|A_{ia}\|_\Frob=o(1)
\]
uniformly. Differentiating the finite-dimensional quadratic-form identity under the fixed law $P_0$ gives
\[
 \begin{aligned}
 \E_0\mathcal L_{n,h_i aa}
 &=-\partial_a\mathscr I_{ia}
  +\tfrac12\tr\{-M_aA_{ia}+(I_n-M)\partial_aA_{ia}\},\\
 |\E_0\mathcal L_{n,h_i aa}|
 &\le |\partial_a\mathscr I_{ia}|
  +\tfrac12\|M_a\|_\Frob\|A_{ia}\|_\Frob
  +\tfrac12\|I_n-M\|_\Frob\|\partial_aA_{ia}\|_\Frob=o(1).
 \end{aligned}
\]
Again every bound is uniform on $\Gamma_H$.

The two centered parts are $\tfrac12Q(A_{ia})$ and $\tfrac12Q(\partial_aA_{ia})$. To control their suprema, take either $H_i(u)=A_{ia}(u)$ or $H_i(u)=\partial_aA_{ia}(u)$, and set $\epsilon_n=(\log n)/r_n$. For fixed $p\ge1$, choose $m>\max\{p,k_d+1\}$. Lemma~\ref{lem:gaussian} and the Frobenius and Lipschitz bounds in Lemma~\ref{lem:derivatives}\ref{it:covariance-nuisance} give
\[
 \begin{aligned}
 \sup_{u\in\Gamma_H}\|Q(H_i(u))\|_{L^m}
       &\le C_{H,m}\epsilon_n,\\
 \|Q(H_i(u))-Q(H_i(v))\|_{L^m}
       &\le C_{H,m}\epsilon_n\|u-v\|.
 \end{aligned}
\]
The set $\Gamma_H$ is fixed, compact, and convex. Lemma~\ref{lem:compact-supremum}
therefore gives
\[
 \left\|\sup_{u\in\Gamma_H}|Q(H_i(u))|\right\|_{L^p}
 \le\left\|\sup_{u\in\Gamma_H}|Q(H_i(u))|\right\|_{L^m}
 \le C_{H,m}\epsilon_n=o(1).
\]
Combining the mean and centered bounds proves both mixed observed assertions, in fact in every fixed $L^p$. This proves part~\ref{it:studentization-mixed}.

\noindent\textbf{Part (iii): the full inverse and its identifiable block.}
At a tight global fit, arrange the coordinates as $(\vartheta,a)$ in low dimension and write
\[
 \begin{pmatrix}E_n&0\\0&1\end{pmatrix}^{\!\top}
 \widehat{\mathcal A}_n
 \begin{pmatrix}E_n&0\\0&1\end{pmatrix}
 =\begin{pmatrix}J_n+o_P(1)&z_n\\z_n^\top&\iota_n\end{pmatrix}.
\]
The local bounds give $z_n=o_P(1)$. For expected curvature, $\iota_n\ge c$; for observed curvature, Lemma~\ref{lem:range-curvature} gives $|\iota_n|^{-1}=O_P(1)$.
Thus
\[
 E_n^\top\operatorname{Sc}_a(\widehat{\mathcal A}_n)E_n
 =J_n+o_P(1)-z_n\iota_n^{-1}z_n^\top=J_n+o_P(1).
\]
The Schur complement is positive definite and the range entry is nonzero with probability tending to one. The full matrix is consequently invertible, and the block-inverse identity derived in Appendix~\ref{app:full-inverse-block} gives
\[
 E_n^{-1}[\widehat{\mathcal A}_n^{-1}]_{\vartheta\vartheta}E_n^{-\top}
 =\{E_n^\top\operatorname{Sc}_a(\widehat{\mathcal A}_n)E_n\}^{-1}
 =J_n^{-1}+o_P(1).
\]
The full observed matrix itself need not be positive definite at a boundary range fit; the identifiable covariance block has the required positive definiteness. Direct identifiable inversion and full inversion in high dimension follow from the local curvature comparison alone.

Since $E_n^{-1}e_t=\sqrt n\,e_t$, the cumulant term in \eqref{eq:feasible-sandwich} becomes $\widetilde\kappa_{4,n}e_te_t^\top$ after normalization. Lemmas \ref{lem:kurtosis} and~\ref{lem:noise-diagonals} prove \eqref{eq:covariance-substitution} and its positive definiteness.

\smallskip
\noindent\emph{Other nuisance coordinates.}
Let $a=\alpha(q)$ as in Lemma~\ref{lem:range-curvature}. The expected curvature transforms by the block diagonal Jacobian, so its range Schur complement is unchanged. For the observed negative Hessian,
\[
 \mathcal H_{qq}=(\alpha')^2\mathcal H_{aa}
                      -\alpha''\partial_a\ell_n,
 \qquad \mathcal H_{\vartheta q}=\alpha'\mathcal H_{\vartheta a}.
\]
The normalized mixed block still tends to zero. At an interior fit the score term vanishes; at an endpoint it must be retained. Lemma~\ref{lem:range-curvature} gives
$|\mathcal H_{qq}|^{-1}=O_P(1)$ in either case, so the same Schur argument proves \eqref{eq:covariance-substitution} using the full matrix in $(\vartheta,q)$. In particular, $q$ can be the $a,\log a,l$ or $\log l$ actually used in the selected row, without changing the fitted model. Identifiable-coordinate transformations are handled
directly from this spectral covariance in the proof of Corollary~\ref{cor:uncertainty}.

All constants and probability bounds are uniform over $\mathcal K\times\mathcal P$. Restricting smoothness takes the corresponding principal blocks and removes the associated mixed derivatives. The allowed range sections are the restrictions of $I(\vartheta)$ to the fixed-smoothness slice and retain the same endpoint regularity. \qed

\subsection{Proof of Lemma~\ref{lem:range-curvature}}
\label{app:range-curvature-proof}
We check the hypotheses of Lemma~\ref{lem:nuisance-curvature}. Put $I=[a_-,a_+]$, let $a=\alpha(q)$ be the stated change of range coordinate from an interval $J$ onto $I$, and write
\[
 \mathcal L_n(h,a)=\ell_n(\vartheta_0+E_nh,a),\qquad
 G_n(a)=\mathcal L_n(0,a),\qquad
 \widetilde G_n(q)=G_n(\alpha(q)).
\]
The criterion is defined on the whole enclosing interval $I=I_*$, but the global fit maximizes over the actual section:
\[
 \widehat a_n\in\argmax_{a\in I_n(\widehat h_n)}
                         \mathcal L_n(\widehat h_n,a),\qquad
 \widehat h_n=O_P(1).
\]
The second assertion is Lemma~\ref{lem:tails}. These identities are needed only on $\{\widehat\vartheta_n\in U\}$, whose probability tends to one uniformly. Formally, on its complement replace the fitted indices by $h=0$ and a measurable maximizer of $G_n$ on $I(\vartheta_0)$; this changes any conclusion in probability by $o(1)$. Each fixed $h$-ball has nonempty feasible sections eventually, and \eqref{eq:feasible-endpoint-motion} verifies the endpoint-displacement condition of Lemma~\ref{lem:nuisance-curvature}. No consistency of $\widehat a_n$ is used. All covariance matrices below are evaluated at $(\vartheta_0,a)$, and
\[
 M(a)=C_0^{-1/2}C_n(\vartheta_0,a)C_0^{-1/2},\qquad R=M^{-1},
 \qquad P_0=C_0^{-1/2}S_0C_0^{-1/2}.
\]
The smooth Mat\'ern covariance and the positive nugget supply the differentiability of the criterion. We verify the four probabilistic hypotheses in turn, with constants uniform in the truth and noise law.

\noindent\textbf{(i) Bounded score densities.}
Since the field is Gaussian and independent of the errors,
\[
 X=C_0^{-1/2}F_n+C_0^{-1/2}\varepsilon,
 \qquad X\mid\varepsilon\sim N(z,P_0),\qquad
 z=C_0^{-1/2}\varepsilon.
\]
Equivalently, conditionally on the errors we may write $X=P_0^{1/2}G+z$, with $G\sim N(0,I_n)$ independent of them.

Put $B_a=-RM_aR$. Lemma~\ref{lem:derivatives}\ref{it:covariance-signal} gives $0\preceq B_a\preceq CP_0$ and shows that, for all sufficiently large $n$, $P_0^{1/2}B_aP_0^{1/2}$ has two eigenvalues bounded below by a common positive constant. Take $n_0$ large enough for this bound to hold. The score has the form
\[
 G_n'(a)=-\tfrac12X^\T B_aX-\tfrac12\tr(RM_a).
\]
Condition on $\varepsilon$, so $z$ is fixed. Choose an orthogonal matrix $O$ such that
\[
 P_0^{1/2}B_aP_0^{1/2}
 =O\operatorname{diag}(\lambda_1,\ldots,\lambda_n)O^\T,
 \qquad \lambda_1,\lambda_2\ge c,
\]
and put $\xi=O^\T G$ and $v=O^\T P_0^{1/2}B_a z$. The coordinates $\xi_j$ are independent standard normal variables, also conditionally on $\varepsilon$. Expanding the quadratic form gives
\[
 \begin{aligned}
 X^\T B_aX
 &=G^\T P_0^{1/2}B_aP_0^{1/2}G
       +2G^\T P_0^{1/2}B_a z+z^\T B_a z\\
 &=\sum_{j=1}^n\lambda_j\xi_j^2
       +2\sum_{j=1}^n v_j\xi_j+z^\T B_a z.
 \end{aligned}
\]
For $j=1,2$, define $b_j=v_j/\lambda_j$. Then
\[
 \lambda_j\xi_j^2+2v_j\xi_j
 =\lambda_j\left(\xi_j+\frac{v_j}{\lambda_j}\right)^2
       -\frac{v_j^2}{\lambda_j}
 =\lambda_j(\xi_j+b_j)^2-\frac{v_j^2}{\lambda_j}.
\]
Consequently,
\[
 \begin{aligned}
 X^\T B_aX
 &=\lambda_1(\xi_1+b_1)^2+\lambda_2(\xi_2+b_2)^2+W_\varepsilon,\\
 W_\varepsilon
 &=\sum_{j=3}^n(\lambda_j\xi_j^2+2v_j\xi_j)
       +z^\T B_a z-\sum_{j=1}^2\frac{v_j^2}{\lambda_j}.
 \end{aligned}
\]
Conditionally on $\varepsilon$, the shifts $b_1,b_2$ are fixed, and $W_\varepsilon$ is independent of $(\xi_1,\xi_2)$ because it depends only on the remaining Gaussian coordinates. Under $y_j=\sqrt{\lambda_j}(\xi_j+b_j)$, the inverse map is $\xi_j=y_j/\sqrt{\lambda_j}-b_j$. Its absolute Jacobian determinant is $(\lambda_1\lambda_2)^{-1/2}$. Conditionally on the errors, independence of $\xi_1,\xi_2$ gives the joint density
\[
 p(y_1,y_2)
 =\frac{\exp\!\left[-\tfrac12\sum_{j=1}^2
           (y_j/\sqrt{\lambda_j}-b_j)^2\right]}
        {2\pi\sqrt{\lambda_1\lambda_2}}
 \le\frac1{2\pi\sqrt{\lambda_1\lambda_2}}.
\]
Polar coordinates therefore give, for $s>0$,
\[
 p_{y_1^2+y_2^2}(s)
 =\frac12\int_0^{2\pi}
      p_{(y_1,y_2)}(\sqrt{s}\cos\phi,\sqrt{s}\sin\phi)\,d\phi
 \le\frac1{2\sqrt{\lambda_1\lambda_2}}.
\]
Convolution with the contribution of the other coordinates, scaling by $-1/2$, and averaging over the errors preserve a uniform bound. Thus
\begin{equation*}
  \sup_{n\ge n_0}\sup_{a\in I,\,\theta_0\in\mathcal K,\,P\in\mathcal P,\,x\in\R}
       p_{G_n'(a)}(x)\le C,
\end{equation*}
which verifies part~\ref{it:nuisance-curvature-density}.

\noindent\textbf{(ii) Tight derivative bounds.}
Equation~\eqref{eq:nuisance-baseline-moments} gives, for every fixed
$1\le p<\infty$,
\[
 \max_{1\le j\le3}
 \left\|\sup_{a\in I_*}|G_n^{(j)}(a)|\right\|_{L^p}\le C_p.
\]
This verifies part~\ref{it:nuisance-curvature-third}, including the first two derivative bounds needed when an endpoint moves.

\noindent\textbf{(iii) Endpoint anti-concentration.}
At either deterministic baseline endpoint
$q_e\in\partial\{\alpha^{-1}(I(\vartheta_0))\}$, put
\[
 A_q=RM_{qq}R-2RM_qRM_qR,\qquad
 \iota_q=\tfrac12\tr(RM_qRM_q),\qquad
 W_q=\widetilde G_n''(q_e).
\]
Since $M_q=\alpha'M_a$ and $|\alpha'|$ is bounded away from zero, Lemma~\ref{lem:derivatives}\ref{it:covariance-nuisance} gives $\iota_q\ge c$. The Hessian identity and the quadratic-form variance bound give
\[
 \E W_q=-\iota_q+\tfrac12\tr\{(I_n-M)A_q\},\qquad
 \Var W_q\le C\|A_q\|_\Frob^2,
 \qquad \|M-I_n\|_\Frob\le C.
\]
Fix a small $\eta>0$. If $\|A_q\|_\Frob\le\eta$, then
\[
 |\E W_q|\ge\iota_q-\tfrac12\|I_n-M\|_\Frob\|A_q\|_\Frob
             \ge c/2
\]
when $\eta$ is sufficiently small. For $0<\epsilon<c/4$, Chebyshev's inequality therefore gives
\[
 P(|W_q|\le\epsilon)
 \le P(|W_q-\E W_q|\ge c/4)\le C\eta^2.
\]
If instead $\|A_q\|_\Frob>\eta$, Lemma~\ref{lem:derivatives}\ref{it:covariance-signal}, equivalently \eqref{eq:range-gaussian-part}, gives
\[
 \|P_0^{1/2}A_qP_0^{1/2}\|_\Frob\ge c_\eta>0.
\]
Conditionally on $\varepsilon$, $W_q$ is a Gaussian quadratic polynomial with quadratic coefficient $P_0^{1/2}A_qP_0^{1/2}/2$. Its linear and constant coefficients
may depend on the errors, but Lemma~\ref{lem:quadratic-smallball} is uniform in those coefficients. Hence
\[
 P(|W_q|\le\epsilon\mid\varepsilon)
       \le C\sqrt{\epsilon/c_\eta}.
\]
Combining the two deterministic cases for $\|A_q\|_\Frob$, and then averaging over the errors, gives
\[
 \sup_{n\ge n_0,\,\theta_0,P}\max_{q_e\in\partial\{\alpha^{-1}(I(\vartheta_0))\}}
 P\{|\widetilde G_n''(q_e)|\le\epsilon\}
 \le C\eta^2+C\sqrt{\epsilon/c_\eta}.
\]
Let $\epsilon\downarrow0$ first and then $\eta\downarrow0$. This verifies part~\ref{it:nuisance-curvature-endpoints} for the chosen coordinate $q$, including the term $\alpha''G_n'$ at an endpoint. The estimates hold uniformly at every deterministic $a\in I_*$, so allowing these two endpoints to depend on the truth
does not change the bound.

\noindent\textbf{(iv) Transfer to bounded local sets.}
For $j=1,2$, integration along the segment $sh$ gives
\[
 \partial_a^j\mathcal L_n(h,a)-G_n^{(j)}(a)
   =\int_0^1\{\nabla_h\partial_a^j\mathcal L_n(sh,a)\}^\T h\,ds.
\]
Lemma~\ref{lem:studentization}\ref{it:studentization-mixed} therefore implies, for each fixed $H$,
\begin{equation*}
 \begin{split}
 &\sup_{\|h\|\le H,\ a\in I}\max_{j=1,2}
       |\partial_a^j\mathcal L_n(h,a)-G_n^{(j)}(a)|\le H\sup_{\|h\|\le H,\ a\in I}\max_{j=1,2}
                \|\nabla_h\partial_a^j\mathcal L_n(h,a)\|=o_P(1).
 \end{split}
\end{equation*}
This verifies part~\ref{it:nuisance-curvature-transfer}.

All hypotheses of Lemma~\ref{lem:nuisance-curvature} now hold, so \eqref{eq:nuisance-curvature-general} proves the asserted reciprocal bound. Fixing smoothness removes one identifiable coordinate and leaves these verifications unchanged. For iid sites, the same bounds hold on the good designs of Lemma~\ref{lem:randomgood}, with common constants; these events hold eventually almost surely. The general lemma therefore also gives the claimed conditional result. \qed

\subsection{Proof of Lemma~\ref{lem:logtransform}}
\label{app:logtransformproof}
Use the forward map $\zeta=g(\vartheta)$ throughout. Define
\[
 T_n^\zeta(h)
 =E_n^{-1}Dg(\vartheta_0)^{-1}
     \{g(\vartheta_0+E_nh)-g(\vartheta_0)\}.
\]
Thus $T_n^\zeta(0)=0$, $DT_n^\zeta(0)=I$, and
\[
 DT_n^\zeta(h)
 =E_n^{-1}Dg(\vartheta_0)^{-1}
       Dg(\vartheta_0+E_nh)E_n.
\]
The normalized Hessians are obtained directly from those in the proof of Corollary~\ref{cor:rates}:
\[
 H_{i,n}^\zeta(h)\coloneqq \nabla_h^2T_{n,i}^\zeta(h)
 =\sum_j\left[E_n^{-1}Dg(\vartheta_0)^{-1}\right]_{ij}
      E_n^\T\nabla_\vartheta^2g_j(\vartheta_0+E_nh)E_n.
\]
All matrices below use the order $(h_u,h_v,h_a,h_t)$. Omit the range coordinate when $d\le3$.

\noindent\textbf{Normalized Hessians.}
Write $\vartheta=\vartheta_0+E_nh$, $s=\sigma^2$, $B(\nu)=\partial_\nu\log A_d(\nu)$, and
\[
 q_c=2L_n-B(\nu),\qquad
 q_\sigma=2L_n-B(\nu)-\log a.
\]
A subscript $0$ denotes evaluation at the truth. The coordinates $c,l,s$ have their values at $\vartheta$, where $l=a^{-1/2}$.

For the spectral vector, $T_n^\zeta(h)=h$, so every Hessian is zero. For the three log parameterizations in Table~\ref{tab:parameterizations}, the last three Hessians are
the same:
\[
 \begin{aligned}
 H_{2,n}^\zeta&=\operatorname{diag}
     \left(0,-\frac{\nu_0}{\nu^2r_n},0,0\right),~~
 H_{3,n}^\zeta&=\operatorname{diag}
     \left(0,0,-\frac{a_0}{a^2b_n},0\right),~~
 H_{4,n}^\zeta=\operatorname{diag}
     \left(0,0,0,-\frac{t_0}{t^2\sqrt n}\right).
 \end{aligned}
\]
For the spectral log vector,
\[
 H_{1,n}^\zeta
 =\operatorname{diag}\left(0,\frac{2L_n\nu_0}{\nu^2r_n},0,0\right).
\]
For the microergodic log vector,
\[
 H_{1,n}^\zeta
 =\operatorname{diag}\left(0,
       \frac{\nu_0q_{c,0}/\nu^2-B'(\nu)}{r_n},0,0\right).
\]
The auxiliary chart in \eqref{eq:logcharts} gives exactly the same $T_n^\zeta$ and Hessians: replacing $\log l$ by $\log a=-2\log l$ is a fixed linear coordinate change.

For the variance--lengthscale log vector, the forward Jacobian at the truth is
\begin{equation}\label{eq:rawJacobian}
 Dg(\vartheta_0)=
 \begin{pmatrix}
 1&-B_0-\log a_0&-\nu_0/a_0&0\\
 0&0&-1/(2a_0)&0\\
 0&1/\nu_0&0&0\\
 0&0&0&1/t_0
 \end{pmatrix},
\end{equation}
and its first normalized Hessian is
\[
 H_{1,n}^\zeta=
 \begin{pmatrix}
 0&0&0&0\\
 0&\{\nu_0q_{\sigma,0}/\nu^2-B'(\nu)\}/r_n&-1/(ab_n)&0\\
 0&-1/(ab_n)&r_n(\nu-\nu_0)/(a^2b_n^2)&0\\
 0&0&0&0
 \end{pmatrix}.
\]

For both original-scale vectors $(c,\nu,l,t)$ and $(s,l,\nu,t)$,
\[
 H_{2,n}^\zeta=H_{4,n}^\zeta=0,\qquad
 H_{3,n}^\zeta=\operatorname{diag}
       \left(0,0,-\frac{3l^5}{2l_0^3b_n},0\right).
\]
For the microergodic vector,
\[
 H_{1,n}^\zeta=\frac{c}{c_0r_n}
 \begin{pmatrix}
 1&q_c&0&0\\
 q_c&q_c^2-B'(\nu)&0&0\\
 0&0&0&0\\
 0&0&0&0
 \end{pmatrix}.
\]
For the original-scale variance--lengthscale vector,
\[
 \begin{aligned}
 H_{1,n}^\zeta
 &=\frac{s}{s_0}
 \begin{pmatrix}
 r_n^{-1}&q_\sigma/r_n&-\nu/(ab_n)&0\\
 q_\sigma/r_n&(q_\sigma^2-B'(\nu))/r_n
       &-(1+\nu q_\sigma)/(ab_n)&0\\
 -\nu/(ab_n)&-(1+\nu q_\sigma)/(ab_n)
       &\nu(\nu+1)r_n/(a^2b_n^2)&0\\
 0&0&0&0
 \end{pmatrix}-\operatorname{diag}
       \left(0,0,\frac{3\nu_0lr_n}{2l_0a^2b_n^2},0\right).
 \end{aligned}
\]

\noindent\textbf{Reading off the remainder bounds.}
The resulting orders are recorded in Table~\ref{tab:forward-remainder-orders}. Every displayed entry is bounded by $C\epsilon_n^\zeta$, except for the $(h_a,h_a)$ entry of $H_{1,n}^\zeta$ in the log-scale variance--lengthscale coordinates when $d\ge4$. Here we use $b_n(1+L_n)/r_n\to0$ to simplify the original-scale variance--lengthscale bound. The remaining entry is $h_v/(a^2b_n^2)$, since $\nu-\nu_0=h_v/r_n$; we bound its
contribution after inserting it in the integrals below.

\begin{table}[!htbp]
\centering\small
\renewcommand{\arraystretch}{1.35}
\caption{Orders in every fixed $L^p$ norm for the normalized Taylor remainder and Jacobian increment. Terms involving $b_n$ are omitted when $d\le3$. The variance--lengthscale rows apply only when $d\ge4$; the original-scale row tends to zero when $d>8$.}
\label{tab:forward-remainder-orders}
\begin{tabular}{@{}ll@{}}\toprule
Parameterization & $\epsilon_n^\zeta$\\\midrule
Spectral & $0$\\
Spectral, log scale & $(1+L_n)/r_n+b_n^{-1}+n^{-1/2}$\\
Microergodic & $(1+L_n)^2/r_n+b_n^{-1}$\\
Microergodic, log scale & $(1+L_n)/r_n+b_n^{-1}+n^{-1/2}$\\
Variance--lengthscale, log scale
 & $(1+L_n)/r_n+b_n^{-1}+n^{-1/2}$\\
Variance--lengthscale & $r_n/b_n^2+b_n^{-1}$\\\bottomrule
\end{tabular}
\end{table}
\FloatBarrier

Put $\mathcal R_n^\zeta=T_n^\zeta(\widehat h_n)-\widehat h_n$. Taylor's formula and the fundamental theorem of calculus give
\[
 \begin{aligned}
 [\mathcal R_n^\zeta]_i
 &=\int_0^1(1-u)\widehat h_n^\T
       H_{i,n}^\zeta(u\widehat h_n)\widehat h_n\,du,\\
 [E_n^{-1}Dg(\vartheta_0)^{-1}
       Dg(\widehat\vartheta_n)E_n-I]_{i\cdot}
 &=\int_0^1\widehat h_n^\T
       H_{i,n}^\zeta(u\widehat h_n)\,du.
 \end{aligned}
\]
Theorem~\ref{thm:normalized} gives uniform bounds for every fixed moment of $\widehat h_n$. In particular, all its coordinates are $O_{L^p}(1)$ for each fixed $p\ge1$. For the log-scale variance--lengthscale coordinates, the dependence of $T_{n,1}^\zeta(h)$ on $a$ is through $-h_v\log(a/a_0)$. Along the Taylor segment, its second $h_a$ derivative is
\[
 [H_{1,n}^\zeta(u\widehat h_n)]_{33}
 =\frac{u\widehat h_{n,v}}
        {\{a_0+u(\widehat a_n-a_0)\}^2b_n^2}.
\]
Since $a_0+u(\widehat a_n-a_0)\ge a_-$, this entry contributes at most $C|\widehat h_{n,v}|\widehat h_{n,a}^2/b_n^2$ and $C|\widehat h_{n,v}\widehat h_{n,a}|/b_n^2$ to the two integrals. H\"older's inequality and the moment bounds give
\[
 \begin{aligned}
 \left\|\frac{|\widehat h_{n,v}|\widehat h_{n,a}^2}{b_n^2}
       \right\|_{L^p}
 &\le\frac1{b_n^2}\|\widehat h_{n,v}\|_{L^{3p}}
                       \|\widehat h_{n,a}\|_{L^{3p}}^2
 \le\frac{C_p}{b_n^2},\\
 \left\|\frac{|\widehat h_{n,v}\widehat h_{n,a}|}{b_n^2}
       \right\|_{L^p}
 &\le\frac1{b_n^2}\|\widehat h_{n,v}\|_{L^{2p}}
                       \|\widehat h_{n,a}\|_{L^{2p}}
 \le\frac{C_p}{b_n^2}.
 \end{aligned}
\]
These two integral terms are therefore $O_{L^p}(b_n^{-2})$. For this row of the table, $b_n^{-2}\le\epsilon_n^\zeta$ for all sufficiently large $n$.

Every remaining Hessian entry is bounded by $C\epsilon_n^\zeta$ along the Taylor segment. In the Taylor remainder, its term is multiplied by two coordinates of $\widehat h_n$; in the Jacobian increment, it is multiplied by one. Explicitly, for each such entry,
\[
 \begin{aligned}
 \left|\int_0^1(1-u)\widehat h_{n,j}
       [H_{i,n}^\zeta(u\widehat h_n)]_{jk}
       \widehat h_{n,k}\,du\right|
 &\le C\epsilon_n^\zeta
                   |\widehat h_{n,j}\widehat h_{n,k}|,\\
 \left|\int_0^1\widehat h_{n,j}
       [H_{i,n}^\zeta(u\widehat h_n)]_{jk}\,du\right|
 &\le C\epsilon_n^\zeta|\widehat h_{n,j}|.
 \end{aligned}
\]
Summing over the fixed number of indices gives a bound $C\epsilon_n^\zeta\|\widehat h_n\|^2$ for these terms in the Taylor remainder, and $C\epsilon_n^\zeta\|\widehat h_n\|$ for these terms in the Jacobian increment. Taking $L^p$ norms and using the moment bounds gives, for every row and every fixed $p\ge1$,
\begin{equation}\label{eq:common-chart-remainder-bound}
 \begin{aligned}
 \|\mathcal R_n^\zeta\|_{L^p}
 &\le C_p\epsilon_n^\zeta,\\
 \left\|\left\|E_n^{-1}Dg(\vartheta_0)^{-1}
       Dg(\widehat\vartheta_n)E_n-I\right\|_{\mathrm{op}}
       \right\|_{L^p}
 &\le C_p\epsilon_n^\zeta.
 \end{aligned}
\end{equation}
The coefficients tend to zero for the vectors covered by the lemma; for the original-scale variance--lengthscale vector this uses $r_n/b_n^2=R_n^{4-d/2}\to0$ when $d>8$.
This proves convergence in every fixed positive moment. Continuity of matrix inversion gives the assertion for the inverse Jacobian product.

Finally, the scalar Taylor expansions and remainder orders in Table~\ref{tab:coordinate-expansions} give \eqref{eq:raw-marginal-remainder} directly, including when $4\le d\le8$. \qed

\subsection{Proof of Lemma~\ref{lem:randomgood}}\label{app:randomgoodproof}
Use the cell-count event $\mathcal C_n$ constructed in the proof of Lemma~\ref{lem:sampling}, with the fixed order $s=\nu_-+d/2>d/2$. It depends only on the observation locations, and \eqref{eq:cellcountprob} gives
\[
 P_X(\mathcal C_n^c)\le Ce^{-cn^\delta},
 \qquad \delta=1-d/(2s)>0.
\]
On this event every cell of side $\ell_n\asymp n^{-1/(2s)}$ contains an observation. Thus the polynomial coverage bound holds on $Q_+$, and hence on $Q_-$, with the fixed exponent $\gamma=1/(2s)<1/d$.

On the same event, \eqref{eq:averagedsampling} holds simultaneously for every $f\in H^s(\Omega)$. The uniform Sobolev comparison and variational deduction in the proof of Lemma~\ref{lem:sampling} therefore give \eqref{eq:predictionrisk} with every parameter in the compact rectangle taken as truth, with common constants.

Taking $\mathcal G_n=\mathcal C_n$ proves \eqref{eq:randomfailure} and both design properties, on the common smaller cube $Q_-$ if necessary. The deterministic information, local tightness, and normal-approximation arguments use only these properties with their fixed constants, the geometry, and the compact parameter bounds. Their constants therefore apply to every design in $\mathcal G_n$. \qed

\subsection{Proof of Lemma~\ref{lem:spectral-covariance}}
\label{app:spectral-covariance-proof}
For ordinary test functions, Fubini and the evenness of $f_\theta$ give
\begin{align*}
 \Cov_\theta(Z_{g_1},Z_{g_2})
 &=\int_{\R^d}\!\int_{\R^d}
       g_1(x)\overline{g_2(y)}K_\theta(x-y)\,dx\,dy\\
 &=\int_{\R^d}f_\theta(\omega)
       \left(\int_{\R^d}g_1(x)e^{-2\pi i\omega^\T x}\,dx\right)
       \overline{\left(\int_{\R^d}g_2(y)e^{-2\pi i\omega^\T y}\,dy\right)}
       \,d\omega\\
 &=\int_{\R^d}f_\theta(\omega)\widehat g_1(\omega)
                         \overline{\widehat g_2(\omega)}\,d\omega.
\end{align*}
The functions are integrable because they have compact support, and $f_\theta$ is integrable, so these interchanges are justified by absolute integrability. The same calculation with finite sums in place of spatial integrals covers point masses and their mixed terms. Equation \eqref{eq:sampledspectralcovariance} is the corresponding
spectral covariance representation at the sampled sites.

Apply the first identity with $g_1=g_2=L^{-d/2}\chi e^{-2\pi i\lambda^\T x}$ to obtain
\begin{align*}
 v_0(\lambda)
 &=L^{-d}\int_{\R^d}f_0(\omega)
                  |\widehat\chi(\omega+\lambda)|^2\,d\omega=L^{-d}\int_{\R^d}f_0(\omega)
                  |\widehat\chi(\lambda-\omega)|^2\,d\omega.
\end{align*}
The last equality changes $\omega$ to $-\omega$ and uses
$f_0(-\omega)=f_0(\omega)$. Hence
\[
 \int_{\R^d}f_0(\omega)|\ell_\lambda(\omega)|^2\,d\omega
 =\frac{L^{-d}}{v_0(\lambda)}
       \int_{\R^d}f_0(\omega)|\widehat\chi(\omega+\lambda)|^2
                    \,d\omega
 =1.
\]
For the prediction, use its signal-plus-noise representation
\[
 \widehat U_{\lambda,n}
 =Z_{\check p_{\lambda,n}}+b_{\lambda,n}\varepsilon_n,\qquad
 \check p_{\lambda,n}=\sum_i b_{\lambda,n,i}\delta_{x_{n,i}}.
\]
The two terms are centered and independent. Also $\widehat{\check p_{\lambda,n}}=p_{\lambda,n}$, since the Fourier transform of $\delta_{x_{n,i}}$ is $e^{-2\pi i\omega^\T x_{n,i}}$. Applying \eqref{eq:spectralcovariance} to this finite sum gives directly
\begin{align*}
 \E_0|\widehat U_{\lambda,n}|^2
 &=\E_0|Z_{\check p_{\lambda,n}}|^2
       +\E_0|b_{\lambda,n}\varepsilon_n|^2\\
 &=\int_{\R^d}f_0(\omega)
       |\widehat{\check p_{\lambda,n}}(\omega)|^2\,d\omega
       +t_0\sum_i|b_{\lambda,n,i}|^2\\
 &=\int_{\R^d}f_0(\omega)|p_{\lambda,n}(\omega)|^2\,d\omega
       +t_0\sum_i|b_{\lambda,n,i}|^2.
\end{align*}
This proves \eqref{eq:fouriernormalizations}. \qed

\subsection{Proof of Lemma~\ref{lem:cubesampling}}
\label{sec:cubesampling}
This is the local $p=q=\kappa=2$, $l=0$, $r=s$ case of the sampling inequality of Arcang\'eli et al., restated by \citet[Theorem~1.1]{corrigan2009}. The following self-contained proof makes the powers of $h$ and the absence of a separation assumption explicit.

\noindent\textbf{Polynomial approximation on a fixed cube.}
For a cube $V$, write $s=k+\sigma$, with $k$ a nonnegative integer and $0\le\sigma<1$, and use the highest-order seminorm
\[
 [g]_{s,V}^2=
 \begin{cases}
 \displaystyle\sum_{|\alpha|=k}\norm{D^\alpha g}_{L^2(V)}^2,
       &\sigma=0,\\[6pt]
 \displaystyle\sum_{|\alpha|=k}\int_V\int_V
       \frac{|D^\alpha g(x)-D^\alpha g(y)|^2}
            {\norm{x-y}^{d+2\sigma}}\,dx\,dy,
       &0<\sigma<1.
 \end{cases}
\]
Together with the lower-order $L^2$ derivative terms, this gives an equivalent Sobolev norm on a fixed cube.

First work on $U=(0,1)^d$ and $V=(-1,2)^d$. Put $m=\lceil s\rceil-1$, and let $\mathcal P_m$ denote the polynomials of total degree at most $m$. For every $g\in H^s(V)$ there is $p\in\mathcal P_m$ such that
\[
 \norm{g-p}_{H^s(V)}+\norm{g-p}_{L^\infty(V)}
       \le C[g]_{s,V}.
\]
Here is the polynomial approximation step. Choose the coefficients of $p$, in descending degree, so that $\int_V D^\alpha(g-p)=0$ for all $|\alpha|\le m$. For integer $s$, successive applications of Poincar\'e's inequality control all lower derivatives by those of order $s$. For noninteger $s=k+\sigma$, first apply the fractional Poincar\'e estimate to derivatives of order $k$:
\[
 \norm{w-\overline w}_{L^2(V)}^2
 =\frac1{2|V|}\int_V\int_V|w(x)-w(y)|^2\,dx\,dy
 \le C[w]_{\sigma,V}^2,
 \qquad \overline w=|V|^{-1}\int_Vw.
\]
Then apply the ordinary Poincar\'e inequality to the lower derivatives. This proves the $H^s$ bound; the $L^\infty$ bound follows from Sobolev embedding since $s>d/2$.

\noindent\textbf{Stable evaluation of the polynomial.}
Choose the fixed tensor grid
\[
 \{\xi_1,\ldots,\xi_M\}
 =\left\{\left(\frac{j_1+1}{m+2},\ldots,
                    \frac{j_d+1}{m+2}\right):
                    j_1,\ldots,j_d\in\{0,\ldots,m\}\right\}
 \subset U.
\]
Values on this grid determine every polynomial in $\mathcal P_m$. Choose an $L^2(U)$-orthonormal basis $q_1,\ldots,q_N$ of $\mathcal P_m$, where $N=\dim\mathcal P_m$, and define its evaluation matrix by $A_{jr}=q_r(\xi_j)$. For $p=\sum_{r=1}^N c_rq_r$,
\[
 \norm p_{L^2(U)}^2=\norm c^2,
 \qquad
 \sum_{j=1}^M|p(\xi_j)|^2=\norm{Ac}^2.
\]
Since the grid determines $p$, the matrix $A$ has trivial kernel. Compactness of the unit sphere therefore gives
\[
 \sigma_*\coloneqq \min_{\norm c=1}\norm{Ac}>0,
 \qquad
 \norm c^2\le\sigma_*^{-2}\norm{Ac}^2.
\]
By continuity of the evaluation matrix in its evaluation points, this bound holds with a uniform constant for $\widetilde A_{jr}=q_r(\widetilde\xi_j)$ whenever the points
are sufficiently close to the original grid. Thus, for some fixed $\varepsilon>0$,
\[
 \norm p_{L^2(U)}^2\le C\sum_{j=1}^M|p(\widetilde\xi_j)|^2,
 \qquad
 \max_j\norm{\widetilde\xi_j-\xi_j}\le\varepsilon,
 \quad p\in\mathcal P_m.
\]
Choose $\varepsilon$ small enough that these perturbed points lie in $V$ and are distinct. Apply this bound to the polynomial $p$ above and use the $L^\infty$ error bound at its nodes to get
\[
 \norm g_{L^2(U)}^2
       \le C\left\{[g]_{s,V}^2
                    +\sum_{j=1}^M|g(\widetilde\xi_j)|^2\right\}.
\]

\noindent\textbf{Rescaling and summing over the coverage cells.}
Now fix $A_*$ so large that $A_*^{-1}\le\varepsilon$, put $\ell=A_*h$, and take the grid cells $E=b+\ell U$, with $b\in\ell\mathbb Z^d$, that meet $Q_-$. Thus $b$ is the
lower corner of each cell. For $h_*$ small enough, every enlargement $E^+=b+\ell V$ lies in $Q_+$. Coverage supplies a site $x_{i(E,j)}$ within distance $h$ of $b+\ell\xi_j$. Thus its scaled position $\widetilde\xi_j=(x_{i(E,j)}-b)/\ell$ is within $A_*^{-1}$ of $\xi_j$. Apply the preceding estimate to $g(y)=f(b+\ell y)$ and use $[g]_{s,V}^2=\ell^{2s-d}[f]_{s,E^+}^2$ to obtain
\[
 \norm f_{L^2(E)}^2
 \le C\left\{\ell^{2s}[f]_{s,E^+}^2
              +\ell^d\sum_{j=1}^M|f(x_{i(E,j)})|^2\right\}.
\]

The enlargements $E^+$ have bounded overlap. Therefore
\[
 \sum_E[f]_{s,E^+}^2\le C[f]_{s,\Omega}^2
       \le C\norm f_{H^s(\Omega)}^2.
\]
For fractional $s$, the first inequality follows by counting the bounded number of enlargements containing each pair $(x,y)$. Each observation index is selected at most a bounded number of times, because selection for $E$ requires $x_i\in E^+$. Summing over the cells covering $Q_-$, and using $\ell=A_*h$ with $A_*$ fixed, proves \eqref{eq:sobolevsampling}.
\qed

\subsection{Proof of Lemma~\ref{lem:gaussian-hessian}}
Put $Q(A)=X^\T AX-\tr A$. Differentiate $RM=I$ to obtain $\partial_iR=-RM_iR$. The exact Hessian identity is
\begin{equation*}
 \partial_{ij}\mathcal L_n
 =\tfrac12X^\T A_{ij}X+\tfrac12\tr(RM_iRM_j-RM_{ij}).
\end{equation*}
Its expectation depends only on $\E XX^\T=I$. In particular,
\[
 \E\partial_{ij}\mathcal L_n+\tfrac12\tr(RM_iRM_j)
 =\tfrac12\tr\{(I-M)A_{ij}\},\qquad
 \partial_{ij}\mathcal L_n-\E\partial_{ij}\mathcal L_n=\tfrac12Q(A_{ij}).
\]
The deterministic trace correction is $o(1)$, since $\|M-I\|_\Frob\le C$ and $\|A_{ij}\|_\Frob=o(1)$. To compare the remaining trace directly, note that $R-I=-R(M-I)$ gives $\|R-I\|_\Frob\le C$. Hence
\[
 \begin{aligned}
 |\tr(RM_iRM_j)-\tr(M_iM_j)|
 &\le C\|R-I\|_\Frob\bigl(
      \|M_i\|_\op\|M_j\|_\Frob
       +\|M_j\|_\op\|M_i\|_\Frob\bigr)=o(1),\\
 |\tr(M_iM_j)-\tr(B_iB_j)|
 &\le\|M_i-B_i\|_\Frob\|M_j\|_\Frob
       +\|B_i\|_\Frob\|M_j-B_j\|_\Frob=o(1).
 \end{aligned}
\]
Thus $\tfrac12\tr(RM_iRM_j)=J_{n,ij}+o(1)$. This proves the deterministic information comparison in \eqref{eq:commonHessianconclusion}.

\noindent\textbf{Centered Hessian via the compact-set bound.}
Fix $p\ge1$ and choose a fixed $m>\max\{p,k\}$. For one Hessian entry, define the continuous process
\[
 \xi_{ij,n}(u)=Q(A_{ij}(u)),\qquad u\in\Gamma.
\]
The quadratic-form moment bound in Lemma~\ref{lem:gaussian}, with $\|A\|_\op\le\|A\|_\Frob$, gives
\[
 \sup_{u\in\Gamma}\|\xi_{ij,n}(u)\|_{L^m}
       \le C_m\epsilon_n.
\]
Linearity of $Q$ gives the increment bound from the same lemma:
\[
 \begin{aligned}
 \|\xi_{ij,n}(u)-\xi_{ij,n}(v)\|_{L^m}
 &=\|Q(A_{ij}(u)-A_{ij}(v))\|_{L^m}\\
 &\le C_m\|A_{ij}(u)-A_{ij}(v)\|_\Frob
 \le C_mL\|u-v\|.
 \end{aligned}
\]
Apply the Ibragimov--Has'minskii compact-set bound in Lemma~\ref{lem:compact-supremum}, with pointwise bound $C_m\epsilon_n$ and Lipschitz constant $C_mL$. The diameter of $\Gamma$ is fixed, so
\[
 \E\sup_{u\in\Gamma}|\xi_{ij,n}(u)|^m
 \le C_{\Gamma,m}\{\epsilon_n^m+L^k\epsilon_n^{m-k}\}.
\]
Since $m>k$, both powers of $\epsilon_n$ are positive. Since $m>p$, monotonicity of probability-space norms now gives
\[
 \begin{aligned}
 \left\|\sup_{u\in\Gamma}|Q(A_{ij}(u))|\right\|_{L^p}
 &\le\left\|\sup_{u\in\Gamma}|\xi_{ij,n}(u)|\right\|_{L^m}\\
 &\le C_{\Gamma,m}
      \{\epsilon_n^m+L^k\epsilon_n^{m-k}\}^{1/m}=o(1).
 \end{aligned}
\]
There are only finitely many Hessian entries. Combining this with the deterministic mean bounds proves the Hessian conclusion in \eqref{eq:commonHessianconclusion} in $L^p$, hence also in probability.
\qed

\subsection{Proof of Lemma~\ref{lem:score-comparison}}
The score changes with $a$ through a mixed derivative of the criterion. We first bound that derivative in $L^p$ at each $a$, then integrate over the fixed nuisance interval.

Differentiate $RM=I_n$ to obtain $\partial_aR=-RM_aR$. The score and its derivative are
\[
 \begin{aligned}
 \partial_{h_i}\mathcal L_n
   &=\tfrac12X_n^\T RM_iRX_n-\tfrac12\tr(RM_i),\\
 \partial_{h_i a}\mathcal L_n
   &=\tfrac12X_n^\T T_{ia}X_n
       +\tfrac12\tr(RM_aRM_i-RM_{ia}),\\
 T_{ia}&=RM_{ia}R-RM_aRM_iR-RM_iRM_aR.
 \end{aligned}
\]
All matrices here are evaluated at $(0,a)$. Cyclicity of the trace and $MR=I_n$ give
\[
 \tr(MT_{ia})=\tr(RM_{ia})-2\tr(RM_iRM_a).
\]
Consequently the preceding identity can be written as
\[
 \partial_{h_i a}\mathcal L_n
 =\tfrac12Q(T_{ia})-\tfrac12\tr(RM_iRM_a)
       +\tfrac12\tr\{(I_n-M)T_{ia}\}.
\]
This separates the centered quadratic fluctuation from the two deterministic terms. The inverse and nuisance-derivative bounds give
\[
 \|T_{ia}\|_\Frob
 \le C\{\|M_{ia}\|_\Frob+\|M_i\|_\op\|M_a\|_\Frob\}
 \le C\epsilon_n.
\]
The assumed quadratic-form moment bound therefore gives $\|Q(T_{ia})\|_{L^p}\le C_p\epsilon_n$. The mixed trace is at most $\rho_n$, and the last trace is at most
$\|M-I_n\|_\Frob\|T_{ia}\|_\Frob\le C\epsilon_n$. Thus
\[
 \sup_{a\in I}\max_i\|\partial_{h_i a}\mathcal L_n(0,a)\|_{L^p}
 \le C_p(\epsilon_n+\rho_n).
\]
For every $a\in I$, the fundamental theorem of calculus gives
\[
 s_{n,i}(a)-s_{n,i}(a_0)
   =\int_{a_0}^a\partial_{h_i b}\mathcal L_n(0,b)\,db.
\]
Taking the supremum pathwise and then applying Minkowski's inequality yields
\[
 \left\|\sup_{a\in I}|s_{n,i}(a)-s_{n,i}(a_0)|\right\|_{L^p}
 \le\int_I\|\partial_{h_i b}\mathcal L_n(0,b)\|_{L^p}\,db
 \le C_{p,I}(\epsilon_n+\rho_n).
\]
Summing over the fixed number of score coordinates proves \eqref{eq:score-comparison-general}. Finally, trace-norm duality gives
\[
 |\tr(RM_iRM_a)|
 \le\|RM_iR\|_\op\|M_a\|_1
 \le C\|M_i\|_\op\|M_a\|_1,
\]
which proves the stated sufficient condition. All constants come from the uniform bounds in the statement, so the same proof gives the asserted uniformity. \qed

\subsection{Proof of Lemma~\ref{lem:nuisance-curvature}}
An interior maximizer has zero score. If its curvature were also small, the score would be nearly tangent to zero. The first part of the proof is a quantitative form of Bulinskaya's lemma; see \citet[Proposition~2.3]{azais2001regularity} for the classical
no-tangency statement. We give the grid argument because it also provides uniformity in $n$ and accommodates the approximation in part~\ref{it:nuisance-curvature-transfer}.

\noindent\textbf{Nearly flat stationary points.}
Fix $0<\epsilon\le1$ and $K>0$. On the event $T_n\le K$, suppose some $a\in I$ satisfies
\[
 |G_n'(a)|\le\epsilon,\qquad |G_n''(a)|\le2\epsilon.
\]
Choose a deterministic grid of spacing at most $\sqrt\epsilon$, including the endpoints of $I$. It has at most $C_I\epsilon^{-1/2}$ points, and some grid point $b$ satisfies $|b-a|\le\sqrt\epsilon$. Taylor's theorem gives
\[
 \begin{aligned}
 |G_n'(b)|
 &\le |G_n'(a)|+|G_n''(a)|\,|b-a|
               +\tfrac12K|b-a|^2\\
 &\le\epsilon+2\epsilon^{3/2}+\tfrac12K\epsilon
 \le(3+K/2)\epsilon.
 \end{aligned}
\]
At each fixed grid point the density bound implies
\[
 P\{|G_n'(b)|\le(3+K/2)\epsilon\}
 \le 2C(3+K/2)\epsilon.
\]
A union bound therefore yields
\begin{equation}\label{eq:range-flat-stationary}
 \begin{split}
 P\{\exists a\in I:
          |G_n'(a)|\le\epsilon,\ |G_n''(a)|\le2\epsilon\}
 \le P(T_n>K)+C_1(1+K)\sqrt\epsilon,
 \end{split}
\end{equation}
where $C_1$ depends only on the density bound and the length of $I$. Thus a score cannot be nearly zero and nearly flat at a common point with appreciable probability, uniformly in $n$.

\noindent\textbf{Interior and endpoint fits.}
Let $c_\alpha=\inf_J|\alpha'|>0$, and fix $H<\infty$. On the event $\{\|\widehat h_n\|\le H,\ d_{n,H}\le\delta\}$, first suppose that $\widehat a_n$ is interior to
$I_n(\widehat h_n)$. Maximization in $a$ gives $\partial_a\mathcal L_n(\widehat h_n,\widehat a_n)=0$, so
\[
 \partial_q^2\widetilde{\mathcal L}_n(\widehat h_n,\widehat q_n)
 =\{\alpha'(\widehat q_n)\}^2
                      \partial_a^2\mathcal L_n(\widehat h_n,\widehat a_n).
\]
Consequently, if its absolute value is at most $\delta$, then
\[
 |G_n'(\widehat a_n)|\le\delta,\qquad
 |G_n''(\widehat a_n)|\le(1+c_\alpha^{-2})\delta.
\]
These inequalities imply the event in \eqref{eq:range-flat-stationary} with $\epsilon=C_\alpha\delta$ for a fixed sufficiently large $C_\alpha$.

If $\widehat a_n$ is an endpoint of $I_n(\widehat h_n)$, let $q_e$ be the corresponding endpoint of $\mathcal J_n(0)$. The inverse of $\alpha$ is Lipschitz, so on $\{\|\widehat h_n\|\le H\}$,
\[
 |\widehat q_n-q_e|\le c_\alpha^{-1}\rho_{n,H}.
\]
The endpoints need not be fixed as $h$ changes. We therefore control the change in the baseline second derivative between these two points. Let $\omega_{\alpha''}(r)$ be the deterministic modulus of continuity of $\alpha''$ on $J$. The identity
\[
 \widetilde G_n''(q)
  =(\alpha'(q))^2G_n''(\alpha(q))
                  +\alpha''(q)G_n'(\alpha(q))
\]
and the derivative bound $T_n$ give, for $q,q'\in J$,
\[
 |\widetilde G_n''(q)-\widetilde G_n''(q')|
 \le C_\alpha T_n
        \{|q-q'|+\omega_{\alpha''}(|q-q'|)\}.
\]
Indeed the changes in $G_n'$ and $G_n''$ are bounded by $T_n$ times the change in their argument, the change in $(\alpha')^2$ is Lipschitz, and the remaining coefficient uses the displayed modulus. Since $\rho_{n,H}\to0$ and $T_n=O_P(1)$,
\[
 \omega_{n,H}\coloneqq 
 \sup_{|q-q'|\le c_\alpha^{-1}\rho_{n,H}}
       |\widetilde G_n''(q)-\widetilde G_n''(q')|=o_P(1).
\]
The transfer hypothesis and the chain rule now yield
\[
 |\partial_q^2\widetilde{\mathcal L}_n(\widehat h_n,\widehat q_n)
                         -\widetilde G_n''(q_e)|
 \le C_\alpha d_{n,H}+\omega_{n,H}.
\]
Thus, on $\{d_{n,H}\le\delta,\ \omega_{n,H}\le\delta\}$, small fitted curvature implies
$|\widetilde G_n''(q_e)|\le C_\alpha\delta$. The endpoint score term $\alpha''G_n'$ has been retained throughout.

Combining the two cases gives, for sufficiently small $\delta$,
\[
 \begin{split}
 &~~~~P\{|\partial_q^2\widetilde{\mathcal L}_n(\widehat h_n,\widehat q_n)|
                                                    \le\delta\}\\
 &\le P(\|\widehat h_n\|>H)+P(d_{n,H}>\delta)
                +P(\omega_{n,H}>\delta)+P(T_n>K)\\
 &\quad~+C_2(1+K)\sqrt\delta
       +\sum_{q_e\in\partial \mathcal J_n(0)}
          P\{|\widetilde G_n''(q_e)|\le C_\alpha\delta\}.
 \end{split}
\]
Take $\limsup_n$, then let $\delta\downarrow0$. The transfer error, the grid term, and the endpoint probabilities disappear. Letting $K\to\infty$ and $H\to\infty$ removes the two tightness terms. This proves \eqref{eq:nuisance-curvature-general} and hence the reciprocal bound. The same inequalities with a supremum over the family of experiments prove the uniform assertion. \qed

\subsection{Proof of Lemma~\ref{lem:quadratic-smallball}}
This is the degree-two case of the Carbery--Wright inequality \citep[Theorem~8, with $q=4$]{carbery2001distributional} for Gaussian measure. With $P(G)=G^\T AG+b^\T G+c-z$, that inequality gives
\[
 P\{|P(G)|\le\epsilon\}
 \le C\epsilon^{1/2}\{\E|P(G)|^2\}^{-1/4}.
\]
The Gaussian quadratic and linear terms are uncorrelated, so
\[
 \E|P(G)|^2\ge\Var(P(G))=2\|A\|_\Frob^2+\|b\|^2\ge2\tau^2.
\]
This proves the assertion uniformly in $z,b,c,N$. Conditioning on independent random coefficients preserves the same bound. \qed

\subsection{Proof of Lemma~\ref{lem:gaussian-affinity}}
Write $P=C_1^{-1}$ and $Q=C_2^{-1}$. Multiplying the square roots of the two Gaussian densities and integrating gives
\begin{align*}
 \Aff\{N(0,C_1),N(0,C_2)\}
 &=(\det C_1\det C_2)^{-1/4}(2\pi)^{-n/2}
       \int\exp\{-\tfrac14x^\T(P+Q)x\}\,dx\\*
 &=(\det C_1\det C_2)^{-1/4}
                       \det\{(P+Q)/2\}^{-1/2}\\*
 &=\frac{\det(C_1)^{1/4}\det(C_2)^{1/4}}
                    {\det\{(C_1+C_2)/2\}^{1/2}}.
\end{align*}
The last equality follows from
$P+Q=C_1^{-1}(C_1+C_2)C_2^{-1}$ and multiplicativity of the determinant. Factoring out $C_1$ and diagonalizing $C_1^{-1/2}C_2C_1^{-1/2}$ then yields
\[
 \Aff\{N(0,C_1),N(0,C_2)\}
 =\prod_{k=1}^n
        \left(\frac{2\sqrt{\lambda_k}}{1+\lambda_k}\right)^{1/2}
 =\prod_{k=1}^n
        \cosh\!\left(\tfrac12\log\lambda_k\right)^{-1/2}.
\]
Taking negative logarithms proves \eqref{eq:gaussian-affinity-formula}.

To verify \eqref{eq:logcosh-comparison}, let $g(x)=\log\cosh x$. It is even, $g(0)=g'(0)=0$, and $g''(x)=1/\cosh^2x$. For $0\le x\le1$, integration twice gives
\[
 g(x)=\int_0^x(x-u)g''(u)\,du,
 \qquad \frac{x^2}{2\cosh^2 1}\le g(x)\le\frac{x^2}{2}.
\]
For $x\ge1$, convexity and $g(0)=0$ give $g(x)\ge xg(1)$, while $\cosh x\le e^x$ gives $g(x)\le x$. Evenness handles negative $x$. These estimates prove the claimed
comparison on all of $\R$, including the two different growth regimes.

For the precision bound, use $F(P)=-\log\det P$. The determinant formula just proved can be written as
\[
 -\log\Aff\{N(0,C_1),N(0,C_2)\}
 =\frac12\left\{\frac{F(P)+F(Q)}2-F\!\left(\frac{P+Q}2\right)\right\}.
\]
The assumption gives $P,Q\preceq t_-^{-1}I$, and this inequality holds throughout their line segment. For a symmetric increment $U$,
\[
 D^2F(P)[U,U]=\tr(P^{-1}UP^{-1}U)
 =\|P^{-1/2}UP^{-1/2}\|_\Frob^2
 \ge t_-^2\|U\|_\Frob^2.
\]
Thus $F(P)-(t_-^2/2)\|P\|_\Frob^2$ is convex on this set. Applying its midpoint inequality gives
\begin{align*}
 \frac{F(P)+F(Q)}2-F\!\left(\frac{P+Q}2\right)
 &\ge\frac{t_-^2}{2}
       \left\{\frac{\|P\|_\Frob^2+\|Q\|_\Frob^2}{2}
                    -\left\|\frac{P+Q}2\right\|_\Frob^2\right\}\\*
 &=\frac{t_-^2}{8}\|P-Q\|_\Frob^2.
\end{align*}
Multiplying by $1/2$ proves \eqref{eq:precision-affinity}. \qed

\subsection{Proof of Lemma~\ref{lem:affinity-processing}} This is the data-processing inequality for R\'enyi divergence of order $1/2$, since $D_{1/2} P\|Q)=-2\log\Aff(P,Q)$; see \citet[Theorem~9 and Example~2]{vanerven2014renyi}.
Here is the direct conditional Cauchy--Schwarz proof.

Let $R=(P+Q)/2$, and write $p=dP/dR$, $q=dQ/dR$. Generate $X\sim R$ and then an output $Y$ with conditional law $\mathsf K(X,\cdot)$. When an expectation involves both $X$ and $Y$, $\E_R$ denotes their joint law $R(dx)\mathsf K(x,dy)$; define $\E_P$ similarly. Let $\widetilde p(Y)=\E_R[p(X)\mid Y]$. For any bounded measurable
function $f$ of the output, we have
\begin{align*}
 \E_{P\mathsf K}[f(Y)]
 &=\E_P[f(Y)]\\
 &=\int\!\int f(y)\,\mathsf K(x,dy)\,p(x)\,R(dx)\\
 &=\E_R[f(Y)p(X)]\\
 &=\E_R\!\left[f(Y)\E_R[p(X)\mid Y]\right]\\
 &=\E_{R\mathsf K}[f(Y)\widetilde p(Y)].
\end{align*}
The penultimate equality is the defining conditional-expectation identity, because $f(Y)$ is measurable with respect to $Y$. The last equality uses that the integrand depends only on $Y$, whose law under the joint reference measure is $R\mathsf K$.
Since this holds for every such $f$, $\widetilde p$ is the density of $P\mathsf K$ with respect to $R\mathsf K$. Repeating the calculation with $q$ gives
\[
 \frac{d(P\mathsf K)}{d(R\mathsf K)}(Y)=\E_R[p(X)\mid Y],
 \qquad
 \frac{d(Q\mathsf K)}{d(R\mathsf K)}(Y)=\E_R[q(X)\mid Y].
\]
Conditional Cauchy--Schwarz gives
\[
 \E_R[\sqrt{p(X)q(X)}\mid Y]
 \le\sqrt{\E_R[p(X)\mid Y]\E_R[q(X)\mid Y]}.
\]
Taking expectations yields
\begin{align*}
 \Aff(P,Q)
 &=\E_R\sqrt{p(X)q(X)}\\*
 &\le\E_R\sqrt{\E_R[p(X)\mid Y]\E_R[q(X)\mid Y]}
   =\Aff(P\mathsf K,Q\mathsf K).
\end{align*}
The definition of affinity can use any common dominating measure, including $R$; changing that measure cancels inside the square root. This proves \eqref{eq:affinity-processing}. Taking $Y=T(X)$ gives the statistic assertion. Taking $Y=X+Z$, with independent $Z\sim N(0,H)$ under both input laws, proves \eqref{eq:affinity-common-noise}. This also allows singular $H$: no Lebesgue density is required. \qed

\subsection{Proof of Lemma~\ref{lem:local-fractional}}
\noindent\textbf{A small likelihood power has a decaying expectation.}
Put $A=M^{-1}-I_n$ and $Q(A)=X^\T AX-\tr A$. Then
\[
 \mathcal L_X(M)=-D(M)-\tfrac12Q(A),\qquad
 D(M)=\tfrac12\{\tr M^{-1}-n+\log\det M\}.
\]
To compare the deterministic term with the covariance distance, write $\lambda_1,\ldots,\lambda_n\in[b,B]$ for the eigenvalues of $M$. The scalar function
$\phi(\lambda)=\tfrac12(\lambda^{-1}-1+\log\lambda)$ is positive except at $1$. Taylor expansion at $1$ gives $\phi(\lambda)=(\lambda-1)^2/4+O(|\lambda-1|^3)$. Thus $\phi(\lambda)/(\lambda-1)^2$, extended at $1$ by $1/4$, has positive lower and finite upper bounds on $[b,B]$. Consequently,
\begin{equation}\label{eq:local-contrast-comparison}
 D(M)=\sum_j\phi(\lambda_j)\asymp\sum_j(\lambda_j-1)^2
       =\|M-I_n\|_\Frob^2.
\end{equation}
Also
\[
 \|A\|_\Frob^2
   =\sum_j\frac{(\lambda_j-1)^2}{\lambda_j^2}
   \asymp\|M-I_n\|_\Frob^2,
 \qquad \|A\|_\op\le 1+b^{-1}.
\]
The exponential-moment bound in Lemma~\ref{lem:gaussian}, applied with $t=-s/2$, therefore gives
\[
 \begin{aligned}
 \log \E e^{s\mathcal L_X(M)}
 &=-sD(M)+\log \E e^{-(s/2)Q(A)}\\
 &\le-sD(M)+Cs^2\|A\|_\Frob^2\\
 &\le(-c_2s+C_2s^2)\|M-I_n\|_\Frob^2.
 \end{aligned}
\]
Choose $s_*>0$ so that the exponential-moment bound applies for $0\le s\le s_*$ and $C_2s_*\le c_2/2$. This proves the first assertion with $c=c_2/2$. Choosing a small fixed power makes the negative contrast dominate the centered quadratic fluctuation.

\noindent\textbf{The likelihood power is Lipschitz in $L^m$.}
Fix $M,N$ in the stated covariance interval. Interpolate the matrices:
\[
 \Delta=N-M,\qquad M_t=M+t\Delta,\qquad R_t=M_t^{-1},
 \qquad 0\le t\le1.
\]
Every $M_t$ still satisfies $bI_n\preceq M_t\preceq BI_n$. This interpolation only uses positive definite matrices, so it requires no differentiability of a parameterization of the covariance family. Differentiating the inverse and log determinant gives
\[
 \begin{aligned}
 \frac{d}{dt}\mathcal L_X(M_t)
 &=\tfrac12X^\T R_t\Delta R_tX-\tfrac12\tr(R_t\Delta)\\
 &=\tfrac12Q(R_t\Delta R_t)
      +\tfrac12\tr\{(R_t^2-R_t)\Delta\}.
 \end{aligned}
\]
Let $d_t=\|M_t-I_n\|_\Frob$. Since $\|R_t\|_\op\le b^{-1}$,
\[
 \|R_t\Delta R_t\|_\Frob\le b^{-2}\|\Delta\|_\Frob,
 \qquad
 \|R_t^2-R_t\|_\Frob
 =\|R_t(I_n-M_t)R_t\|_\Frob\le b^{-2}d_t.
\]
The quadratic-form moment bound in Lemma~\ref{lem:gaussian} and Frobenius Cauchy--Schwarz therefore imply
\[
 \left\|\frac{d}{dt}\mathcal L_X(M_t)\right\|_{L^{2m}}
 \le C_m(1+d_t)\|\Delta\|_\Frob.
\]
Now put $\xi_t=e^{s\mathcal L_X(M_t)/m}$. H\"older's inequality and the first assertion, applied with $2s\le s_*$, give
\[
 \begin{aligned}
 \|\dot\xi_t\|_{L^m}
 &\le\frac{s}{m}
      (\E e^{2s\mathcal L_X(M_t)})^{1/(2m)}
      \left\|\frac{d}{dt}\mathcal L_X(M_t)\right\|_{L^{2m}}\\
 &\le C_{m,s}e^{-csd_t^2/m}(1+d_t)\|\Delta\|_\Frob
 \le C'_{m,s}\|\Delta\|_\Frob.
 \end{aligned}
\]
The last step uses the finite supremum of $(1+x)e^{-csx^2/m}$ over $x\ge0$. Thus the decay of the likelihood power also controls the deterministic part of its derivative.
The fundamental theorem of calculus, followed by Minkowski, now yields
\[
 \|\xi(N)-\xi(M)\|_{L^m}
 \le\int_0^1\|\dot\xi_t\|_{L^m}\,dt
 \le C'_{m,s}\|N-M\|_\Frob.
\]
Finally, substitute the covariance separation and Lipschitz bounds:
\[
 \E|\xi(u)|^m
 \le e^{-cs\|M(u)-I_n\|_\Frob^2}
 \le e^{csC_0}e^{-csc_0\|h\|^2},
\]
and apply the preceding increment estimate to $M(u),M(v)$ 
Continuity follows from continuity of the inverse and determinant on positive definite matrices. This proves all the assertions. \qed

\subsection{Proof of Lemma~\ref{lem:compact-supremum}}
The proof balances the maximum on a finite mesh against a modulus of continuity bound.

If $M=0$, the process vanishes almost surely on a countable dense subset of $\Gamma$, and hence everywhere by continuity. Assume $M>0$. Fix $a\in\Gamma$, and let $\pi$ denote nearest-point projection onto $\Gamma$. It is $1$-Lipschitz: if $u=\pi(x)$ and $v=\pi(y)$, minimizing squared distance along segments in $\Gamma$ gives $\langle x-u,v-u\rangle\le0$ and $\langle y-v,u-v\rangle\le0$. Adding the resulting inequalities
in the direction $u-v$ yields
\[
 \|u-v\|^2\le\langle u-v,x-y\rangle
               \le\|u-v\|\|x-y\|.
\]
Existence follows from compactness, and this inequality also gives uniqueness. Since $\Gamma\subset a+[-\delta,\delta]^k$, the process
\[
 \eta(x)=\xi\bigl(\pi(a-\delta\mathbf1+2\delta x)\bigr),
 \qquad x\in[0,1]^k,
\]
has the same supremum as $\xi$ on $\Gamma$. Its paths are continuous, its pointwise $L^m$ norms are at most $M$, and
\[
 \|\eta(x)-\eta(y)\|_{L^m}\le L'\|x-y\|,
 \qquad L'=2L\delta.
\]

\noindent\textbf{The modulus of continuity on the cube.}
If $L'=0$, continuity and the increment bound make $\eta$ constant almost surely, so the following modulus bound is immediate. Otherwise apply \citet[Section~1.3, equation~(1.10), p.~4]{talagrand2021upper} to $\eta/L'$. His index dimension is our $k$, and his moment and increment exponents are both our $m$: indeed,
\[
 \E\left|\frac{\eta(x)-\eta(y)}{L'}\right|^m
       \le\|x-y\|^m,\qquad m>k.
\]
His estimate on the dense dyadic grid extends to the whole cube by continuity and gives
\[
 \left\|\sup_{\|x-y\|\le2^{-J}}|\eta(x)-\eta(y)|\right\|_{L^m}
 \le C_{k,m}L'2^{-J(1-k/m)}.
\]
For $0<r\le1$, choose $2^{-J-1}<r\le2^{-J}$. Then
\[
 \left\|\sup_{\|x-y\|\le r}|\eta(x)-\eta(y)|\right\|_{L^m}
 \le C_{k,m}L\delta\,r^{1-k/m}.
\]

\noindent\textbf{Choosing the mesh.}
Cover the unit cube by at most $(2/r)^k$ cubes of side at most $r$, choosing one point $x_B$ in each. Then
\[
 \sup_x|\eta(x)|\le\max_B|\eta(x_B)|
             +\sup_{\|x-y\|\le\sqrt{k}r}|\eta(x)-\eta(y)|.
\]
The first term has $L^m$ norm at most $(2/r)^{k/m}M$. For the second, divide each segment into at most $\lceil\sqrt{k}\rceil$ pieces of length at most $r$ and use the
modulus bound. Thus
\[
 \left\|\sup_{u\in\Gamma}|\xi(u)|\right\|_{L^m}
 \le C_{k,m}\{Mr^{-k/m}+L\delta\,r^{1-k/m}\}.
\]
If $M\ge L\delta$, take $r=1$, obtaining $C_{k,m}M$. Otherwise take $r=M/(L\delta)$, obtaining $C_{k,m}M^{1-k/m}(L\delta)^{k/m}$. Raising these bounds to the power $m$ proves the result. All suprema are measurable by continuity on compact sets. \qed

\subsection{Proof of Lemma~\ref{lem:quadratic-envelope}}
Cover the $h$ coordinates by closed unit cubes $j+[0,1]^p$, $j\in\Z^p$. If $a$ is present, cover $I$ by at most $2+\operatorname{diam}(I)$ closed intervals of length at most $1$ and take products. Each resulting cube or rectangle $Q$ has diameter at most $\sqrt{k}$. The nonempty intersections $\Gamma_Q=Q\cap\mathcal U$ are compact and convex. Put
\[
 r_Q=\inf_{(h,a)\in\Gamma_Q}\norm h,\qquad
 M_Q=A^{1/m}e^{-c_0r_Q^2/m},
\]
with $a$ omitted in the case without a nuisance. Lemma~\ref{lem:compact-supremum} gives
\begin{align*}
 \E\sup_{u\in\Gamma_Q}|\xi(u)|^m
 &\le C_{k,m}\{M_Q^m+L^kM_Q^{m-k}\}\\
 &\le C_{k,m}\{A+L^kA^{1-k/m}\}e^{-\beta r_Q^2},
 \qquad \beta=c_0(1-k/m)>0.
\end{align*}
This is where the strict inequality $m>k$ preserves exponential decay after taking a supremum within a cube.

Only cubes meeting $\{\norm h>H\}\cap\mathcal U$ need be retained. If $j$ is the $h$-corner of such a cube, then
\[
 \norm j>H-\sqrt p,\qquad
 r_Q\ge(\norm j-\sqrt p)_+,\qquad
 r_Q^2\ge\tfrac12\norm j^2-p.
\]
There are only a bounded number of nuisance intervals for each $j$. The supremum over a union is bounded by the sum of the nonnegative suprema over its pieces. Consequently, for $H\ge2\sqrt p$,
\begin{align*}
 \E\sup_{\substack{u\in\mathcal U\\\norm h>H}}|\xi(u)|^m
 &\le C\{A+L^kA^{1-k/m}\}
       \sum_{\substack{j\in\Z^p\\\norm j>H-\sqrt p}}
                 e^{-(\beta/2)\norm j^2}\\
 &\le C\{A+L^kA^{1-k/m}\}
       e^{-\beta(H-\sqrt p)^2/4}
       \sum_{j\in\Z^p}e^{-(\beta/4)\norm j^2}\\
 &\le C'\{A+L^kA^{1-k/m}\}e^{-\beta H^2/16}.
\end{align*}
The lattice sum is finite: it is the $p$th power of $\sum_{r\in\Z}e^{-(\beta/4)r^2}$.
For $0\le H<2\sqrt p$, sum over all cubes and enlarge the constant. This proves the claim, for example with $c_1=\beta/16$.
\qed

\section{Auxiliary calculations and interpretations}\label{app:auxiliary}
\subsection{Expected curvature and sampling covariance}
Under non-Gaussian errors the Gaussian criterion still has the correct expected covariance contrast. At the truth, its expected negative Hessian is the Gaussian curvature:
\[
 \E_{\theta_0,P}\{-\nabla_\theta^2\ell_n(\theta_0)\}
       =\mathcal I_n^{\mathrm{full}}(\theta_0).
\]
Away from the truth, the expected Hessian under the fixed true law also contains the covariance-mismatch term displayed in Lemma~\ref{lem:gaussian-hessian}. At the truth, the difference between the score covariance and the expected negative Hessian is exactly the fourth-cumulant term in \eqref{eq:quadratic-cumulant}, applied to the score matrices. Lemma~\ref{lem:noise-diagonals} gives its asymptotic form in normalized coordinates. Thus curvature alone need not give the sampling covariance required for Wald inference. If $\E_P\xi^4=3$, the score covariance and expected negative Hessian agree exactly at the truth for every $n$, even when the errors are non-Gaussian. However, this equality does not make the Gaussian criterion the true
log-likelihood.

The covariance structure in Theorem~\ref{thm:normalized} also has an interpretation for the fitted errors. Its moment comparison gives $\Cov(\widehat h_n)=V_n+o(1)$ at fixed locations, and conditionally on almost every random-design sequence. Together with \eqref{eq:covariance-blocks}, this shows that the normalized nugget error is asymptotically uncorrelated with every signal error. At $d=4$, \eqref{eq:critical-covariance-block} gives the same conclusion for the normalized range error and the centered amplitude and smoothness errors. Their variances stay bounded away from zero,
so the corresponding correlations also tend to zero.

\subsection{Examples and the nugget nondegeneracy condition}
If every standardized error law is supported on the same interval $[-M,M]$, the exponential-moment bound in Assumption~\ref{ass:noise} holds with $K_\varepsilon=M/\sqrt{\log 2}$. The standard Gaussian may be included by increasing $K_\varepsilon$. More generally, the support need not be bounded: a common sub-Gaussian
bound suffices for this part of the assumption, and mixtures of laws satisfying the same exponential-moment bound inherit it by averaging. The condition $\Var(\xi^2)=\E\xi^4-1\ge v_*>0$ is separate: for symmetric errors taking only $-1$ and $1$, $\xi^2=1$ almost surely. Such errors still satisfy upper concentration and moment inequalities, but the leading root-$n$ nugget variance is degenerate and they are excluded from the uniformly nondegenerate joint normal theorem.

\subsection{Full inversion and identifiable-block inversion}
\label{app:full-inverse-block}
In low dimension, order a full parameter vector as $(\zeta,\rho)$, where $\zeta$ is identifiable and $\rho$ is the range nuisance. At the same parameter and locations, write the positive definite full Gaussian Fisher information as
\[
 \mathcal I=\begin{pmatrix}A&B\\B^\top&D\end{pmatrix},
 \qquad S=A-BD^{-1}B^\top.
\]
The Schur complement $S$ is the information remaining for $\zeta$ after accounting for the nuisance. Since $\mathcal I\succ0$, its principal block $D$ is positive definite, and, for every $x\ne0$,
\[
 x^\top Sx
 =\begin{pmatrix}x\\-D^{-1}B^\top x\end{pmatrix}^{\!\top}
   \mathcal I
   \begin{pmatrix}x\\-D^{-1}B^\top x\end{pmatrix}>0.
\]
The block inversion formula therefore gives
\[
 \mathcal I^{-1}
 =\begin{pmatrix}
 S^{-1}&-S^{-1}BD^{-1}\\
 -D^{-1}B^\top S^{-1}
 &D^{-1}+D^{-1}B^\top S^{-1}BD^{-1}
 \end{pmatrix}.
\]
In particular,
\[
 [\mathcal I^{-1}]_{\zeta\zeta}
 =(A-BD^{-1}B^\top)^{-1}=S^{-1}.
\]
Since $BD^{-1}B^\top\succeq0$, we have $0\prec S\preceq A$. Congruencing by $A^{-1/2}$ and inverting gives
\[
 A^{-1/2}SA^{-1/2}\preceq I
 \quad\Longrightarrow\quad
 A^{1/2}S^{-1}A^{1/2}\succeq I
 \quad\Longrightarrow\quad
 [\mathcal I^{-1}]_{\zeta\zeta}\succeq A^{-1}.
\]
Thus full inversion gives marginal information-based variances at least as large as direct inversion of the identifiable block. The latter corresponds to treating the nuisance as known. This is an ordering of the two information-based covariance calculations; it is not an identity for the finite-sample covariance of the MLE.
Adding the same fourth-cumulant correction to both calculations preserves their difference and hence this ordering.

The normalized Schur correction is negligible by Lemma~\ref{lem:studentization}, although its finite-sample effect can be substantial, as in Appendix~\ref{apdx:sim-range}. The same algebraic ordering holds for a positive definite observed negative Hessian. A constrained fit need not have such a Hessian.

\subsection{Failure of joint linear normality on the variance--lengthscale scale}
\label{app:raw-public-failure}
The exception in Corollary~\ref{cor:log} is a failure of the claimed joint limit, not just of a particular remainder bound. Fix a deterministic design satisfying Assumption~\ref{ass:sharp} and an interior truth with $\nu_0\ne1/2$ and $4\le d\le8$. Put $\zeta=g(\vartheta)=(s,l,\nu,t)$, where $s=\sigma^2$, using the forward map from the core spectral coordinates. We examine the entire linearly transformed error vector
\[
 \mathbf Z_n
 =E_n^{-1}Dg(\vartheta_0)^{-1}(\widehat\zeta_n-\zeta_0).
\]
Under the linear approximation this vector is $\widehat h_n$. The question is whether the nonlinear remainder remains negligible after applying this normalization.

The coefficients can be read directly from the derivative rows $b_{g,n}$ in Table~\ref{tab:coordinate-expansions}. Stacking the rows for $s,l,\nu,t$, and writing
$q_{\sigma,0}=2L_n-B_0-\log a_0$, where $B_0=\partial_\nu\log A_d(\nu_0)$, gives
\[
 Dg(\vartheta_0)E_n
 =\begin{pmatrix}
 s_0/r_n&s_0q_{\sigma,0}/r_n&-s_0\nu_0/(a_0b_n)&0\\
 0&0&-l_0/(2a_0b_n)&0\\
 0&1/r_n&0&0\\
 0&0&0&1/\sqrt n
 \end{pmatrix}.
\]
Inverting this matrix gives
\[
 E_n^{-1}Dg(\vartheta_0)^{-1}
 =\begin{pmatrix}
 r_n/s_0&-2\nu_0r_n/l_0&-r_nq_{\sigma,0}&0\\
 0&0&r_n&0\\
 0&-2a_0b_n/l_0&0&0\\
 0&0&0&\sqrt n
 \end{pmatrix}.
\]
In particular, the first row is obtained by solving the $b_{s,n}$ row for $h_u$, using the $b_{l,n}$ and $b_{\nu,n}$ rows to eliminate $h_a$ and $h_v$. Since $-q_{\sigma,0}=B_0-2\log l_0-2L_n$, the vector is exactly
\[
 \mathbf Z_n
 =\begin{pmatrix}
 Z_n\\
 \widehat h_{n,v}\\
 -2a_0b_n(\widehat l-l_0)/l_0\\
 \widehat h_{n,t}
 \end{pmatrix},
 \qquad
 Z_n=r_n\left\{
 \frac{\widehat s-s_0}{s_0}
 -2\nu_0\frac{\widehat l-l_0}{l_0}
 -q_{\sigma,0}(\widehat\nu-\nu_0)\right\}.
\]

We now retain the quadratic terms in the coordinate expansions. Let $\mathcal R_{s,n}$ and $\mathcal R_{l,n}$ be the Taylor remainders defined in \eqref{eq:coordinate-taylor-recipe}. The Hessians displayed in the proof of Corollary~\ref{cor:rates}, together with $\widehat h_n=O_P(1)$, give
\[
 \begin{aligned}
 \frac{\mathcal R_{s,n}}{s_0}
 &=\frac{\nu_0(\nu_0+1)}{2a_0^2b_n^2}\widehat h_{n,a}^2
   +O_P\left(\frac{L_n}{r_nb_n}
              +\frac{L_n^2}{r_n^2}+\frac1{b_n^3}\right),\\
 \frac{\mathcal R_{l,n}}{l_0}
 &=\frac{3}{8a_0^2b_n^2}\widehat h_{n,a}^2+O_P(b_n^{-3}).
 \end{aligned}
\]
Here the $aa$ entries of the Hessians give the displayed terms. The mixed amplitude--range and smoothness--range entries contribute $O_P\{L_n/(r_nb_n)\}$, and the amplitude--smoothness block contributes $O_P(L_n^2/r_n^2)$. Evaluating the $aa$ coefficient along the Taylor segment instead of at the truth adds terms bounded
by these orders and $O_P(b_n^{-3})$. Since $\nu$ and $t$ are linear coordinates, their remainders are zero. Consequently
\[
 \begin{aligned}
 Z_n
 &=\widehat h_{n,u}
   +r_n\left\{\frac{\mathcal R_{s,n}}{s_0}
              -2\nu_0\frac{\mathcal R_{l,n}}{l_0}\right\}\\
 &=\widehat h_{n,u}
   +\gamma_0\frac{r_n}{b_n^2}\widehat h_{n,a}^2
   +O_P\left(\frac{L_n}{b_n}
             +\frac{L_n^2}{r_n}+\frac{r_n}{b_n^3}\right),\\
 -2a_0b_n\frac{\widehat l-l_0}{l_0}
 &=\widehat h_{n,a}
   -\frac{3}{4a_0b_n}\widehat h_{n,a}^2+O_P(b_n^{-2}),
 \qquad
 \gamma_0=\frac{\nu_0(2\nu_0-1)}{4a_0^2}\ne0.
 \end{aligned}
\]
Write $c_n=r_n/b_n^2$. Since $b_n\to\infty$ and $L_nb_n/r_n\to0$, the first remainder is $o_P(1+c_n)$. Thus the full-vector expansion is
\[
 \mathbf Z_n
 =\widehat h_n+\gamma_0c_n\widehat h_{n,a}^2e_u
   +o_P(1+c_n)e_u+o_P(1)e_a.
\]
Only the centered-amplitude coordinate has a potentially nonvanishing error; the other three coordinates agree with $\widehat h_n$ up to $o_P(1)$.

This vector directly tests joint linear standardization. Indeed,
\[
 \Sigma_n^\zeta
 =Dg(\vartheta_0)E_nV_nE_n^\top Dg(\vartheta_0)^\top,
\]
so $V_n^{-1/2}\mathbf Z_n$ is the reporting error multiplied by a valid covariance-standardizing factor. Any other such factor differs by an orthogonal matrix. Since $V_n\asymp I$, tightness of the standardized error requires tightness of $\mathbf Z_n$.

For $d=4$, $c_n=R_n^2/\log n\to\infty$; for $5\le d\le7$, $c_n=R_n^{4-d/2}\to\infty$.
Along any subsequence with $V_n\to V$, Theorem~\ref{thm:normalized} gives $\widehat h_n\Rightarrow H\sim N(0,V)$ with $V\succ0$. The full-vector expansion therefore gives
\[
 \frac{\mathbf Z_n}{c_n}\Rightarrow\gamma_0H_a^2e_u,
 \qquad \frac{Z_n}{c_n}\Rightarrow\gamma_0H_a^2.
\]
Because $V_{aa}>0$ and $\gamma_0\ne0$, this scalar limit has no atom at zero. For any fixed $M>0$ and any $\delta>0$, eventually $M/c_n\le\delta$, so the closed-set bound for convergence in distribution gives
\[
 \limsup_n\PP(|Z_n|\le M)
 \le\limsup_n\PP(|Z_n/c_n|\le\delta)
 \le\PP(|\gamma_0H_a^2|\le\delta).
\]
Letting $\delta\downarrow0$ makes the right side zero. Hence $\PP(|Z_n|>M)\to1$ for every fixed $M$: $Z_n$ is not $O_P(1)$, and $\|\mathbf Z_n\|\to\infty$ in probability.
These conclusions hold on the full sequence, because every subsequence has a further subsequence on which the bounded, uniformly positive definite matrices $V_n$ converge.
Thus joint linear standardization cannot have the limit in \eqref{eq:joint-reporting-normal} in dimensions four through seven.

At $d=8$, $c_n=1$, and the whole vector instead has the subsequential limit
\[
 \mathbf Z_n\Rightarrow H+\gamma_0H_a^2e_u.
\]
Its mean is $\gamma_0V_{aa}e_u\ne0$, whereas joint linear normality would require the limit $N(0,V)$. The limiting vector has finite moments, so these two laws cannot agree. Thus the claimed joint limit is generally false in dimension eight as well, even
with Gaussian errors.

Using the fitted curvature does not remove the quadratic term. To see why its sign changes, put $\delta_n=\widehat\vartheta_n-\vartheta_0=E_n\widehat h_n$. Taylor expansion backwards from the fit gives the exact identity
\[
 \widehat\zeta_n-\zeta_0
 =Dg(\widehat\vartheta_n)\delta_n
  -\int_0^1 u\,D^2g(\vartheta_0+u\delta_n)
                    [\delta_n,\delta_n]\,du.
\]
Multiplying by $E_n^{-1}Dg(\widehat\vartheta_n)^{-1}$ gives $\widehat h_n$ from the linear term. In the first coordinate of the integral, the $aa$ entries of the variance and lengthscale Hessians give a coefficient of $c_n\widehat h_{n,a}^2$ equal to
\[
 \int_0^1 u\left\{
       \frac{\nu_0(\nu_0+1)}{a_0^2}
       -2\nu_0\frac{3}{4a_0^2}\right\}\,du+o_P(1)
 =\gamma_0+o_P(1).
\]
Here the fitted parameter and every point on the Taylor segment converge to $\vartheta_0$, so the coefficient can be evaluated at the truth up to $o_P(1)$. The mixed terms and the amplitude--smoothness terms have the same respective bounds $O_P(L_n/b_n)$ and $O_P(L_n^2/r_n)$ as in the forward expansion. The minus sign in front of the integral therefore gives
\[
 E_n^{-1}Dg(\widehat\vartheta_n)^{-1}(\widehat\zeta_n-\zeta_0)
 =\widehat h_n-\gamma_0c_n\widehat h_{n,a}^2e_u
   +o_P(1+c_n)e_u+o_P(1)e_a.
\]
In these dimensions the fit is interior with probability tending to one. On that event the fitted curvature transforms exactly by congruence, as shown in the proof of Corollary~\ref{cor:uncertainty}. Together with the covariance approximation in the core coordinates, \eqref{eq:covariance-substitution}, this expresses the feasible Wald statistic as the quadratic form of the vector just displayed with inverse covariance $(V_n+o_P(1))^{-1}$. It diverges in probability for $4\le d\le7$. At $d=8$ its subsequential limit is
\[
 (H-\gamma_0H_a^2e_u)^\top V^{-1}(H-\gamma_0H_a^2e_u),
\]
whose expectation is $k_d+3\gamma_0^2(V^{-1})_{uu}V_{aa}^2>k_d$, using the vanishing third moments of a centered Gaussian vector. It is therefore not $\chi^2_{k_d}$. Marginal Wald inference is unaffected, since its normalization uses only each coordinate's own rate. For $d>8$, the quadratic term vanishes. The transformation bounds in Lemma~\ref{lem:logtransform} then give the joint normal limit by Corollary~\ref{cor:log} and relative consistency of the fitted covariance by Corollary~\ref{cor:uncertainty}, also in the original variance--lengthscale coordinates.

\section{Simulation methods and additional results}\label{apdx:sim}
This appendix gives the computation and numerical checks behind Section~\ref{sec:simulation}, followed by detailed performance summaries and supplementary figures.
\renewcommand{\thefigure}{E\arabic{figure}}
\renewcommand{\theHfigure}{E.\arabic{figure}}
\setcounter{figure}{0}
\renewcommand{\thetable}{E\arabic{table}}
\renewcommand{\theHtable}{E.\arabic{table}}
\setcounter{table}{0}

\subsection{Exact covariance and information calculations}
\subsubsection{Exact evaluation of the likelihood}\label{apdx:sim-lik}

The sample sizes in \Cref{sec:simulation} are too large for direct $O(n^3)$ evaluation; at $n=10^5$ a dense covariance matrix alone requires $80$ gigabytes. The isotropic Mat\'ern covariance does not factorize across grid coordinates as the RBF covariance does. We instead use the symmetry of the regular grid, which yields the same
Gaussian likelihood \eqref{eq:likelihood} up to floating-point rounding.

\noindent\textbf{Symmetry of the grid covariance.}
Index the midpoint grid by $i=(i_1,\ldots,i_d)\in\{0,\ldots,G-1\}^d$. The entry of $C_n(\theta)$ for $i,i'$ depends only on $\sum_k(i_k-i'_k)^2$. Hence $C_n(\theta)$ commutes, for every $\theta$, with the permutation matrices of the hyperoctahedral group generated by the reflections $i_k\mapsto G-1-i_k$ of each axis and the permutations of the axes. We construct a fixed orthogonal matrix $Q$, independent of
$\theta$, such that $Q^\T C_n(\theta)Q$ is block diagonal, in two steps.\begin{enumerate}[label=(\roman*)]
\item \emph{Reflections.} On each axis, the even and odd vectors $(e_j\pm e_{G-1-j})/\sqrt2$, $0\le j<\lfloor G/2\rfloor$, together with the central vector $e_{(G-1)/2}$ when $G$ is odd, form an orthonormal basis. Their tensor products split $C_n$ into $2^d$ parity blocks, indexed by $\epsilon\in\{0,1\}^d$, with $\epsilon_q=1$
denoting an odd axis.

To specify the folded entries, let $k_\theta(D)$ be the latent covariance at physical distance $\sqrt D/G$. For representative multi-indices $i,j$ in a parity block, its latent covariance entry is
\[
 w_iw_j\sum_{c\in\{0,1\}^d}
 (-1)^{\sum_q\epsilon_qc_q}
 k_\theta\!\left(\sum_{q=1}^d\delta_q^2\right),
 \qquad
 \delta_q=\begin{cases}
 i_q-j_q,&c_q=0,\\
 G-1-i_q-j_q,&c_q=1,
 \end{cases}
\]
where $w_i=2^{-m(i)/2}$ and $m(i)$ counts the central coordinates of $i$ on an odd grid; these occur only on even axes. For even $G$, $m(i)=0$. The nugget adds $tI$ to each parity block. In one dimension, the paired even and odd entries are $k_\theta((i-j)^2)\pm k_\theta((G-1-i-j)^2)$. An even pair against the central site $j_*=(G-1)/2$ gives $\sqrt2\,k_\theta((j_*-i)^2)$, and the central site against itself gives $k_\theta(0)$. Thus each latent parity-block entry is a weighted signed sum of at most $2^d$ distinct kernel values. The axes remain coupled through the total squared distance.
\item \emph{Permutations.} Patterns with the same number $k$ of odd axes give identical blocks up to a permutation of indices, so only $d+1$ blocks are distinct. Within the block with the first $k$ axes odd, the stabilizer $S_k\times S_{d-k}$ permutes multi-indices. Using its irreducible representations, constructed from Young's orthogonal form, a compatible orthonormal basis splits this block further into one block $B_\rho(\theta)$ for each irreducible representation $\rho$ that occurs, repeated $\binom dk\dim\rho$ times. The entries of these further blocks are linear combinations of the parity-block entries above.
\end{enumerate}
Each basis vector of $Q$ is supported on a single orbit of the group, which keeps $Q$ sparse. Because $C_n(\theta)$ commutes with the group, every block $B_\rho(\theta)$ can be assembled from the rows of the parity block at one representative of each orbit, without forming the parity block itself. \Cref{tab:sim-blocks} lists the resulting block sizes.

\begin{table}[htbp]
\centering\small
\caption{Symmetry blocks of the grid covariance at the largest sample sizes. The multiplicity is the number of parity patterns times the dimension of the irreducible representation. Distinct copies have independent data under Gaussian errors. Under non-Gaussian errors the transformed responses are uncorrelated but need not be independent.}
\label{tab:sim-blocks}
\vspace{3pt}
\begin{tabular}{rrrrr}
\toprule
$d$ & $G$ & $n$ & number of distinct blocks & largest block\\
\midrule
$3$ & $40$ & $64{,}000$ & $10$ & $4200$\\
$3$ & $46$ & $97{,}336$ & $10$ & $6348$\\
$4$ & $18$ & $104{,}976$ & $20$ & $2160$\\
$5$ & $10$ & $100{,}000$ & $36$ & $600$\\
\bottomrule
\end{tabular}
\end{table}

\noindent\textbf{Blockwise likelihood.}
Write $Z_\rho$ for the matrix whose columns are the components of $Q^\T Y_n$ belonging to the copies of block $\rho$, the log-likelihood is
\[
 \ell_n(\theta)
 =
 -\frac12\sum_\rho
 \Bigl\{
 m_\rho\log\det B_\rho(\theta)
 +\tr\bigl(B_\rho(\theta)^{-1}Z_\rho Z_\rho^\T\bigr)
 \Bigr\},
\]
up to the constant $-n\log(2\pi)/2$, where $m_\rho$ is the multiplicity of the block. Each evaluation requires one Cholesky factorization per distinct block. The score uses the derivatives of the kernel in the coordinates $(\eta,\nu,\log a,\log t)$ of \eqref{eq:chart}. With $z=r\sqrt a$ and $M_\nu(z)=2^{1-\nu}z^\nu\mathsf K_\nu(z)/\Gamma(\nu)$,
\[
 z\,\partial_zM_\nu(z)
 =
 -M_\nu(z)\,\frac{z\,\mathsf K_{\nu-1}(z)}{\mathsf K_\nu(z)},
 \qquad
 \partial_\nu\log\mathsf K_\nu(z)
 =
 \frac{\int_0^\infty u\sinh(\nu u)\,e^{-z\cosh u}\,du}
      {\int_0^\infty \cosh(\nu u)\,e^{-z\cosh u}\,du}.
\]
The integrals are computed by a stabilized trapezoidal rule; its settings and accuracy checks are given in Appendix~\ref{apdx:sim-validation}. Kernel values are evaluated only at the at most $d(G-1)^2+1$ distinct squared lag lengths of the grid.

\noindent\textbf{Cost.}
On one core of an AMD EPYC 7713 processor, a single evaluation of the log-likelihood and its gradient takes about $2.3$ seconds for $d=5$ ($n=10^5$), $9.3$ seconds for $d=4$ ($n=104{,}976$), and $45$ seconds for $d=3$ ($n=97{,}336$); with four threads it takes $0.74$, $2.8$, and $13.8$ seconds. For $d=3$ the time is dominated by the Cholesky factorization and inversion of the $6348$-row blocks. The block assembly is compiled with Numba and runs in parallel over rows.

\subsubsection{Computation of the Fisher information}
\label{apdx:sim-information}

Write $\mathcal I_n^{\mathrm{full}}$ for the full four-parameter Gaussian information in $w=(\eta,\nu,\log a,\log t)$. This distinguishes it from the three-coordinate information used to state the theorem when $d\le3$. It is computed by applying the trace formula \eqref{eq:info} to all four coordinates, using the same blocks, with $\partial_j=\partial/\partial w_j$:
\[
 [\calI_n^{\mathrm{full}}]_{jk}
 =
 \frac12\sum_\rho m_\rho
 \tr\bigl(B_\rho^{-1}\partial_jB_\rho\,B_\rho^{-1}\partial_kB_\rho\bigr),
\]
and is mapped to the reported parameters by the delta method.

Under a non-Gaussian error law the covariance of the score is also exact at finite $n$. Write $A_j=C_n^{-1}(\partial_jC_n)C_n^{-1}$ at $\theta_0$, so that $\partial_j\ell_n(\theta_0)=\frac12\{Y_n^\T A_jY_n-\tr(C_nA_j)\}$. Then \eqref{eq:quadratic-cumulant} gives
\[
 [\Omega_n]_{jk}
 =\Cov_{\theta_0,P}\{\partial_j\ell_n(\theta_0),\partial_k\ell_n(\theta_0)\}
 =[\calI_n^{\mathrm{full}}]_{jk}+\frac{\kappa_4(P,t_0)}4\sum_{i=1}^n(A_j)_{ii}(A_k)_{ii}.
\]
Order the columns of $Q$ by block, then by copy, then by position within the block. The diagonal reconstruction uses
\[
 Q^\T C_nQ=\bigoplus_\rho(I_{m_\rho}\otimes B_\rho).
\]
Thus
\[
 A_j=Q\left\{\bigoplus_\rho
 \left[I_{m_\rho}\otimes
 \{B_\rho^{-1}(\partial_jB_\rho)B_\rho^{-1}\}\right]\right\}Q^\T.
\]
Each distinct block is inverted once, and its contribution is added for all $m_\rho$ copies. The implementation checks $\sum_\rho m_\rho\dim B_\rho=n$ whenever a grid is constructed; this identity was verified on all $31$ simulated and information-only grids. The diagonal reconstruction is also checked against dense calculations on odd and even grids in Appendix~\ref{apdx:sim-validation}. For the four-parameter fit the benchmark in $w$ coordinates is $\mathcal I^{-1}\Omega_n\mathcal I^{-1}$, where $\mathcal I=\calI_n^{\mathrm{full}}$. When $d=3$, let $I=\{1,2,4\}$ index $(\eta,\nu,\log t)$ in $w$. We compare two benchmarks for these three coordinates:
\[
 \underbrace{[\mathcal I^{-1}\Omega_n\mathcal I^{-1}]_{II}}
   _{\text{range estimated}},
 \qquad
 \underbrace{(\mathcal I_{II})^{-1}(\Omega_n)_{II}(\mathcal I_{II})^{-1}}
   _{\text{range known}}.
\]
Both are evaluated at the same truth, including $a_0$, and the same locations. In the first, range is a fitted coordinate; in the second, the range row and column are removed before inversion. We transform each benchmark to the reported coordinates using the corresponding Jacobian as in \Cref{sec:sim-setup} and compare both with the same four-parameter fits.

Because these quantities require no optimization, they are evaluated beyond the Monte Carlo range, up to $n=97{,}336$ for $d=3$, $n=234{,}256$ for $d=4$, and $n=248{,}832$ for $d=5$; this is the extension into the grey regions of the figures for $d=4,5$. The information grids also include intermediate sample sizes at which no Monte Carlo fits were computed.

\subsection{Simulation design and fitting}
\subsubsection{Model, grids, and error distributions}\label{apdx:sim-design}
Section~\ref{sec:sim-setup} gives the full simulation setup, including the grids, four error laws, parameter bounds, shared-field design, and performance measures. The following subsections explain the choice of truth, data generation, and numerical fitting in more detail.

\subsubsection{Choice of the true parameter}\label{apdx:sim-truth}

The true parameter was chosen from the exact Fisher information alone, before the Monte Carlo study, by computing the asymptotic standard deviations for several candidate values of $(l_0,\tau_0^2)$ with $\sigma_0^2=1$. The relevant quantity is how many of the $n$ directions are dominated by the nugget, which is governed by the ratio of the high-frequency signal spectrum to the noise level. Short lengthscales (large $c=\sigma^2l^{-2\nu}$) and small nuggets push the asymptotic regime to much larger $n$. For example, for $d=3$, $\nu_0=1/2$ and $n=32{,}768$, use the identifiable information block with range held at $a_0$. The ratio of the resulting standard deviation of $\widehat\tau^2$ to its limit $\sqrt2\tau_0^2n^{-1/2}$ is $174$ for $(l_0,\tau_0^2)=(0.1,0.01)$, $83$ for $(0.2,0.01)$, $11$ for $(0.2,0.1)$, $6.6$ for $(0.4,0.1)$, and $2.6$ for $(0.4,0.5)$, and the local slopes for $\nu$ are correspondingly farther from their limits. We use $(l_0,\tau_0^2)=(0.4,0.25)$, a moderate signal-to-noise ratio $\sigma_0^2/\tau_0^2=4$ for which the asymptotic regime is approached within $n\le10^5$, at least for $\nu_0=3/2$.

\subsubsection{Data generation}\label{apdx:sim-generation}

Since $Q$ is orthogonal and does not depend on $\theta$, the criterion depends on $Y_n$ only through $Q^\T Y_n$. We therefore generate the latent field directly in the symmetry-adapted coordinates,
\[
 F_\rho=U_\rho\Lambda_\rho^{1/2}\,\Xi_\rho,
 \qquad
 U_\rho\Lambda_\rho U_\rho^\T=B_\rho(\theta_0)-t_0I,
\]
with independent standard normal matrices $\Xi_\rho$. Here $B_\rho(\theta_0)-t_0I$ is the block of the signal covariance $C_n(\theta_0)-t_0I_n$, and $U_\rho\Lambda_\rho U_\rho^\T$ is its eigendecomposition. The eigendecompositions are computed in double precision. The code replaces a negative computed eigenvalue by zero; it imposes no positive floor, jitter, or further spectral truncation. In an audit on an AMD processor with four threads, no computed eigenvalue was negative in any of the $42$ grid--smoothness settings, so this replacement was never applied in that audit. An additional check on the six largest settings found no negative eigenvalues on AMD, AMD AVX-512, or Intel AVX-512 processor configurations.

Let $S$ be the assembled signal covariance and $\widetilde S$ the covariance represented by the computed factors. On the six largest grids, one for each $(d,\nu_0)$, the same audit measured
\[
 W=C_n(\theta_0)^{-1/2}(\widetilde S-S)C_n(\theta_0)^{-1/2},
 \qquad
 \|W\|_{\mathrm{op}}\le2.2\times10^{-9},\qquad
 \|W\|_{\mathrm F}\le1.7\times10^{-8}.
\]
The operator norm is the maximum over the symmetry blocks, and the squared Frobenius norm sums their squared norms with multiplicity. These values are double-precision numerical checks, rather than certified error bounds. They support agreement at the level of numerical rounding on the audited processor configuration.

The measurement errors are generated on the grid, $\varepsilon=\tau_0\xi$ with $\xi_i$ iid from the error law, and then rotated, so that the data blocks are $Z_\rho=F_\rho+[Q^\T\varepsilon]_\rho$. Thus the block construction introduces no model approximation beyond numerical rounding. The errors must be generated on the grid: for a non-Gaussian law the components of $Q^\T\varepsilon$ are uncorrelated but not independent. The responses on the grid, which the fourth-cumulant estimator uses, are $Y_n=QF+\varepsilon$.

Replicate $r$ of $(d,\nu_0,G)$ draws $\Xi$ from one random stream and the errors for each law from another, independent stream. The field is therefore shared by the four error laws within a replicate (common random numbers), and every replicate is new.

\subsubsection{Computation of the MLE}\label{apdx:sim-mle}

The target parameter box $\Theta$ is specified in Section~\ref{sec:sim-setup}, in the coordinates $(\eta,\nu,a,t)$. We use L-BFGS-B with the exact Gaussian criterion and gradient, relative function tolerance $10^{-15}$, projected-gradient tolerance $10^{-6}$, at most $2000$ iterations and $4000$ function evaluations, and at most $50$ line-search steps per iteration. The differentiation coordinates are
\[
 w=(\eta,\nu,\log a,\log t)^\top.
\]
To improve conditioning, the optimizer uses
\[
 \phi=(\eta-2L_n\nu,\nu,\log a,\log t)^\top,
 \qquad L_n=\frac{\log n}{2\nu_0+d},
\]
and rescales each coordinate by the square root of its diagonal Fisher information at the truth. Thus the optimizer coordinates $x$ satisfy an affine relation $w=w_0+T_{\mathrm{opt}}x$. Both this preconditioning and the production initialization use the truth. The optimizer coordinate $u=\eta-2L_n\nu$ differs from the centered error in \eqref{eq:shear} by the constant $\eta_0-2L_n\nu_0$.

The exact transformed amplitude constraint would be
\[
 \eta_-\le u+2L_n\nu\le\eta_+.
\]
The implementation imposed the $\nu$, $a$, and $t$ bounds, but used the enclosing interval $\eta_--2L_n\nu_+\le u\le\eta_+-2L_n\nu_-$. Consequently the set actually searched was
\[
 \begin{split}
 \widetilde\Theta_n=\{(\eta,\nu,a,t):\;&
   \nu_-\le\nu\le\nu_+,\ a_-\le a\le a_+,\ t_-\le t\le t_+,\\
 &\eta_- -2L_n(\nu_+-\nu)\le\eta
       \le\eta_+ +2L_n(\nu-\nu_-)\}.
 \end{split}
\]
Here the endpoints are the simulation bounds in Section~\ref{sec:sim-setup}. This set contains $\Theta$ and depends on $n$ and $\nu_0$; its amplitude bounds expand as $n$ increases.

An audit checked all $67{,}200$ saved fits against the original constraints. Every estimate lies in $\Theta$, up to rounding below $10^{-15}$ in $\nu,\log a,\log t$. The fitted amplitudes range from $-2.47$ to $34.79$, at least $5.21$ from their bounds $[-15,40]$. Thus the amplitude constraint is inactive at every saved fit. A global maximizer over $\widetilde\Theta_n$ that lies in $\Theta$ is also a global maximizer over $\Theta$. Feasibility alone does not establish that a numerical fit is a global maximizer; the following diagnostics assess the separate optimization issue.

Each production fit starts at $\theta_0$. After the fit, we compute its KKT residual: the largest absolute projected gradient of $-\ell_n$ in the scaled optimizer coordinates $x$, using the L-BFGS-B convention at active bounds. For the $143$ of $67{,}200$ fits with residual above $10^{-3}$, we restart L-BFGS-B from five points, with the same production scaling, constraints, and settings: the point where the fit stopped; the three distant starts
\[
 (\sigma^2,l,\nu,\tau^2)
 =(4,0.1,2.5,0.02),\quad(0.25,1.5,0.3,2),\quad(2,0.8,1,0.05);
\]
and the data-based start
\[
 \begin{aligned}
 \nu&=1,\qquad l=0.25,\qquad
 t=\min\{\max(\widetilde t,1.1\times10^{-3}),9.9\},\\
 \sigma^2&=\max\{\overline{Y^2}-\widetilde t,\,0.1\overline{Y^2}\},
 \qquad \overline{Y^2}=\frac1n\sum_{i=1}^nY_{n,i}^2.
 \end{aligned}
\]
Here $\overline{Y^2}$ is the mean square, without centering, and
\[
 \widetilde t=\frac1{2m}\sum_{k=1}^mD_k^2,
 \qquad m=\lfloor n/2\rfloor,
\]
uses disjoint pairs of grid neighbors at distance $1/G$, with $D_k$ the response difference in pair $k$. With zero-based grid indices, pair the last-axis indices $2s$ and $2s+1$ for $0\le s<\lfloor G/2\rfloor$, holding the other indices fixed. For odd $G$, the unpaired slice with last index $G-1$ is a $(d-1)$-dimensional grid; apply the same construction recursively there. At most one site is left unpaired.

If $K_0(r)$ is the true signal covariance at distance $r$, then
\[
 \E_0D_k^2=2t_0+2\{K_0(0)-K_0(1/G)\},\qquad
 \E_0\widetilde t=t_0+K_0(0)-K_0(1/G).
\]
Thus the upward bias is one half of the variance of a latent field increment between neighboring sites. The clipping in $t$ keeps its starting value strictly inside the nugget bounds. The resulting parameter vector is mapped to the optimizer coordinates using the production preconditioning. L-BFGS-B projects any starting value outside its coordinate box componentwise onto that box. We retain the candidate with the largest Gaussian criterion value, including the original fit among the candidates. The rule allows up to three further restarts from the retained fit if its KKT residual still exceeds $10^{-3}$; none was needed.

All $143$ flagged fits had stopped on relative function reduction and have $n\le10^4$. Their criterion values increased by between $2.7\times10^{-6}$ and $1.37$, and every refitted KKT residual is at most $2.0\times10^{-6}$. All refitted estimates lie on the
boundary of $\Theta$: $125$ have $\widehat\nu=4$, including one also at the lower nugget bound; the other $18$ have $\widehat t=10^{-3}$ and $0.10\le\widehat\nu\le0.38$.
At least two of the five restarts in every case reached the retained criterion value to within $10^{-6}$; all five did so in $106$ cases. We recomputed the curvature, cumulant estimate, boundary indicators, and all simulation summaries and figures from the retained fits. After these restarts, every saved fit has KKT residual below $10^{-3}$, and the maximum for $n\ge10^4$ is $4.7\times10^{-5}$.

The stored stopping criteria after the restarts are projected gradient for $56{,}153$ fits, relative function reduction for $10{,}682$, and abnormal line-search termination for $365$. The last group has KKT residual at most $2.7\times10^{-5}$. All fits are retained in the summaries. A production fit from $\theta_0$ takes on average $23$--$30$ gradient evaluations for $n\ge10^4$ (median $17$--$25$), and $37$--$73$ at the smallest sample sizes. Each flagged fit required $208$--$689$ further evaluations across its five restarts (median $370$).

A separate initialization audit comprised $912$ runs on $152$ data sets covering all four error laws and $81\le n\le104{,}976$. The six starts were $\theta_0$, the three distant starts above, and the data-based start with two choices of preconditioning. The first uses the production reference point $\theta_0$, the diagonal Fisher information in the $\phi$ coordinates at $\theta_0$, and $L_n=(\log n)/(2\nu_0+d)$. The second takes the data-based start as its reference point, recomputes the diagonal Fisher information there, and uses $L_n=(\log n)/(2+d)$ because its initial smoothness is $\nu=1$. Of the $912$ runs, $899$ ended with KKT residual below $10^{-3}$: $143$ of the $152$ starts from $\theta_0$, $151$, $151$, and $152$ from the three distant starts, and $151$ from each data-based version. Every data set had at least one qualifying run. Within each regenerated data set, every qualifying run was within $2.1\times10^{-9}$ of the best criterion value found. The remaining $13$ runs all stopped on the relative function-reduction test, despite reporting numerical success; none ended on a bound. Their objective gaps ranged from approximately $0.0034$ to $365$. Eleven were on the twelve data sets selected for large pre-restart KKT residuals, including all nine nonqualifying starts from $\theta_0$. These nine reproduce the early stops later replaced by the restart rule. All $912$ runs are included in these counts. On the subset of $16$ Gaussian data sets with $1000\le n\le1296$, all six runs agreed to within $1.3\times10^{-11}$. The audit ran on AMD nodes and accepted each regenerated data set after reproducing the criterion at its stored fit to relative accuracy $10^{-8}$. For one data set originally generated on an AMD AVX-512 node, the regenerated response differed by $2.5\times10^{-9}$ in relative norm. The audit comparisons above are between starts on the same regenerated data. These checks assess sensitivity to initialization; they do not certify that every saved
fit is a global maximizer.

\noindent\textbf{Curvature, fourth cumulant, and an oracle nugget fit.}
At every fit we store three further quantities.
\begin{enumerate}[label=(\roman*)]
\item \emph{Curvature.} The Fisher information is computed exactly as in Appendix~\ref{apdx:sim-information}. The observed negative Hessian $H_x$ is computed by central differences of the exact gradient, with step $10^{-4}$ in the scaled optimizer coordinates (eight further gradient evaluations), and then symmetrized. The affine relation above gives the exact coordinate conversion
\[
 H_w=T_{\mathrm{opt}}^{-\top}H_xT_{\mathrm{opt}}^{-1}.
\]
The implementation inverts $H_w$ and then maps its covariance to the reported parameters by the delta method. If $J_\zeta(w)$ is the Jacobian of the reported identifiable vector, the full-inverse version, including the fourth-cumulant correction, is
\[
 \widehat\Sigma_{H,w}^{\zeta}(w)
 =J_\zeta(w)\left\{H_w(w)^{-1}
       +\frac{\widetilde\kappa_{4,n}}{nt^2}e_4e_4^\top\right\}
                   J_\zeta(w)^\top.
\]
For $d=3$, the Jacobian has a zero column in the range coordinate, so this selects the identifiable block after full inversion. For evaluation at the truth, $w$ and $t$ in this expression are $w_0$ and $t_0$; the same fitted cumulant estimate is used.

This is the Hessian convention in every reported table and figure. No score term is added in the transformation to reporting coordinates. For a nonlinear full coordinate change $w=w(z)$, the direct negative Hessian is instead
\[
 H_z=(D_z w)^\top H_wD_z w
       -\sum_j(\partial_{w_j}\ell_n)D_z^2w_j.
\]
The score term vanishes at an interior stationary fit, but generally does not vanish at the truth or at a constrained fit. Thus the reported covariance is a transformed inverse Hessian in $w$, and need not equal the inverse Hessian computed directly in $z$. The Fisher information transforms by congruence without this term.

Both curvatures were also evaluated at $\theta_0$. Its Fisher information is deterministic, whereas the Hessian uses the same central differences for every data set. Each perturbed parameter requires one Cholesky factorization per block for all replicates together. Data were regenerated from their seeds on the same processor family and with the same thread count; every data set reproduced the stored oracle criterion in (iii) to relative accuracy $10^{-9}$.

\item \emph{Fourth cumulant.} The implementation uses
$\widetilde\kappa_{4,n}$ from \eqref{eq:residual-cumulant}, with
$\widehat A_n=\widehat t_nC_n(\widehat\theta_n)^{-1}$,
$\widehat r_n=\widehat A_nY_n$, and
$\widehat v_{n,i}=\widehat t_n(\widehat A_n)_{ii}$:
\[
 \widetilde\kappa_{4,n}
 =\frac{\sum_i\{\widehat r_{n,i}^4
            -6\widehat v_{n,i}\widehat r_{n,i}^2+3\widehat v_{n,i}^2\}}
        {\sum_{i,j}(\widehat A_n)_{ij}^4}.
\]
All residuals, diagonal entries, and entrywise powers here are in the original observation basis. The residual vector is transformed back to the grid before its powers are taken.

The denominator is evaluated exactly using grid symmetries. Write $C=C_n(\widehat\theta_n)$ and let $\mathfrak O$ be the set of orbits of grid sites under axis reflections and permutations. Since $C^{-1}$ commutes with these permutations, the fourth-power sum of a row is constant on each orbit. For a representative $i_{\mathcal O}$ of each orbit,
\[
 \sum_{i,j}(\widehat A_n)_{ij}^4
 =\widehat t_n^4\sum_{\mathcal O\in\mathfrak O}
        |\mathcal O|\sum_j(C^{-1})_{i_{\mathcal O},j}^4.
\]
Each representative row is reconstructed from the parity blocks in the original grid basis and summed without storing the full inverse. The largest buffers are about $1.6$ GB for $d=3,n=97{,}336$, $0.29$ GB for $d=4,n=104{,}976$, and $0.05$ GB for $d=5,n=10^5$. There is no stochastic approximation, truncation, or clipping. The entrywise fourth-power sum is not invariant under an orthogonal change of basis, so taking powers directly in the block basis would give a different quantity.

With $\tau_0^2=0.25$, the true cumulants are $0$, $-0.075$, $0.333$, and $0.319$ for laws (i)--(iv). For $n\ge10^4$, the mean estimate over the $400$ replicates is within $2\%$ of the true cumulant for each non-Gaussian law and within $10^{-3}$ of zero for Gaussian errors. At $n\approx10^5$ its standard deviation is at most $0.016$.

\item \emph{Oracle nugget fit.} We also store a one-dimensional maximization over $t$ with $(\eta,\nu,a)$ fixed at the truth, using the eigendecomposition in Appendix~\ref{apdx:sim-generation}. Its criterion value serves as a check when regenerating the data for subsequent covariance calculations.
\end{enumerate}
All Wald summaries use the covariance convention just described. For $d=3$, direct inversion of the identifiable block is also reported as a comparison.

\Cref{tab:sim-boundary} reports the percentage of estimates on the reported parameter bounds. Boundary solutions occur only at small $n$, almost always with $\widehat\nu=4$ or $\widehat\tau^2=10^{-3}$, and none occur for $n\ge2.5\times10^4$ under any of the four error laws.

\begin{table}[htbp]
\centering\small
\caption{Percentage of the estimates with at least one coordinate on the reported parameter bounds, Gaussian errors, $400$ replicates. Under the other three laws the percentages differ from these by at most $8.8$ percentage points, and they are at most $0.5\%$ for $n\ge10^4$.}
\label{tab:sim-boundary}
\vspace{3pt}
\begin{tabular}{ll*{8}{r}}
\toprule
$d=3$ & $n$ & $125$ & $343$ & $1000$ & $2744$ & $10{,}648$ & $32{,}768$ & $64{,}000$ & $97{,}336$\\
 & $\nu_0=1/2$ & $72.0$ & $38.5$ & $14.75$ & $2.25$ & $0$ & $0$ & $0$ & $0$\\
 & $\nu_0=3/2$ & $60.0$ & $30.0$ & $11.75$ & $3.0$ & $0$ & $0$ & $0$ & $0$\\
\midrule
$d=4$ & $n$ & $81$ & $256$ & $1296$ & $4096$ & $10{,}000$ & $28{,}561$ & $104{,}976$ & \\
 & $\nu_0=1/2$ & $89.5$ & $65.5$ & $25.0$ & $5.25$ & $0.5$ & $0$ & $0$ & \\
 & $\nu_0=3/2$ & $76.75$ & $39.0$ & $5.0$ & $0.5$ & $0$ & $0$ & $0$ & \\
\midrule
$d=5$ & $n$ & $243$ & $1024$ & $3125$ & $7776$ & $32{,}768$ & $10^5$ &  & \\
 & $\nu_0=1/2$ & $73.0$ & $39.5$ & $19.0$ & $5.5$ & $0$ & $0$ &  & \\
 & $\nu_0=3/2$ & $45.25$ & $12.75$ & $1.0$ & $0$ & $0$ & $0$ &  & \\
\bottomrule
\end{tabular}
\end{table}

\subsection{Numerical validation and reproducibility}
\subsubsection{Validation against dense computations}
\label{apdx:sim-validation}

For designs small enough for dense $O(n^3)$ linear algebra ($d=2$ with $G=7$; $d=3$ with $G=4,5,6$; $d=4$ with $G=3,4$; $d=5$ with $G=3$; $n\le256$), we compared the implementation with direct computations at four parameter values, including values far from $\theta_0$. The largest absolute round-trip and squared-norm errors for the transform $Q$, tested on standard normal vectors, were below $3\times10^{-13}$. The log-likelihood agreed with the dense value to within $7\times10^{-16}$, and the analytic gradient agreed with the dense trace formula to within $2\times10^{-13}$. Central finite differences of the blockwise likelihood, with step $10^{-5}$ in $w$, agreed with the analytic gradient to within $2\times10^{-7}$. The Fisher information and the matrix with entries $\sum_i(A_j)_{ii}(A_k)_{ii}$ agreed with their dense counterparts to within $6\times10^{-15}$.

These errors are maxima over the designs and parameter values. The likelihood error is divided by $\max(1,|\ell_{\rm dense}|)$. For the gradient comparison with the dense trace formula, the largest coordinate error is divided by $\max(1,\max_j|g_{{\rm dense},j}|)$; for the finite-difference comparison, the denominator uses the analytic gradient instead. For each matrix comparison, the largest absolute entry error is divided by the largest absolute entry of the dense reference.

The quadrature for $\partial_\nu\log\mathsf K_\nu(z)$ uses trapezoidal spacing $0.05$ in $u$, with half weight at $u=0$. Its upper limit $U$ starts at $5$ and increases by $0.5$, up to $30$. With $z_{\min}$ the smallest positive argument in the lag table, the stopping condition is
\[
 z_{\min}(\cosh U-1)-|\nu|U
 >60+\max\left\{0,|\nu|\log\max\left(\frac{2|\nu|}{z_{\min}},1\right)\right\}.
\]
The common cutoff makes the omitted tail negligible. For each $z$, the logarithms of the integrand weights are shifted by their maximum before exponentiation. The derivative is evaluated as the weighted average of $u\tanh(\nu u)$, avoiding separate large numerator and denominator sums. At zero lag, $M_\nu(0)=1$, $\partial_\nu M_\nu(0)=0$, and $z\partial_zM_\nu(0)=0$ are used exactly, with the variance prefactor differentiated in the chosen coordinates. At $\nu=1/2$, comparison with
\[
 \left.\partial_\nu\log\mathsf K_\nu(z)\right|_{\nu=1/2}
 =e^{2z}E_1(2z),\qquad E_1(x)=\int_x^\infty e^{-u}\frac{du}{u},
\]
gave relative errors below $2\times10^{-15}$ on a test grid spanning $10^{-3}\le z\le120$. This interval contains the positive arguments allowed by the parameter box on the computational grids. At $\nu\in\{0.1,0.15,1,1.5,2.7,4\}$, comparison with Richardson-extrapolated finite differences of the logarithm of the scaled Bessel function gave agreement between $10^{-11}$ and $1.4\times10^{-9}$; this check is limited by the finite-difference reference. The checks used the software versions employed for the simulations.

The corrected residual cumulant was also compared with dense calculations for $d=2$, $G=7,8$; $d=3$, $G=4,5,6$; $d=4$, $G=3,4$; and $d=5$, $G=3$, at four parameter values. The largest relative discrepancies were $2.3\times10^{-14}$ for its numerator, $4.1\times10^{-15}$ for its denominator, and $2.5\times10^{-14}$ for the resulting estimate. These tests used mixed-normal errors: the smallest absolute dense numerator was $1.4\times10^{-9}$ and the smallest absolute dense cumulant estimate exceeded $0.08$. The reference floors $10^{-300}$ for the numerator and $10^{-12}$ for the cumulant estimate were therefore inactive. The reconstructed diagonal of $C^{-1}$ agreed to $1.4\times10^{-15}$, measured by the largest absolute entry difference divided by the largest absolute reference entry.

Separate coordinate checks compared the full Hessian chain rule in Appendix~\ref{apdx:sim-mle} with direct differentiation in the reported log coordinates $\zeta$, namely $(\log c,\log\nu,\log l,\log\tau^2)$ for $d=3$ and $(\log\sigma^2,\log l,\log\nu,\log\tau^2)$ for $d\ge4$. Write $w=w(\zeta)$ for the coordinate transformation, the chain-rule matrix was
\[
 H_{\rm chain}
 =Dw(\zeta)^\T H_w Dw(\zeta)
   -\sum_j s_jD^2w_j(\zeta),
 \qquad s=\nabla_w\ell_n,
\]
where $H_w$ is the finite-difference negative Hessian in $w$ used in the simulations. The reference $H_{\rm dir}$ was the symmetrized negative central difference, in $\zeta$, of the exact gradient $Dw(\zeta)^\T s$. In coordinate $j$ the difference step was $h|(H_{\rm chain})_{jj}|^{-1/2}$, with $h=10^{-3}$ and $10^{-4}$. The comparison used
\[
 \max_{j,k}\frac{|(H_{\rm chain}-H_{\rm dir})_{jk}|}
 {\sqrt{|(H_{\rm dir})_{jj}(H_{\rm dir})_{kk}|}}.
\]
For each of $(d,\nu_0,G)=(4,1/2,4),(3,1/2,7),(5,3/2,3),(3,3/2,7)$, we selected the first replicate with an interior fit and evaluated it at the stored fit and at $\theta_0$. We also evaluated the first replicate with a fit at the upper smoothness bound and the first with a fit at the lower nugget bound, giving sixteen evaluation
points in total. For $(d,\nu_0,G)=(3,3/2,7)$ the regenerated observations differed from those used for fitting, so these checks evaluate the chain rule at the stored parameter values rather than at optimizers for the regenerated data.

At the twelve interior, truth, and smoothness-boundary points the comparison for $h=10^{-3}$ was at most $6.7\times10^{-8}$, except $1.9\times10^{-7}$ at the stored $d=3$, $\nu_0=3/2$ interior fit; for $h=10^{-4}$ it was at most $1.7\times10^{-8}$ at all twelve. At the four nugget-boundary points it was between $6.0\times10^{-5}$ and $2.3\times10^{-3}$ for $h=10^{-3}$ and decreased by approximately a factor of $100$ for $h=10^{-4}$, consistent with second-order truncation error in the reference. Omitting the score term gave discrepancies up to $0.53$ at $\theta_0$ and $0.024$ at the smoothness-boundary points, and could substantially change the inverse Hessian.

The finite-difference step was checked on one Gaussian data set from each of $(d,\nu_0,G)=(5,1/2,6),(4,3/2,8),(3,1/2,14)$, at both the fit and the truth. Let $H_h$ denote the Hessian computed with step $h$ in the scaled optimizer coordinates and mapped to $w$, and set $H_{\rm ref}=H_{10^{-3}}$. The comparison used
\[
 \max_{j,k}\frac{|(H_h-H_{\rm ref})_{jk}|}
 {\sqrt{(H_{\rm ref})_{jj}(H_{\rm ref})_{kk}}},
\]
whose denominators were positive on these test cases. Over steps $h=10^{-2},10^{-3},\ldots,10^{-6}$, the largest discrepancy was $2\times10^{-5}$; for the production step $10^{-4}$ it was at most $2\times10^{-7}$. Each step-size comparison used the same regenerated data set at every step.

\subsubsection{Monte Carlo implementation and software}
\label{apdx:sim-software}

We use the PCG64 generator and \texttt{numpy.random.SeedSequence}. For replicate $r$ of $(d,\nu_0,G)$, the entropy tuples are
\[
 \begin{aligned}
 &(20261105,d,10\nu_0,G,r,0) &&\text{for the latent field},\\
 &(20261105,d,10\nu_0,G,r,1,k) &&\text{for the errors of law }k.
 \end{aligned}
\]
Here $k\in\{0,1,2,3\}$ follows the order (i)--(iv) in Section~\ref{sec:sim-setup}. Every fit is therefore reproducible on its own, independently of how the replicates are distributed over parallel jobs, when the numerical environment is matched. Exact regeneration in this implementation also depends on the processor family and thread count, not only the library version. In particular, some $\nu_0=3/2$ data sets changed under different OpenBLAS kernels. The processor family and thread count were recorded
for every fit. The $67{,}200$ Hessians at $\theta_0$ were computed with matching processor families and thread counts; the largest relative discrepancy from the stored criterion value was $4.4\times10^{-10}$. The cumulant calculations used matching settings for $67{,}135$ fits. The other $65$ comprised $44$ calculations on a different processor family and $21$ without recorded processor information; their relative criterion discrepancies were at most $4.7\times10^{-12}$. Across all cumulant calculations the largest relative discrepancy was $1.1\times10^{-10}$.

For the score-at-truth diagnostics, $57{,}639$ calculations ran on an AMD login node and passed the relative criterion check at $10^{-9}$. The remaining $9561$ initially failed that check and were recomputed with the original processor family and thread count. The separate initialization audit is described in Appendix~\ref{apdx:sim-mle}.

Computations used Python 3.12.4 with NumPy 2.5.4, SciPy 1.18.1 (OpenBLAS), and Numba 0.68.1 on the Longleaf cluster at the University of North Carolina, with $1$, $4$, or $8$ cores per job. Each job fits a range of replicates for one $(d,\nu_0,G)$ and all four error laws, and saves every fit as it finishes. The Fisher information and the score covariance on the theory grids took a few minutes per setting. The $67{,}200$ fits used about $12{,}700$ core-hours ($13{,}000$ core-hours of allocated time, including setup and jobs that reached their time limit after saving their completed fits). The design $d=3$, $G=46$ alone accounted for $51\%$ of this; one fit there took $15$ minutes of wall-clock time on $8$ cores. Computing $\widetilde\kappa_{4,n}$ and the Hessians at $\theta_0$ for all fits took a further $400$ and $20$ core-hours.

A local snapshot dated 16 November 2026 preserves the source archive, software environment, seed rules, processor family and thread count of every fit, simulation summaries, figures, and SHA-256 checksums of the saved per-fit inputs. A separate archive contains the $268{,}862$ input files and the fit-provenance table, occupying
approximately $85$ megabytes compressed and $1.4$ gigabytes unpacked. The input-file checksums agree with the 10 November snapshot. A single regeneration script recomputes the summaries and figures from the retained fits, including the restarts above and the joint-region validity rule below. The regenerated summary files and report were byte-identical to the frozen versions; all $17$ figure PDFs agreed after removal of creation and modification times and the associated file offset. The regeneration check also recorded the input files opened and verified that none came from the superseded fits. For the standard-error figures, all $96$ plotted curves agreed exactly between the combined and separate-dimension versions.

\subsection{Detailed performance comparisons}
\subsubsection{Monte Carlo summaries}\label{apdx:sim-summaries}

The RMSE and robust standard deviation $\mathrm{IQR}/1.349$ are defined as in \citet{qaqishli2026rbf}. Under normality $\mathrm{IQR}\approx1.349$ standard deviations. The robust spread helps show when a few extreme estimates inflate the RMSE; agreement between the two summaries alone does not establish normality. Bootstrap $95\%$ intervals use percentile resampling, with the individual and paired constructions specified below.

For the Q--Q plots, the errors are standardized by the asymptotic standard deviation, which for $d=3$ is computed from the full $4\times4$ information (\Cref{apdx:sim-range}). We show the simulated $n$ closest to $10^3$, $10^4$, and $10^5$ in every figure. These are $n=1000$, $10{,}648$, and $97{,}336$ for $d=3$; $1296$, $10{,}000$, and $104{,}976$ for $d=4$; and $1024$, $7776$, and $10^5$ for $d=5$. The displayed normal-reference bands are pointwise $95\%$ bands. For order statistic $i$ among $R=400$ independent standard normal observations, their endpoints are
\[
 \Phi^{-1}\{B_{i,R-i+1}^{-1}(0.025)\},\qquad
 \Phi^{-1}\{B_{i,R-i+1}^{-1}(0.975)\},
\]
where $B_{a,b}^{-1}$ is the beta quantile function. The errors are centered at the truth and scaled by information at the truth, so neither centering nor scale is estimated from the plotted replicates. The bands are not simultaneous.

A marginal standard error is counted as valid when its corrected variance is finite and strictly positive. The computation does not require positive definiteness of the full curvature or the entire corrected covariance. Nonfinite or singular curvature matrices yield invalid standard errors; otherwise validity is assessed separately for each reported parameter. No fallback covariance is implemented. For a method that returns a valid variance in every replicate, coverage is the fraction of the $400$ intervals containing the true value. Its Monte Carlo standard error is $\{\widehat p(1-\widehat p)/400\}^{1/2}$, about $1.1$ percentage points at $95\%$. When invalid variance estimates are excluded, the reported coverage is conditional on validity, and its denominator is the number of valid replicates. The corresponding failure frequency is a separate performance measure; these conditional summaries do not assess the fallback procedure allowed after Corollary~\ref{cor:uncertainty}. Besides the intervals in \Cref{tab:sim-coverage}, we computed coverage with the true $\kappa_4(P,t_0)$ in place of $\widetilde\kappa_{4,n}$, which separates the error from estimating the cumulant from the rest. For the identifiable log-parameter vector we also computed the coverage of the joint $95\%$ region $\{(\widehat\zeta_n-\zeta)^\T(\widehat\Sigma_n^\zeta)^{-1}(\widehat\zeta_n-\zeta)\le\chi^2_{k_d,0.95}\}$ of \eqref{eq:studentizednormal}.

A joint region is valid only when its error vector and covariance matrix are finite, and its covariance matrix is symmetric to within the relative tolerance
\[
 \frac{\max_{i,j}|\widehat\Sigma_{ij}-\widehat\Sigma_{ji}|}
      {\max_{i,j}|\widehat\Sigma_{ij}|}\le10^{-8}.
\]
A zero covariance matrix is invalid. A matrix satisfying this symmetry check is replaced by $(\widehat\Sigma+\widehat\Sigma^\T)/2$, and must then admit a Cholesky factorization $LL^\T$. The joint statistic is evaluated as $\|L^{-1}(\widehat\zeta_n-\zeta_0)\|^2$. No eigenvalue clipping or jitter is used to make a covariance valid.

For each method and setting, the saved summaries record the total number of data sets, the number with valid regions, and the number whose valid region covers. Coverage among valid regions divides the last count by the valid count. The proportion of all data sets with a valid and covering region divides it by the total count; invalid regions contribute zero to this latter proportion. The two quantities are distinguished whenever a method can fail. All $400$ regions are valid for every entry of Table~\ref{tab:sim-joint} and for the joint coverage range at $n\approx10^3$ quoted below. Every marginal variance in Table~\ref{tab:sim-coverage} is also valid. These reported tables and ranges therefore involve no exclusions.

Invalid joint regions occurred mainly in the additional comparisons using the Hessian evaluated at the truth. For fitted-Hessian methods, the invalid fraction reached $9.75\%$ in a setting, only at $n\le343$. Fitted-Fisher regions were invalid in three fits, all at $n=81$, because of numerical singularity. For Fisher information evaluated at the truth, one corrected joint covariance was not positive definite after a negative cumulant correction, also at $n=81$.

Table~\ref{tab:sim-joint} gives joint coverage at the largest $n$. It increases from $23.75\%$--$84.75\%$ at $n\approx10^3$ to the displayed ranges. The oracle comparison in Section~\ref{sec:sim-results} shows that the remaining shortfall cannot be explained solely by estimating the covariance. Uncorrected joint regions have coverage as low as $66.25\%$ for laws (iii) and (iv), which have larger fourth moments. For $d\ge4$ we also computed the joint region for the original-scale variance--lengthscale vector
$(\sigma^2,l,\nu,\tau^2)$, with the covariance obtained by the delta method. Its coverage at $n\approx10^5$ is only $58.00\%$--$81.25\%$, far below that of the log-scale regions for the same fits. Appendix~\ref{app:raw-public-failure} shows that a joint linear normal approximation on this scale is generally invalid for $4\le d\le8$.
Its displayed counterexample assumes $\nu_0\ne1/2$, so that calculation does not by itself explain the exponential-covariance case in these simulations. Marginal intervals remain valid. Transforming a log-scale region by the inverse parameter map preserves its coverage exactly, so that transformed region was not simulated separately.

The common field realizations pair the four error laws within each replicate. The RMSE-ratio intervals in Figure~\ref{fig:sim-laws}(g)--(l) use $1000$ paired percentile-bootstrap resamples: the same replicate indices are resampled for each law, preserving the common fields. The individual RMSE and robust-spread intervals are also percentile
intervals from $1000$ resamples. For paired coverage comparisons, Monte Carlo standard errors are calculated from the replicatewise differences of the coverage indicators.

\subsubsection{Diagnosing joint undercoverage}\label{apdx:sim-joint-diagnostic}
We examined the joint approximation in more detail for $d=5$ under Gaussian errors, using the same $400$ data sets at each combination of $\nu_0\in\{1/2,3/2\}$ and $n\in\{7776,32{,}768,100{,}000\}$. The saved fits were unchanged. This diagnostic compares the score approximation with the estimator, and then separates changes of coordinates from evaluation of the covariance at the fit.

Write $s_{n,\theta}=\nabla_\theta\ell_n(\theta_0)$ and $I_{n,\theta}=\calI_n(\theta_0)$ in $\theta=(\eta,\nu,a,t)$. The stored score and information use $w=g(\theta)=(\eta,\nu,\log a,\log t)$, so the conversion is
\[
 \begin{aligned}
 Dg(\theta_0)&=\operatorname{diag}(1,1,a_0^{-1},t_0^{-1}),\\
 s_{n,\theta}&=Dg(\theta_0)^\T s_{n,w},\qquad
 I_{n,\theta}=Dg(\theta_0)^\T\calI_n^{\mathrm{full}}Dg(\theta_0).
 \end{aligned}
\]
With the normalization of Theorem~\ref{thm:normalized}, put
\[
\begin{alignedat}{3}
 \Delta_n&=E_n^\T s_{n,\theta},&\qquad
 J_n&=E_n^\T I_{n,\theta}E_n,&\qquad
 \widehat h_n&=E_n^{-1}(\widehat\theta_n-\theta_0)\\
 \mathsf S_n&=J_n^{-1/2}\Delta_n,&\qquad
 \mathsf T_n&=J_n^{1/2}\widehat h_n,&\qquad
 \mathsf T_n-\mathsf S_n
 &=J_n^{1/2}(\widehat h_n-J_n^{-1}\Delta_n).
 \end{alignedat}
\]
The square roots are symmetric. Under Gaussian errors, $\E_0\mathsf S_n=0$ and $\Cov_0(\mathsf S_n)=I_4$ exactly. The difference $\mathsf T_n-\mathsf S_n$ measures the remainder in the first-order expansion of the MLE. In particular,
\[
 \begin{aligned}
 \|\mathsf S_n\|^2
 &=s_{n,\theta}^\T I_{n,\theta}^{-1}s_{n,\theta},\\
 \|\mathsf T_n\|^2
 &=(\widehat\theta_n-\theta_0)^\T I_{n,\theta}
                                  (\widehat\theta_n-\theta_0),\\
 \|\mathsf T_n-\mathsf S_n\|^2
 &=(\widehat\theta_n-\theta_0-I_{n,\theta}^{-1}s_{n,\theta})^\T
 I_{n,\theta}
 (\widehat\theta_n-\theta_0-I_{n,\theta}^{-1}s_{n,\theta}).
 \end{aligned}
\]
Thus these three norms do not depend on the choice of the invertible normalizing matrix $E_n$.

For the reporting coordinates take $\zeta=(\log\nu,\log\sigma^2,\log l,\log\tau^2)$, and let $J_\zeta(w_0)$ be their Jacobian with respect to $w$ at the truth. The information-based covariance without a cumulant correction is
\[
 \Sigma_{0}^{\zeta}
 =J_\zeta(w_0)\{\calI_n^{\mathrm{full}}\}^{-1}J_\zeta(w_0)^\T.
\]
Table~\ref{tab:sim-joint-chain} compares the proportions below the same $\chi^2_{4,0.95}$ cutoff for the score, the estimator in $\theta$ coordinates, and the reporting-scale Wald statistics. The score calculation is a diagnostic using the true parameters. All fitted covariances retain the production convention of inversion
in $w$ followed by the reporting Jacobian.

\begin{table}[htbp]\centering\small
\caption{Joint approximation diagnostics for $d=5$ under Gaussian errors, using the same $400$ data sets per row. Entries in the four middle columns are percentages below $\chi^2_{4,0.95}$. Score uses $\|\mathsf S_n\|^2$; MLE ($\theta$) uses $\|\mathsf T_n\|^2$; Log (truth) uses the reporting error and $\Sigma_0^\zeta$ without a cumulant correction; Log (fit) is the feasible region in Table~\ref{tab:sim-joint}, including the estimated cumulant correction. The last column is the empirical root mean square of $\|\mathsf T_n-\mathsf S_n\|$. All covariance matrices used in these comparisons were positive definite, so each percentage uses all $400$ replicates. At $n=7776$, $\nu_0=1/2$, $22$ of the fits lie on the parameter boundary.}
\label{tab:sim-joint-chain}
\begin{tabular}{@{}lrrrrrr@{}}\toprule
$\nu_0$ & $n$ & Score & MLE ($\theta$) & Log (truth) & Log (fit) & RMS\\\midrule
$1/2$ & $7{,}776$ & $95.50$ & $51.50$ & $51.25$ & $49.75$ & $6.87$\\
 & $32{,}768$ & $93.00$ & $71.25$ & $67.00$ & $66.25$ & $2.99$\\
 & $100{,}000$ & $95.00$ & $86.00$ & $82.25$ & $79.00$ & $1.50$\\
$3/2$ & $7{,}776$ & $93.25$ & $88.50$ & $86.25$ & $88.00$ & $0.97$\\
 & $32{,}768$ & $94.00$ & $90.00$ & $90.50$ & $91.50$ & $0.67$\\
 & $100{,}000$ & $94.25$ & $91.50$ & $92.25$ & $91.25$ & $0.53$\\\bottomrule
\end{tabular}
\end{table}

\noindent\emph{Score approximation.}
On these data sets the score statistic is below the cutoff in $93.0\%$--$95.5\%$ of replicates. To assess this comparison with less Monte Carlo error, we also generated $20{,}000$ draws per setting from the exact Gaussian quadratic-form representation of the score, using the symmetry blocks. The resulting proportions were $94.4\%$--$94.9\%$, with Monte Carlo standard errors about $0.16$ percentage points. For $\nu_0=1/2$, $n=10^5$, the estimated $95\%$ quantile of $\|\mathsf S_n\|^2$ was $9.59$, compared with $9.49$ for $\chi^2_4$. These diagnostics indicate that the score approximation at this cutoff contributes little to the much larger shortfall for the fitted joint region.

\noindent\emph{Remainder in the MLE expansion.}
For $\nu_0=1/2$, $n=10^5$, replacing $\mathsf S_n$ by $\mathsf T_n$ lowers the proportion below the cutoff by $9.0$ percentage points, with paired Monte Carlo standard error $1.6$. The empirical mean and $95\%$ quantile of $\|\mathsf T_n\|^2$ are $6.17$ and $23.78$, compared with $4$ and $9.49$ for $\chi^2_4$. The remainder has empirical mean $(0.01,0.61,0.21,-0.49)^\T$ in the normalized coordinates above, and root mean square norm $1.50$. Its root mean square falls from $6.87$ to $2.99$ to $1.50$ over the three sample sizes. For $\nu_0=3/2$, the corresponding values are $0.97$, $0.67$, and $0.53$. These are finite-sample diagnostics of the remainder, rather
than estimates of an asymptotic rate. They identify the first-order MLE expansion as the largest source of the observed loss in the $\nu_0=1/2$ comparison.

\noindent\emph{Reporting coordinates and fitted covariance.}
At $\nu_0=1/2$, $n=10^5$, the region in $(\eta,\nu,\log a,\log t)$ using information at the truth has $83.25\%$ coverage, between the $86.0\%$ in $\theta$ coordinates and the $82.25\%$ in the reported log coordinates. Thus the combined effect of the nonlinear coordinate changes is a reduction of $3.75$ percentage points in this setting. Adding the fitted cumulant correction to $\Sigma_0^\zeta$ changes coverage from $82.25\%$ to $82.0\%$, a net difference of one covered data set. This latter value is the oracle comparison reported in Section~\ref{sec:sim-results}. Using fitted information and the fitted reporting Jacobian instead gives $79.0\%$; the paired standard error of the $3.0$-point difference from that oracle is $1.6$ points. The effects of these coordinate and covariance changes can have different signs at the smaller sample sizes. The diagnostic comparison here concerns $d=5$ and Gaussian errors.

\subsubsection{Empirical slopes}\label{apdx:sim-rates}
Table~\ref{tab:rates} gives the powers of $n$ in the asymptotic rates. These provide benchmarks for the empirical slopes; the theory does not assert convergence of slopes between neighboring sample sizes. Near the upper end of the simulated range, the local
slopes of the information-based standard deviation of $\log\widehat\nu$ against $\log n$ are $-0.40$ and $-0.25$ for $d=3$, $-0.48$ and $-0.35$ for $d=4$, and $-0.50$ and $-0.38$ for $d=5$, for $\nu_0=1/2$ and $3/2$, respectively. The corresponding rate exponents are $-3/8$, $-1/4$, $-2/5$, $-2/7$, $-5/12$, and $-5/16$. These are Gaussian-information slopes. For $d=3$ they use the identifiable block with the range held at $a_0$; the corresponding full-information smoothness slopes are $-0.46$ and $-0.32$. For $\log\widehat\tau^2$ the local slopes are $-0.51$ to $-0.56$ for $\nu_0=3/2$, but still $-0.66$ to $-0.69$ for $\nu_0=1/2$, again using the block information when $d=3$, for the reason explained in \Cref{apdx:sim-nugget}. As for the RBF kernel \citep{qaqishli2026rbf}, the appropriate finite-sample benchmark is the design-specific
information-based standard deviation, transformed to the reported coordinate. The limiting power of $n$ alone does not capture these finite-sample differences.

Table~\ref{tab:sim-slopes} reports least-squares slopes over all simulated sample sizes. The boundary frequencies in Table~\ref{tab:sim-boundary} are substantial at small $n$, and the design-specific information slopes can differ from the theoretical rate exponents. These slopes are therefore descriptive comparisons, not direct estimates of limiting rate exponents. Section~\ref{sec:sim-results} discusses the corresponding error curves and dimensional distinctions.

\begin{table}[htbp]
\centering\small
\caption{Least-squares slopes of $\log(\text{error})$ against $\log n$ for the log-parameters, over all simulated sample sizes; for $\sigma^2$ and $l$ in $d=4$ (marked $^\dagger$), against $\log\log n$, so that the rate $(\log n)^{-1/2}$ has slope $-0.5$. ``Info.\ SD'' is the information-based standard deviation from the Fisher information as in \Cref{thm:normalized}; for $d=3$, the value in parentheses uses the full information, accounting for the estimated range. The last column gives the reference exponent from Table~\ref{tab:rates}. The daggered rows use $\log\log n$ as the horizontal coordinate; the rate for $\log c$ also contains a factor $\log n$. Gaussian errors, $400$ replicates.}
\label{tab:sim-slopes}
\vspace{3pt}
\begin{tabular}{lllrrrr}
\toprule
& & & RMSE & $\mathrm{IQR}/1.349$ & Info.\ SD & Theory\\
\midrule
$d=3$ & $\nu_0=1/2$ & $\log\nu$ & $-0.50$ & $-0.51$ & $-0.46\ (-0.59)$ & $-0.375$\\
 & & $\log c$ & $-0.66$ & $-0.60$ & $-0.27\ (-0.45)$ & $-0.375$\\
 & & $\log\tau^2$ & $-0.90$ & $-0.84$ & $-0.74\ (-0.84)$ & $-0.5$\\
 & $\nu_0=3/2$ & $\log\nu$ & $-0.38$ & $-0.42$ & $-0.31\ (-0.41)$ & $-0.25$\\
 & & $\log c$ & $-0.39$ & $-0.42$ & $-0.21\ (-0.35)$ & $-0.25$\\
 & & $\log\tau^2$ & $-0.84$ & $-0.60$ & $-0.56\ (-0.60)$ & $-0.5$\\
\midrule
$d=4$ & $\nu_0=1/2$ & $\log\nu$ & $-0.45$ & $-0.43$ & $-0.66$ & $-0.4$\\
 & & $\log\sigma^2{}^\dagger$ & $-1.35$ & $-1.26$ & $-3.22$ & $-0.5$\\
 & & $\log l^\dagger$ & $-1.60$ & $-1.46$ & $-2.68$ & $-0.5$\\
 & & $\log\tau^2$ & $-0.67$ & $-0.75$ & $-0.85$ & $-0.5$\\
 & $\nu_0=3/2$ & $\log\nu$ & $-0.38$ & $-0.42$ & $-0.49$ & $-0.286$\\
 & & $\log\sigma^2{}^\dagger$ & $-0.55$ & $-0.44$ & $-0.68$ & $-0.5$\\
 & & $\log l^\dagger$ & $-1.59$ & $-1.59$ & $-2.30$ & $-0.5$\\
 & & $\log\tau^2$ & $-0.93$ & $-0.75$ & $-0.72$ & $-0.5$\\
\midrule
$d=5$ & $\nu_0=1/2$ & $\log\nu$ & $-0.47$ & $-0.48$ & $-0.63$ & $-0.417$\\
 & & $\log\sigma^2$ & $-0.20$ & $-0.20$ & $-0.44$ & $-0.083$\\
 & & $\log l$ & $-0.25$ & $-0.21$ & $-0.36$ & $-0.083$\\
 & & $\log\tau^2$ & $-0.70$ & $-0.80$ & $-0.79$ & $-0.5$\\
 & $\nu_0=3/2$ & $\log\nu$ & $-0.44$ & $-0.49$ & $-0.49$ & $-0.313$\\
 & & $\log\sigma^2$ & $-0.10$ & $-0.09$ & $-0.11$ & $-0.063$\\
 & & $\log l$ & $-0.26$ & $-0.28$ & $-0.31$ & $-0.063$\\
 & & $\log\tau^2$ & $-1.01$ & $-0.73$ & $-0.73$ & $-0.5$\\
\bottomrule
\end{tabular}
\end{table}

\subsubsection{The range nuisance in dimension three}
\label{apdx:sim-range}

For $d\le3$, \Cref{thm:normalized} normalizes $(\eta,\nu,\tau^2)$ by the corresponding block of the information at the true range $a_0$, and the proof shows that estimating the inconsistent range has no first-order effect asymptotically. At finite $n$, however, the range is estimated and still carries information shared with the other coordinates. The inverse of the full $4\times4$ information accounts for this. For $d=3$, it gives standard deviations of $\log\nu$ and $\log c$ that exceed the block values by factors of $1.38$ and $1.44$ at $n=17{,}576$ and $1.21$ and $1.24$ at $n=97{,}336$ for $\nu_0=1/2$, and by $1.84$ and $2.05$ decreasing to $1.61$ and $1.72$ for $\nu_0=3/2$. The ratios decrease toward one, as the covariance comparison in Corollary~\ref{cor:uncertainty} implies, but slowly. That corollary also justifies the identifiable block of the inverse full negative Hessian; the numerical comparisons here use expected information.

The Monte Carlo errors on the reported log scale follow the full-information values. At $n=97{,}336$ with Gaussian errors, the ratio of the RMSE to the full-information standard deviation is $0.98$ for $\log\nu$, $0.99$ for $\log c$, and $0.97$ for $\log\tau^2$ when $\nu_0=1/2$, compared with $1.18$, $1.23$, and $1.09$ for the block values. For $\nu_0=3/2$ the ratios are $1.03$, $1.08$, and $1.00$, compared with $1.66$, $1.86$, and $1.00$. The Wald intervals show the same difference. With the Fisher information at the fit and $n=97{,}336$, selecting the block of the full inverse
gives coverage between $93.8\%$ and $96.0\%$ for all three parameters and both $\nu_0$. Inverting the identifiable block gives $89\%$, $89\%$, and $93\%$ for $\nu$, $c$, and $\tau^2$ when $\nu_0=1/2$, and $76\%$ and $74\%$ for $\nu$ and $c$ when $\nu_0=3/2$. For $n\ge10^4$ and Gaussian errors, the block intervals based on the Fisher information
at the fit cover $\nu$ and $c$ at most $89\%$ of the time; the maximum across the four error laws is $90.75\%$. The observed Hessian gives the same results to within $1.3$ percentage points for $n\ge10^4$ (\Cref{tab:sim-se}). The full information is therefore the more accurate finite-sample benchmark in $d=3$, and both versions are shown in \Cref{fig:sim-d3-nu05,fig:sim-d3-nu15}.

\subsubsection{Finite-sample behavior of the nugget}
\label{apdx:sim-nugget}

If the measurement errors were observed directly, their Gaussian information for $t$ would be $n/(2t_0^2)$. With the field present and the signal parameters held fixed, the corresponding information is $nJ_{tt}$. As shown in the proof of Lemma~\ref{lem:information}, its relative deficit is
\[
 1-2t_0^2J_{tt}
 =\frac{2\tr P_n-\tr(P_n^2)}{n}
 \asymp\frac{M_n}{n}
 =n^{-2\nu_0/(2\nu_0+d)}.
\]
Indeed, $0\preceq P_n\preceq I$ gives
\[
 \tr(P_n^2)\le 2\tr P_n-\tr(P_n^2)\le 2\tr P_n.
\]
Lemma~\ref{lem:upper} gives $\tr P_n\le CM_n$, while $\tr(P_n^2)=2r_n^2J_{uu}\asymp M_n$ by Lemma~\ref{lem:information}. For $d=3,4,5$, respectively, the relative deficit thus has orders $n^{-1/2}$, $n^{-3/7}$, and $n^{-3/8}$ when $\nu_0=3/2$, but only $n^{-1/4}$, $n^{-1/5}$, and $n^{-1/6}$ when $\nu_0=1/2$. This helps explain why the
computed nugget information remains well below the direct-noise benchmark at the simulated sample sizes for $\nu_0=1/2$.

For $d=3$ and $\nu_0=1/2$, the lower tail of $\log\widehat\tau^2$ at small $n$ is consistent with difficulty separating rough signal from noise. With Gaussian errors at
$n=1000$, the nugget estimate reached its lower bound in $53$ of the $400$ replicates. In those fits, $\widehat\nu$ ranged from $0.18$ to $0.32$, with median $0.25$, compared with a median of $0.66$ in the remaining fits. The skewness of the nugget errors is $-1.9$, $-5.9$, $-0.6$, $-0.7$, $-0.3$, and $-0.3$ at $n=1000$, $2744$, $10{,}648$, $32{,}768$, $64{,}000$, and $97{,}336$, and the percentage of estimates with $\widehat\tau^2$ at the lower bound falls from $13.3\%$ at $n=1000$ to $1\%$ at $n=2744$ and $0$ for larger $n$. For $\nu_0=3/2$ the nugget errors have skewness between $-0.19$ and $0.10$ for $n\ge2744$.

The error law hardly changes these shapes. Under the other three laws, the percentage at the lower bound at $n=1000$ is $10.0\%$--$14.0\%$ for $d=3$, $\nu_0=1/2$, and for $\nu_0=3/2$ and $n\ge2744$ the skewness lies between $-0.22$ and $0.15$ for all four laws, including the skewed Bernoulli law. The skewness of $\xi$ does not enter the first-order behavior of $\log\widehat\tau^2$, which depends on the errors through $\xi_i^2$; the skewness of the average of the $\xi_i^2$ is of order $n^{-1/2}$.

In the score covariance $\Omega_n$, the error law enters only through $\kappa_4(P,t_0)$. To first order this adds $(\E\xi^4-3)/n$ to the variance of $\log\widehat\tau^2$ (\Cref{cor:log}). This correction is negative for the uniform law. For $\nu_0=1/2$ and $n\le10^5$, the law-specific standard deviations in \Cref{sec:sim-results} differ by at most $13\%$ across the four error laws. For $\nu_0=3/2$, at the largest sample size the law-specific standard deviations are within $3.5\%$ of their limits $\{(\E\xi^4-1)/n\}^{1/2}$ for $d=3$, but exceed them by up to $13\%$ for $d=4$ and $39\%$ for $d=5$. Section~\ref{sec:sim-results} reports the information-based and Monte Carlo nugget errors.

\subsection{Supplementary graphical results}
\subsubsection{Estimated standard errors}
\begin{figure}[!htbp]
\centering
\includegraphics[width=0.85\textwidth]{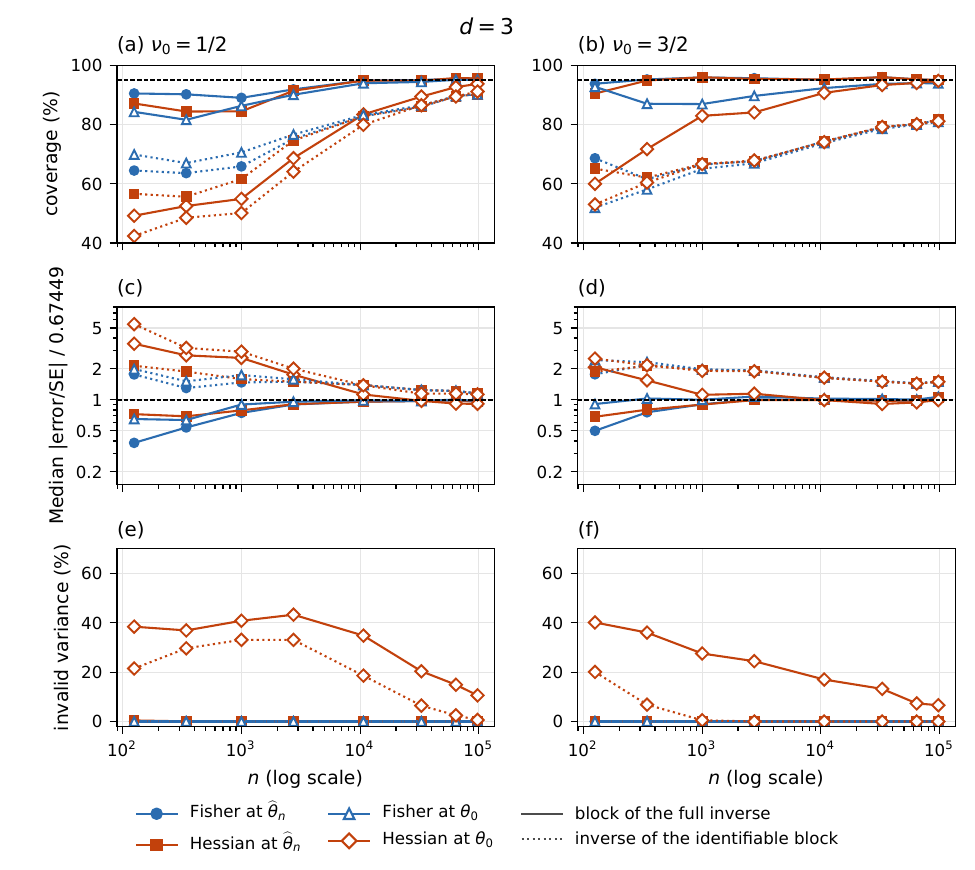}
\caption{Accuracy of the eight standard errors as a function of $n$ for $d=3$, one column for each $\nu_0$, averaged over the identifiable log-parameters $\log\nu$, $\log c$, and $\log\tau^2$ and the four error laws, with $400$ replicates per $n$ and law. Circles and squares: Fisher information and observed negative Hessian at $\widehat\theta_n$; triangles and diamonds: the same at $\theta_0$. The Hessian is computed in $w=(\eta,\nu,\log a,\log t)$ and its inverse is transformed by the delta method. Solid lines: identifiable block of the inverse of the full curvature; dotted lines: inverse of the identifiable block. (a)--(b): coverage of nominal $95\%$ Wald intervals, averaged arithmetically over parameters and laws after conditioning on a valid variance within each cell. (c)--(d): for each parameter, law, and $n$, form $R=\operatorname{median}_r |(\widehat\zeta_r-\zeta_0)/\widehat{\mathrm{se}}_r|/
\Phi^{-1}(0.75)$ over valid replicates, with $\Phi^{-1}(0.75)\approx0.67449$. The plotted value is the geometric mean of these ratios, with equal weights over parameters and laws, on a logarithmic axis. Errors are measured from the truth. The ratio is calibrated to one for standard normal standardized errors; larger values can reflect underestimated standard errors, bias, or nonnormality. Agreement near one concerns the center of the distribution and does not establish correct tail coverage. (e)--(f): percentage of invalid marginal variance estimates, pooled over replicates, parameters, and laws. Validity means a finite, strictly positive corrected marginal variance; invalid estimates are excluded from (a)--(d). Dashed lines mark $95\%$ in (a)--(b) and one in (c)--(d).}
\label{fig:sim-se}
\end{figure}

\begin{figure}[!htbp]
\centering
\includegraphics[width=\textwidth]{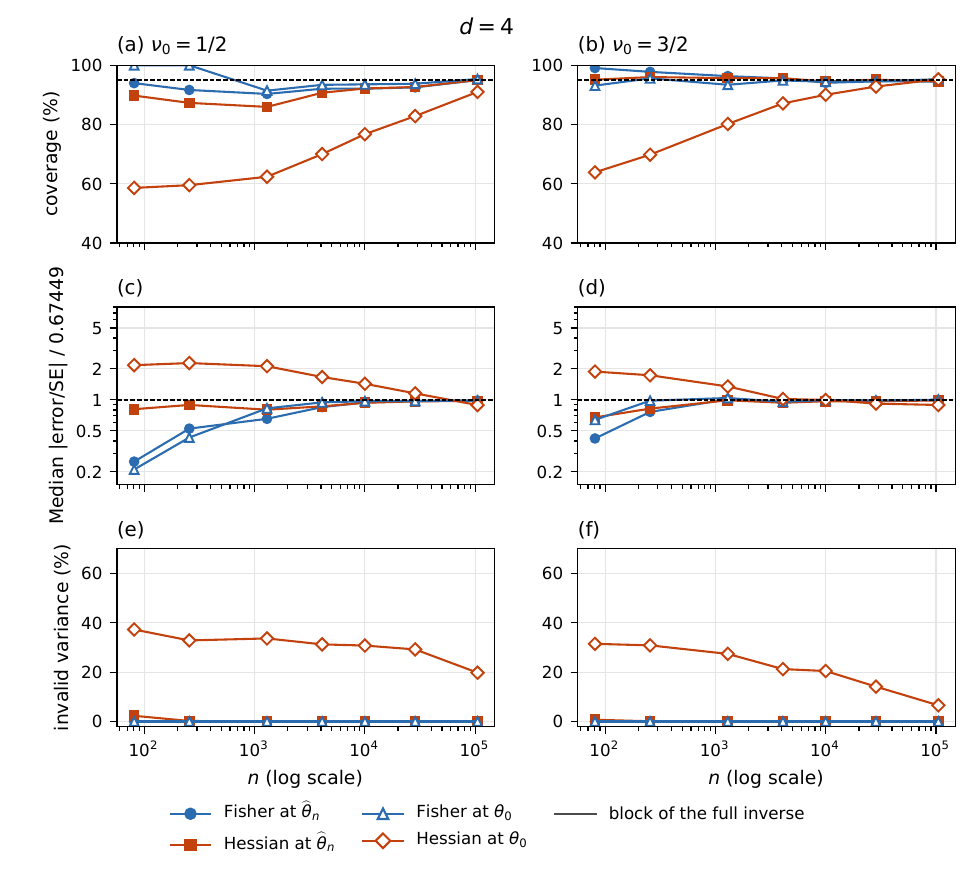}
\caption{As Figure~\ref{fig:sim-se}, for $d=4$. The averages are over the identifiable log-parameters $\log\nu$, $\log\sigma^2$, $\log l$, and $\log\tau^2$ and the four error laws. For $d\ge4$ all four parameters are identifiable, so the two block constructions coincide and only one line is shown for each curvature.}
\label{fig:sim-se-d4}
\end{figure}

\begin{figure}[!htbp]
\centering
\includegraphics[width=\textwidth]{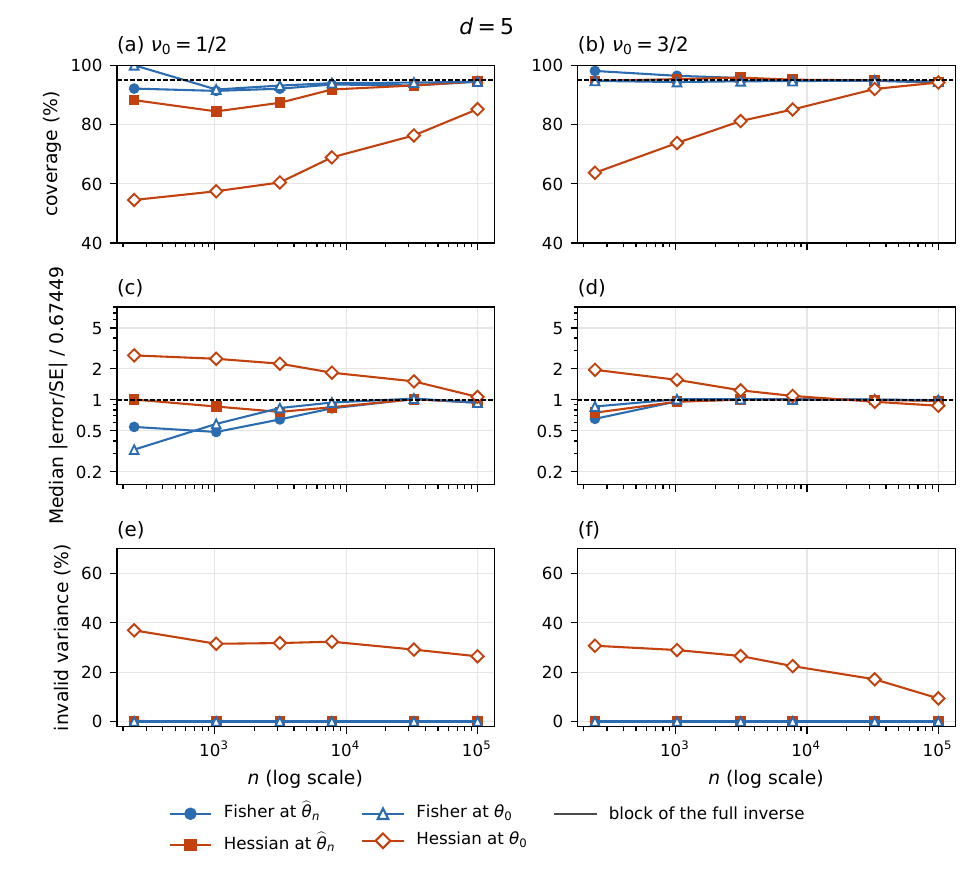}
\caption{As Figure~\ref{fig:sim-se}, for $d=5$. The averages are over the identifiable log-parameters $\log\nu$, $\log\sigma^2$, $\log l$, and $\log\tau^2$ and the four error laws; only one line is shown for each curvature, as in Figure~\ref{fig:sim-se-d4}.}
\label{fig:sim-se-d5}
\end{figure}

\Cref{tab:sim-se} gives the coverage of the eight standard errors at the largest $n$, and \Cref{fig:sim-se,fig:sim-se-d4,fig:sim-se-d5} follow their accuracy over $n$ for $d=3$, $4$, and $5$. The largest differences concern the treatment of the range nuisance and the instability of the Hessian evaluated at the truth.

For $d=3$, the treatment of the range nuisance has a large effect. At $n=97{,}336$ with Gaussian errors, the block methods in Table~\ref{tab:sim-se} cover $\nu$ and $c$ in $87.75\%$--$91.00\%$ of replicates when $\nu_0=1/2$, and $71.00\%$--$76.50\%$ when $\nu_0=3/2$. Appendix~\ref{apdx:sim-range} gives the corresponding full-inverse coverages, standard-deviation ratios, and RMSE comparisons.

With the full inverse, the Fisher information and the observed Hessian at the fit give similar results for $n\gtrsim10^4$ in these experiments. Under Gaussian errors, at the largest $n$ their median standard errors differ by at most $2\%$ and their coverages by at most $0.8$ percentage points. Across all four error laws, the largest coverage difference is $1.0$ percentage point, with paired Monte Carlo standard error $0.61$ points at that comparison; the paired standard errors are at most $0.66$ points across the settings. At smaller $n$ the Fisher information tends to give conservative intervals, and the Hessian shorter ones that can undercover: for $d=5$, $\nu_0=1/2$, $n=1024$, and Gaussian errors, the log-scale interval for $\nu$ from the full inverse covers in all $400$ data sets with the Fisher information at the fit and in $320$ of $400$ ($80\%$) with the Hessian at the fit; under the four error laws the ranges are $99.75\%$--$100\%$ and $78.25\%$--$80.25\%$, and every variance is valid. At the smallest $n$, the corrected variance from the Hessian at the fit is occasionally invalid. For a single setting and error law, the largest share of the $400$ data sets with an invalid corrected marginal variance is $1.5\%$ ($6$ data sets; $d=3$, $\nu_0=1/2$, $n=125$, uniform errors) for $d=3$, $7.5\%$ ($30$; $d=4$, $\nu_0=1/2$, $n=81$, uniform errors) for $d=4$, and $0.5\%$ for $d=5$. In the displayed $d=3$ and $d=4$ settings, the nugget variance is invalid in every such data set. For each dimension, the largest share for any single parameter is also the percentage just given.
Most invalid variances concern the nugget, when a negative cumulant correction exceeds its uncorrected variance at the lower nugget bound. Across the study, this correction turns $98$ previously positive nugget variances negative. The full Hessian itself is not
positive definite in $53$ of $67{,}200$ fits, all with $n\le343$; the largest share for a setting and error law is $2.25\%$ ($9$ of $400$; $d=4$, $\nu_0=1/2$, $n=81$, Gaussian errors). Without the cumulant correction, the largest invalid-variance shares, on the same basis, are $0.75\%$, $2.0\%$, and $0.5\%$ for $d=3,4,5$, respectively. The following counts distinguish failures of the uncorrected curvature, marginal variances, and joint covariance. For the Hessian at the fit, $148$ data sets have at least one invalid corrected marginal variance: $50$ already had an invalid marginal variance before correction, and $98$ became invalid after correction. The corrected joint covariance fails its validity checks in $199$ data sets: $197$ fail positive definiteness and two fail the symmetry check. Of these $199$, $51$ have every corrected marginal variance valid. All these failures occur at $n\le343$. For Fisher information at the fit, three corrected joint covariances fail the symmetry check because the fitted information is numerically singular. They all occur at $d=4$, $\nu_0=1/2$, $n=81$, with the fitted range at or near its upper bound. The cumulant correction does not create any additional failures for this fitted-Fisher method.

Oracle evaluation at the truth gives no uniform finite-sample improvement in these settings. The Fisher information at $\theta_0$ is the same for every replicate, whereas the standard errors from the Fisher information at $\widehat\theta_n$ have coefficients of variation up to $23\%$ at the largest $n$. Nevertheless, for $n>1.5\times10^4$ the mean absolute deviation of the coverage from $95\%$, over all identifiable parameters and laws, is $1.0$ percentage points at $\widehat\theta_n$ and $1.2$ at $\theta_0$ for $d=3$, and $1.3$ and $0.9$ for $d\ge4$. For $d=3$ and smaller $n$, the version at the fit is better: $6.6$ against $11.6$ points for $n<1500$, and $1.8$ against $3.7$ for $1500\le n\le1.5\times10^4$. The standard errors at the fit adapt to the realized data: under Gaussian errors and with the Fisher information at the fit, the standard errors of $\log\nu$ and $\log c$ for $d=3$ have Spearman correlations with the absolute errors between $0.17$ and $0.91$ for $n\le2744$.

The covariance obtained from the inverse Hessian in $w$ at $\theta_0$ is particularly unstable. All numbers here use the coordinate convention in Appendix~\ref{apdx:sim-mle}. Even at $n\approx10^5$ it gives a negative variance estimate for up to $36.5\%$ of
the data sets with Gaussian errors, its standard errors have coefficients of variation between $0.4$ and $1.2$ for the signal parameters, and among the data sets with a positive estimate the coverage falls to $85\%$ for $d=5$ and $\nu_0=1/2$. The Hessian at
$\theta_0$ contains the full random fluctuation of the second derivative of the criterion. In the slowly estimated directions this fluctuation is not small compared with the information. For $d=3$, inverting only the identifiable block of the same matrix gives a positive variance for at least $98.5\%$ of the data sets at the largest $n$, so the instability comes mainly from the range direction. At the fit, the Hessian is instead close to the Fisher information at the fit, as for the variance of a Gaussian sample, where the observed and expected information coincide at the maximum likelihood estimate.

At the largest $n$, the median standard errors of the three full versions (the two curvatures at $\widehat\theta_n$ and the Fisher information at $\theta_0$) are within $4\%$ of each other: over the $88$ combinations of dimension, smoothness, error law, and identifiable parameter, the ratio of the largest to the smallest of the three medians is at most $1.034$. Its maximum is approximately $1.0335$, for $\log c$ with $d=3$, $\nu_0=3/2$, and Gaussian errors.

\begin{table}[!htbp]
\centering\footnotesize
\setlength{\tabcolsep}{4pt}
\caption{Coverage (\%) of nominal $95\%$ Wald intervals from the eight standard errors at the largest $n$ ($97{,}336$, $104{,}976$, and $10^5$ for $d=3,4,5$), Gaussian errors, $400$ replicates. The curvature is the Fisher information (F) or the observed negative Hessian (H), evaluated at the fit $\widehat\theta_n$ or at the truth $\theta_0$. H is computed in $w=(\eta,\nu,\log a,\log t)$ and its inverse is mapped to the reported coordinates by the delta method, as in Appendix~\ref{apdx:sim-mle}. ``Full'': identifiable block of the inverse of the full curvature, as in \eqref{eq:estimatedlogcov}; ``Block'': inverse of the identifiable block of the curvature ($d=3$ only). All $400$ estimates are valid in every column except full H at the truth, whose coverage is conditional on a finite, strictly positive variance estimate. ``Inv.'' is the percentage of the $400$ data sets excluded from that column.
Entries are exact ratios of counts rounded to one decimal place, with ties rounded to even. The Monte Carlo standard error is about $1.1$ percentage points at $95\%$ with $400$ valid estimates; it is larger when fewer estimates are valid.}
\label{tab:sim-se}
\vspace{3pt}
\begin{tabular}{lllrrrrrrrrr}
\toprule
 & & & \multicolumn{5}{c}{Full} & \multicolumn{4}{c}{Block}\\
\cmidrule(lr){4-8}\cmidrule(l){9-12}
 & & & F, $\widehat\theta_n$ & H, $\widehat\theta_n$ & F, $\theta_0$ & H, $\theta_0$ & Inv.
 & F, $\widehat\theta_n$ & H, $\widehat\theta_n$ & F, $\theta_0$ & H, $\theta_0$\\
\midrule
$d=3$ & $\nu_0=1/2$ & $\log\nu$ & $94.8$ & $95.0$ & $94.5$ & $94.2$ & $10.0$ & $89.0$ & $89.2$ & $90.0$ & $91.0$\\
 &  & $\log c$ & $95.5$ & $95.0$ & $94.5$ & $93.3$ & $10.5$ & $89.0$ & $89.0$ & $87.8$ & $90.2$\\
 &  & $\log\tau^2$ & $95.8$ & $95.5$ & $96.2$ & $95.9$ & $8.0$ & $93.0$ & $93.0$ & $93.5$ & $95.5$\\
 & $\nu_0=3/2$ & $\log\nu$ & $96.0$ & $96.0$ & $93.8$ & $95.6$ & $9.2$ & $76.2$ & $76.2$ & $76.2$ & $76.5$\\
 &  & $\log c$ & $96.0$ & $95.8$ & $93.5$ & $93.9$ & $9.2$ & $73.8$ & $73.2$ & $71.2$ & $71.0$\\
 &  & $\log\tau^2$ & $93.8$ & $93.8$ & $93.8$ & $93.8$ & $0.0$ & $93.8$ & $93.8$ & $93.8$ & $93.5$\\
\midrule
$d=4$ & $\nu_0=1/2$ & $\log\nu$ & $94.8$ & $94.8$ & $94.8$ & $93.4$ & $20.8$ &  &  &  & \\
 &  & $\log\sigma^2$ & $94.5$ & $94.5$ & $96.2$ & $90.2$ & $11.0$ &  &  &  & \\
 &  & $\log l$ & $95.5$ & $95.2$ & $95.5$ & $90.8$ & $18.5$ &  &  &  & \\
 &  & $\log\tau^2$ & $94.0$ & $94.5$ & $94.0$ & $91.9$ & $19.8$ &  &  &  & \\
 & $\nu_0=3/2$ & $\log\nu$ & $94.8$ & $95.2$ & $94.2$ & $96.7$ & $9.0$ &  &  &  & \\
 &  & $\log\sigma^2$ & $92.2$ & $92.2$ & $95.0$ & $95.3$ & $9.0$ &  &  &  & \\
 &  & $\log l$ & $92.8$ & $92.5$ & $94.8$ & $96.2$ & $9.0$ &  &  &  & \\
 &  & $\log\tau^2$ & $95.8$ & $95.8$ & $95.8$ & $96.0$ & $0.5$ &  &  &  & \\
\midrule
$d=5$ & $\nu_0=1/2$ & $\log\nu$ & $93.8$ & $93.8$ & $93.5$ & $86.6$ & $36.5$ &  &  &  & \\
 &  & $\log\sigma^2$ & $94.8$ & $94.8$ & $95.0$ & $91.8$ & $12.0$ &  &  &  & \\
 &  & $\log l$ & $93.5$ & $93.8$ & $95.0$ & $84.8$ & $24.2$ &  &  &  & \\
 &  & $\log\tau^2$ & $94.5$ & $94.8$ & $95.5$ & $85.7$ & $35.2$ &  &  &  & \\
 & $\nu_0=3/2$ & $\log\nu$ & $96.8$ & $96.8$ & $94.8$ & $95.7$ & $12.8$ &  &  &  & \\
 &  & $\log\sigma^2$ & $91.2$ & $91.2$ & $93.2$ & $93.4$ & $12.8$ &  &  &  & \\
 &  & $\log l$ & $93.0$ & $93.8$ & $94.8$ & $94.3$ & $12.8$ &  &  &  & \\
 &  & $\log\tau^2$ & $96.5$ & $96.5$ & $96.0$ & $95.3$ & $3.5$ &  &  &  & \\
\bottomrule
\end{tabular}
\end{table}

\subsubsection{Results on the original scale}
\label{apdx:sim-original}

\begin{figure}[!htbp]
\centering
\includegraphics[width=0.96\textwidth]{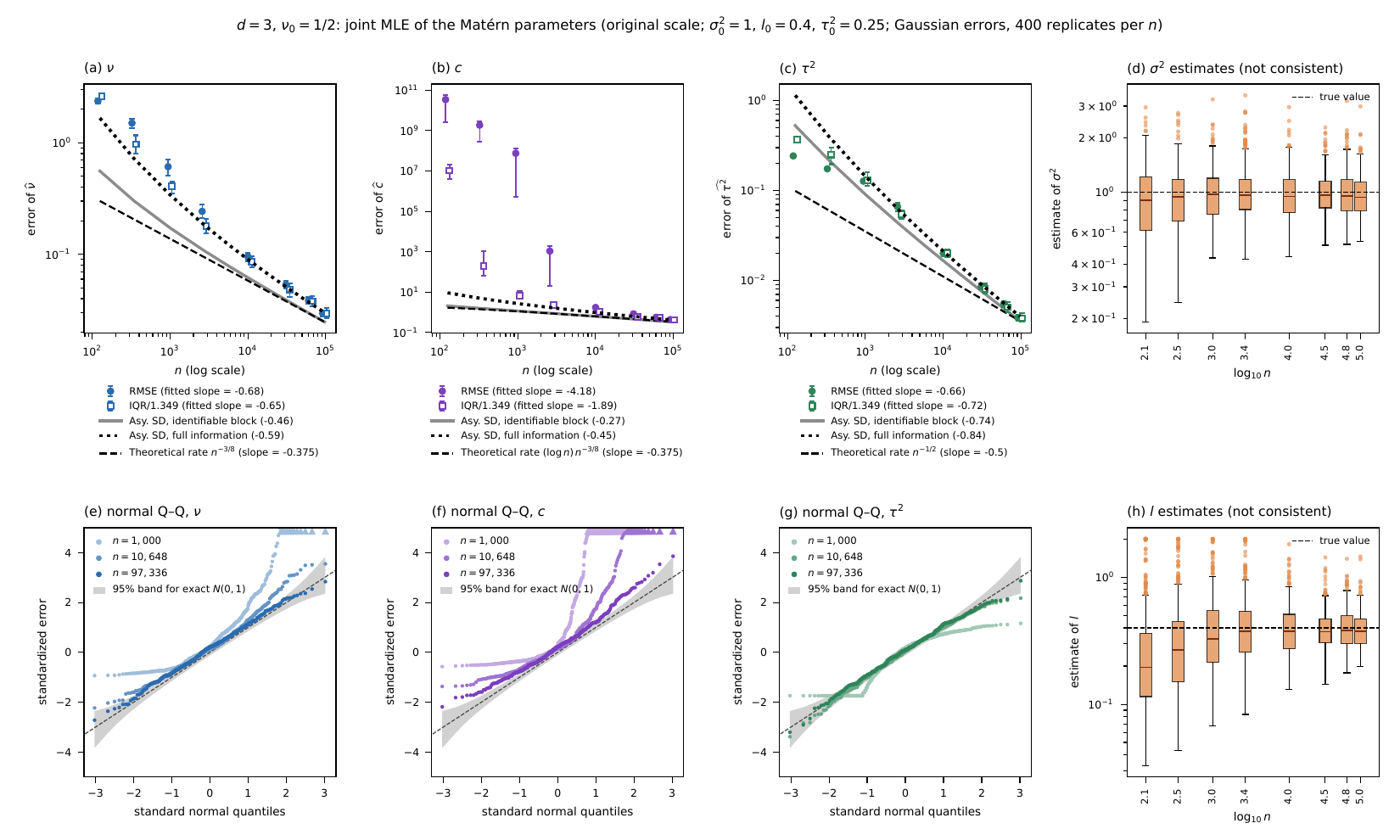}
\caption{Simulation results for $d=3$, $\nu_0=1/2$ on the original parameter scale, with the layout and graphical elements as in \Cref{fig:sim-d3-nu05}.}
\label{fig:sim-raw-d3-nu05}
\end{figure}

\begin{figure}[!htbp]
\centering
\includegraphics[width=0.96\textwidth]{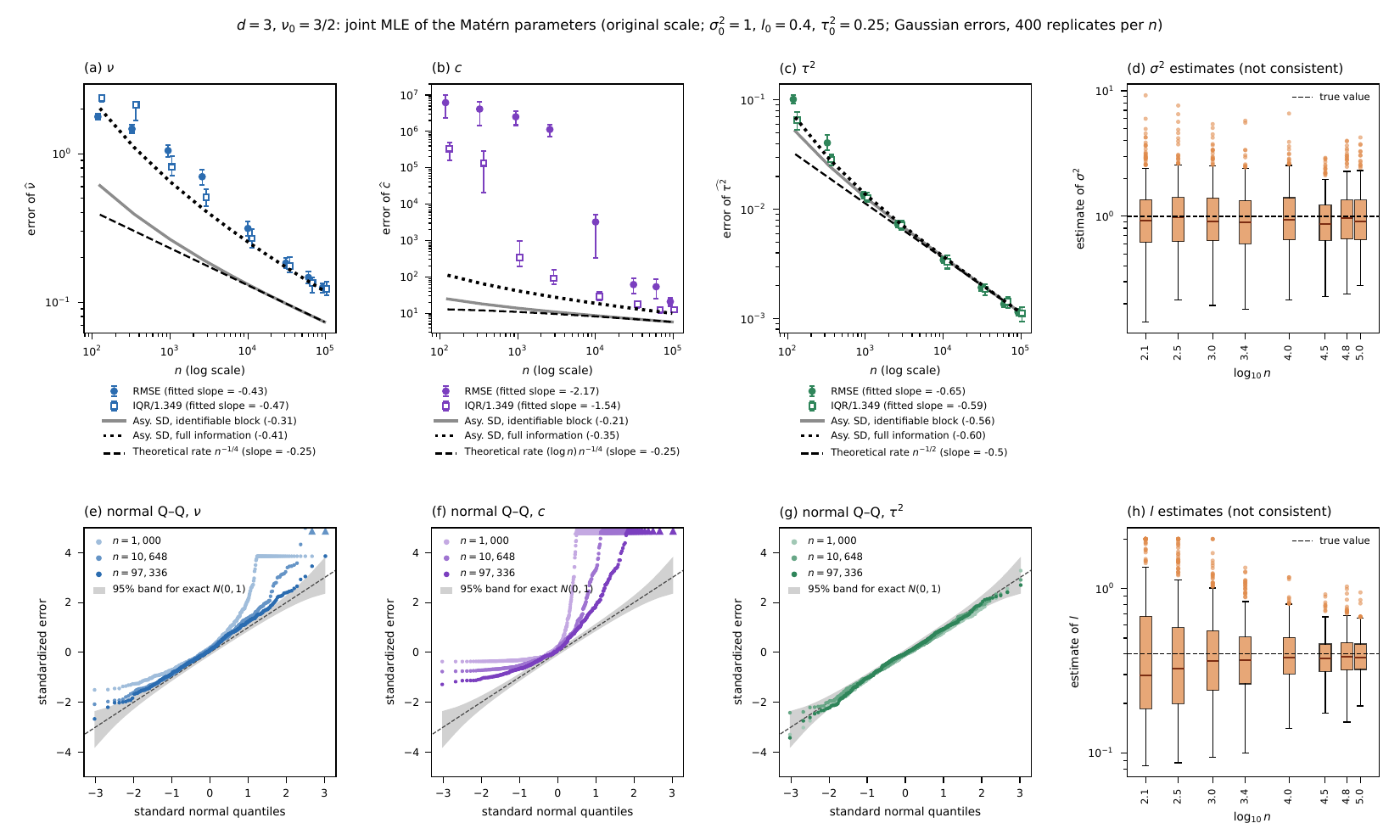}
\caption{Simulation results for $d=3$, $\nu_0=3/2$ on the original parameter scale, with the layout and graphical elements as in \Cref{fig:sim-d3-nu05}.}
\label{fig:sim-raw-d3-nu15}
\end{figure}

\begin{figure}[!htbp]
\centering
\includegraphics[width=0.96\textwidth]{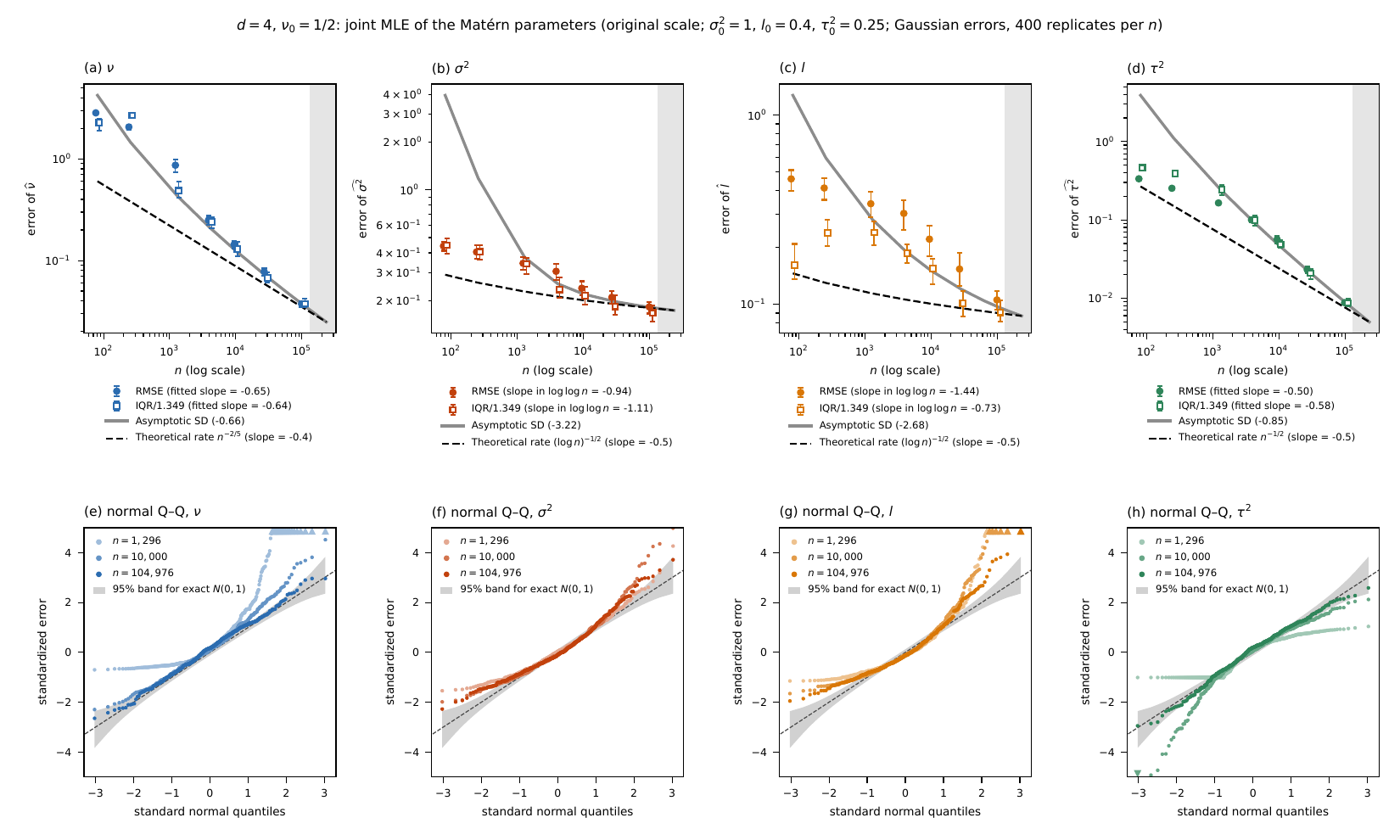}
\caption{Simulation results for $d=4$, $\nu_0=1/2$ on the original parameter scale, with the layout and graphical elements as in \Cref{fig:sim-d4-nu05}.}
\label{fig:sim-raw-d4-nu05}
\end{figure}

\begin{figure}[!htbp]
\centering
\includegraphics[width=0.96\textwidth]{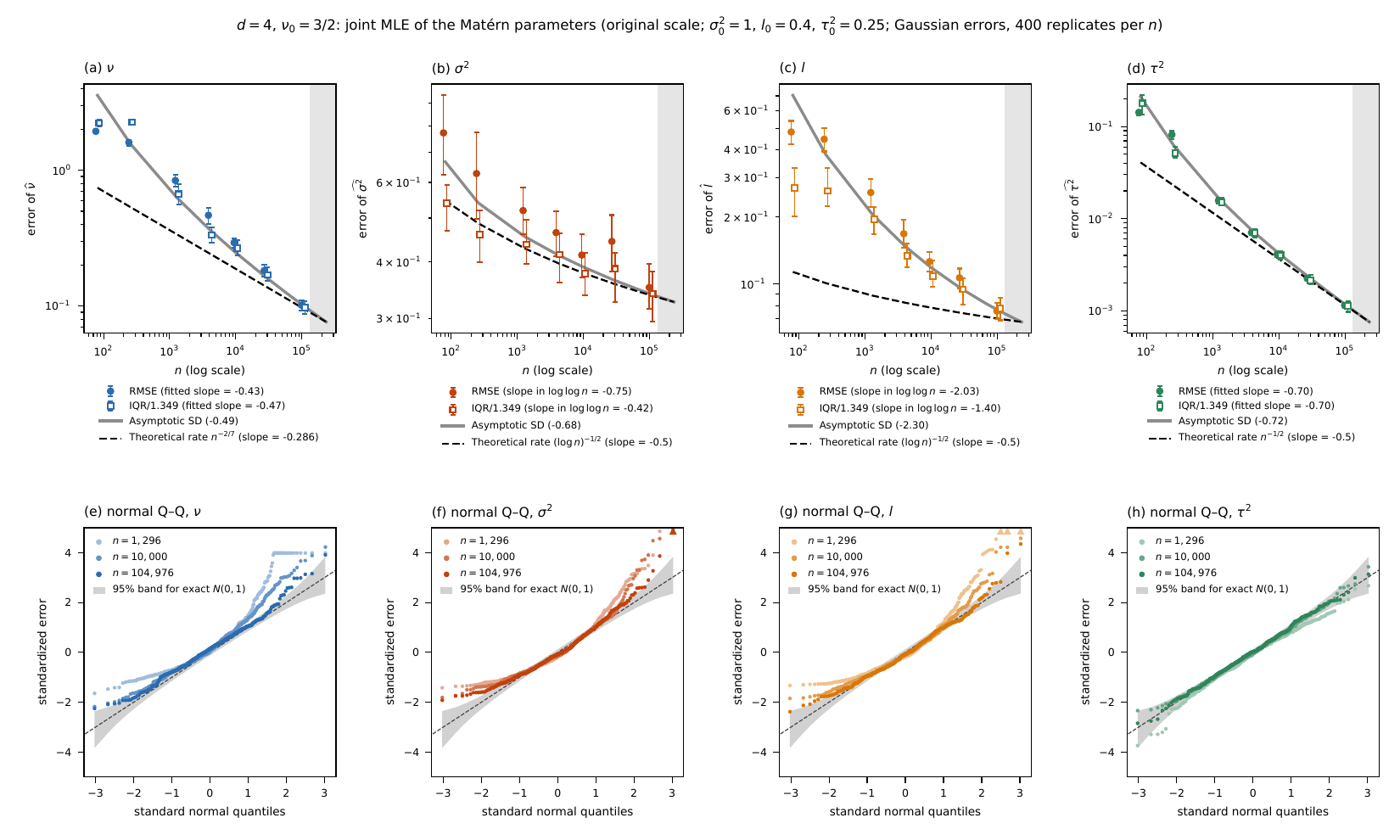}
\caption{Simulation results for $d=4$, $\nu_0=3/2$ on the original parameter scale, with the layout and graphical elements as in \Cref{fig:sim-d4-nu05}.}
\label{fig:sim-raw-d4-nu15}
\end{figure}

\begin{figure}[!htbp]
\centering
\includegraphics[width=0.96\textwidth]{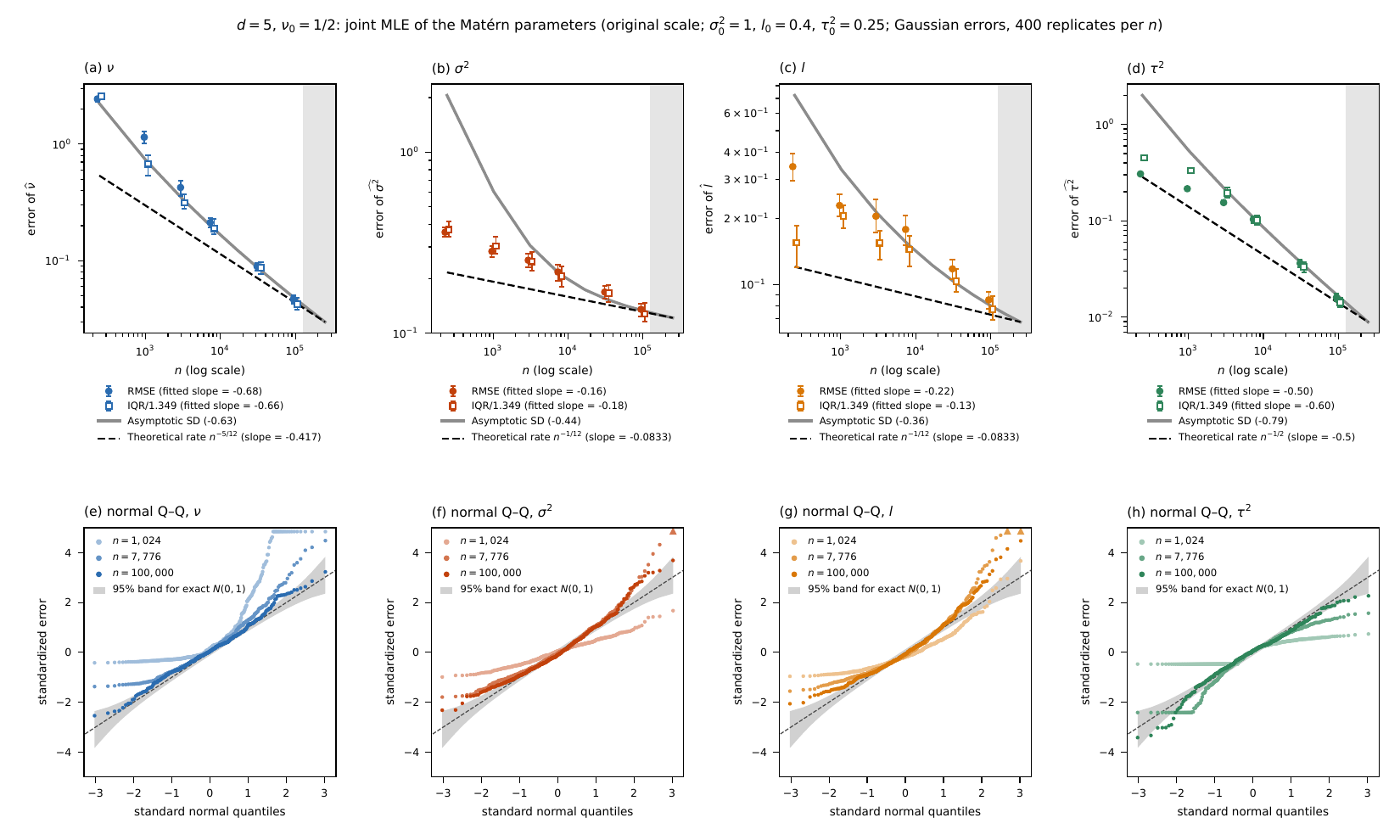}
\caption{Simulation results for $d=5$, $\nu_0=1/2$ on the original parameter scale, with the layout and graphical elements as in \Cref{fig:sim-d4-nu05}.}
\label{fig:sim-raw-d5-nu05}
\end{figure}

\begin{figure}[!htbp]
\centering
\includegraphics[width=0.96\textwidth]{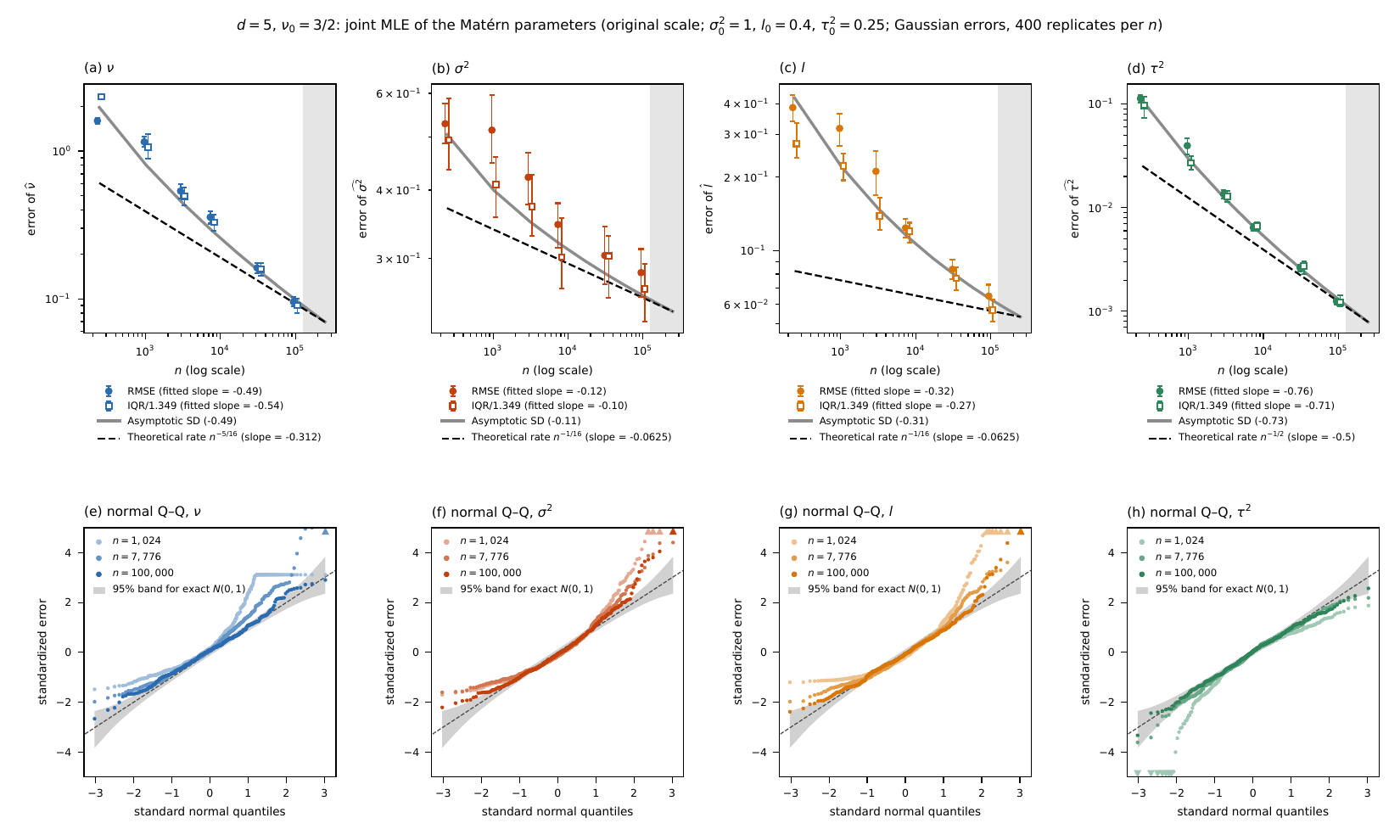}
\caption{Simulation results for $d=5$, $\nu_0=3/2$ on the original parameter scale, with the layout and graphical elements as in \Cref{fig:sim-d4-nu05}.}
\label{fig:sim-raw-d5-nu15}
\end{figure}

\Cref{fig:sim-raw-d3-nu05,fig:sim-raw-d3-nu15,fig:sim-raw-d4-nu05,fig:sim-raw-d4-nu15,fig:sim-raw-d5-nu05,fig:sim-raw-d5-nu15} repeat the Gaussian results in Section~\ref{sec:simulation} on the original parameter scale. The reference rates are unchanged. The inverse information $\{\calI_n^{\mathrm{full}}(\theta_0)\}^{-1}$ is transformed by the reporting Jacobian on both sides; the square roots of the resulting
diagonal entries give the asymptotic standard deviations. The conclusions about the rates are the same. The main visible difference is a heavy right tail of $\widehat c$ at small $n$: when $\widehat\nu$ is near its upper bound, $\widehat c=\widehat\sigma^2\widehat l^{-2\widehat\nu}$ is very large, which inflates the RMSE relative to the robust standard deviation. The log transformation removes most of this effect.

\end{document}